\documentclass[12pt]{report}

\usepackage{tamuconfig}

\usepackage[T1]{fontenc}

\usepackage{graphicx}

\DeclareGraphicsExtensions{.png}
\DeclareGraphicsExtensions{.pdf}

\graphicspath{ {./graphic/} }

\usepackage[hidelinks]{hyperref}

\usepackage{graphicx, bm, amsmath, siunitx, longtable, tabularx, float, mathtools, booktabs, eqnarray, color, xcolor, amsfonts, amssymb, appendix, tensor, subcaption, dsfont, makecell, soul, changepage}
\usepackage[version=4]{mhchem}
\usepackage{lscape}
\usepackage[export]{adjustbox}

\usepackage{calligra}

\makeatletter
\def\blfootnote{\xdef\@thefnmark{}\@footnotetext}
\makeatother

\newcommand{\B}[1]{{\bm #1}}
\newcommand{\T}{^{\mbox{\tiny T}}}
\newcommand{\ds}{\displaystyle}
\newcommand{\dd}{\; \text{d}}
\newcommand{\bdot}[1]{\dot{\B{#1}}}

\newcommand{\ces}{constrained expressions}
\renewcommand{\ce}{constrained expression}
\newcommand{\p}[2]{\prescript{(#1)}{}{#2}}
\newcommand{\C}[1]{\tensor*[^{}_{}]{\mathfrak{C}}{^{}_{#1}}}
\newcommand{\R}{\mathbb{R}}
\newcommand{\andd}{\quad \text{and} \quad}
\newcommand{\phiz}[2]{\prescript{#1}{}{#2}}

\newcommand{\Hs}{\mathds{1}}

\usepackage{placeins, tocloft}
 
\newcommand{\listExamplesname}{List of Examples}
\newlistof[chapter]{Examples}{ex}{\listExamplesname}
\renewcommand{\theExamples}{\thechapter.\arabic{Examples}}
\newcommand{\Examples}[1]{%
    \par\textbf{Example \theExamples : #1}
    \addcontentsline{ex}{Examples}{\protect\numberline{\theExamples}#1}\par}
    
\definecolor{titleColor}{rgb}{0.19, 0.55, 0.91}

\definecolor{backColor}{rgb}{0.19, 0.55, 0.91}

\usepackage{tcolorbox}
\tcbuselibrary{many}

\newtcolorbox{mybox}{colback=white, colframe=titleColor,arc=5pt,width=6.5in,height=9in,halign=center,valign=center}

\newlength\defaultparindent
\AtBeginDocument{\setlength\defaultparindent{\parindent}}
\newtcolorbox{exBox}[1][]{
    colback=backColor!20,
    colframe=titleColor,
    fonttitle=\bfseries,
    colbacktitle=titleColor,
    coltitle=white,
    enhanced,
    breakable,
    arc=5pt,
    attach boxed title to top right ={yshift=-\tcboxedtitleheight/2,yshifttext=-\tcboxedtitleheight/2,xshift=-0.5cm},
    title=#1,
    before upper={\parindent\defaultparindent},
}

\newenvironment{example}[1]
    {%
    \refstepcounter{Examples}
    \begingroup \allowdisplaybreaks\vspace{0.4cm}\begin{exBox}[\Examples{#1}]\relax}
    {\end{exBox}\endgroup}
    
\newtcolorbox{whiteExBox}[1][]{
    colback=white,
    colframe=titleColor,
    fonttitle=\bfseries,
    colbacktitle=titleColor,
    coltitle=white,
    enhanced,
    breakable,
    arc=5pt,
    attach boxed title to top right ={yshift=-\tcboxedtitleheight/2,yshifttext=-\tcboxedtitleheight/2,xshift=-0.5cm},
    title=#1,
    before upper={\parindent\defaultparindent},
}

\newenvironment{whiteExample}[1]
    {%
    \refstepcounter{Examples}
    \begingroup \allowdisplaybreaks\vspace{0.4cm}\begin{whiteExBox}[\Examples{#1}]\relax}
    {\end{whiteExBox}\endgroup}
    
\newtcolorbox{blankExBox}[1][]{
    colback=backColor!20,
    colframe=titleColor,
    fonttitle=\bfseries,
    colbacktitle=titleColor,
    coltitle=white,
    enhanced,
    breakable,
    arc=5pt,
    attach boxed title to top right ={yshift=-\tcboxedtitleheight/2,yshifttext=-\tcboxedtitleheight/2,xshift=-0.5cm},
    title=#1,
    before upper={\parindent\defaultparindent},
}

\newenvironment{blankBox}[1]
    {%
    \begingroup \allowdisplaybreaks\vspace{0.4cm}\begin{blankExBox}[{#1}]\relax}
    {\end{blankExBox}\endgroup}
    
\newtcolorbox{wblankExBox}[1][]{
    colback=white,
    colframe=titleColor,
    fonttitle=\bfseries,
    colbacktitle=titleColor,
    coltitle=white,
    enhanced,
    breakable,
    arc=5pt,
    attach boxed title to top right ={yshift=-\tcboxedtitleheight/2,yshifttext=-\tcboxedtitleheight/2,xshift=-0.5cm},
    title=#1,
    before upper={\parindent\defaultparindent},
}

\newenvironment{wblankBox}[1]
    {%
    \begingroup \allowdisplaybreaks\vspace{0.4cm}\begin{wblankExBox}[{#1}]\relax}
    {\end{wblankExBox}\endgroup}
    
\newcounter{Theorems}
\newcommand{\Theorems}{%
    \par\textbf{Claim \theTheorems}\par
}
    
\newtcolorbox{thrmBox}[1][]{
    colback=backColor!20,
    colframe=titleColor,
    fonttitle=\bfseries,
    colbacktitle=titleColor,
    coltitle=white,
    enhanced,
    breakable,
    arc=5pt,
    attach boxed title to top right ={yshift=-\tcboxedtitleheight/2,yshifttext=-\tcboxedtitleheight/2,xshift=0.5cm},
    title=#1,
    before upper={\noindent\parindent\defaultparindent},
}

\newenvironment{theorem}
    {%
    \refstepcounter{Theorems}
    \begingroup \allowdisplaybreaks\vspace{0.4cm}\begin{thrmBox}[\Theorems]\relax}
    {$\blacksquare$\end{thrmBox}\endgroup}
    
\newcounter{Properties}
\newcommand{\Properties}{%
    \par\textbf{Property \theProperties}\par
}
    
\newtcolorbox{propDefBox}[1][]{
    colback=white,
    colframe=titleColor,
    fonttitle=\bfseries,
    colbacktitle=titleColor,
    coltitle=white,
    enhanced,
    breakable,
    arc=5pt,
    attach boxed title to top right ={yshift=-\tcboxedtitleheight/2,yshifttext=-\tcboxedtitleheight/2,xshift=0.5cm},
    title=#1,
}

\newenvironment{property}
    {%
    \refstepcounter{Properties}
    \begingroup \allowdisplaybreaks\vspace{0.4cm}\begin{propDefBox}[\Properties]\relax}
    {\end{propDefBox}\endgroup}

\newcounter{Definitions}
\newcommand{\Definitions}{%
    \par\textbf{Definition \theDefinitions}\par
}

\newenvironment{definition}
    {%
    \refstepcounter{Definitions}
    \begingroup \allowdisplaybreaks\vspace{0.4cm}\begin{propDefBox}[\Definitions]\relax}
    {\end{propDefBox}\endgroup}

\begin{document}

\renewcommand{\tamumanuscripttitle}{The Theory of Functional Connections \linebreak A journey from theory to application}

\renewcommand{\tamupapertype}{Dissertation}

\renewcommand{\tamufullname}{Hunter Reed Johnston}

\renewcommand{\tamudegree}{Doctor of Philosophy}
\renewcommand{\tamuchairone}{Daniele Mortari}

\renewcommand{\tamumemberone}{John E. Hurtado}
\newcommand{\tamumembertwo}{Srinivas Vadali}
\newcommand{\tamumemberthree}{Yalchin Efendiev}
\renewcommand{\tamudepthead}{Srinivas Vadali}

\renewcommand{\tamugradmonth}{August}
\renewcommand{\tamugradyear}{2021}
\renewcommand{\tamudepartment}{Aerospace Engineering}

%
%
%
%


\providecommand{\tabularnewline}{\\}

\begin{titlepage}
\begin{center}
\MakeUppercase{\tamumanuscripttitle}
\vspace{4em}

A \tamupapertype

by

\MakeUppercase{\tamufullname}

\vspace{4em}

\begin{singlespace}

Submitted to the Office of Graduate and Professional Studies of \\
Texas A\&M University \\

in partial fulfillment of the requirements for the degree of \\
\end{singlespace}

\MakeUppercase{\tamudegree}
\par\end{center}
\vspace{2em}
\begin{singlespace}
\begin{tabular}{ll}
 & \tabularnewline
& \cr
Chair of Committee, & \tamuchairone\tabularnewline
Committee Members, & \tamumemberone\tabularnewline
 & \tamumembertwo\tabularnewline
 & \tamumemberthree\tabularnewline
Head of Department, & \tamudepthead\tabularnewline

\end{tabular}
\end{singlespace}
\vspace{3em}

\begin{center}
\tamugradmonth \hspace{2pt} \tamugradyear

\vspace{3em}

Major Subject: \tamudepartment \par
\vspace{3em}
Copyright \tamugradyear \hspace{.5em}\tamufullname 
\par\end{center}
\end{titlepage}
\pagebreak{}

%

\chapter*{ABSTRACT}
\addcontentsline{toc}{chapter}{ABSTRACT} 

\pagestyle{plain} 
\pagenumbering{roman} 
\setcounter{page}{2}

\indent The Theory of Functional Connections (TFC) is a general methodology for functional interpolation that can embed a set of user-specified linear constraints. The functionals derived from this method, called \emph{constrained expressions}, analytically satisfy the imposed constraints and can be leveraged to transform constrained optimization problems to unconstrained ones. By simplifying the optimization problem, this technique has been shown to produce a numerical scheme that is faster, more accurate, and robust to poor initialization. The content of this dissertation details the complete development of the Theory of Functional Connections. First, the seminal paper on the Theory of Functional Connections is discussed and motivates the discovery of a more general formulation of the constrained expressions. Leveraging this formulation, a rigorous structure of the constrained expression is produced with associated mathematical definitions, claims, and proofs. Furthermore, the second part of this dissertation explains how this technique can be used to solve ordinary differential equations providing a wide variety of examples compared to the state-of-the-art. The final part of this work focuses on unitizing the techniques and algorithms produced in the prior sections to explore the feasibility of using the Theory of Functional Connections to solve real-time optimal control problems, namely optimal landing problems. 

\pagebreak{}
%
%
%
%

\chapter*{DEDICATION}
\addcontentsline{toc}{chapter}{DEDICATION}  

\begin{center}
\vspace*{\fill}
To my mother and father.\\ And to the friends (C, L, \& M) who have been there from the beginning,\\ and those who I've met along the way.
\end{center}
\vspace{1.5in}
\begin{center}
All things inevitably come to an end.\\
Some day the machine stops running.\\
We can share paths for a while, but\\
ultimately we all have our own\\
separate destinations.\\
\vspace{0.1in}
--- Unravel, \emph{ColdWood Interactive}
\vspace*{\fill}
\end{center}

\pagebreak{}

%
%
%
%

\chapter*{ACKNOWLEDGMENTS}
\addcontentsline{toc}{chapter}{ACKNOWLEDGMENTS}  

The path to completing this document involved not just numbers and equations but loving and caring human beings --- family, friends, teachers, and mentors. Although I encountered many roadblocks, dead ends, and unfavorable terrain, you, knowingly or unknowingly, have propelled me. While I could easily fill this page with names, I restrain over the fear of forgetting just one. However, to those to who I am referring, you know who You are ... \\\phantom{a}
\begin{center}
    \includegraphics[width=.25\linewidth]{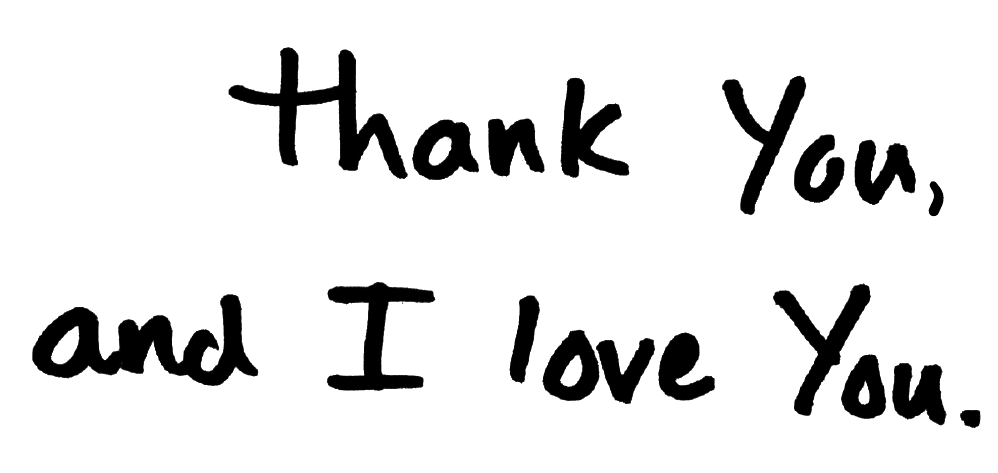}
\end{center}

Regardless, a few people were fundamental to my education and the completion of this document, and I would like to identify them by name specifically.

First, Dr. Daniele Mortari, my advisor and friend. Thank you for taking a chance to bring me in as one of your graduate students. My four years at Texas A\&M were memorable, to say the least, and I will cherish the brainstorming session we've had, ALL of the meals we've shared, and our conversations about literature, life, and philosophy.

Second, my labmates and true friends, (soon to be Dr.) Carl Leake and Dr. Stoian Borissov. You both have given me unmatched support in dealing with the rough terrain of graduate school and graduate student life. Thanks not only for challenging me and providing unmatched feedback, but for also pulling me away from graduate life and distracting me with good food, good company, loud drums, and above all, many MANY ``coffees.''

Next, my colleagues from the University of Arizona, (also, soon to be Dr.) Enrico Schiassi and Dr. Roberto Furfaro. Thank you for your amazing collaboration on many projects and for welcoming me into your research group during my month-long visit to Arizona. Specifically, thank you Enrico for video chatting with me to watch F1 races throughout this crazy year of COVID.

Additionally, my committee members Drs. John E. Hurtado, Rao Vadali, and Yalchin Efendiev. Of the two I've been fortunate enough to take classes with, I would like to thank for their inspiration and guidance; your classes are two of my most memorable ones from my time at Texas A\&M. Additionally, I thank Dr. Efendiev for the many Saturday mornings he spent with the TFC research group and his unmatched guidance and feedback.

Lastly, I would like to thank my NASA/NSTRF collaborators, Drs. Chris D'Souza and Martin Lo (who was also my Visiting Technologist Experience host at JPL). Thank you for your support and guidance in both research and my career goals. Additionally, thank you for our long conversations when it seemed the world was falling down around us.

\pagebreak{}
%
%
%
%

\chapter*{CONTRIBUTORS AND FUNDING SOURCES}
\addcontentsline{toc}{chapter}{CONTRIBUTORS AND FUNDING SOURCES}  

\subsection*{Contributors}
This work was supported by a dissertation committee consisting of Daniele Mortari (advisor) and John E. Hurtado and Srinivas Vadali of the Department of Aerospace Engineering, and Yalchin Efendiev of the Department of Mathematics.

The Theory of Functional Connections was collaboratively developed by Daniele Mortari (advisor), Carl Leake (Ph.D. candidate), and Hunter Johnston (author/Ph.D. candidate). To clarify the major contributions of each, the following figure is included.

\begin{center}
    \includegraphics[width=\linewidth]{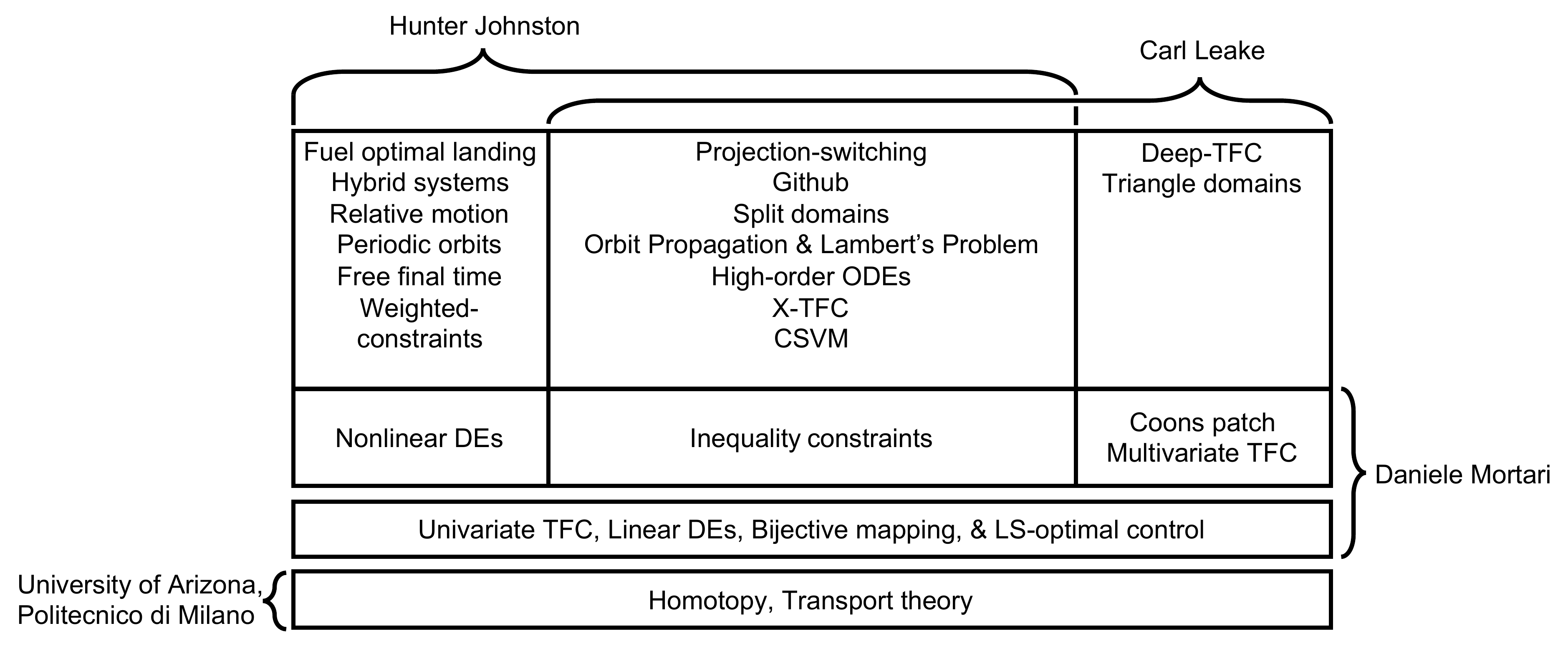}
\end{center}

All other work conducted for the dissertation was completed by the student independently.
\subsection*{Funding Sources}
Graduate study was supported by teaching and research assistantships from Texas A\&M University from August 2017 - August 2019, and by the NASA Space Technology Research Fellowship, Johnston [NSTRF 2019] Grant \#: 80NSSC19K1149, from August 2019 - August 2021.
\pagebreak{}
%
%
%
%


\chapter*{NOMENCLATURE}
\addcontentsline{toc}{chapter}{NOMENCLATURE}  


\hspace*{-1.25in}
\vspace{12pt}
\begin{spacing}{1.0}
	\begin{longtable}[htbp]{@{}p{0.35\textwidth} p{0.62\textwidth}@{}}
	   ELM & Extreme Learning Machine \\ [2ex]
	   FEM & Finite Element Method \\ [2ex]
	   LS-SVM & Least-Squares Support Vector Machine \\ [2ex]
	   NN & Neural Network \\ [2ex]
	   NSTRF  & NASA Space Technology Research Fellowship\\ [2ex]
	   ODE  & Ordinary differential equation\\ [2ex]
	   PDE  & Partial differential equation\\
	   [2ex]
	   PMP & Pontryagin Minimum Principle \\
	   [2ex]
	   SVM & Support Vector Machine \\ [2ex]
	   TFC  & Theory of Function Connections\\ [2ex]
	   TPBVP  & two-point boundary-value problem\\ [2ex]
	   X-TFC & Extreme Theory of Functional Connections \\ [2ex]
	   $c$  & Slope in the linear map for the independent variable that maps the basis function domain to the problem domain. \\ [2ex]
	   $b$  & The square-root of the slope in the linear map for the independent variable that maps the basis function domain to the problem domain. $b^2 = c$ \\ [2ex]
	   $\C{i}$ & Constraint operator for the $i$-th constraint \\ [2ex]
	   $\delta_{ij}$ & Kronecker delta \\ [2ex]
	   $g(x)$ & Free function $\mathbb{R}\mapsto\mathbb{R}$. Note that a superscript may be used to denote the free function for a specific dependent variable, e.g., $g^u(x)$ is the free function for the dependent variable $u$. \\ [2ex]
	   $\mathbb{J}$ & Jacobian matrix of the loss vector function $\mathbb{L}$\\ [2ex]
	   $\kappa_i(x)$ & Portion of the $i$-th constraint of the independent variable that does not contain the dependent variable. \\ [2ex]
	   $\mathbb{L}$ & Loss vector function $\mathbb{R}^m\mapsto\mathbb{R}^n$ \\ [2ex]
	   $\rho_i(x,g(x))$ & Projection functional for the $i$-th constraint on the independent variable. \\ [2ex]
	   $\phi_i(x)$ & Switching function for the $i$-th constraint on the independent variable. \\ [2ex]
	   $\mathbb{R}$ & Field of real numbers \\ [2ex]
	   $\mathbb{S}_{ij}$ & Support matrix \\ [2ex]
	   $\tau$ & Alternative definition of the basis function independent variable. Note, this is used when $z$ is used as an independent variable. \\ [2ex]
	   $\mathbb{Z}^+$ & Set of positive integers \\ [2ex]
	   $z$ & Basis function independent variable. Note, this is replaced with $\tau$ in some cases. \\ [2ex]
	   $\mathds{1}(x,x_1)$ & Heaviside function, $\mathbb{R}\mapsto\mathbb{R}$ \\ [2ex]
	   $\mathds{1}_0(x)$ & Heaviside function where $x_1 = 0$ \\ [2ex]
	\end{longtable}
\end{spacing}

\pagebreak{}

%
%
%
%

\phantomsection
\addcontentsline{toc}{chapter}{TABLE OF CONTENTS}  

\begin{singlespace}
\renewcommand\contentsname{\normalfont} {\centerline{TABLE OF CONTENTS}}

\setcounter{tocdepth}{4} 

\setlength{\cftaftertoctitleskip}{1em}
\renewcommand{\cftaftertoctitle}{%
\hfill{\normalfont {Page}\par}}

\tableofcontents

\end{singlespace}

\pagebreak{}


\phantomsection
\addcontentsline{toc}{chapter}{LIST OF FIGURES}  

\renewcommand{\cftloftitlefont}{\center\normalfont\MakeUppercase}

\setlength{\cftbeforeloftitleskip}{-12pt} 
\renewcommand{\cftafterloftitleskip}{12pt}

\renewcommand{\cftafterloftitle}{%
\\[4em]\mbox{}\hspace{2pt}FIGURE\hfill{\normalfont Page}\vskip\baselineskip}

\begingroup

\begin{center}
\begin{singlespace}
\setlength{\cftbeforechapskip}{0.4cm}
\setlength{\cftbeforesecskip}{0.30cm}
\setlength{\cftbeforesubsecskip}{0.30cm}
\setlength{\cftbeforefigskip}{0.4cm}
\setlength{\cftbeforetabskip}{0.4cm}



\listoffigures

\end{singlespace}
\end{center}

\pagebreak{}

%
\phantomsection
\addcontentsline{toc}{chapter}{LIST OF TABLES}  

\renewcommand{\cftlottitlefont}{\center\normalfont\MakeUppercase}

\setlength{\cftbeforelottitleskip}{-12pt} 

\renewcommand{\cftafterlottitleskip}{1pt}

\renewcommand{\cftafterlottitle}{%
\\[4em]\mbox{}\hspace{2pt}TABLE\hfill{\normalfont Page}\vskip\baselineskip}

\begin{center}
\begin{singlespace}

\setlength{\cftbeforechapskip}{0.4cm}
\setlength{\cftbeforesecskip}{0.30cm}
\setlength{\cftbeforesubsecskip}{0.30cm}
\setlength{\cftbeforefigskip}{0.4cm}
\setlength{\cftbeforetabskip}{0.4cm}

\listoftables 

\end{singlespace}
\end{center}

\pagebreak{}  

%
\phantomsection
\addcontentsline{toc}{chapter}{LIST OF EXAMPLES}  

\renewcommand{\cftextitlefont}{\center\normalfont\MakeUppercase}

\setlength{\cftbeforeextitleskip}{-12pt} 

\renewcommand{\cftafterextitleskip}{1pt}

\renewcommand{\cftafterextitle}{%
\\[4em]\mbox{}\hspace{2pt}EXAMPLE\hfill{\normalfont Page}\vskip\baselineskip}

\begin{center}
\begin{singlespace}

\setlength{\cftbeforechapskip}{0.4cm}
\setlength{\cftbeforesecskip}{0.30cm}
\setlength{\cftbeforesubsecskip}{0.30cm}
\setlength{\cftbeforefigskip}{0.4cm}
\setlength{\cftbeforetabskip}{0.4cm}
\setlength{\cftbeforeExamplesskip}{0.4cm}

\setlength{\cftExamplesindent}{1.5em}
\setlength{\cftExamplesnumwidth}{2.3em}

\listofExamples

\end{singlespace}
\end{center}
\endgroup
\pagebreak{}  

%

\pagestyle{plain} 
\pagenumbering{arabic} 
\setcounter{page}{1}

\chapter{INTRODUCTION}
The topics presented in this dissertation can be split into three distinct areas which flow from the general formulation of the Theory of Functional Connections (TFC) (Chapter \ref{chap:tfc_intro} and Chapter \ref{chap:tfc_general}) to its application to the solution of differential equations (Chapter \ref{chap:ode}) and finally leveraging the method to solve optimal control problems (Chapter \ref{chap:opt_con}), namely the energy-optimal landing (Chapter \ref{chap:eol}) and fuel-optimal landing (Chapter \ref{chap:fol}) problems. Ultimately, the goal of this work is to develop a fast, accurate, and robust numerical system to solve problems relevant in aerospace engineering; however, the development of TFC and its initial application to differential equations are vital stepping stones in this effort since each chapter is heavily reliant on those coming before. 

Since this work covers the full journey from the initial theory first published by Mortari \cite{U-TFC} in 2017 to applications in aerospace engineering, I have opted to provide multiple literature reviews directly before the chapters they pertain to. For example, Chapter \ref{chap:tfc_intro} provides an overview of the mathematical concept of interpolation and how they have been utilized. Similarly, the beginning of Chapter \ref{chap:ode} reviews current numerical techniques available to solve ordinary differential equations, and Chapter \ref{chap:opt_con} provides
background on the techniques to solve optimal control problems. 

The following sections of this chapter provide a summary of the work in this dissertation. This is provided to give the reader insight into the structure of the document and highlight the new contributions made to current literature.

\begin{wblankBox}{Part 1 --- Theory}
\noindent\textbf{Chapter \ref{chap:tfc_intro}: An Introduction to the \emph{Theory of Functional Connections}}

This chapter introduces the reader to the original work on TFC, at that time, published simply as the Theory of Connections \cite{U-TFC}. Through this review, TFC is presented in the broader context of interpolation to show this method is a generalized interpolation scheme enabling functional interpolation. This provides the mathematical framework to generate functionals (functions of functions) that analytically satisfy all imposed linear constraints and represent the real-valued set of functions satisfying the constraints. Additionally, to familiarize the reader with the specific vocabulary of TFC and how the method is used, specific examples are provided with increasing complexity. The scope of these examples are two-fold as they 1) provide the reader with concrete, step-by-step derivations and 2) develop an understanding of the theory such that the general formulation of the univariate framework, provided in Chapter \ref{chap:tfc_general}, is easily understood. After these examples, an ad-hoc approach is developed to handle inequality type constraints. Then the chapter concludes with a section highlighting how the functionals derived through the TFC framework can be over-constrained.

\noindent\textbf{Chapter \ref{chap:tfc_general}: A General Formulation of the Univariate \emph{Theory of Functional Connections}}

Leveraging the intuition of the TFC method provided in Chapter \ref{chap:tfc_intro}, this chapter provides a rigorous definition of TFC, and the terminology used and is an expanded version of the general formulation first published by Leake, Johnston, and Mortari \cite{M-TFC-new}. Whereas Chapter \ref{chap:tfc_intro} highlights the consistent structure of the interpolating functionals, this chapter utilizes this discovery to define the terms, identify their associated mathematical properties, and ultimately provide straightforward proofs on the existence and uniqueness of these functionals. These proofs have further implications when the expressions are used to solve differential equations, which is covered in Chapter \ref{chap:ode}. Moreover, the development in this section facilitates the generalization of TFC to $n$-dimensions.
\end{wblankBox}

\begin{wblankBox}{Part 2 --- Application}
\noindent\textbf{Chapter \ref{chap:ode}: Application to the Solution of Ordinary Differential Equations}

As mentioned earlier, the use of TFC expressions to solve ordinary differential equations is one of the three major pillars of the work presented in this dissertation. Consequently, careful attention is paid to developing the numerical framework and consistent notation throughout to allow ease of implementation. Similar to the examples provided in Chapter \ref{chap:tfc_intro} to derive the interpolating functionals, this chapter provides example solutions of differential equations starting with linear ordinary differential equations and culminates in the solution of systems of coupled, nonlinear ordinary differential equations. The examples presented (i.e., the Lane-Emden equations, perturbed orbit propagation,  perturbed Lambert's problem, etc.) are meant to guide the reader in implementing the method and provide solutions to some relevant equations in the field of science and engineering. Following these examples, two unpublished additions to the numerical application of TFC are introduced. First, the method is adapted for the solution of hybrid systems --- where the dynamics exhibit discrete jumps over the solution domain. Following this, a numerical technique to handle unspecified time, i.e., unknown final time problems, is introduced and highlighted with examples. Lastly, the author provides some numerical applications to problems of over-constrained differential equations.
\end{wblankBox}

\begin{wblankBox}{Part 3 --- Optimal Control}
\noindent\textbf{Chapter \ref{chap:opt_con}: Use for Real-time Optimal Controllers in Aerospace Systems}

This chapter contains an overview of the current techniques to solve optimal control problems, emphasizing real-time implementation. After distinguishing between the direct and indirect methods to solve optimal control problems, the first-order necessary conditions for optimality are derived from first principles using the indirect method. This derivation is used as a background for the reader. It precisely shows where the TFC approach fits into the solution of the resultant system of equations by analytically satisfying a portion of these equations. Additionally, this chapter serves as a high-level literature review for the specific problems presented in Chapter \ref{chap:eol} and Chapter \ref{chap:fol}, where the indirect method and TFC are ultimately used to solve these problems. Finally, further insight is provided with comparisons between TFC, spectral, and collocation methods already studied in the context of optimal control theory.

\noindent\textbf{Chapter \ref{chap:eol}: Energy-Optimal Landing}

In this chapter, the full three degree-of-freedom energy-optimal landing problem is formulated according to two different TFC based schemes, Outer-loop and Single-loop. The Outer-loop relies on an external optimizer to solve for the final time in the problem (i.e., MATLAB \verb"fsolve()"), while the Single-loop incorporates all necessary conditions into a single TFC scheme. First, these schemes are compared to the feedback solution for the constant gravity case of this problem to ensure the method's accuracy. With this said, both TFC schemes are formulated ``blind'' to the feedback form to study the effects of the free-final time on the algorithm. Finally, the two developed approaches are studied through a Monte Carlo simulation for varying initial conditions and compared to a similar implementation using a spectral method.

\noindent\textbf{Chapter \ref{chap:fol}: Fuel-Optimal Landing}

With the increasing interest in human spaceflight operations to the Moon, Mars, and possibly beyond, maximizing the amount of payload that can be landed on these bodies is of utmost importance. Optimizing the landing trajectory and minimizing fuel consumption over this landing sequence is one of the many avenues to achieve this. In all, this fuel-optimal landing problem is still an active area of research in the field of aerospace engineering. Therefore, this chapter is the culmination of the prior chapters, and herein, the three degree-of-freedom fuel-optimal landing problem is formulated and solved using the TFC framework. Similar to Chapter \ref{chap:eol}, all necessary conditions are derived using the indirect method, which poses the problem as a nonlinear system of equations that is solved using the TFC framework. Ultimately, the resulting algorithm is used to solved for these trajectories and compared to current state-of-the-art, commercially available methods.
\end{wblankBox}

\noindent\textbf{Chapter \ref{chap:summary}: Summary and Conclusions}

While this dissertation stretches from the basis of the analytical method to embed constraints (the Theory of Functions Connections) to the numerical solution of optimal control problems, it was infeasible to include everything that has been done with regards to this new theory. Therefore, along with drawing the major takeaways of the work presented in this dissertation, I have also devoted the final section of this work (Section \ref{sec:outro}) to comment on the state of TFC at the date of publication. This includes a comprehensive list of all available publications with a summary of the major contributions and results. Additionally, I have also noted the current work in progress and the key problems moving forward. Therefore, after reading, this section can be leveraged by new researchers as a path to interesting and fruitful topics in the greater field of TFC.

\begin{mybox}
\begin{center}
\begin{huge}
Part 1\\
\vspace{0.2in}
Theory
\end{huge}

\vspace{0.5in}

How beautiful this was, when it was new.\\
And how beautiful it still is, even\\
though time has made it different.\\
\vspace{0.1in}
--- Unravel, \emph{ColdWood Interactive}
\end{center}
\end{mybox}
\pagebreak{}

%

\chapter{AN INTRODUCTION TO THE \emph{THEORY OF FUNCTIONAL CONNECTIONS}\label{chap:tfc_intro}}
Interpolation is the mathematical process of estimating an unknown function's values within the range of $k$ given data points, called constraints, provided by some unknown continuous process. Occasionally, in engineering and science, a function is expressed as data points, whether through sampling or experimentation. These data points represent a finite series (or reconstruction) of the governing process (function) at specific independent variable values. Given this data, it is often desired to estimate the function value at some point in between the given data. In another case, the function might be known but is defined by a complicated set of equations that are computationally inefficient to evaluate. In this context, it may be more desirable to approximate the function using a simpler function (with some associated interpolation error) that is easier to evaluate.

Our first mathematical understanding of interpolation can be traced back to elementary algebra. We were first introduced to interpolation when we looked for the numerical description of the line passing through two points $(x_1, y_1)$ and $(x_2, y_2)$. Recall the equation takes the form,
\begin{equation}\label{eq:interp_1}
    y(x) = y_1 + (y_2 - y_1) \frac{x - x_1}{x_2 - x_1},
\end{equation}
where $x$ is a point along the domain. However, as we look to include more data points, we must substitute this linear interpolation method with other techniques such as polynomial interpolation, where the entire function is described by a polynomial, or spline interpolation, where the function is described by piecewise polynomials between data points. Regardless, these techniques provide us with general interpolation schemes to include a given number of points. As one specific example, a popular technique for polynomial interpolation is Lagrange polynomials\footnote{The author notes that the name ``Lagrange polynomials'' is an academic misnomer since the formula was actually first discovered by Edward Waring \cite{interp_waring} in 1779, then by Leonhard Euler in 1783, and eventually Joseph-Louis Lagrange in 1795.}
\begin{equation*}
    L_k(x) = \sum_{i=0}^k y_j \, \phi_j(x)
\end{equation*}
where the polynomial $L_k(x)$ passes through each set of $k+1$ data points $(x_j, y_j)$, and $\phi_j(x)$ are polynomials based on the equation,
\begin{equation*}
    \phi_j(x) = \prod_{\begin{smallmatrix}0\leq i\leq k\\i\neq j\end{smallmatrix}}  \frac{x - x_i}{x_j - x_i} = \frac{(x - x_0)}{(x_j - x_0)} \hdots \frac{(x - x_{j-1})}{(x_j - x_{j-1})} \frac{(x - x_{j+1})}{(x_j - x_{j+1})} \hdots \frac{(x - x_{k})}{(x_j - x_{k})},
\end{equation*}
where $0 \leq j \leq k$. For example, if two data points are selected ($k=1$), then the formula reduces to our simple description of a line,
\begin{equation*}
    L_2(x) = y_1 \Big(\frac{x - x_2}{x_1-x_2}\Big) + y_2 \Big(\frac{x - x_1}{x_2 - x_1}\Big).
\end{equation*}
Creating the interpolating polynomial in this way makes it easy to see how the constraints of $y_1$ and $y_2$ are satisfied. The $\phi_j(x)$ terms multiplying the constraint terms act as continuous switches that evaluate to $1$ at the constraint they are associated with and $0$ when evaluated at all other constraints. In the case of the polynomial $L_2(x)$, we can see that the term multiplying $y_1$ is $1$ when $x = x_1$ and is $0$ when $x = x_2$. Furthermore, by simple algebraic manipulation, we can see that this equation is identical to Equation \eqref{eq:interp_1}.

At this point, some questions may arise:
\begin{itemize}
    \item What if we have data associated with derivatives as well?
    \item What if we are interested in all possible functions that interpolate these point and derivative values?
    \item What if the function of interest is based on a combination of data measurements?
\end{itemize}

In general, a method that provides answers to these questions is interested in the interpolation of functions rather than just points: in other words, a method for ``functional interpolation.'' Whereas Lagrange polynomials provide the polynomial expression that passes through all given points, the method of interest here is a functional\footnote{Also known as a higher-order function or a function of functions.} that represents all possible functions satisfying some given data set conditions, where these ``conditions'' are not limited to points. The questions mentioned and the search for a functional interpolation framework led to the development of the Theory of Functional Connections (TFC) \footnote{This theory was originally published under the name ``Theory of Connections.'' However, this name conflicted with a specific theory in differential geometry and was not the most accurate description of the functional interpolation method. Therefore, in 2019, this name was changed to the ``Theory of Functional Connections,'' to highlight the tie to functional interpolation and the fact that it provides \emph{all} functions satisfying a set of linear constraints in rectangular domains of $n$-dimensional space.} in the seminal paper by Mortari \cite{U-TFC}. 

The foundation of this work is built on a straightforward method to derive analytical expressions (or functionals), which represent the set of all functions satisfying a specified combination of constraints. In his original paper, Mortari identified three unique ways to build these functionals, including linear, additive, and rational forms. 
\begin{align*}
    y (x, g(x)) &= g (x) (x - x_0) + y_0 \qquad && \qquad{\rm (linear)} \\[5pt]
    y (x, g(x)) &= g (x) + [y_0 - g (x_0)] \qquad && \qquad{\rm (additive)} \\[5pt]
    y (x, g(x)) &= \dfrac{g (x)}{g (x_0)}\, y_0 \qquad && \qquad{\rm (rational)} 
\end{align*}
However, the additive form proved to be the most fruitful and therefore, the name the ``Theory of Functional Connections'' refers to functional interpolation using the additive form.\footnote{Note that linear, additive, and rational forms are equivalent through functional transformations. For example, by performing the logarithm of the rational formulation, an additive formulation is obtained. The additive formulation can also be recovered from the linear formulation by simply setting the function $g(x)$ in the additive formulation as $x \, g(x)$. Therefore, the additive form was adopted as the main formalism because of its simplicity.} In this approach, the resulting functional was coined as a ``constrained expression'' since they constrain the functional to analytically satisfy the imposed constraints. Mortari's original work \cite{U-TFC} provided examples of constraints in $k$ points, constraints in $k$ points and derivatives, and relative constraints. It hinted at the idea of linear constraints, something that this dissertation introduces along with a unified notation and associated claims. In all, the original work produced a generalized interpolation technique, as will soon be demonstrated. In fact, in the cases where only function values are considered, i.e., point constraints, it is easy to see that Lagrange polynomials are a specific case of the more general TFC.

While the idea of functional interpolation is not new, prior methods only existed for a class (or sub-class) of functions and not all of function space \cite{interp_1,interp_2,interp_3,interp_x}. More current techniques also include distributed approximating functions (DAFs) \cite{interp_4,interp_5}, which use Hermite DAFs and Sinc DAFs. However, the theory discovered by Mortari \cite{U-TFC} is the first interpolation technique not restricted to a specific class of functions. In the following section, a summary of the major points in this discovery is provided, along with a step-by-step development of the functional interpolation method called TFC. In all, what was discovered in this seminal paper is leveraged to develop a general technique to handle general linear constraints.

\section{An introduction to constrained expressions}
The idea for TFC started with an attempt to derive an expression for all functions passing through the specific point $(x_0,y_0)$. Using algebra, one can easily define all \emph{straight lines} with the equation, $y(x,m) = m (x - x_0) + y_0$, where $y(x_0) = y_0$ and $m$ represents the constant value of the slope. Yet, the slope could be defined by a function, $m(x): \R \to \R$, where the only restriction on $m(x)$ is it must be defined at $x_0$. By making this modification, the expression now becomes a functional, $y(x,m(x)): \R \to \R$ that represents all functions that evaluate to $y_0$ at $x = x_0$. Although this functional always satisfies the constraints, and is thereby a valid constrained expression,\footnote{A rigorous definition of a constrained expression is provided in Chapter \ref{chap:tfc_general}.} the derivation process did not provide a clear path to add multiple constraints. Therefore, a different approach is desired.

Said approach came from the realization that the \emph{additive} form of the \ce\ describes all functions passing through the point defined earlier. Let $g(x): \R \to \R$, be a user defined function that is defined at $x_0$, then the expression,
\begin{equation}\label{eq:s2a_additive}
y(x,g(x)) = g(x) + (y_0 - g(x_0)),
\end{equation}
produces a similar result to the constrained expression $y(x,m(x)) = m(x) (x - x_0) + y_0$, however, the function $g(x)$ appears linearly, which we will soon find to be invaluable. The next step was to determine the general methodology to derive Equation \eqref{eq:s2a_additive}. Without changing the constrained expression, the latter term could be multiplied with the value 1, or in fact, any function $s(x)$ such that $s(x_0) = 1$. Let us define this function as simply $s(x) = 1$. Adding this to Equation \eqref{eq:s2a_additive} leads to,
\begin{equation*}
    y(x,g(x)) = g(x) + s(x)(y_0 - g(x_0))
\end{equation*}
Analyzing this equation, the term $y_0 - g(x_0)$ is constant for a a given $g(x)$ and is the only term containing information of the constraint point, let us denote this constant by $\eta$, and insert it into the equation and rearrange,
\begin{equation}\label{eq:s2a_CE1}
    y(x,g(x)) =  g(x) + s(x) \eta.
\end{equation}
It becomes clear that in order to determine the coefficient $\eta$ this equation must be evaluated at the constraint point ($x_0,y_0$). This realization was a pivotal moment in the discovery of the constrained expression, and it quickly followed that a general expression to Equation \eqref{eq:s2a_CE1} could be written as,
\begin{equation}\label{eq:s2a_ogTFC}
    y(x,g(x)) = g(x) + \sum_{j=1}^k s_j(x) \eta_j
\end{equation}
where again $g(x): \R \to \R$ is the free function. Additionally, the summation term is a linear combination of the functions, $s_j(x): \R \to \R$, which we will call support functions, and the $\eta_j$ coefficients, which we have already seen capture the constraint information. In fact, from this general expression we can quickly return to Equation \eqref{eq:s2a_additive}. For this problem, the number of constraints $k$ is one, so the expression becomes
\begin{equation*}
     y(x,g(x)) = g(x) + s(x) \eta.
\end{equation*}
Evaluating the expression at the point ($x_0,y_0$), solving for $\eta$, and inserting it back into the expression above yields,
\begin{equation*}
     y(x,g(x)) = g(x) + \frac{s(x)}{s(x_0)}(y_0 - g(x_0)).
\end{equation*}
Defining $s(x) = 1$ this equation reduces to Equation \eqref{eq:s2a_additive}. Finally, Equation \eqref{eq:s2a_CE1} facilitates the derivation of constrained expressions for even more complicated sets of constraints.

\section{Adding a second constraint}
The next logical step is to find the constrained expression passing through two points. While in the previous derivation $s(x)$ was set loosely and without explanation, this example provides insight into how the support function, $s(x)$, must be chosen. Using Equation \eqref{eq:s2a_ogTFC} as a template, let us derive an expression such that $y(x_1) = y_1$ and $y(x_2) = y_2$. 

\begin{example}{Constraints at two points}\label{ex:s2a_twoPoints}
Since there are two constraints, Equation \eqref{eq:s2a_ogTFC} takes the form, 
\begin{equation}\label{eq:s2a_2pointexample}
    y(x,g(x)) = g(x) + s_1(x)\eta_1 + s_2(x)\eta_2.
\end{equation}
Evaluating this expression at the two constraint points (e.g., for the first constraint, this means evaluating the right hand side of the equation at $x_1$ and setting it equal to $y_1$), leads to a system of equations,
\begin{align*}
    y_1 = g(x_1) + s_1(x_1)\eta_1 + s_2(x_1)\eta_2 \\
    y_2 = g(x_2) + s_1(x_2)\eta_1 + s_2(x_2)\eta_2
\end{align*}
where the only unknowns are the $\eta_k$ coefficients. Writing these in vector-matrix form leads to a system of equations for these coefficients,
\begin{equation*}
    \begin{Bmatrix} y_1 -  g(x_1) \\ y_2 - g(x_2)\end{Bmatrix} = \begin{bmatrix} s_1(x_1) & s_2(x_1) \\ s_1(x_2) & s_2(x_2) \end{bmatrix} \begin{Bmatrix} \eta_1 \\ \eta_2\end{Bmatrix}.
\end{equation*}
By inverting the matrix composed of the support functions evaluated at the constraints, we can solve for the unknown coefficients $\eta_1$ and $\eta_2$. This highlights the major restriction on our definition of the support functions since to solve for $\eta$ coefficients, the matrix \emph{must} be invertible. In other words, the columns, and therefore the support functions, \emph{must} be linearly independent.

Continuing with this example, by selecting $s_1(x) = 1$ and $s_2(x) = x$, which are linearly independent, the system of equations becomes,
\begin{equation*}
    \begin{Bmatrix} y_1 -  g(x_1) \\ y_2 - g(x_2)\end{Bmatrix} = \begin{bmatrix} 1 & x_1 \\ 1 & x_2 \end{bmatrix} \begin{Bmatrix} \eta_1 \\ \eta_2\end{Bmatrix}.
\end{equation*}
Solving this system yields the $\eta_1$ and $\eta_2$ values,
\begin{align*}
    \eta_1 &= \frac{1}{x_2 - x_1} \Big( x_2[y_1 - g(x_1)] - x_1[y_2 - g(x_2)]\Big)\\
    \eta_2 &=  \frac{1}{x_2 - x_1} \Big( [y_2 - g(x_2)] - [y_1 - g(x_1)]\Big)
\end{align*}
which can then be substituted into Equation \eqref{eq:s2a_2pointexample} to produce the constrained expression,
\begin{align*}
    y(x,g(x)) = g(x) &+ \frac{1}{x_2 - x_1} \Big( x_2[y_1 - g(x_1)] - x_1[y_2 - g(x_2)]\Big)\\ &+ \frac{x}{x_2 - x_1} \Big( [y_2 - g(x_2)] - [y_1 - g(x_1)]\Big).
\end{align*}
While it may seem there is an excessive use of parenthesis, these are used to highlight that the terms $y_1 - g(x_1)$ and $y_2 - g(x_2)$ show up in the latter two terms, and thus, the equation can be rearranged by collecting on these two terms. Doing this leads to the familiar result obtained in the original derivation in Reference \cite{U-TFC},
\begin{equation}\label{eq:s2a_canon_form_2point}
    y(x,g(x)) = g(x) + \frac{x_2 - x}{x_2 - x_1}\Big(y_1 - g(x_1)\Big) +  \frac{x - x_1}{x_2 - x_1}\Big(y_2 - g(x_2)\Big).
\end{equation}
\end{example}
Using the constrained expression from Equation \eqref{eq:s2a_canon_form_2point}, it is easy to see that if this functional is evaluated at either $x_1$ or $x_2$, the corresponding constraint value of $y_1$ or $y_2$ is obtained regardless of the function $g(x)$. Further analyzing this equation, we might ask, what happens if we select the simplest expression for the free function such that $g(x) = 0$? If $g(x) = 0$, then Equation \eqref{eq:s2a_canon_form_2point} reduces to
\begin{equation*}
    y = \frac{x_2 - x}{x_2 - x_1} \, y_1 +  \frac{x - x_1}{x_2 - x_1} \, y_2,
\end{equation*}
which the reader may recognize as the Lagrange polynomial for two points discussed earlier. This result should come as no surprise since the original goal was to derive a functional that represents all possible functions passing through the given set of constraints, or in this simple case, points. In the context of our constrained expression, the Lagrange polynomial is the simplest interpolating function of the functional $y(x,g(x))$, when $1$ and $x$ are chosen as support functions. While this generalization is insightful, TFC should not be taken as a simple generalization of Lagrange polynomials. The following examples highlight that point constraints are merely the beginning of the theory.

\section{The structure of the constrained expression}
The prior example hinted at an interesting form of the constrained expression but did not give a mechanized method to arrive at the end result. This section explores Equation \eqref{eq:s2a_canon_form_2point}, specifically, and the terms dictating the constraints, to bring to light a structure within the constrained expression that can be utilized to create the aforementioned mechanized method. Moreover, the said method will ultimately reveal itself to be a unified, consistent way to develop constrained expressions for many different types of constraints.

First, notice that the latter two terms in the constrained expression consist of two unique parts, 1) a term composed of only the support functions and their values at the constraint locations and 2) a term composed of the constraint condition and the function $g(x)$ evaluated at this constraint condition. As an example, consider the first of these terms from Equation \eqref{eq:s2a_canon_form_2point},
\begin{equation*}
   \underbrace{\frac{x_2 - x}{x_2 - x_1}}_\text{\normalsize$\phi_1(x)$}\underbrace{(y_1 - g(x_1))}_\text{\normalsize$\rho_1(x,g(x))$}.
\end{equation*}
The first part of this structure we will call the switching function, $\phi_j(x)$. This function is defined such that it is equal to 1 when evaluated at the constraint it is referencing, and equal to 0 when evaluated at all other constraints. In our example, when evaluating the switching function, $\phi_1(x)$, at the constraint it is referencing it is equal to 1 (i.e., $\phi_1(x_1) = 1$), and when it is evaluated at the other constraints it is equal to 0 (i.e., $\phi_2(x_1) = \dfrac{x_1 - x_1}{x_2 - x_1} = 0$). 

The second part of the structure, $\rho_1(x,g(x))$, is called the projection functional. In this case, the projection functional is simply the difference between the constraint value and the free function evaluated at that constraint; however, for more complex constraints this is not always the case. We choose the name projection functional because it ``projects'' the free function onto the set of functions that vanish at the constraint. Continuing with our example, the projection functional, $\rho_1(x,g(x))$, is simply the difference between the constraint $y(x_1)=y_1$ and the free function evaluated at the constraint point, $g(x_1)$. This structure is important, as it shows up in all other constraint types we consider. Additionally, notice what happens to the projection functional if $g(x)$ satisfies the constraint,
\begin{property}{}\label{prop:proj1}
The projection functionals for constraints at a point are always equal to zero if the free function, $g(x)$, is selected such that it satisfies the associated constraint. 
\end{property}
This simply means that if $g(x)$ were defined such that $g(x):= y_1$, the entire term would reduce to $0$. This property will be utilized in mathematical claims later in the dissertation.

Based on this structure, consider an alternative structure to Equation \eqref{eq:s2a_ogTFC}, which leverages the fact that the constrained expression can be built as a sum of switching functions and projection functionals expressed as,
\begin{equation}\label{eq:s2a_uniCeAlt}
    y(x,g(x)) = g(x) + \sum_{j=1}^k \phi_j(x) \rho_j(x,g(x)).
\end{equation}
First, based on their composition, the projection functionals, $\rho_j(x,g(x))$, are trivial to derive, but the switching functions, $\phi_j$, require some attention. From the definition of the switching functions, these functions must go to 1 at their associated constraint and 0 at all other constraints. As a result, the following algorithm can be used to derive the switching functions for a set of $k$ constraints:
\begin{wblankBox}{Algorithm to derive the terms of Equation \eqref{eq:s2a_uniCeAlt}}
\begin{enumerate}
    \item Choose the $k$ linearly independent support functions, $s_k$.
    \item Write each switching function as a linear combination of the support functions with $k$ unknown coefficients.
    \item Based on the switching function definition, write a system of equations to solve for the unknown coefficients. 
\end{enumerate}
\end{wblankBox}
To validate this approach, let us rederive the constrained expression from Example \ref{ex:s2a_twoPoints}.

\begin{example}{Constraints at two points (Alternative derivation)} 
Given two constraints, Equation \eqref{eq:s2a_uniCeAlt} takes the form,
\begin{equation}\label{eq:s2a_twoPointsAlt}
    y(x,g(x)) = g(x) + \phi_1(x)\rho_1(x,g(x)) + \phi_2(x)\rho_2(x,g(x))
\end{equation}
where the switching functions are of the form,
\begin{equation*}
    \phi_1(x) = s_i(x)\alpha_{i1} \quad \text{and} \quad \phi_2(x) = s_i(x)\alpha_{i2}
\end{equation*}
for some as yet unknown coefficients $\alpha_{ij}$; note that in the previous expression, and throughout this book, the Einstein summation convention\footnote{For example, $a_i b_i = \B{a} \T \B{b}$ for the inner product.} is used to improve readability. Additionally, the projection functionals are,
\begin{equation*}
    \rho_1(x,g(x)) = y_1 - g(x_1) \quad \text{and} \quad \rho_2(x,g(x)) = y_2 - g(x_2).
\end{equation*}
Now, the definition of the switching function is used to come up with a set of equations. For example, the first switching function has the two equations,
\begin{equation*}
    \phi_1(x_1) = 1, \quad \phi_1(x_2) = 0.
\end{equation*}
The equations for all the switching functions can be combined into the compact form,
\begin{align*}
    \begin{bmatrix} s_1(x_1) & s_2(x_1) \\ s_1(x_2) & s_2(x_2) \end{bmatrix} \begin{bmatrix} \alpha_{11} & \alpha_{12}  \\ \alpha_{21} & \alpha_{22} \end{bmatrix} &= \begin{bmatrix} \phi_1(x_1) & \phi_2(x_1) \\ \phi_1(x_2) & \phi_2(x_2) \end{bmatrix}.
\end{align*}
This equation offers us our first visible connection to the original technique to derive constrained expressions. Notice that the support function matrix, i.e., the matrix composed of the support functions, is identical to the matrix multiplying the $\eta$ coefficients in our prior example. Therefore, it still holds that the support functions must be linearly independent. Therefore, and in order to mirror Example \ref{ex:s2a_twoPoints}, let us define the support functions as, $s_1(x) = 1$ and $s_2(x) = x$, and the matrix of $\phi_j$ is identity by definition. Solving the system provides the values of the coefficients $\alpha_{ij}$
\begin{align*}
    \begin{bmatrix} 1 & x_1 \\ 1 & x_2 \end{bmatrix} \begin{bmatrix} \alpha_{11} & \alpha_{12}  \\ \alpha_{21} & \alpha_{22} \end{bmatrix} &= \begin{bmatrix} 1 & 0 \\ 0 & 1 \end{bmatrix} \\
     \begin{bmatrix} \alpha_{11} & \alpha_{12} \\ \alpha_{21} & \alpha_{22} \end{bmatrix} &= \begin{bmatrix} 1 & x_1 \\ 1 & x_2\end{bmatrix}^{-1} = \frac{1}{x_2-x_1}\begin{bmatrix} x_2 & -x_1 \\ -1 & 1 \end{bmatrix}.
\end{align*}
Substituting the constants back into the switching functions and simplifying yields,
\begin{equation*}
    \phi_1 = \frac{x_2s_1(x) - s_2(x)}{x_2 - x_1} = \frac{x_2 - x}{x_2 - x_1} \andd \phi_2 = \frac{s_2(x) - x_1s_1(x)}{x_2 - x_1} = \frac{x - x_1}{x_2 - x_1}.
\end{equation*}
Lastly, by substituting the switching functions along with the associated projection functionals back into Equation \eqref{eq:s2a_twoPointsAlt}, the constrained expression becomes,
\begin{equation*} 
    y(x,g(x)) = g(x) + \underbrace{\frac{x_2 - x}{x_2 - x_1}}_{\text{\normalsize$\phi_1(x)$}}\underbrace{(y_1 - g(x_1))}_\text{\normalsize$\rho_1(x,g(x))$} +  \underbrace{\frac{x - x_1}{x_2 - x_1}}_\text{\normalsize$\phi_2(x)$}\underbrace{(y_2 - g(x_2))}_\text{\normalsize$\rho_2(x,g(x))$}. 
\end{equation*}
\end{example}
The result is identical to Equation \eqref{eq:s2a_canon_form_2point} and should come as no surprise as it is simply an exploitation of the structure of the constrained expression. At this point, it may be unclear the benefit of using Equation \eqref{eq:s2a_uniCeAlt} to construct constrained expressions; however, the following section provides in-depth examples building up to general, linear-type constraints where the true power of the switching-projection notation will become obvious.

\section{Examples using the switching-projection form of the
constrained expression}
While our motivating example in the prior section was vital to our understanding of the \ce\ and its underlying structure, it is limited to the application of constraints at a point. However, the insight and methodology built up in this example can be applied to various linear constraints. The following sections provide specific examples of the application of Equation \eqref{eq:s2a_uniCeAlt}. Admittedly, one could derive all of the following examples using the original form of the constrained expression, Equation \eqref{eq:s2a_uniCeAlt}, albeit with more difficulty.

\subsection{Point and derivative constraints}
In our first example, we take a small step by including derivative constraints into the constrained expression. The reader will see that this does not add any complexity when using the TFC approach.
\begin{example}{Point and derivative constraints}
Consider the following set of point and derivative constraints defined by,
\begin{equation*}
    y(0) = 1, \quad y_x(1) = 2, \quad y(2) = 3,
\end{equation*}
where the notation $y_x := \frac{\dd y}{\dd x}$ is used for the derivative of the function $y(x)$ with respect to $x$. The projection functionals are immediate and can be written as,
\begin{equation*}
    \rho_1(x,g(x)) = 1 - g(0), \quad \rho_2(x,g(x)) = 2 - g_x(1), \quad \rho_3(x,g(x)) = 3 - g(2).
\end{equation*}
Now, the only terms that remain are the switching functions. Recall that our definition of the switching functions in terms of the support functions, $s_i(x)$, and the unknown coefficients, $\alpha_{ij}$, is $\phi_j(x) = s_i(x)\alpha_{ij}$, and the expressions for the three switching functions are, $\phi_1 = s_i(x)\alpha_{i1}$, $\phi_2(x) = s_i(x) \alpha_{i2}$, and $\phi_3(x) = s_i \alpha_{i3}$. Now, this definition of the switching function is used to come up with a set of equations. For example, the first switching function has the three equations,
\begin{align*}
    \phi_1(0) &= s_1(0)\alpha_{11} + s_2(0)\alpha_{21} + s_3(0)\alpha_{31} = 1 \\ \frac{\partial \phi_1}{\partial x}(1) &= s_{1_x}(1)\alpha_{11} + s_{2_x}(1)\alpha_{21} + s_{3_x}(1)\alpha_{31} = 0 \\ \phi_1(2) &= s_1(2)\alpha_{11} + s_2(2)\alpha_{21} + s_3(2)\alpha_{31} =0.
\end{align*}
where the reader should notice that the second equation involves the derivative of the switching function and is associated with the derivative constraint $y_x(1) = 2$. It is convenient to represent these equations in matrix form,
\begin{equation*}
\begin{bmatrix} s_1(0) & s_2(0) & s_3(0) \\ s_{1_x}(1) & s_{2_x}(1) & s_{3_x}(1) \\ s_1(2) & s_2(2) & s_3(2)\end{bmatrix} \begin{bmatrix} \alpha_{11} \\ \alpha_{21} \\ \alpha_{31} \end{bmatrix} = \begin{Bmatrix} 1 \\ 0 \\ 0 \end{Bmatrix}.
\end{equation*}
Adding the expressions of the other two switching functions, the set of equations becomes,
\begin{equation*}
    \begin{bmatrix} s_1(0) & s_2(0) & s_3(0) \\ s_{1_x}(1) & s_{2_x}(1) & s_{3_x}(1) \\ s_1(2) & s_2(2) & s_3(2)\end{bmatrix} \begin{bmatrix} \alpha_{11} & \alpha_{12} & \alpha_{13} \\ \alpha_{21} & \alpha_{22} & \alpha_{23} \\ \alpha_{31} & \alpha_{32} & \alpha_{33} \end{bmatrix} = \begin{bmatrix} 1 & 0 & 0 \\ 0 & 1 & 0 \\ 0 & 0 & 1\end{bmatrix}.
\end{equation*}
Now, we can determine a valid expression of support functions that ensures the support matrix is non-singular. For example, if the support functions were chosen as $s_i(x) = (1, x, x^2)$, the 2$^{\text{nd}}$ and 3$^{\text{rd}}$ columns of the support matrix would be linearly dependent; hence, this is an invalid set. The simplest set of monomials that satisfies the requirement is $s_i(x) = (1, x^2, x^3)$. Using the defined support functions, the $\alpha_{ij}$ coefficients can be derived as follows,
\begin{align*}
    \begin{bmatrix} 1 & 0 & 0 \\ 0 & 2 & 3 \\ 1 & 4 & 8\end{bmatrix} \begin{bmatrix} \alpha_{11} & \alpha_{12} & \alpha_{13} \\ \alpha_{21} & \alpha_{22} & \alpha_{23} \\ \alpha_{31} & \alpha_{32} & \alpha_{33} \end{bmatrix} &= \begin{bmatrix} 1 & 0 & 0 \\ 0 & 1 & 0 \\ 0 & 0 & 1\end{bmatrix} \\
    \begin{bmatrix} \alpha_{11} & \alpha_{12} & \alpha_{13} \\ \alpha_{21} & \alpha_{22} & \alpha_{23} \\ \alpha_{31} & \alpha_{32} & \alpha_{33} \end{bmatrix} &= \begin{bmatrix} 1 & 0 & 0 \\ 0 & 2 & 3 \\ 1 & 4 & 8\end{bmatrix}^{-1} = \begin{bmatrix} 1 & 0 & 0 \\ \frac{3}{4} & 2 & -\frac{3}{4} \\ -\frac{1}{2} & -1 & \frac{1}{2}\end{bmatrix}.
\end{align*}
Substituting the constants back into the switching functions and simplifying yields,
\begin{equation*}
    \phi_1(x) = \frac{-2x^3+3x^2+4}{4}, \quad \phi_2(x) = -x^3+2x^2, \quad \phi_3(x) =  \frac{2x^3-3x^2}{4}.
\end{equation*}
Finally, substituting the switching functions and projection functionals back into the constrained expression yields,
\begin{align}\label{eq:uniEx1Soln}
    y(x,g(x)) = g(x) &+ \frac{-2x^3+3x^2+4}{4}\Big(1-g(0)\Big) \\ &+ \Big(-x^3+2x^2\Big)\Big(2-g_x(1)\Big) + \frac{2x^3-3x^2}{4}\Big(3-g(2)\Big), \nonumber
\end{align}
\end{example}
It is simple to verify that regardless of how $g(x)$ is chosen, provided $g(x)$ is defined at the constraint points, Equation \eqref{eq:uniEx1Soln} always satisfies the given constraints.

\subsection{Integral constraints}
Moving forward, another constraint type of interest and one that can be easily incorporated using the TFC approach are integral constraints that include an integral over all or part of the domain. While the idea was first presented in Johnston and Mortari \cite{constraints}, this work relied on the original formulation. With the discovery of the switching-projection form, integral constraints become easier to embed. 
\begin{example}{Integral constraints}
Consider the function $y(x)$ subject to,
\begin{equation*}
    \int_0^3 y(x)\dd{x} = 0 \quad \text{and} \quad \int_1^2 y(x)\dd{x} = 2.
\end{equation*}
Following the same process as the prior example, first the projection functionals are determined. For this problem, the projection functions are merely the difference between the constraint value and free function evaluated over the integral. For this example,
\begin{equation*}
    \rho_1(x,g(x)) = -\int_0^3 g(\zeta)\dd{\zeta} \quad \text{and} \quad \rho_2(x,g(x)) = 2-\int_1^2 g(\zeta) \dd{\zeta}.
\end{equation*}
where $\zeta$ is a ``dummy'' variable for the integration of the function, $g(x)$. As before, the switching functions are defined such that they are equal to 1 when evaluated at their associated integral constraint, and equal to 0 when evaluated at all other constraints. For this example,
\begin{equation*}
    \int_0^3 \phi_1(x) \dd{x} = 1, \quad \int_1^2 \phi_1(x) \dd{x} = 0,
\end{equation*}
for the first switching function, and
\begin{equation*}
    \int_0^3 \phi_2(x) \dd{x} = 0, \quad \int_1^2 \phi_2(x) \dd{x} = 1,
\end{equation*}
for the second switching function. Similar to the previous examples, the switching functions are chosen to be a linear combination of support functions. For the first switching function, this form yields,
\begin{align*}
     \int_0^3 \phi_1(x) \dd{x} &= \int_0^3 \Big(s_1(x)\alpha_{11} + s_2(x)\alpha_{21} \Big) \dd{x}\\ &=  \alpha_{11} \int_0^3 s_1(x) \dd{x} + \alpha_{21} \int_0^3 s_2(x) \dd{x} = 1
\end{align*}
where we can see that the unknown $\alpha_{ij}$ terms still appear linearly. The final step is to define the specific support functions, and evaluate them at the constraint conditions to populate the support matrix. For this example, let's choose the support functions $s_1(x) = 1$ and $s_2(x) = x^2$. Expressing the support functions in this way yields,
\begin{align*}
    \begin{bmatrix} 3 & 9 \\ 1 & \frac{7}{3} \end{bmatrix} \begin{bmatrix} \alpha_{11} & \alpha_{12}\\  \alpha_{21} & \alpha_{22} \end{bmatrix} &= \begin{bmatrix} 1 & 0 \\ 0 & 1 \end{bmatrix}\\
    \begin{bmatrix} \alpha_{11} & \alpha_{12}\\  \alpha_{21} & \alpha_{22} \end{bmatrix} &= \begin{bmatrix} 3 & \frac{9}{2} \\ 1 & \frac{3}{2} \end{bmatrix}^{-1} = \begin{bmatrix} -\frac{7}{6} & \frac{9}{2} \\ \frac{1}{2} & -\frac{3}{2}\end{bmatrix},
\end{align*}
The solution of this system yields the following switching functions,
\begin{equation*}
    \phi_1(x) = \frac{3x^2-7}{6} \quad \text{and} \quad \phi_2(x) = \frac{-3x^2+9}{2}.
\end{equation*}
Finally, substituting the switching functions and projection functionals back into the constrained expression given in Equation \eqref{eq:s2a_uniCeAlt} produces,
\begin{equation*}
    y(x,g(x)) = g(x)-\frac{3x^2-7}{6}\int_0^3 g(\zeta)\dd{\zeta}+\frac{-3x^2+9}{2}\Big(2-\int_1^2 g(\zeta) \dd{\zeta}\Big).
\end{equation*}
\end{example}

Again, it is easy to check this constrained expression to ensure that the constraints are met regardless of the value of $g(x)$. The inclusion of integral constraints leads to another property of projection functionals.

\begin{property}{}\label{prop:proj2}
The projection functions for integral constraints are always equal to zero if the free function is selected such that it satisfies the integral constraint. 
\end{property}

For example, if $g(x)$ is selected such that $\int_1^2 g(\zeta)\dd{\zeta}=2$, then the second projection function in this example becomes $\rho_2(x,g(x)) = 2 - \int_1^2 g(\zeta) \dd{\zeta} = 0$.

\subsection{Linear constraints} 
Taking our discussion on the derivation of constrained expressions a step further, the culmination of all prior examples is the linear constraint case. It is noted that by this definition, relative constraints such as $y(0) = y(1)$ are just a specific case of linear constraints. As mentioned earlier, the idea of embedding a general set of linear constraints is not new and was first teased in the seminal TFC paper \cite{U-TFC}; however, the original form proved cumbersome when deriving constrained expressions of this type. In the following example, we highlight that these linear constraints can be embedded in the same way as the prior examples in the new generalized formulation.

\begin{example}{Linear constraints}
For this example, let us consider the linear constraints,
\begin{equation*}
    y(0) = y(1) \quad \text{and} \quad 3 = \int_0^1 y(x) \dd{x} + \pi y_x(0).
\end{equation*}
To generate a constrained expression, first the constraints are arranged such that the constants are collection on one side; for example,
\begin{equation*}
    0 = y(1)-y(0) \quad \text{and} \quad 3 = \int_0^1 y(x) \dd{x} + \pi y_x(0).
\end{equation*}
By organizing the constraints in this manner, the projection functionals, again, are immediate. However, the author notes one extra step must be taken for the general linear constraints. The projection functionals take the form,
\begin{equation*}
    \rho_1(x,g(x)) = g(0)-g(1) \quad \text{and} \quad \rho_2(x,g(x)) = 3 - \int_0^1 g(\zeta) \dd{\zeta} - \pi g_x(0),
\end{equation*}
where again $\zeta$ is the ``dummy'' variable for the integration of the free function. 

The switching functions are again such that they are equal to 1 when evaluated with their associated constraint and equal to 0 when evaluated at all other constraints. However, the word ``evaluation'' in the previous sentence requires clarification. Here, evaluation means to replace the function, $y(x)$ in this case, with the switching function and remove any terms not multiplied by the switching function. For this example, this leads to
\begin{equation*}
    \phi_1(1)-\phi_1(0) = 1, \quad  \int_0^1 \phi_1(x)\dd{x}+\pi\frac{\partial \phi_1}{\partial x}(0) = 0,
\end{equation*}
for the first switching function, and 
\begin{equation*}
    \phi_2(1)-\phi_2(0) = 0, \quad  \int_0^1 \phi_2(x)\dd{x}+\pi\frac{\partial \phi_2}{\partial x}(0) = 1,
\end{equation*}
for the second switching function. As in all prior examples, the switching functions are defined as a linear combination of support functions with unknown coefficients. Again, this can be written compactly in matrix form. For this example, let's choose the support functions $s_1(x) = 1$ and $s_2(x) = x$. Then the set of equations becomes,
\begin{align*}
    \begin{bmatrix} 0 & 1 \\ 1 & \frac{1}{2}+\pi \end{bmatrix} \begin{bmatrix} \alpha_{11} & \alpha_{12}\\  \alpha_{21} & \alpha_{22} \end{bmatrix} &= \begin{bmatrix} 1 & 0 \\ 0 & 1 \end{bmatrix}\\
    \begin{bmatrix} \alpha_{11} & \alpha_{12}\\  \alpha_{21} & \alpha_{22} \end{bmatrix} &= \begin{bmatrix} 0 & 1 \\ 1 & \frac{1}{2}+\pi \end{bmatrix}^{-1} = \begin{bmatrix} -\frac{1}{2}-\pi & 1 \\ 1 & 0 \end{bmatrix}.
\end{align*}
These coefficients are, as always, used to define the switching functions,
\begin{equation*}
    \phi_1(x) = -\frac{1}{2}-\pi + x, \quad \phi_2(x) = 1.
\end{equation*}
Lastly, substituting the switching functions and projection functionals back into the constrained expression form given in Equation \eqref{eq:s2a_uniCeAlt} yields,
\begin{equation*}
    y(x,g(x)) = g(x) + \Big(-\frac{1}{2}-\pi + x\Big)\Big(g(0)-g(1)\Big) + \Big(3 - \int_0^1 g(\zeta) \dd{\zeta} - \pi g_x(0)\Big).
\end{equation*}
\end{example}
By substituting this expression for $y(x)$ back into the constraints, one can verify that this constraint expression satisfies the constraints regardless of the choice of the free function $g(x)$. Therefore, we are lead to a similar property as those observed before.

\begin{property}{}\label{prop:proj3}
The projection functionals for linear constraints are always equal to zero if the free function is selected such that it satisfies the associated constraint.
\end{property}

It should be clear that Property \ref{prop:proj3} extends Property \ref{prop:proj1} and Property \ref{prop:proj2} to any linear constraints. For example, if $g(x)$ is selected such that $g(1) = g(0)$, then the first projection functional in this example becomes $\rho_1(x,g(x)) = g(1)-g(0) = 0$. Thus far, all examples have been for \emph{scalar} univariate equations. In the following examples we will look into vector univariate equations where another interesting constraint case arises: component constraints.

\subsection{Component constraints}\label{sec:s2_comp}
Component constraints involve constraints across dependent variables. Mortari and Furfaro \cite{comp-TFC} first looked at these constraints and their application to solving systems of ordinary differential equations. The following example is used to highlight that the new, generalized, constrained expression with the switching-projection form can easily embed any set of linear component constraints.

\begin{example}{Component linear constraints}
As with the prior constraint types, it is easiest to explore this constraint type through an example. Therefore, consider the vector function where the dependent variables $x$, $y$, and $z$ are all functions of the independent variable $t$ and are constrained by the following,
\begin{equation*}
    x (0) = 2 y (0) + \int_{-1}^{+1} z (t) \dd t \andd \dot{y} (0) = 2 x (1) - z(1).
\end{equation*}
When handling component constraints, one must decide which dependent variable’s constrained expression the component constraint will be embedded. Regardless of which dependent variable is chosen, a valid constrained expression will be produced. For this example, let us choose to embed all constraints into the $x$-component (note: this could have also been done for the $y$-component or $z$-component). Doing this leads to the following constrained expressions,
\begin{alignat}{2}\label{eq:ex4Ce}
    x(t,g^x(t),g^y(t),g^z(t)) &= g^x(t) &+ \phi_1(t)\rho_1(t,g^x(t),g^y(t),g^z(t)) \nonumber\\& &+ \phi_2(t)\rho_2(t,g^x(t),g^y(t),g^z(t)) \nonumber\\
    y(t,g^y(t)) &= g^y(t) \\
    z(t,g^z(t)) &= g^z(t) \nonumber,
\end{alignat}
Now, the definition of the projection functionals become,
\begin{align*}
    \rho_1(t,g^x(t),g^y(t),g^z(t)) &= g^x(0)-2 y(0,g^y(t))-\int_{-1}^{+1} z (\zeta,g^z(\zeta)) \dd \zeta \\ \rho_2(t,g^x(t),g^y(t),g^z(t))  &= \dot{y}(0,g^y(t))-2g^x(1)+z(1,g^z(t)),
\end{align*}
where we can see that $g^{x}(t)$, which represents the free function used for the $x(t)$ constrained expression, is the only free function that shows up in the expressions. Additionally, since the vector equation is a function of the independent variable $t$ the dot operator is used to signify the derivative such that $\dot{y} := \frac{\dd y}{\dd t}$.

Similar to previous examples, the number of switching functions is equal to the number of constraints. The switching functions are derived by evaluating the conditions based on the applied constraints,
\begin{equation*}
\begin{gathered}
    -\phi_1(0) = 1, \quad 2\phi_1(1) = 0\\
    -\phi_2(0) = 0, \quad 2\phi_2(1) = 1.
\end{gathered}
\end{equation*}
The negative sign will be explained in greater detail in Chapter \ref{chap:tfc_general} and is based on the structure of constraints and projection functionals.

As in previous examples, the switching functions are chosen to be a linear combination of support functions. Let the support functions for this example be $s_1(t) = 1$ and $s_2(t) = t$. Then,
\begin{align*}
    \begin{bmatrix} -1 & 0 \\ 2 & 2 \end{bmatrix} \begin{bmatrix} \alpha_{11} & \alpha_{12}\\  \alpha_{21} & \alpha_{22} \end{bmatrix} &= \begin{bmatrix} 1 & 0 \\ 0 & 1 \end{bmatrix}\\
    \begin{bmatrix} \alpha_{11} & \alpha_{12}\\  \alpha_{21} & \alpha_{22} \end{bmatrix} &= \begin{bmatrix}-1 & 0 \\ 2 & 2 \end{bmatrix}^{-1} = \begin{bmatrix} -1 & 0 \\ 1 & \frac{1}{2} \end{bmatrix}.
\end{align*}
where $\phi_1 = \alpha_11+\alpha_21t$ and $\phi_2 = \alpha_12 + \alpha_22t$. Substituting these values into the constrained expressions shown in Equation \eqref{eq:ex4Ce} yields,
\begin{alignat*}{2}
    x(t,g^x(t),g^y(t),g^z(t)) &= g^x(t) &+ (t-1)\Big(g^x(0)-2 y(0,g^y(t)) -\int_{-1}^{+1} z (\zeta,g^z(\zeta)) \dd \zeta \Big) \\& &+\frac{t}{2}\Big( \dot{y}(0,g^y(t))-2g^x(1)+z(1,g^z(t)) \Big)\\
    y(t,g^y(t)) &= g^y(t)\\
    z(t,g^z(t)) &= g^z(t).
\end{alignat*}
\end{example}
As with all prior examples, notice that regardless of how the free functions are chosen, these constrained expressions will always satisfy the constraints. In fact, Property \ref{prop:proj3} can be extended to component constraints.

\begin{property}{}\label{prop:proj4}
The projection functions for component constraints are always equal to zero if the free functions are selected such that they satisfy the component constraints. 
\end{property}

For example, if $g^x(t)$, $g^y(t)$, and $g^z(t)$ are selected such that $\dot{g}^y (0) = 2 g^x (1) - g^z(1)$, then the second projection function in this example becomes $\rho_2(t,g(t))  = \dot{g}^y(0)-2g^x(1)+g^z(1) = 0$.

\subsection{Mixed constraints}
The methods for building constrained expressions shown in the previous four examples can be combined. However, special care must be taken when combining component constraints with the other types of constraints discussed earlier. The nuances of doing so are highlighted in this example. 

\begin{example}{Mixed constraints}
Consider the vector function where the dependent variables $x$ and $y$ are both functions of the independent variable $t$ and are constrained by the following equations,
\begin{equation*}
    x(0) = 0, \quad y(0) = 0, \quad y(1) = y(2), \andd 4 = 2 y(1) - \int_0^3 x(t) \dd{t}. 
\end{equation*}
Based on the previous examples, the four projection functions are defined,
\begin{align*}
    \rho_1(t,g^x(t))  &= -g^x(0), \quad  \rho_3(t,g^y(t))  = g^y(1)-g^y(2), \\ \rho_2(t,g^y(t)) &= -g^y(0),  \quad \rho_4(t,g^x(t),g^y(t))  = 4-2y(1,g^y(t))+ \int_0^3 g^x(\zeta) \dd{\zeta}.
\end{align*}
As there are four constraints, there must also be four switching functions. Based on the constraints, the first must be associated with the $x$ independent variable, and the second and third must be associated with the $y$ independent variable. However, just as in the previous example, with the component constraint, there is freedom to choose where the constraint goes. How we have written $\rho_4(t,g^x(t),g^y(t))$, the constraint will be applied to the $x$-component, but it could have easily been applied to the $y$-component. The resulting constrained expressions are defined as,
\begin{align*}
    x(t,g^x(t),g^y(t)) &= g^x(t) + \phi^x_1(t)\rho_1(t,g^x(t))+\phi^x_2(t)\rho_4(t,g^x(t),g^y(t))\\
    y(t,g^y(t)) &= g^y(t)+\phi^y_1(t)\rho_2(t,g^y(t))+\phi^y_2(t)\rho_3(t,g^y(t)),
\end{align*}
where the switching function equations are,
\begin{equation*}
\begin{gathered}
    \phi^x_1(0) = 1, \quad -\int_0^3\phi^x_1(t)\dd{t} = 0\\
    \phi^x_2(0) = 0, \quad -\int_0^3\phi^x_2(t)\dd{t} = 1\\
    \phi^y_1(0) = 1, \quad \phi^y_1(2)-\phi^y_1(1) = 0\\
    \phi^y_2(0) = 0, \quad \phi^y_2(2)-\phi^y_2(1) = 1.
\end{gathered}
\end{equation*}
Each switching function is again chosen to be a linear combination of support functions, where in this case the support functions are chosen as $s^x_1(t) = 1$ and $s^x_2(t) = t$ for $x(t)$ switching functions and $s^y_1(t) = 1$ and $s^y_2(t) = t$ for $y(t)$ switching functions. Thus, the switching function can be concisely written as,
\begin{alignat*}{2}
    \begin{bmatrix} 1 & 0 \\ -3 & -\frac{9}{2}\end{bmatrix}\begin{bmatrix} \alpha_1 & \alpha_3 \\ \alpha_2 & \alpha_4 \end{bmatrix} &= \begin{bmatrix}1 & 0 \\ 0 & 1\end{bmatrix} &\quad \begin{bmatrix} 1 & 0 \\ 0 & 1\end{bmatrix}\begin{bmatrix} \alpha_5 & \alpha_7 \\ \alpha_6 & \alpha_8 \end{bmatrix} &= \begin{bmatrix}1 & 0 \\ 0 & 1 \end{bmatrix}\\
    \begin{bmatrix} \alpha_1 & \alpha_3 \\ \alpha_2 & \alpha_4 \end{bmatrix} &= \begin{bmatrix}1 & 0 \\ -\frac{2}{3} & -\frac{2}{9}\end{bmatrix} & \begin{bmatrix} \alpha_5 & \alpha_7 \\ \alpha_6 & \alpha_8 \end{bmatrix} &= \begin{bmatrix} 1 & 0\\ 0 & 1 \end{bmatrix}
\end{alignat*}
where $\phi^x_1(t) = \alpha_1+\alpha_2 t$, $\phi^x_2(t) = \alpha_3+\alpha_4 t$, $\phi^y_1(t) = \alpha_5+\alpha_6 t$, and $\phi^y_2(t) = \alpha_7+\alpha_8 t$. Substituting these values into the constrained expressions yields,
\begin{align*}
    x(t,g^x(t),g^y(t)) &= g^x(t) - \Big(1-\frac{2}{3}t\Big)g^x(0)-\frac{2}{9}t\Big(4-2y(1,g^y(t))+ \int_0^3 g^x(\zeta) \dd{\zeta}\Big)\\
    y(t,g^y(t)) &= g^y(t)-g^y(0) + t \Big(g^y(1)-g^y(2)\Big).
\end{align*}
\end{example}
As in all previous examples, notice that regardless of how the free functions are chosen, the constraints will be satisfied exactly. 

\subsection{Infinite constraints}

The derivation of constrained expression with infinite constraints was first solved by Johnston and Mortari \cite{constraints} and requires greater attention to the selection of support functions. To understand this, first, consider a single infinite constraint on the value of the function as it approaches infinity,
\begin{equation*}
    \ds\lim_{x\to\infty} y (x) = y_{\infty}.
\end{equation*}
When dealing with this single constraint, it should be straightforward to determine a simple constrained expression satisfying this constraint as,
\begin{equation*}
    y (x) = g (x) + \phi(x) \left(y_{\infty} - g(\infty)\right).
\end{equation*}
Here, the switching function can be simply defined as a constant value, $\phi(x) := 1$. As with all other types of constraints, the free function must be defined at the constraint. Therefore, $g(x)$ must be finite as $x\to\infty$. Additionally, as shown in the following example, the support functions must all be defined and finite at infinity.

\begin{example}{Infinite constraints}
Consider a mixture of finite and infinite constraints as defined in the Falkner-Skan boundary layer equation \cite{falk_skan},
\begin{equation*}
    y(0) = 0, \quad y_x(0) = 0, \quad \text{and} \quad y_x(\infty) = 1.
\end{equation*}
It follows that the projection functionals are,
\begin{equation*}
    \rho_1(x,g(x)) = -g(0), \quad \rho_2(x,g(x)) = -g_x(0), \quad \text{and} \quad \rho_3(x,g(x)) = 1 - g_x(\infty).
\end{equation*}
Let the support functions be,
\begin{equation*}
s_1(x) = 1, \quad s_2(x) = x, \quad \text{and} \quad s_3(x) = \frac{x-1}{x+1} 
\end{equation*}
Here, the selection of $s_3(x)$ is not arbitrary and is selected such that the last row of the support matrix is not zero and is therefore invertible. This leads to the system of equations.
\begin{equation*}
    \begin{bmatrix} s_1(0) & s_2(0) & s_3(0) \\ s_{1_x}(0) & s_{2_x}(0) & s_{3_x}(0) \\ s_{1_x}(\infty) & s_{2_x}(\infty) & s_{3_x}(\infty)\end{bmatrix} \begin{bmatrix} \alpha_{11} & \alpha_{12} & \alpha_{13} \\ \alpha_{21} & \alpha_{22} & \alpha_{23} \\ \alpha_{31} & \alpha_{32} & \alpha_{33} \end{bmatrix} = \begin{bmatrix} 1 & 0 & 0 \\ 0 & 1 & 0 \\ 0 & 0 & 1\end{bmatrix}
\end{equation*}
which through matrix inversion leads to the solution of the $\alpha_{ij}$ coefficients 
\begin{align*}
    \begin{bmatrix} 1 & 0 & -1 \\ 0 & 1 & 2 \\ 0 & 1 & 0\end{bmatrix} \begin{bmatrix} \alpha_{11} & \alpha_{12} & \alpha_{13} \\ \alpha_{21} & \alpha_{22} & \alpha_{23} \\ \alpha_{31} & \alpha_{32} & \alpha_{33} \end{bmatrix} &= \begin{bmatrix} 1 & 0 & 0 \\ 0 & 1 & 0 \\ 0 & 0 & 1\end{bmatrix} \\
    \begin{bmatrix} \alpha_{11} & \alpha_{12} & \alpha_{13} \\ \alpha_{21} & \alpha_{22} & \alpha_{23} \\ \alpha_{31} & \alpha_{32} & \alpha_{33} \end{bmatrix} &= \begin{bmatrix} 1 & 0 & -1 \\ 0 & 1 & 2 \\ 0 & 1 & 0\end{bmatrix}^{-1} = \begin{bmatrix} 1 & \frac{1}{2} & -\frac{1}{2} \\ 0 & 0 & 1 \\ 0 & \frac{1}{2} & -\frac{1}{2}\end{bmatrix}.
\end{align*}
From this solution, the switching functions become,
\begin{equation*}
    \phi_1(x) = 1, \quad 
    \phi_2(x) = \frac{1}{2} + \frac{x-1}{2(x+1)}, \quad 
    \text{and} \quad 
    \phi_3(x) = -\frac{1}{2} + x - \frac{x-1}{2(x+1)}
\end{equation*}
and the full constrained expression is ,
\begin{equation*}
    y(x,g(x)) = g(x) -g(0) + \Big(\frac{1}{2} + \frac{x-1}{2(x+1)}\Big)\Big(- g_x(0)\Big) 
    + \Big(-\frac{1}{2} + x -  \frac{x-1}{2(x+1)}\Big)\Big(1 - g_x(\infty)\Big)
\end{equation*}
\end{example}

With this example, we conclude our exploration of the implications and capabilities of the reformulation of the TFC approach spurred by the switching-projection form. These simply applied the techniques and loosely defined such terms as constrained expression, switching function, projection functional, etc., without much mathematical rigor. Chapter \ref{chap:tfc_general} looks to explicitly define all terms used; however, before doing this, it is important to highlight two other constrained types (inequality constraints and weighted-constraints), which are simply an extension of the constrained expression produced above.

\section{Extension to inequality constraints}
This section is referred to as an extension to inequality constraints since the following theory relies on the earlier sections. Inequality type constraints were first explored in Johnston, Leake, Efendiev, and Mortari \cite{TFC-Selected} and Johnston, Leake, and Mortari \cite{inequality}; however, this dissertation provides invaluable updates from these two works. 

To begin, let us consider a simple case with only one, continuous upper-bound inequality constraint defined on the domain $x\in[a,b]$. Let that constraint be given by the function $f_u(x)$ such that a function $y(x)$ satisfies this constraint if,
\begin{equation*}
    y(x) \leq f_u (x), \qquad \forall\ x\in[a,b].
\end{equation*}
For any given function $g(x)$, we can subtract off the sections of $g(x)$ that are larger than the inequality constraint $f_u(x)$ by using the Heaviside step function,
\begin{equation*} 
   \Hs(z_1,z_2) = \begin{cases} 0 &\quad \text{if } z_1 < 0  \\ z_2 &\quad \text{if } z_1 = 0 \\ 1 &\quad \text{if } z_1 > 0 \end{cases}
\end{equation*}
where the derivative of the Heaviside step function is exactly zero for all $z_1$. Furthermore, the Heaviside step function reduces to a simple step function if $z_2 = 0$, and in those cases will be defined as $\Hs_0(z_1) := \Hs(z_1,0)$. The Heaviside step function can be thought of as the functional form of a gate or switch, and can be used to subtract off the difference between $f_u(x)$ and $g(x)$ when $g(x)>f_u(x)$, but does not affect $g(x)$ when $g(x)\leq f_u(x)$. Mathematically, this can be written as,
\begin{equation}\label{eq:fUpper}
    y(x,g(x)) = g(x) + [f_u(x)-g(x)]\Hs_0\Big(g(x)-f_u(x)\Big),
\end{equation}
where $y(x,g(x))$ now represents the family of all possible functions that satisfy the inequality constraint. Another term can be added to Equation \eqref{eq:fUpper} to accommodate a lower bound inequality constraint as well, $f_\ell(x)$. This is shown in Equation (\ref{eq:fGen}).
\begin{equation}\label{eq:fGen}
    y(x,g(x)) = g(x) + [f_u(x)-g(x)]\Hs_0\Big(g(x)-f_u(x)\Big) + [f_\ell(x)-g(x)]\Hs_0\Big(f_\ell(x)-g(x)\Big)
\end{equation}

\subsection{Combining inequality and equality constraints}
The technique to embed equality and inequality constraints builds on the formulation given in the earlier sections on the TFC approach to equality constraints. For a problem subject to equality and inequality constraints, let the TFC constrained expression for just the equality constraints be given by $\hat{y} (x,g(x))$. As per the univariate TFC, $\hat{y} (x,g(x))$ will represent the family of all possible functions that satisfy the equality constraints. Then, we exchange $g(x)$ in Equation (\ref{eq:fGen}) with $\hat{y} (x,g(x))$, as shown in Equation (\ref{eq:fFull}), to project $\hat{y} (x,g(x))$ onto the set of functions that satisfy the inequality constraints. It must be noted that this approach is limited to point equality constraints --- derivative, integral, or component constraints cannot be combined with inequality constraints.
\begin{align}\label{eq:fFull}
    y (x,g(x)) = \hat{y} (x,g(x)) &+ \left[f_u (x) - \hat{y} (x,g(x))\right] \Hs_0 \Big(\hat{y} (x,g(x)) - f_u (x)\Big) \nonumber\\ &+ \left[f_\ell (x) - \hat{y} (x,g(x))\right] \Hs_0 \Big(f_\ell (x) - \hat{y} (x,g(x))\Big)
\end{align}
The resultant functional, $y(x,g(x))$, is now the TFC constrained expression representing all possible functions that satisfy both the equality constraints and inequality constraints of the problem.

\begin{example}{Numerical example of inequality constraints}
Now, to analyze the expressions provided by Equations (\ref{eq:fGen}) and (\ref{eq:fFull}), a numerical test was constructed where the free function $g(x)$ and the inequality constraints were randomly generated through a linear expansion of $m$ Chebyshev polynomials such that,
\begin{equation}\label{eq:basis}
    g (x) = \ds\sum_{i = 0}^{m-1} a_i \, T_i (x),
\end{equation}
where $a_i$ are random coefficients $a_i \thicksim N(0,1)$ and $T_i(x)$ are the individual terms of the Chebyshev polynomials. Figure \ref{fig:ineqConst} shows Equation (\ref{eq:fGen}) subject to random inequality bounds and random values of $g(x)$. Furthermore, Figure \ref{fig:ineqAndEqConst} shows the application Equation (\ref{eq:fFull}) to both inequality and randomly generated equality point constraints. In both plots, the inequality constraints are shown as dotted black lines, the functions are shown as colored lines, and the three equality constraints are shown as black points.
\begin{figure}[H]
\centering
\begin{minipage}[t]{0.49\linewidth}
    \centering\includegraphics[width=\linewidth]{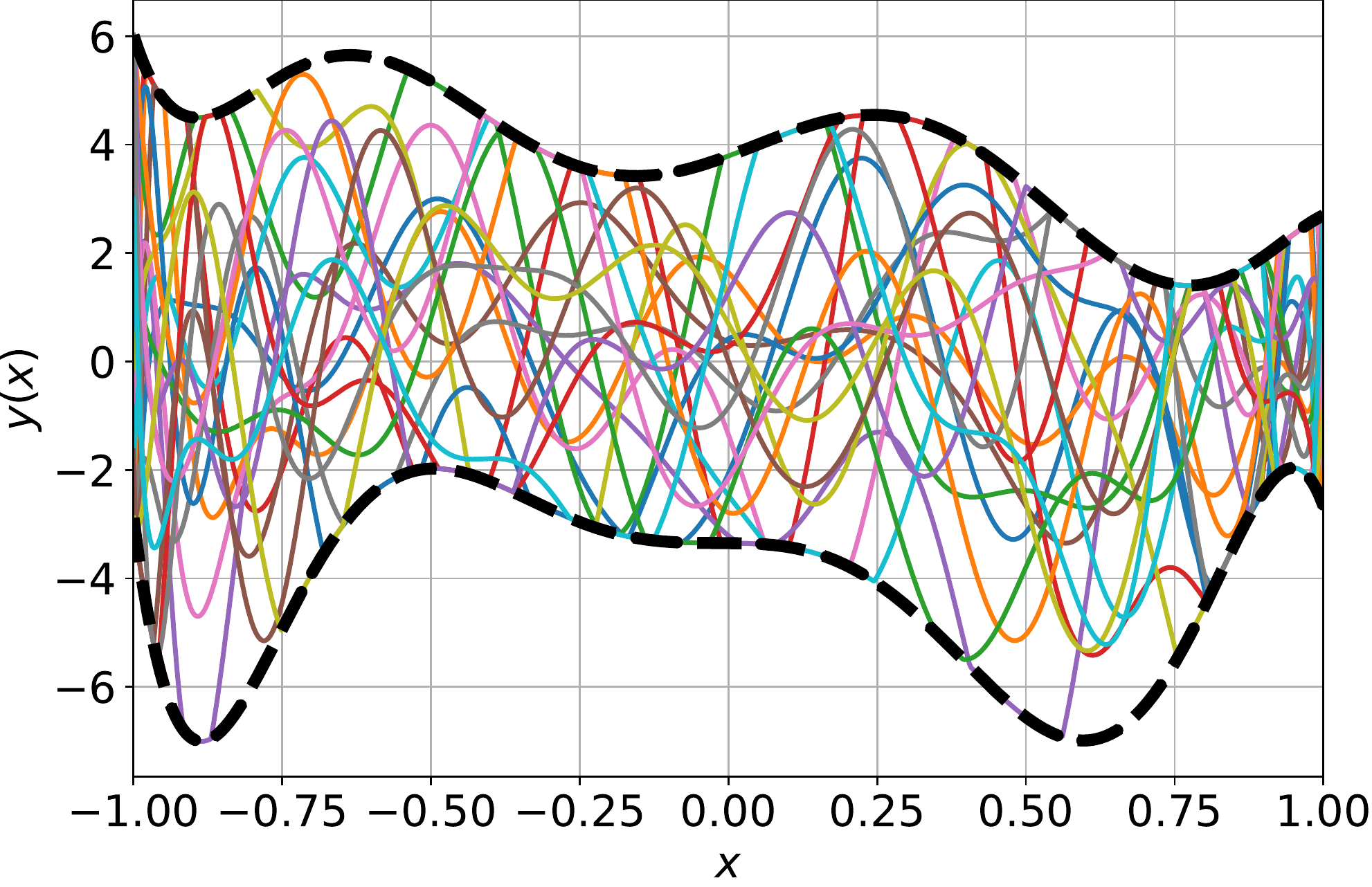}
    \caption{TFC constrained expression for inequality constraints only.}
    \label{fig:ineqConst}
\end{minipage}
    \begin{minipage}[t]{0.49\linewidth}
    \centering\includegraphics[width=\linewidth]{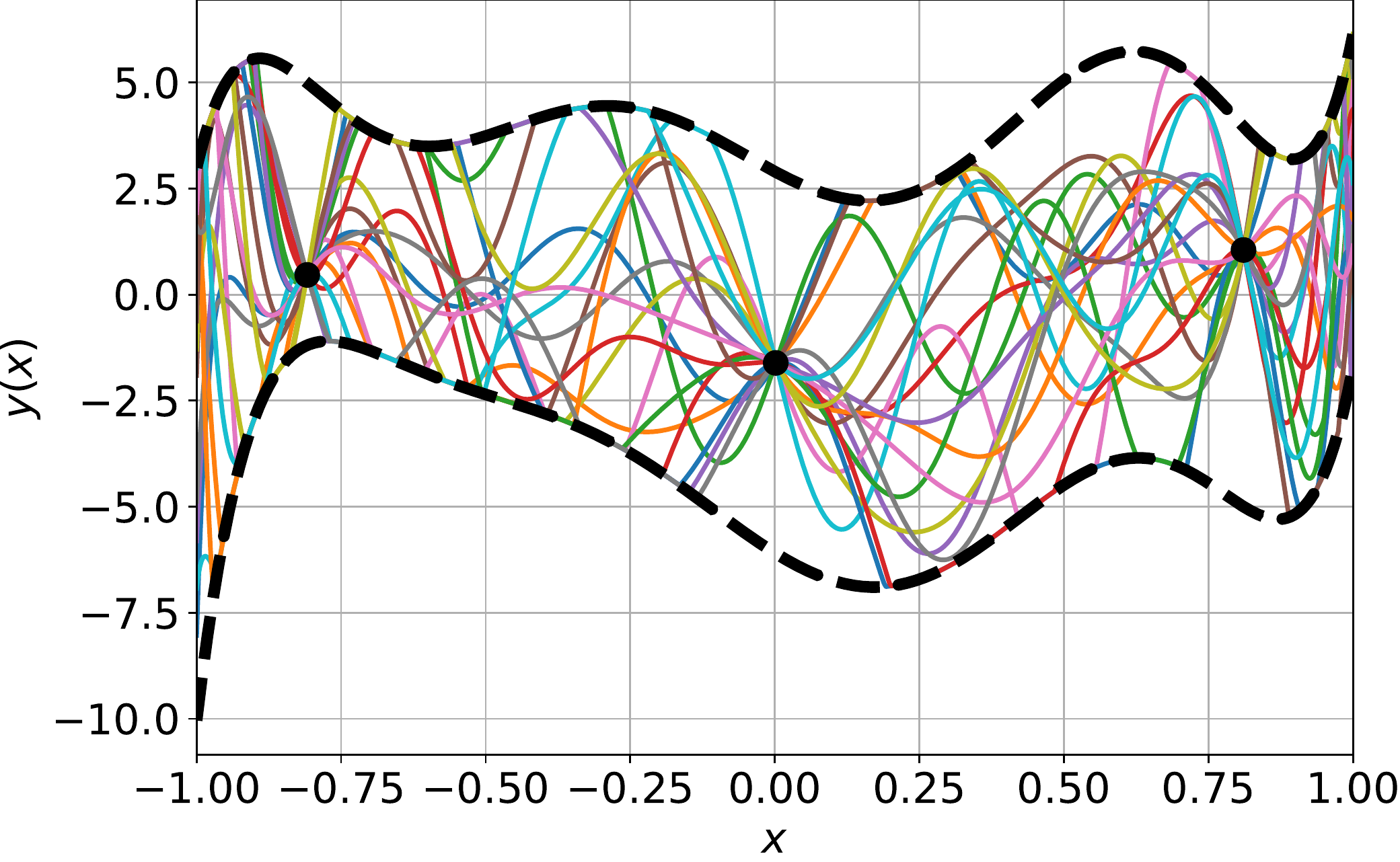}
    \caption{TFC constrained expression for equality and inequality constraints.}
    \label{fig:ineqAndEqConst}
\end{minipage}
\end{figure}
\end{example}

\subsection{Keep-out zones}

Using the univariate formulation of the TFC method subject to inequality constraints, a technique can be constructed for keep-out zones by augmenting Equation (\ref{eq:fGen}) or (\ref{eq:fFull}). This approach requires the constrained expression to be split into multiple constrained expressions for each possible path. 

\begin{example}{Keep-out zones}
As a simple example, let us consider solving all possible trajectories subject to upper and lower inequality constraints such that $f_u(x) = 1$, $f_\ell(x) = -1$ and avoiding an interior box defined by the coordinates, $A(-0.25,-0.25)$, $B(-0.25,0.25)$, $C(0.25,0.25)$, and $D(0.25,-0.25)$ as detailed in Figure \ref{fig:box}.
\begin{figure}[H]
    \centering\includegraphics[width=0.5\linewidth]{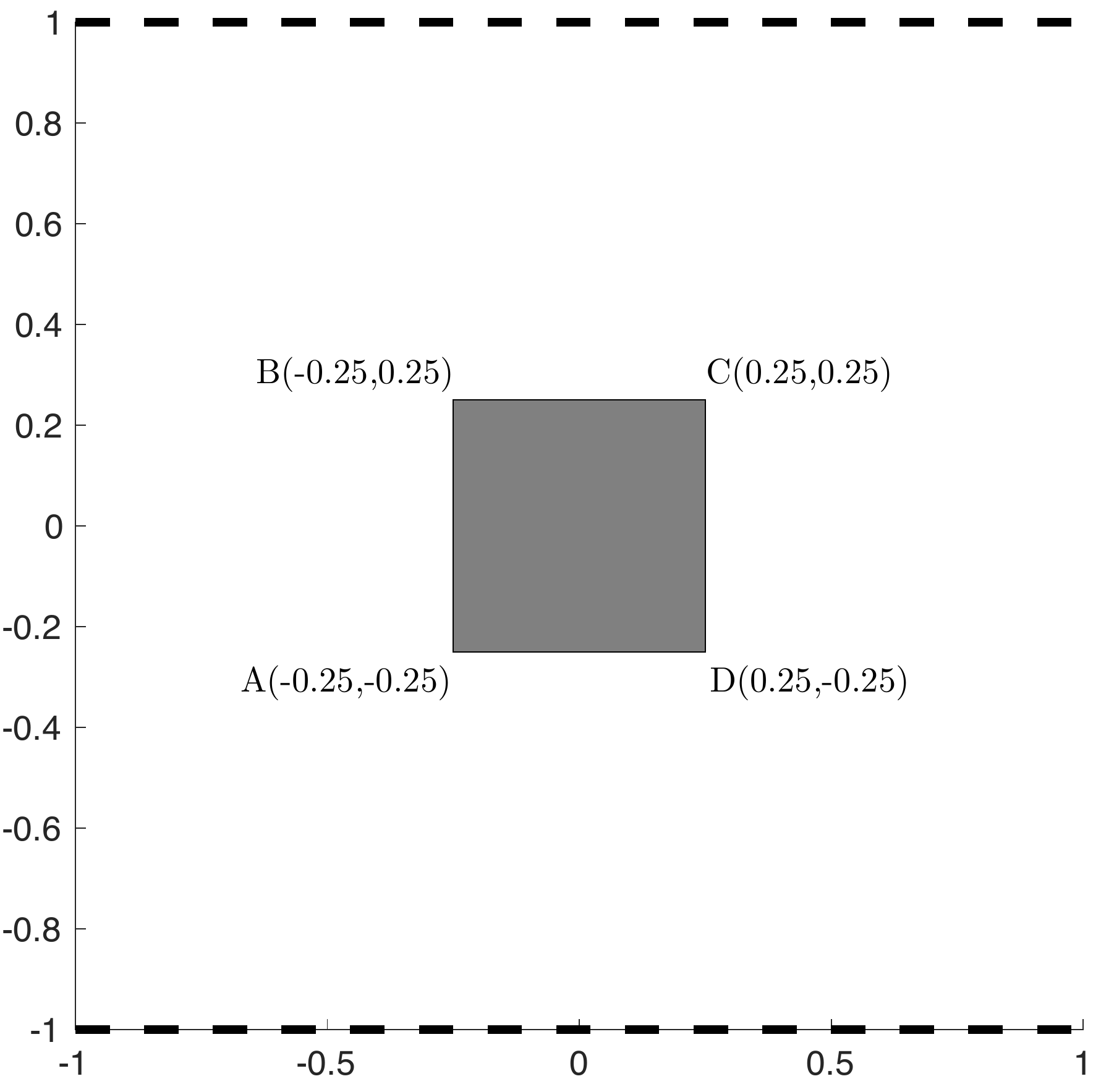}
    \caption{Keep-out box example.}
    \label{fig:box}
\end{figure}
\begin{figure}[H]
\centering
\begin{minipage}[t]{0.49\linewidth}
    \centering\includegraphics[width=0.95\linewidth]{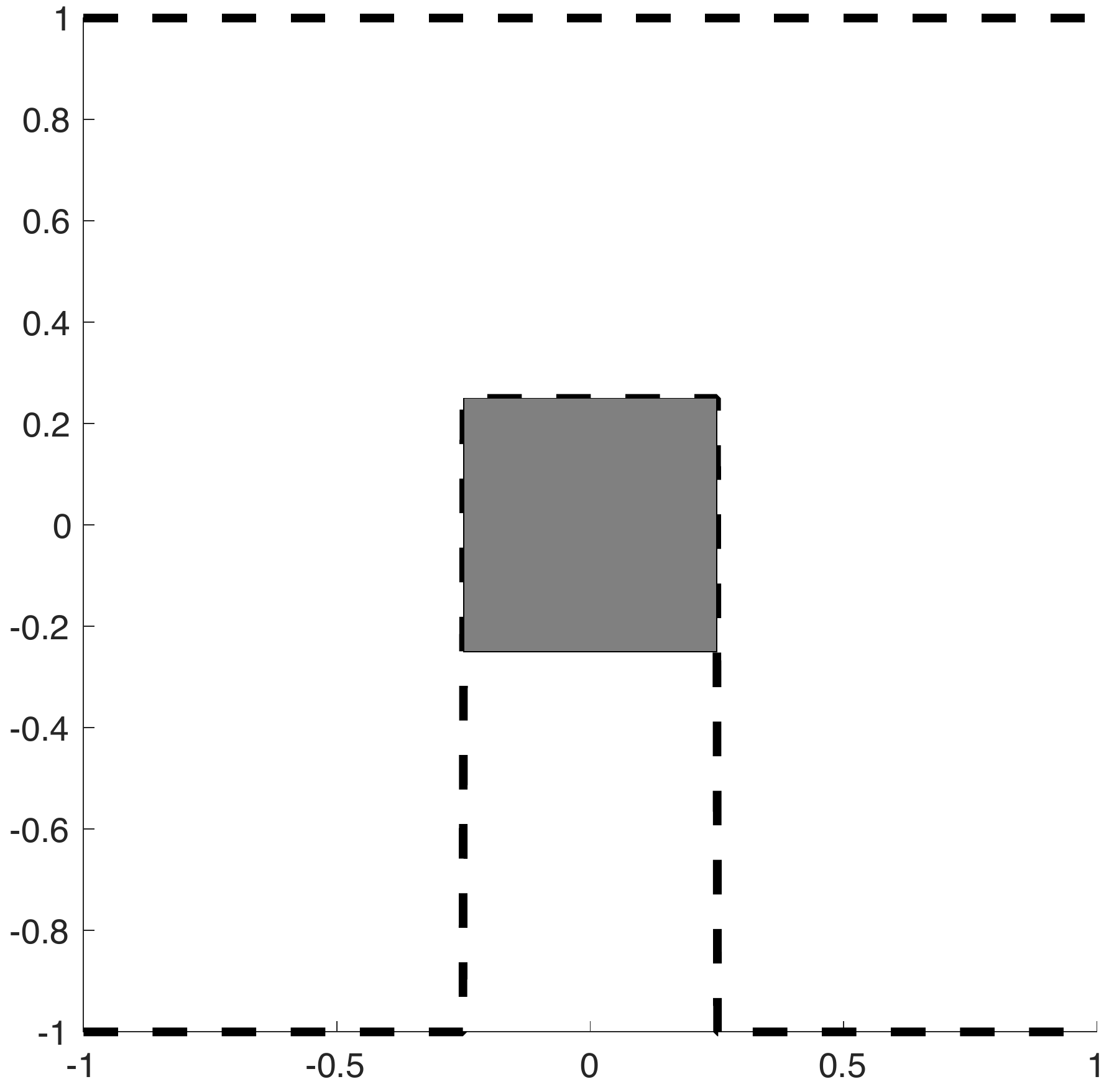}
    \caption{Upper path of keep-out box example.}
    \label{fig:upper_path}
\end{minipage}
\begin{minipage}[t]{0.49\linewidth}
    \centering\includegraphics[width=0.95\linewidth]{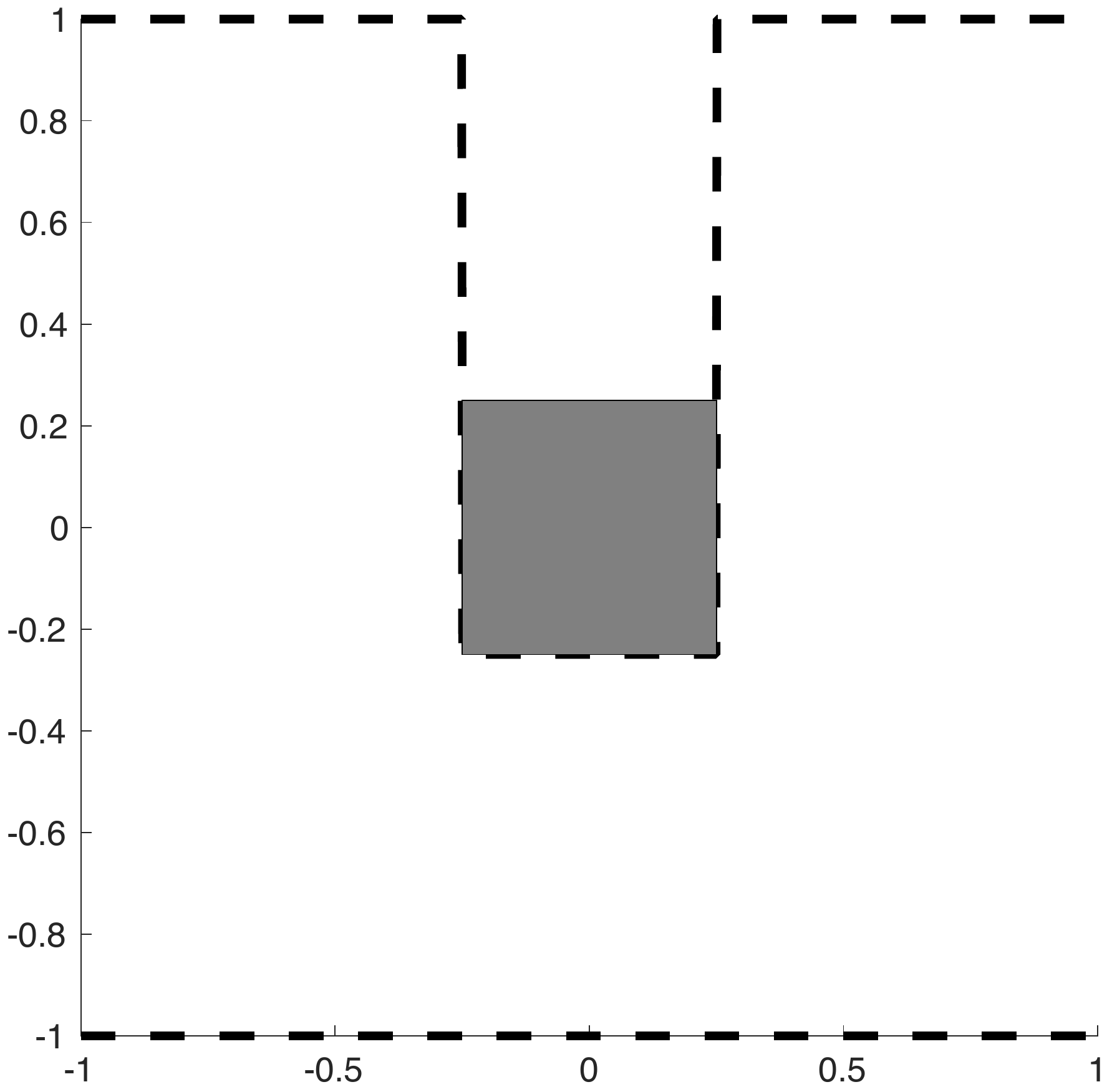}
    \caption{Lower path of keep-out box example.}
    \label{fig:lower_path}
\end{minipage}
\end{figure}
In order to accommodate these constraints, we can split the problem into two individual problems that follows the formulation of the earlier sections. Therefore, we consider the upper path defined in Figure \ref{fig:upper_path} and the lower path defined in Figure \ref{fig:lower_path}. First, to satisfy the upper path, the lower boundary needs to be augmented by the function defining the box which takes the form,
\begin{equation*}
    f_\ell(x) = \begin{cases} -1, \quad &\text{if}\quad x < -0.25 \\0.25, \quad &\text{if}\quad -0.25 \leq x \leq 0.25 \\ -1, \quad &\text{if}\quad x > 0.25\end{cases}
\end{equation*}
It then follows that for the lower path, the upper boundary is augmented by the function of the lower portion of the box such that,
\begin{equation*}
    f_u(x) = \begin{cases} 1, \quad &\text{if}\quad x < -0.25 \\-0.25, \quad &\text{if}\quad -0.25 \leq x \leq 0.25 \\ 1, \quad &\text{if}\quad x > 0.25\end{cases}.
\end{equation*}
Searching over both constrained expressions produces all the possible trajectories around the object. This method is analyzed by expressing $g(x)$ by Equation (\ref{eq:basis}) and again defining $\B{\xi} \thicksim \mathcal{N}(0,I_{m \times m})$. In addition to the single box example (Figure \ref{fig:box}), Figures \ref{fig:one_box} - \ref{fig:saturn} experiment with differing geometries and configurations of constraints.
\begin{figure}[H]
\centering
\begin{minipage}[t]{0.49\linewidth}
    \centering\includegraphics[width=0.95\linewidth]{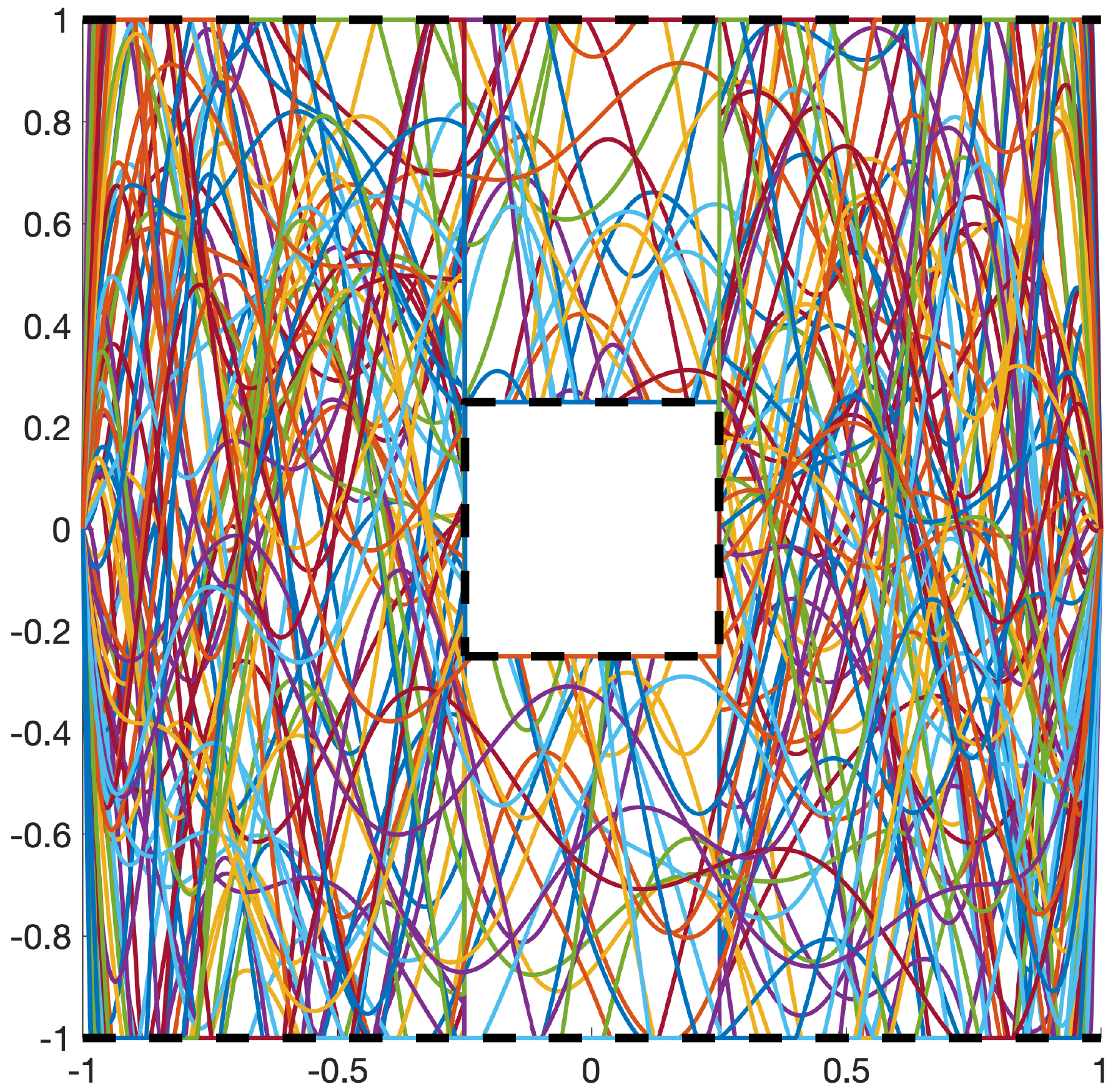}
    \caption{Single box.}
    \label{fig:one_box}
\end{minipage}
\begin{minipage}[t]{0.49\linewidth}
    \centering\includegraphics[width=0.94\linewidth]{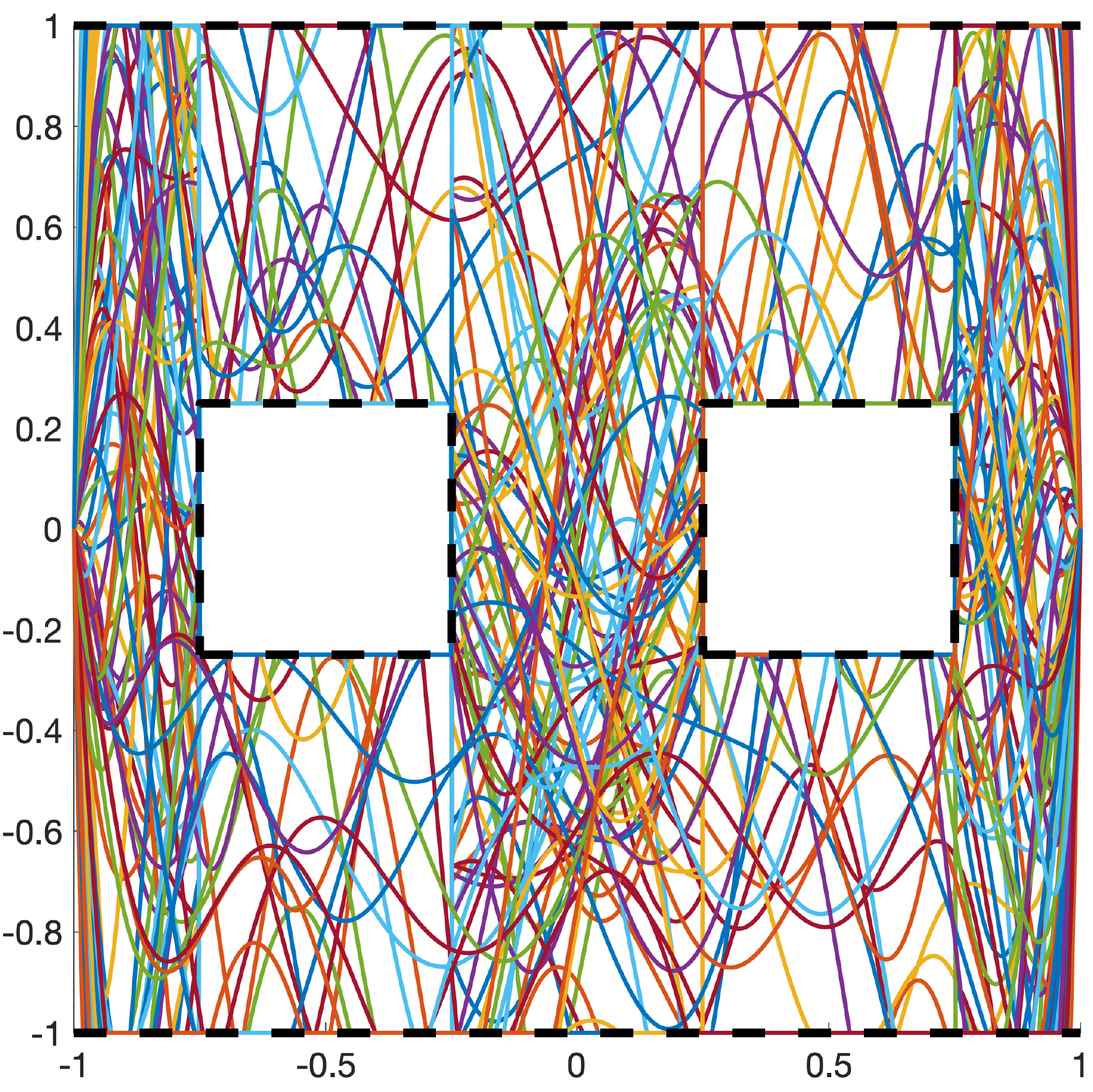}
    \caption{Two boxes horizontally arranged.}
\end{minipage}
\end{figure}
\begin{figure}[H]
\centering
\begin{minipage}[t]{0.49\linewidth}
    \centering\includegraphics[width=0.94\linewidth]{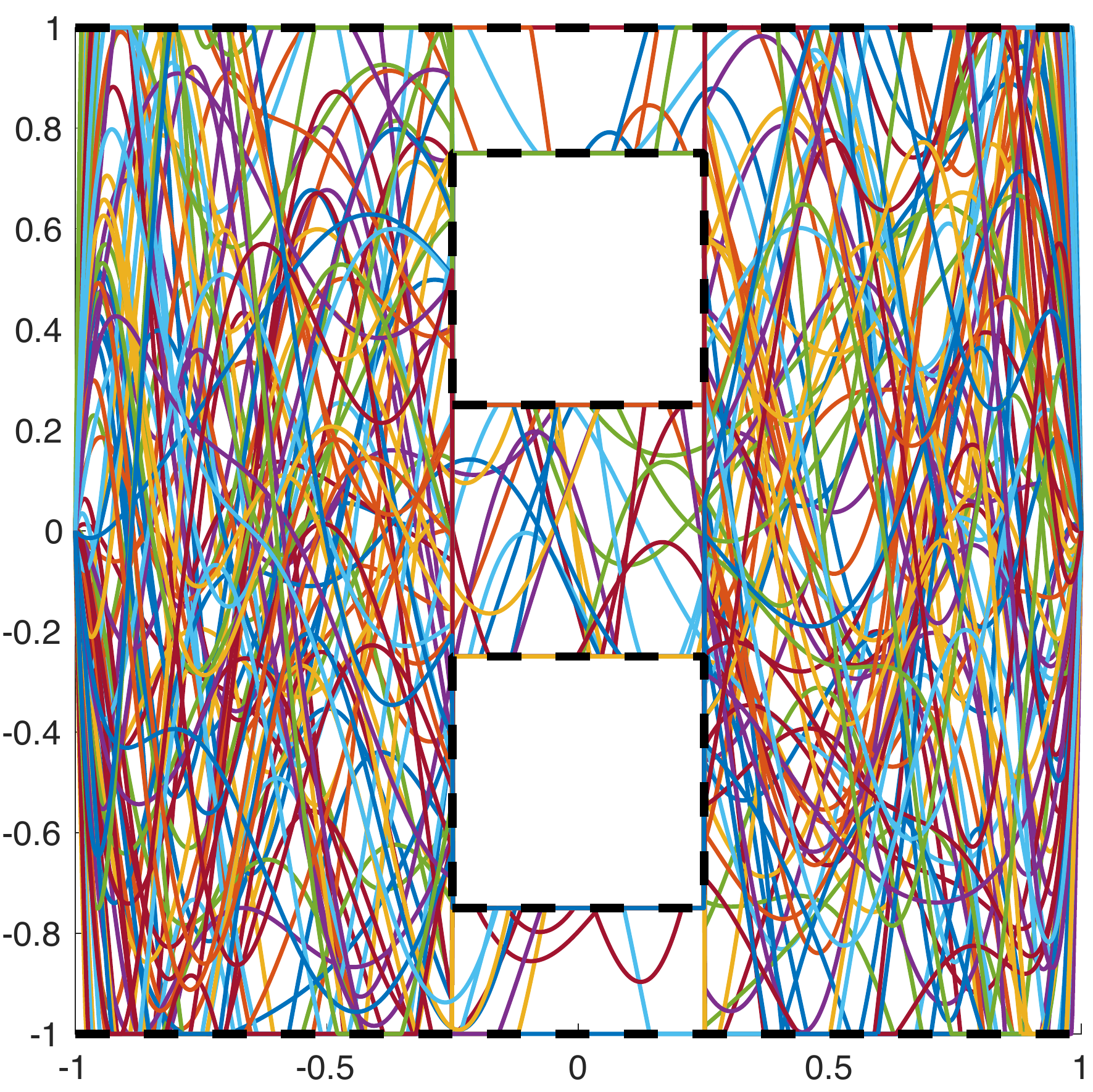}
    \caption{Two boxes vertically arranged.}
\end{minipage}
\begin{minipage}[t]{0.49\linewidth}
    \centering\includegraphics[width=0.95\linewidth]{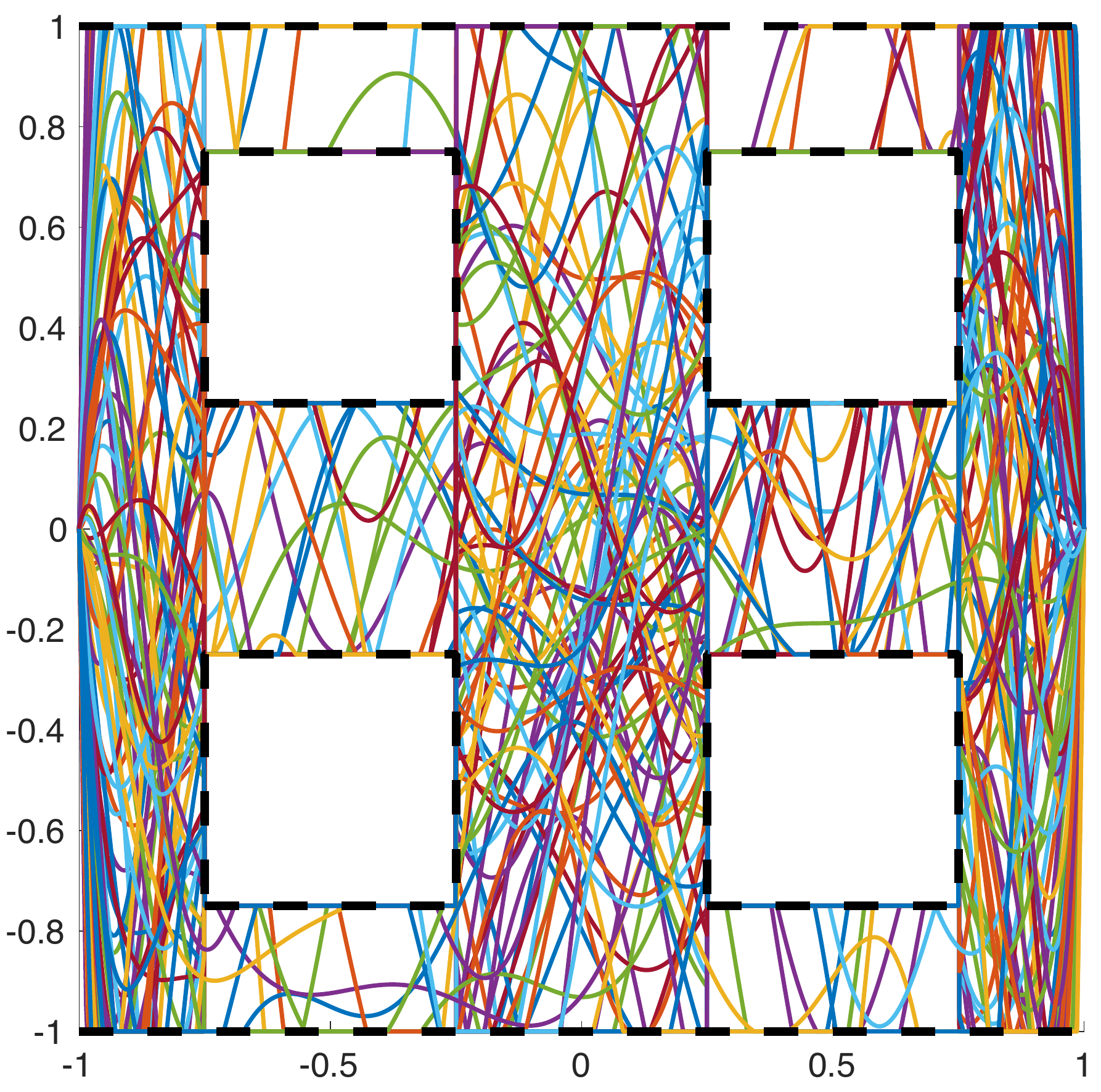}
    \caption{Four boxes.}
    \label{fig:four_box}
\end{minipage}
\end{figure}
Figures \ref{fig:one_box}-\ref{fig:four_box} provide multiple different keep-out box structures, including two horizontally arranged, two vertically arranged, and four equally spaced boxes.
\begin{figure}[H]
\centering
\begin{minipage}[t]{0.49\linewidth}
    \centering\includegraphics[width=0.95\linewidth]{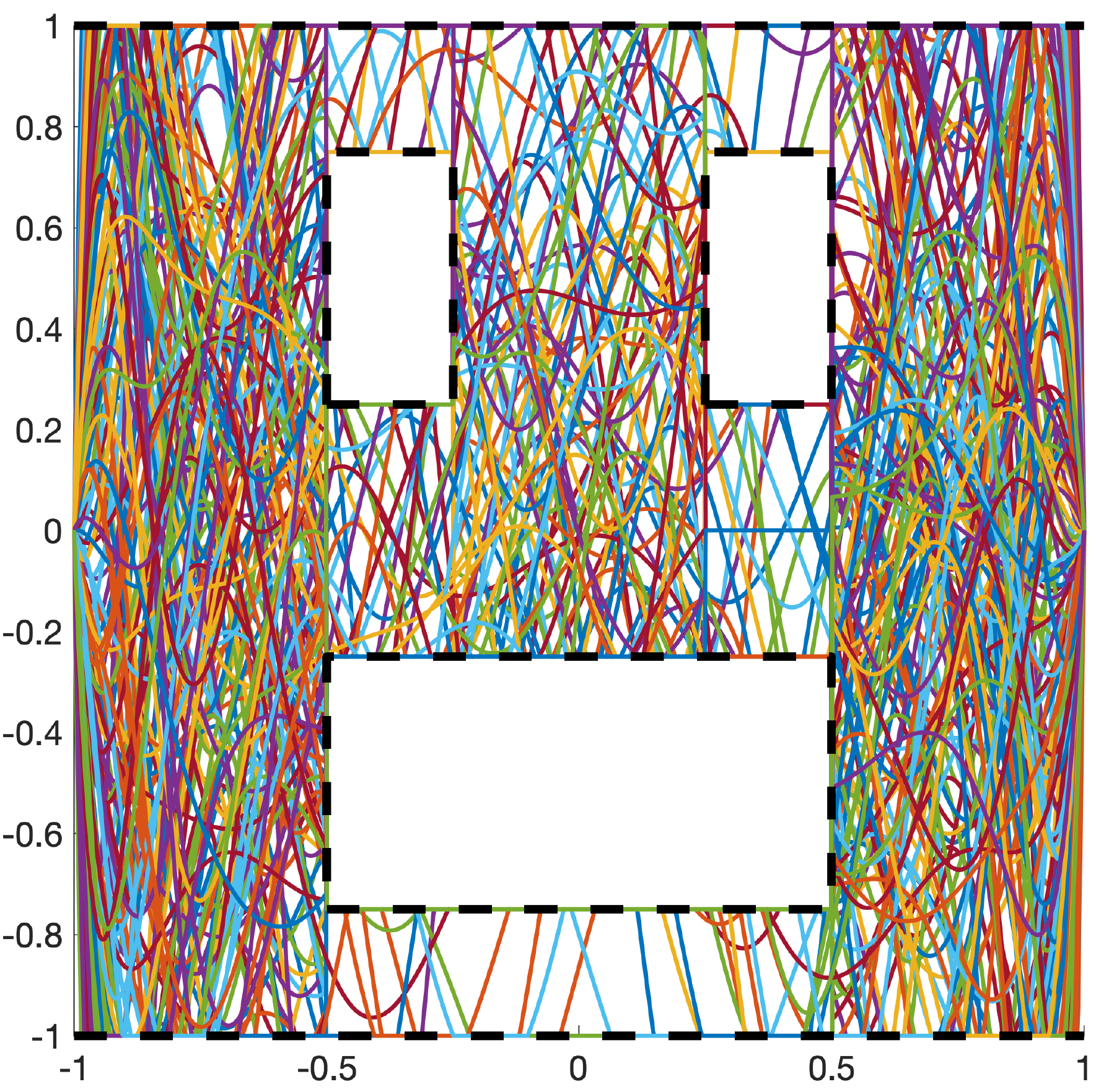}
    \caption{Different rectangles.}
    \label{fig:rectangle}
\end{minipage}
\begin{minipage}[t]{0.49\linewidth}
    \centering\includegraphics[width=0.95\linewidth]{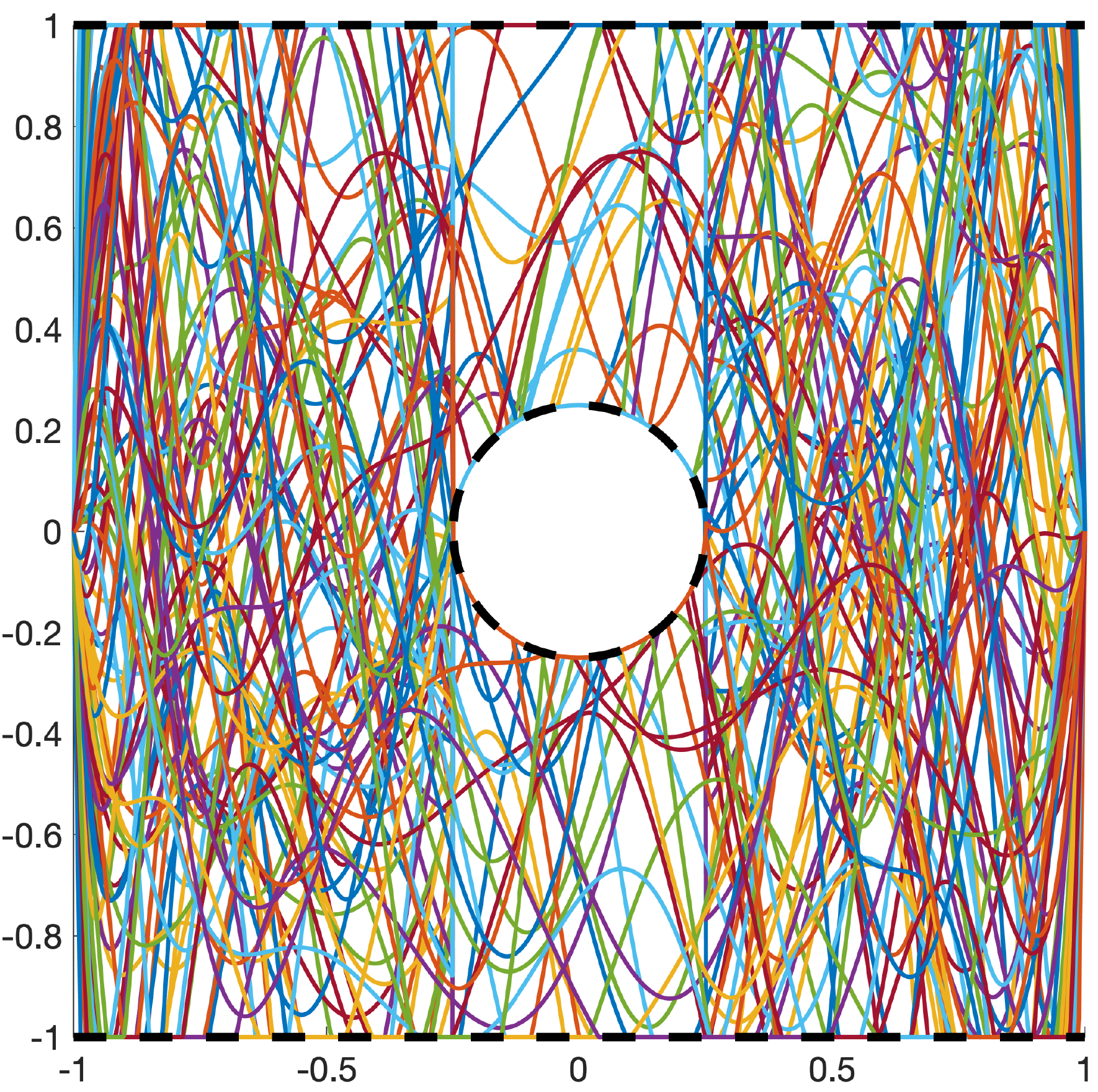}
    \caption{One circle.}
\end{minipage}
\end{figure}
Figures \ref{fig:rectangle} through \ref{fig:saturn} look to push this method to unequally spaced rectangles, circles, and keep-out zones defined by an image. Two important things must be noted about this technique: 1) although the function defining the absolute lower and upper bounds in the example were constrained to $[-1,+1]$, these can be defined by any function similar to the bounds provided in Figures \ref{fig:ineqConst} and \ref{fig:ineqAndEqConst}; 2) the major drawback of this method is that the search space scales with the number of possible trajectories, and therefore, the number of constrained expressions also increases. This implies that any optimization technique using this structure would produce the optimal trajectory for each path. As the number of paths increases drastically, this could become computationally expensive. Additionally, regardless of this method's flexibility, since it is only one-dimensional $y(x,g(x))$, it cannot be used in path planning problems. The next sections explores a two-dimensional, parametric space formulation where $x(t)$ and $y(t)$.
\begin{figure}[H]
\centering
\begin{minipage}[t]{0.49\linewidth}
    \centering\includegraphics[width=0.94\linewidth]{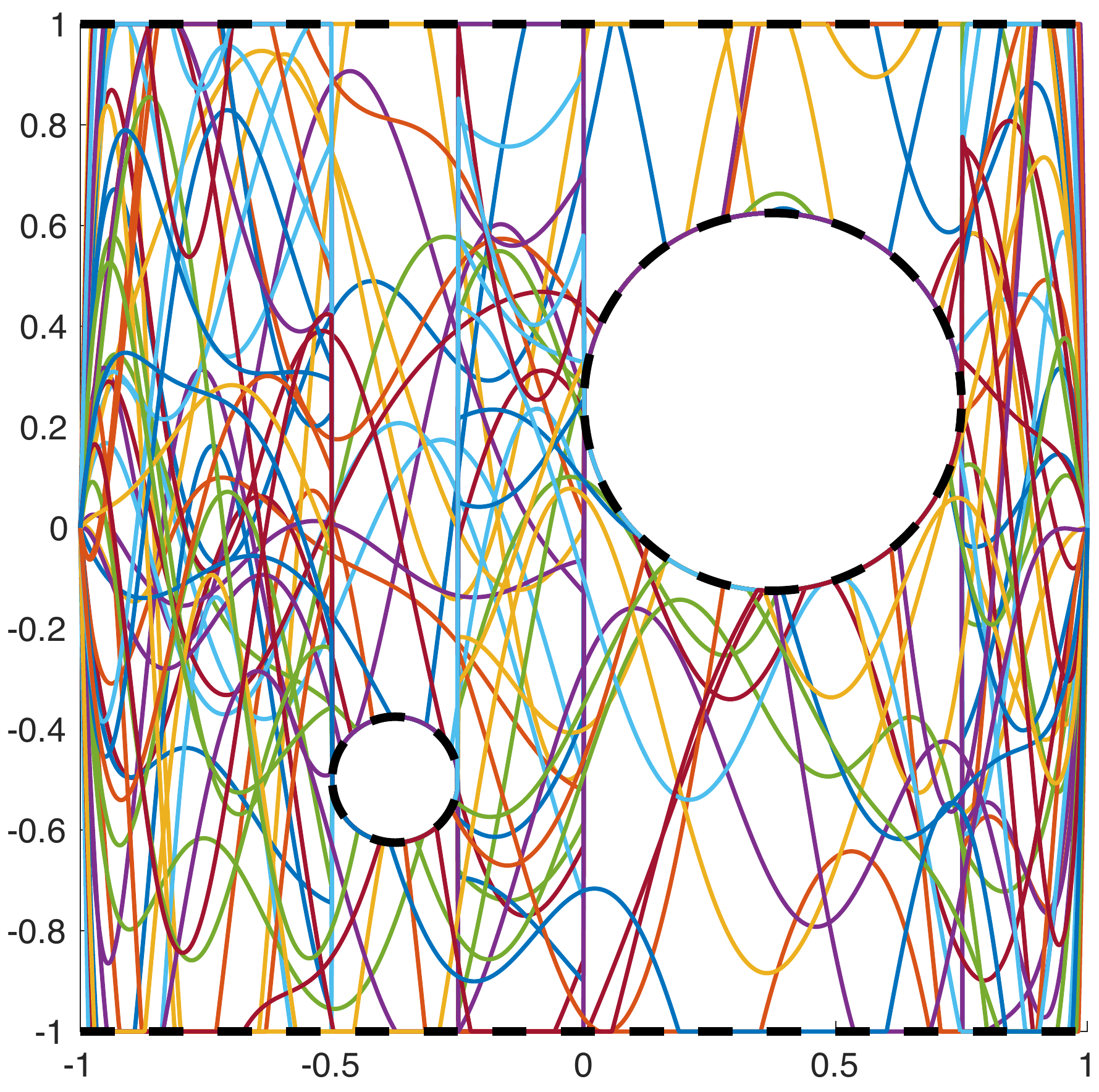}
    \caption{Two circles.}
\end{minipage}
\begin{minipage}[t]{0.49\linewidth}
    \centering\includegraphics[width=0.95\linewidth]{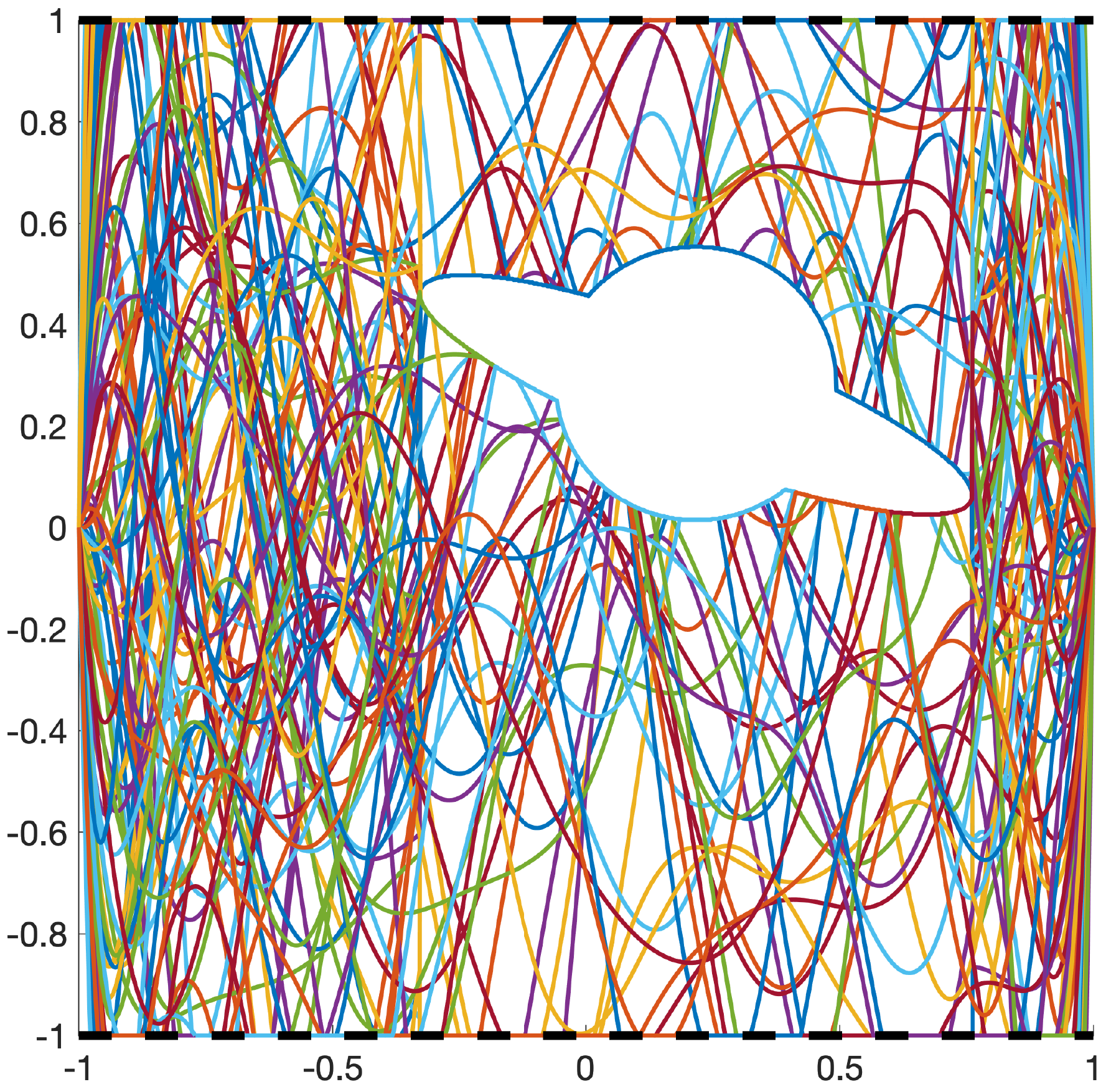}
    \caption{Random object.}
    \label{fig:saturn}
\end{minipage}
\end{figure}
\end{example}

\subsection{Toward 2D inequality constraints}
\begin{figure}[ht]
    \centering\includegraphics[width=0.5\linewidth]{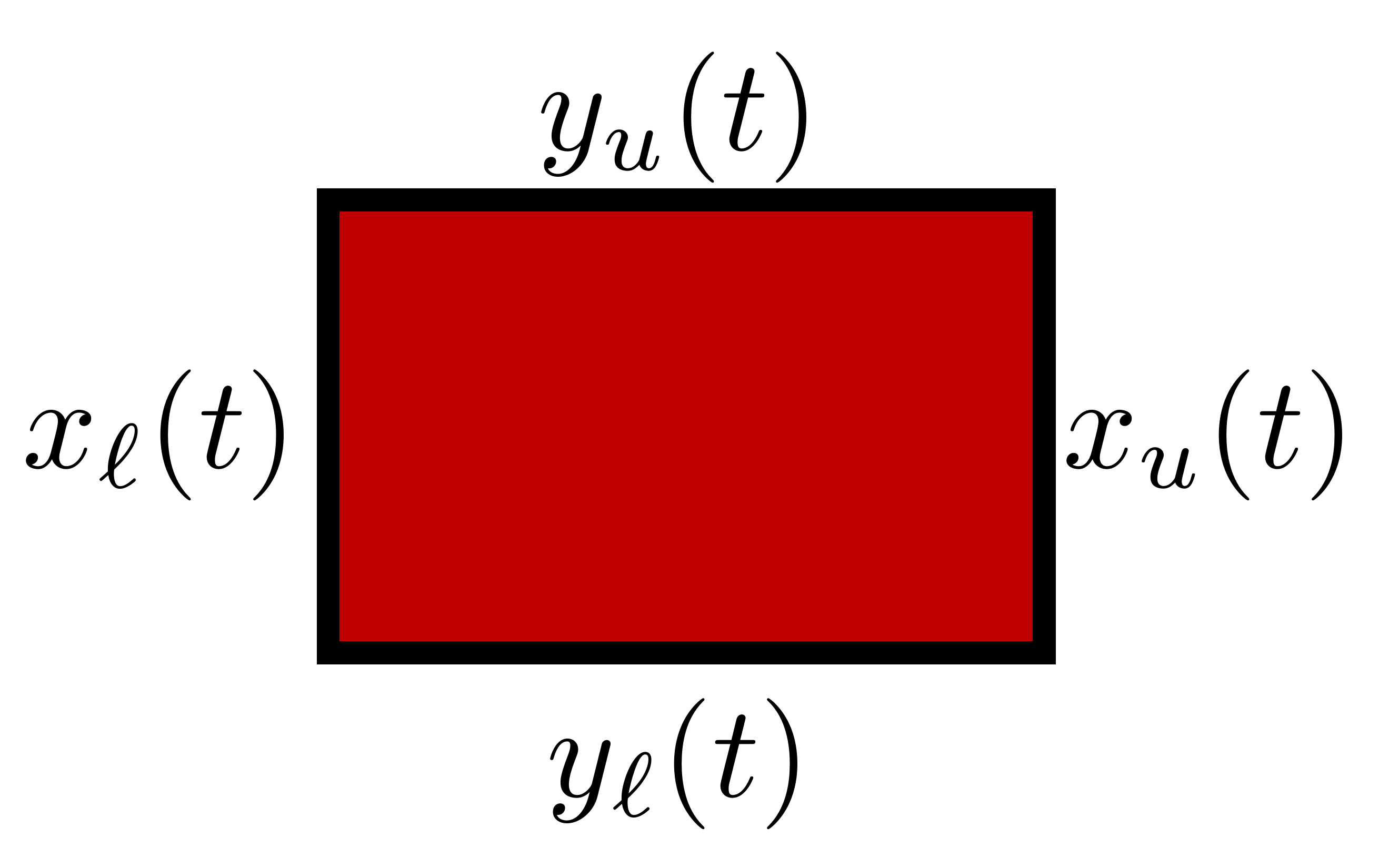}
    \caption{Conceptual keep-out box.}
    \label{fig:keepOut}
\end{figure}

For this theory to be extended for path planning, the constrained expressions must be defined parametrically. For simplicity, let us consider a keep-out box defined by $x_{\ell}(t)$, $x_{u}(t)$, $y_{\ell}(t)$, and $y_{u}(t)$ as defined in Figure \ref{fig:keepOut}. In general, the keep-out zone could be dynamic; however, for now, let us consider the simple example of a static rectangular keep-out region. 
For this formulation let's define the path in terms of parametric variable $t$, using the functionals $x(t,g(t))$ and $y(t,g(t))$ for $x$-position and $y$-position respectively. Associated with these two functions we use the free functions $g^x(t)$ and $g^y(t)$ in defining the constrained expression. Using the TFC method the constrained expression for $x(t,g^x(t))$ and $y(t,g^y(t))$ are as follows,
\begin{equation}\label{eq:2DIneqX}
\begin{aligned}
    x(t,g^x(t)) = g^x(t) &+\big[x_{\ell}(t)-g^x(t)\big]\bigg[\Hs_0( \varphi_1^y)\Hs_0(\varphi_2^y)
    \Hs_0(\varphi_3^x)\Hs_0(\varphi_2^x)\bigg] \\ &+
     \big[x_{u}(t)-g^x(t)\big]\bigg[\Hs_0( \varphi_1^y)\Hs_0(\varphi_2^y)
    \Hs_0(-\varphi_3^x)\Hs_0(\varphi_1^x)\bigg],
\end{aligned}
\end{equation}
\begin{equation}\label{eq:2DIneqY}
\begin{aligned}
    y(t,g^y(t)) = g^y(t) &+\big[y_{\ell}(t)-g^y(t)\big]\bigg[\Hs_0( \varphi_1^x)\Hs_0(\varphi_2^x)
    \Hs_0(\varphi_3^y)\Hs_0(\varphi_2^y)\bigg] \\ &+
     \big[y_{u}(t)-g^y(t)\big]\bigg[\Hs_0( \varphi_1^x)\Hs_0(\varphi_2^x)
    \Hs_0(-\varphi_3^y)\Hs_0(\varphi_1^y)\bigg].
\end{aligned}
\end{equation}
The functions of $\varphi$ (referred to as pseudo-switching functions due to their similarity with the true switching functions defined in the prior sections) are defined in Table \ref{tab:switch}, where $c$ is replaced with either the component $x$ or $y$.
\begin{table}[ht]
    \caption{Pseudo-switching functions for Heaviside functions.}
    \begin{center}
    \begin{tabular}{|c|c|c|c|}
    \hline
    $\varphi_1^c(t)$ & $\varphi_2^c(t)$ & $\varphi^c_3(t)$ \\ \hline\hline
    $c_u(t) - g^c(t)$ & $g^c (t) - c_{\ell} (t)$ & $\dfrac{c_u (t) + c_{\ell}(t)}{2} - g^c(t)$ \\ \hline\hline
    \end{tabular}
    \end{center}
    \label{tab:switch}
\end{table}
The constrained expressions in Equations \eqref{eq:2DIneqX} and \eqref{eq:2DIneqY} are populated by three specific terms, and the interpretation for the $x$-component constrained expression is detailed below (note, the $y$-component constrained expression is of the same structure):
\begin{itemize}
\item The first term is the free functions for the $x$-component.
\item The second term deals with the projection of the lower boundary and has four sigmoid functions as inputs to a 4-way AND gate that is true if and only if the following conditions are met:
    \begin{itemize}
        \item $\varphi_1^y(t)$: the current path's $y$-position is less than $y_{u}(t)$
        \item $\varphi_2^y(t)$: the current path's $y$-position is greater than $y_{\ell}(t)$
        \item $\varphi_3^x(t)$: the current path's $x$-position is less than the average value of $x_u(t)$ and $x_{\ell}(t)$
        \item $\varphi_2^x(t)$: the current path's $x$-position is greater than $x_{\ell}(t)$
    \end{itemize}
    If these four conditions are true, than the current path is inside of the box and closer to the $x_{\ell}(t)$ line than the $x_u(t)$ line. In this case, the line is projected onto the $x_{\ell}(t)$ line by adding the difference between $x_{\ell}(t)$ and $g^x(t)$ to $g^x(t)$.
\item The third term functions in a similar way to the second term, except in this case it deals with the projection of the upper boundary and has four sigmoid functions as inputs to a 4-way AND gate that is true if and only if the following conditions are met:
    \begin{itemize}
        \item $\varphi_1^y(t)$: the current path's $y$-position is less than $y_{u}(t)$
        \item $\varphi_2^y(t)$: the current path's $y$-position is greater than $y_{\ell}(t)$
        \item $-\varphi_3^x(t)$: the current path's $x$-position is greater than the average value of $x_u(t)$ and $x_{\ell}(t)$
        \item $\varphi_1^x(t)$: the current path's $x$-position is greater than $x_{u}(t)$
    \end{itemize}
    If these four conditions are true, than the current path is inside of the box and closer to the $x_{u}(t)$ line than the $x_{\ell}(t)$ line. In this case, the line is projected onto the $x_{u}(t)$ line by adding the difference between $x_{u}(t)$ and $g^x(t)$ to $g^x(t)$.
\end{itemize}

\begin{example}{2D inequality constraints}
In order to analyze this technique, a numerical test was constructed using the constrained expressions given by Equations \eqref{eq:2DIneqX} and \eqref{eq:2DIneqY} where the terms $g^x(t)$ and $g^y(t)$ were defined according to Equation \eqref{eq:basis} with $\B{\xi} \thicksim \mathcal{N}(0,\sigma^2 I_{m \times m})$ where $\sigma = 0.1$. The following tests show an example of single box path avoidance and multiple object path avoidance with boundary conditions such that $(x(t_0),y(t_0)) = (-1,-1)$ and $x(t_f),y(t_f)) = (+1,+1)$. The trajectories of both tests are shown in Figures \ref{fig:box2d} and \ref{fig:multiplebox2d}. It can be seen in both tests the trajectories avoid the boundary displayed by the dashed black line. Additionally, the initial and final constraints on position are always met exactly. Going further, Figure \ref{fig:smoothbox2d} shows specifically the ``smooth'' trajectories that avoid the keep-out zones.
\begin{figure}[H]
    \centering\includegraphics[width=0.55\linewidth]{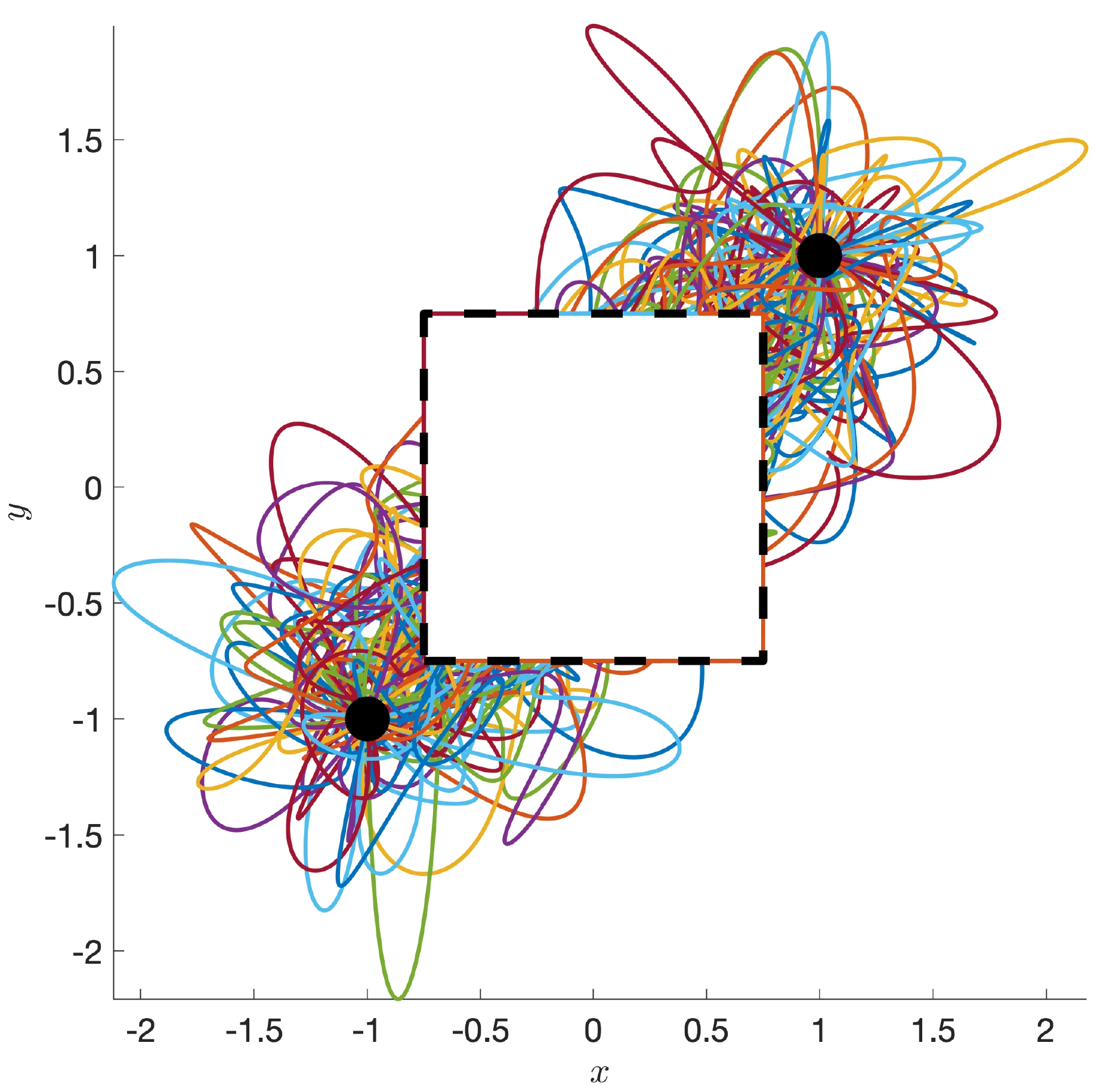}
    \caption{Keep-out box in parametric space.}
    \label{fig:box2d}
\end{figure}
\begin{figure}[H]
    \centering\includegraphics[width=0.55\linewidth]{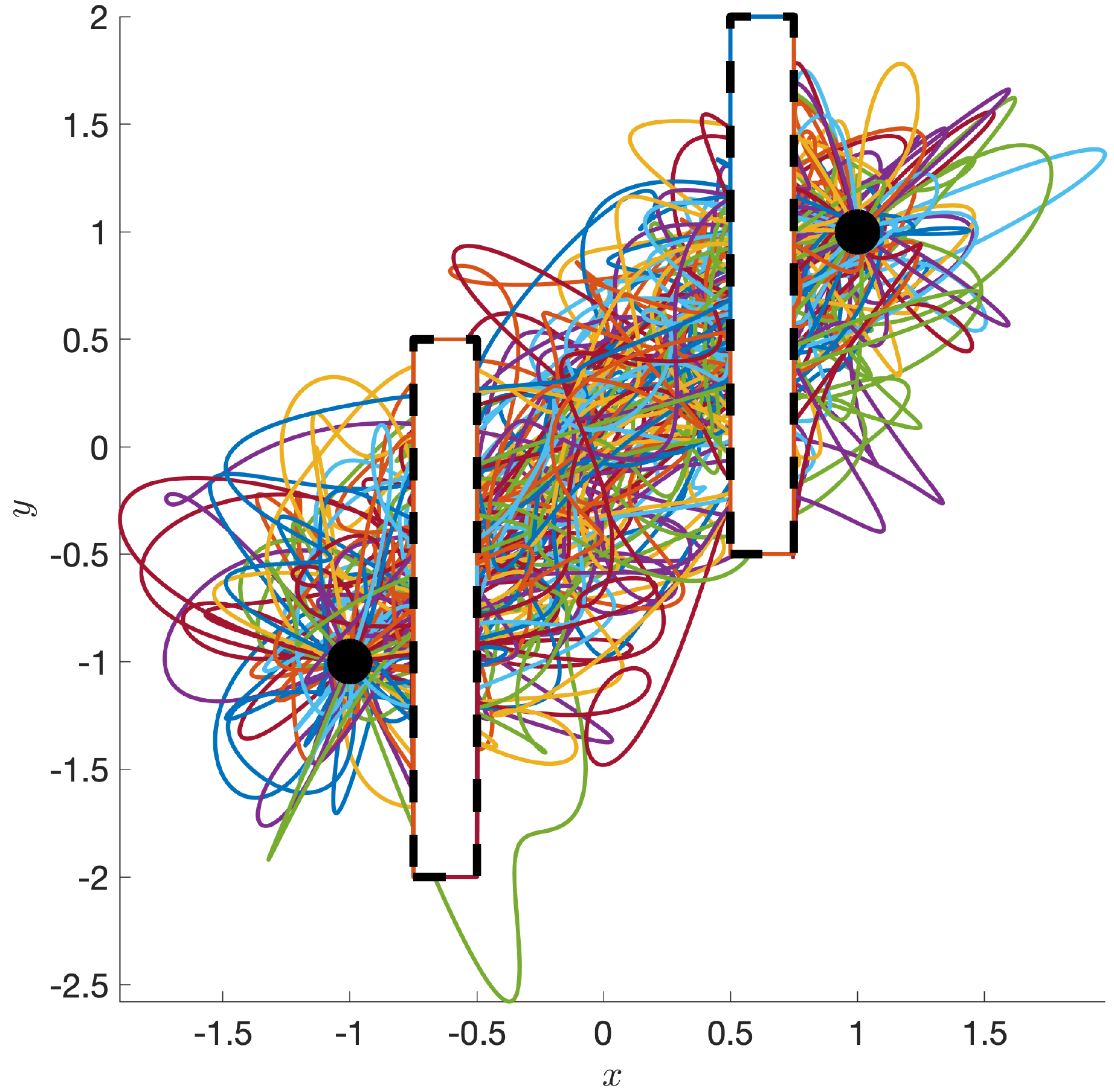}
    \caption{Multiple keep-out zones for parametric formulation.}
    \label{fig:multiplebox2d}
\end{figure}
\begin{figure}[H]
    \centering\includegraphics[width=0.55\linewidth]{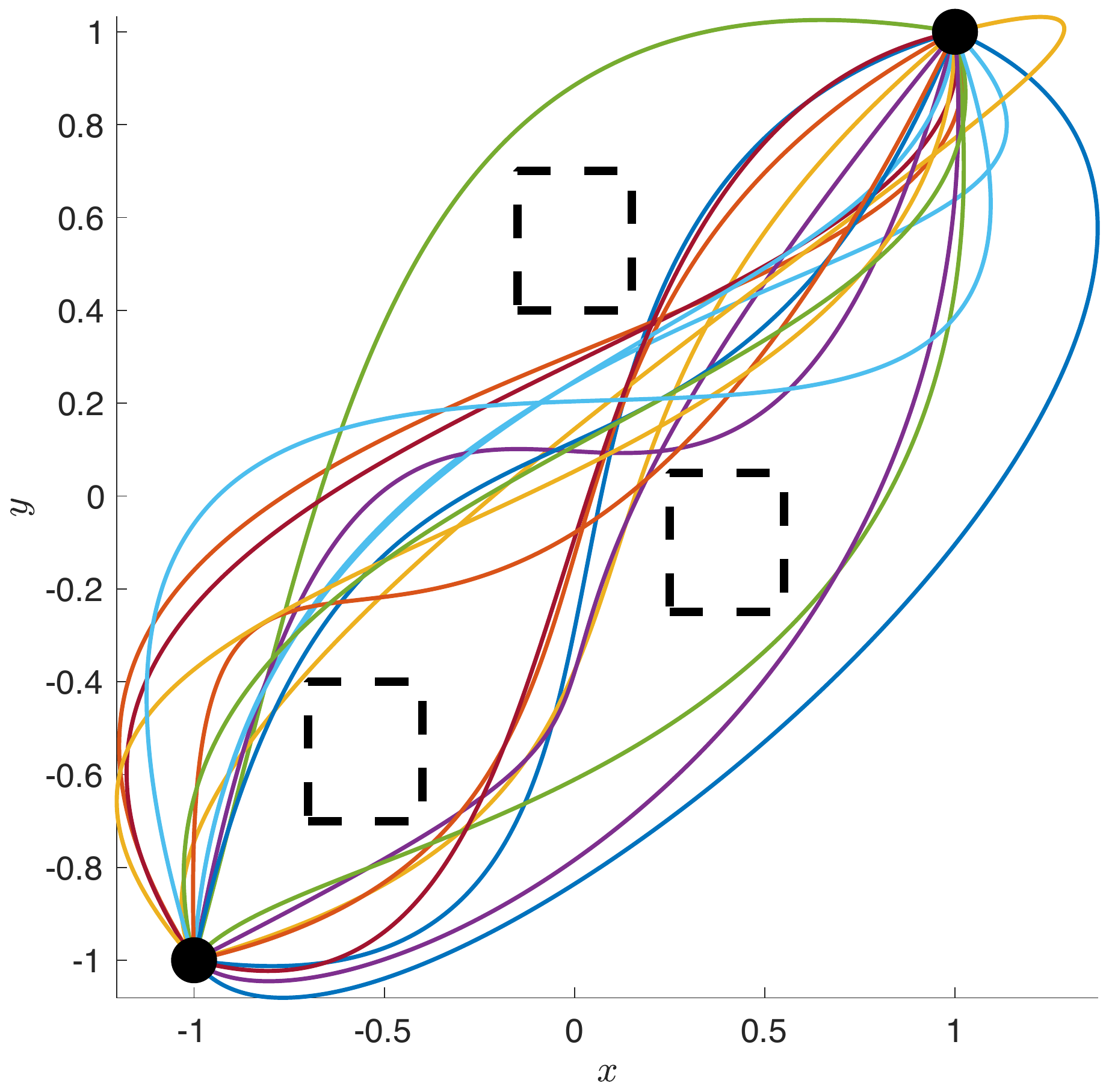}
    \caption{``Smooth'' trajectories avoiding three box keep-out zones.}
    \label{fig:smoothbox2d}
\end{figure}
\end{example}

Although the test shows favorable results, there are potential issues when applying this formulation to optimization problems, namely path planning problems. For example, this method will try to project lines inside the box towards one of the corners. Since the Heaviside functions act as the switches in this problem, there could be cases where lines ``snap" to the corners.

\section{Over-constrained problems}\label{sec:s2a_overCON}
\begin{figure}[ht]
	\centering\includegraphics[width=.8\linewidth]{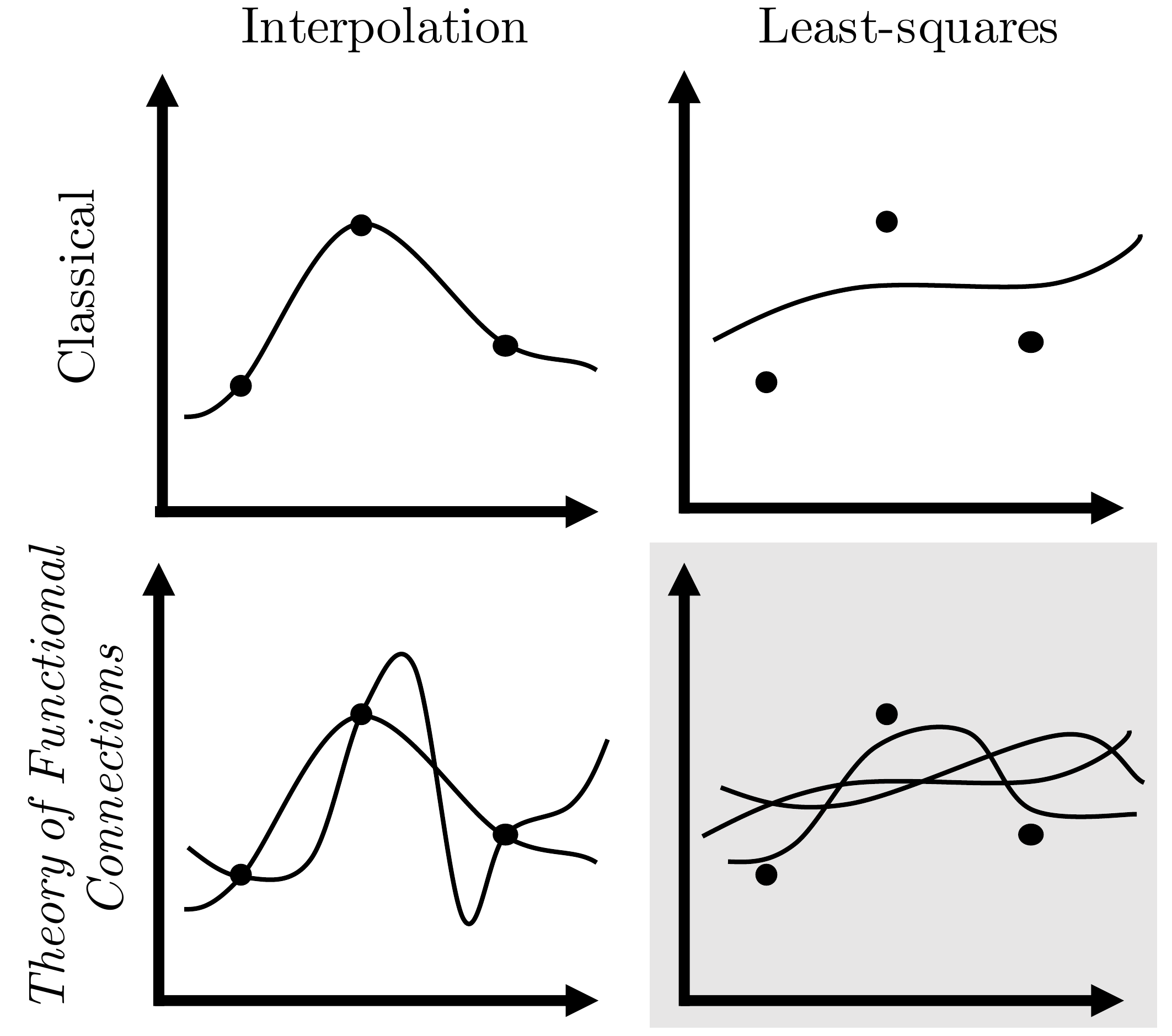}
    \caption{General illustration of classic and TFC approaches for interpolation and least-squares.}
    \label{fig:motivation}
\end{figure}

As we have seen, the TFC framework can incorporate any linear constraints like those developed in the previous sections. Figure \ref{fig:motivation} provides an outline that distinguishes the TFC approach from classical methods in interpolation and least-squares. In the prior development of TFC, the number of the support function was equal to the number of constraints incorporated. It was shown to be a general interpolation approach that described \emph{all} functions passing through $k$ constraints. This section's theory, highlighted in the grey box in Figure \ref{fig:motivation}, combines this general interpolation method with a weighted least-squares technique for the constraints. Doing this allows for a constrained expression to be derived where the number of constraints is greater than the number of support functions, $s_i(x)$, producing a \emph{weighted} constrained expression. Using this expression produces a family of functions minimizing the weighted sum of squares of the constraints. This extension then provides the framework for the solution of over-constrained differential equations (a topic discussed in Chapter \ref{chap:ode}).

\subsection{Two constraints in one degree of freedom}
Consider a constrained expression such that,
\begin{equation*} 
    y(x,g(x))  = g(x)  + \varphi_1(x) \rho_1(x,g(x)) + \varphi_2(x) \rho_2(x,g(x)),
\end{equation*}
Note, in this expression, the $\varphi_i(x)$ is used for the switching functions because for the over-constrained cases these functions do not act like the switching functions, $\phi_i(x)$, discussed earlier. However, as will be seen later, the function $\varphi_i(x)$ can collapse to $\phi_i(x)$. Moving forward, it is desired that this function be subject to two constraints such that,
\begin{equation*}
    \begin{cases} y^{(i)} (x_1) = y^{(i)}_1 \\ y^{(j)} (x_2) = y^{(j)}_2\end{cases} \qquad \text{where} \qquad i, j \in \mathbb{Z}.
\end{equation*}
First, consider the support function as $s(x)$ which will be evaluated at both constraint locations. Applying TFC produces an over-constrained system since there are two constraints but only one support function,
\begin{equation*}
    \begin{Bmatrix} s^{(i)}(x_1) \\ s^{(j)}(x_2)\end{Bmatrix} \begin{Bmatrix} \alpha_1 & \alpha_2 \end{Bmatrix}= \begin{bmatrix}
    1 & 0 \\ 0 & 1
    \end{bmatrix}
\end{equation*}
Therefore, this system can be solved by a weighted least-squares technique where $W$ represents a diagonal matrix of the relative weights.
\begin{equation*}
    W \mathbb{S} \begin{Bmatrix} \alpha_1 & \alpha_2 \end{Bmatrix} =  W
\end{equation*}
This system is then solved for the $\alpha$ coefficients just like in the traditional TFC approach,
\begin{equation*}
     \begin{Bmatrix} \alpha_1 & \alpha_2 \end{Bmatrix} =  \Big(\mathbb{S}\T W \mathbb{S}\Big)^{-1} \mathbb{S}\T W
\end{equation*}
which leads to the expressions
\begin{align*}
    \varphi_1(x) = s(x) \alpha_1  = \dfrac{s(x) s^{(i)}(x_1)w_1}{(s^{(i)}(x_1))^2w_1 + (s^{(j)}(x_2))^2w_2}\\ 
    \varphi_2 = s(x) \alpha_2  = \dfrac{s(x) s^{(j)}(x_2)w_2}{(s^{(i)}(x_1))^2w_1 + (s^{(j)}(x_2))^2w_2}\\
\end{align*}
The final constrained expression is realized as,
\begin{align}
    y(x,g(x)) = g(x) \, &+ \Bigg[\dfrac{s(x) s^{(i)}(x_1)w_1}{(s^{(i)}(x_1))^2w_1 + (s^{(j)}(x_2))^2w_2} \Bigg] \Big(y^{(i)}_1 - g^{(i)}(x_1)\Big) \nonumber \\ &+ \Bigg[\dfrac{s(x) s^{(j)}(x_2)w_2}{(s^{(i)}(x_1))^2w_1 + (s^{(j)}(x_2))^2w_2}\Bigg]\Big(y^{(j)}_2 - g^{(j)}(x_2)\Big). \label{eq:ce_2_solved}
\end{align}
Yet, there remains some conditions on the functions of $g(x)$ and $s(x)$. First, the function $g(x)$ must be differentiable up to the $i$-th and $j$-th derivative. Additionally, by analyzing terms $\varphi_1(x)$ and $\varphi_2(x)$, it can be seen that information on the constraints is lost when $\varphi_1(x)$ or $\varphi_2(x)$ becomes zero. Therefore, the support function $s(x)$ must be selected such that $s^{(i)}(x_1) \ne 0$ and $s^{(j)}(x_2) \ne 0$. Let us now consider the weighting scheme $w_1 = 1$ and $w_2 = 0$. In this case, Equation \eqref{eq:ce_2_solved} reduces to a familiar form,
\begin{equation*}
    y(x,g(x)) = g(x) + \frac{s(x)}{s^{(i)}(x_1)} (y^{(i)}_1 - g^{(i)}(x_1)),
\end{equation*}
which represents one constraint at one point. With this is example in mind, the following sections explore the characteristics of the \emph{weighted} constrained expression for multiple points.

\subsection{Weighted constraints at two points}

\begin{figure}[ht!]
	\centering
    \includegraphics[width=.8\linewidth]{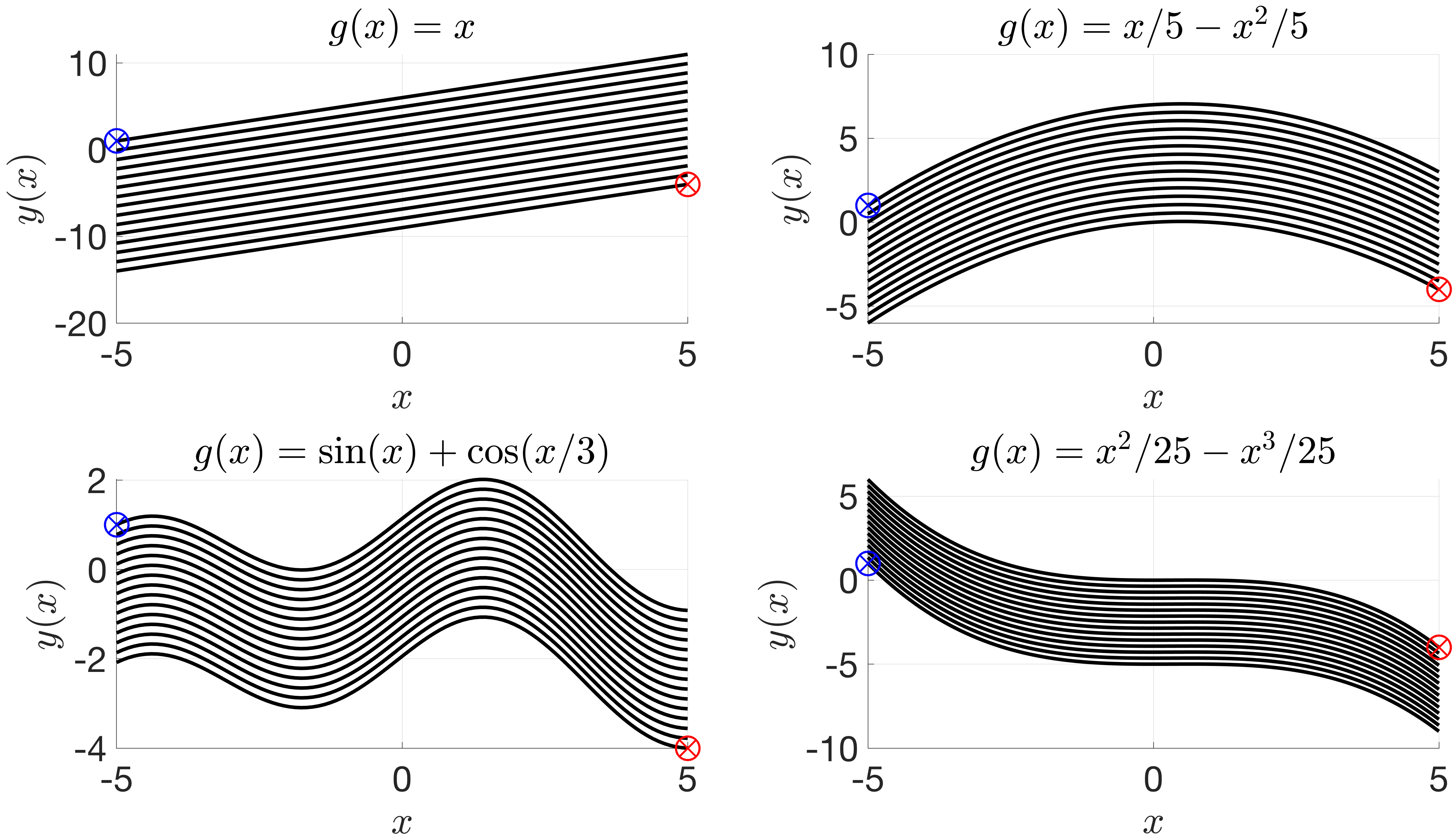}
    \caption{Analysis of Equation (\ref{eq:pound}) for varying values of $g(x)$. It follows that as $\gamma$ increases from 0 to 1, the function translates between the constraint conditions.}
    \label{fig:two_points_multiple_functions}
\end{figure}

\begin{figure}[ht!]
	\centering\includegraphics[width=.8\linewidth]{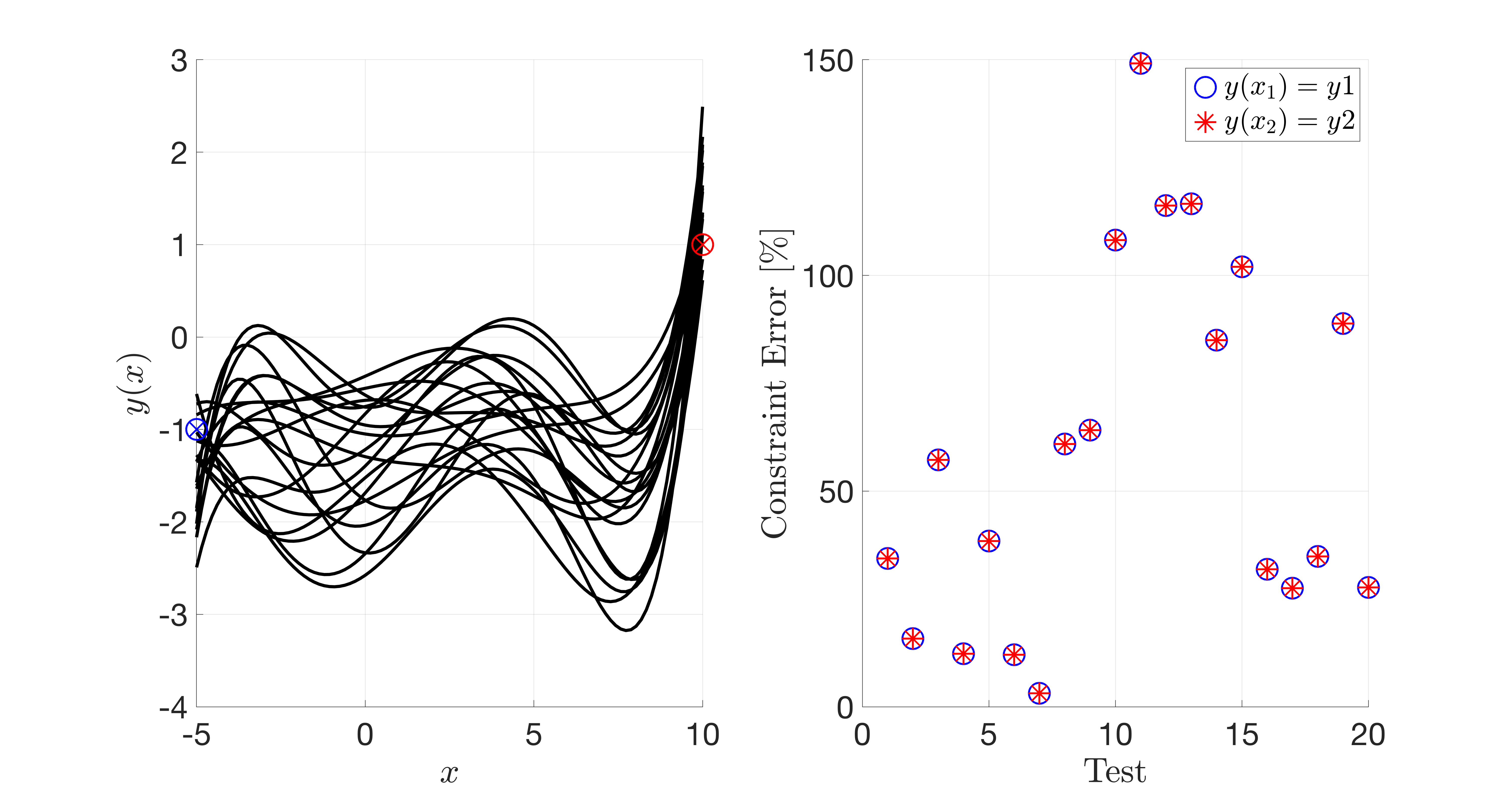}
    \caption{Analysis of Equation (\ref{eq:double_star}) for 20 randomly selected $g(x)$'s. The relative error between constraints is the same for every test.}
    \label{fig:two_points_random_g}
\end{figure}

\begin{example}{Weight-constrained expression for two points}
As an example, let us consider constraints at two points such that,
\begin{equation*}
        y (x_1) = y_1 \quad \text{and} \quad y (x_2) = y_2,
\end{equation*}
which implies that $i = j = 0$. For these constraints, Equation (\ref{eq:ce_2_solved}) reduces to,
\begin{align*}
    y(x,g(x)) = g(x) &+ \left[\dfrac{s(x) \, s(x_1) \, w_1}{s(x_1)^2 \, w_1 + s(x_2)^2 \, w_2}\right] \Big(y_1 - g(x_1)\Big) \\ &+ \left[\dfrac{s(x) \, s(x_2) \, w_2}{s(x_1)^2 \, w_1 + s(x_2)^2 \, w_2}\right] \Big(y_2 - g(x_2)\Big).
\end{align*}
The simplest definition of $s(x)$ such that $s(x_1) \ne 0$ and $s(x_2) \ne 0$ is $s(x) = 1$, leading to the equation,
\begin{equation}\label{eq:ce_twopoints}
    y(x,g(x))  = g(x)  + \left(\dfrac{w_1}{w_1 + w_2}\right) \Big(y_1 - g(x_1)\Big) + \left(\dfrac{w_2}{w_1 + w_2}\right) \Big(y_2 - g(x_2)\Big).
\end{equation}
\end{example}
Analyzing this function, it can be seen that $g(x)$ is the only non-constant term in the equation and all other terms represent the relative weights of the prescribed constraints. Moreover, this equation represents every function that when evaluated at the constraint locations satisfies them relative to the prescribed weights $w_1$ and $w_2$. By setting $w_1 = 1$ and $w_2 = 0$, Equation \eqref{eq:ce_twopoints} reduces to a constrained expression for one point,
\begin{equation*}
    y(x,g(x)) = g(x) + (y_1 - g(x_1)).
\end{equation*}
If $w_1 = 0$ and $w_2 = 1$ is selected, an equation satisfying $y (x_2) = y_2$ is obtained. This gives reason to believe that the weighted least-squares solution occupies the set of functions between these two absolute constraints. Keying in on this notion, let us explore the parametric weight scheme,
\begin{equation}\label{eq:alpha_weights}
    W (\gamma) = \begin{bmatrix} 1 - \gamma & 0 \\ 0 & \gamma\end{bmatrix}, \qquad \text{where} \qquad \gamma \in [0, 1].
\end{equation}
Using these weights, Equation (\ref{eq:ce_twopoints}) becomes,
\begin{equation}\label{eq:pound}
    y(x,g(x)) = g(x) + (1 - \gamma)\Big(y_1 - g(x_1)\Big) + \gamma \Big(y_2 - g(x_2)\Big).
\end{equation}
Equation (\ref{eq:pound}) was analyzed for multiple values of $g(x)$ over the range $x \in [-5,+5]$. The results in Figure \ref{fig:two_points_multiple_functions} show that for each function, varying $\gamma$ corresponds to translating the free function between the two prescribed constraints.
Using  $w_1 = w_2 = w$ (constraints equally weighted), Equation (\ref{eq:ce_twopoints}) becomes,
\begin{equation}\label{eq:double_star}
    y(x,g(x)) = g(x) + \dfrac{1}{2} (y_1 - g(x_1)) + \dfrac{1}{2} (y_2 - g(x_2)).
\end{equation}
Analyzing this equation, it is expected that the constraint will be met with the same \emph{relative} error for any function chosen for $g(x)$. Figure \ref{fig:two_points_random_g} shows the results of $20$ randomly generated functions (left plot) and the constraint errors (right plot).
In Figure \ref{fig:two_points_random_g}, this ``constraint error'' is simply showing that since the projection functionals are equally weighted, the error from their imposed value (either $y_1$ or $y_2$) is the same.

\subsection{Constraints on a function and its derivative}

\begin{figure}[ht]
	\centering\includegraphics[width=.8\linewidth]{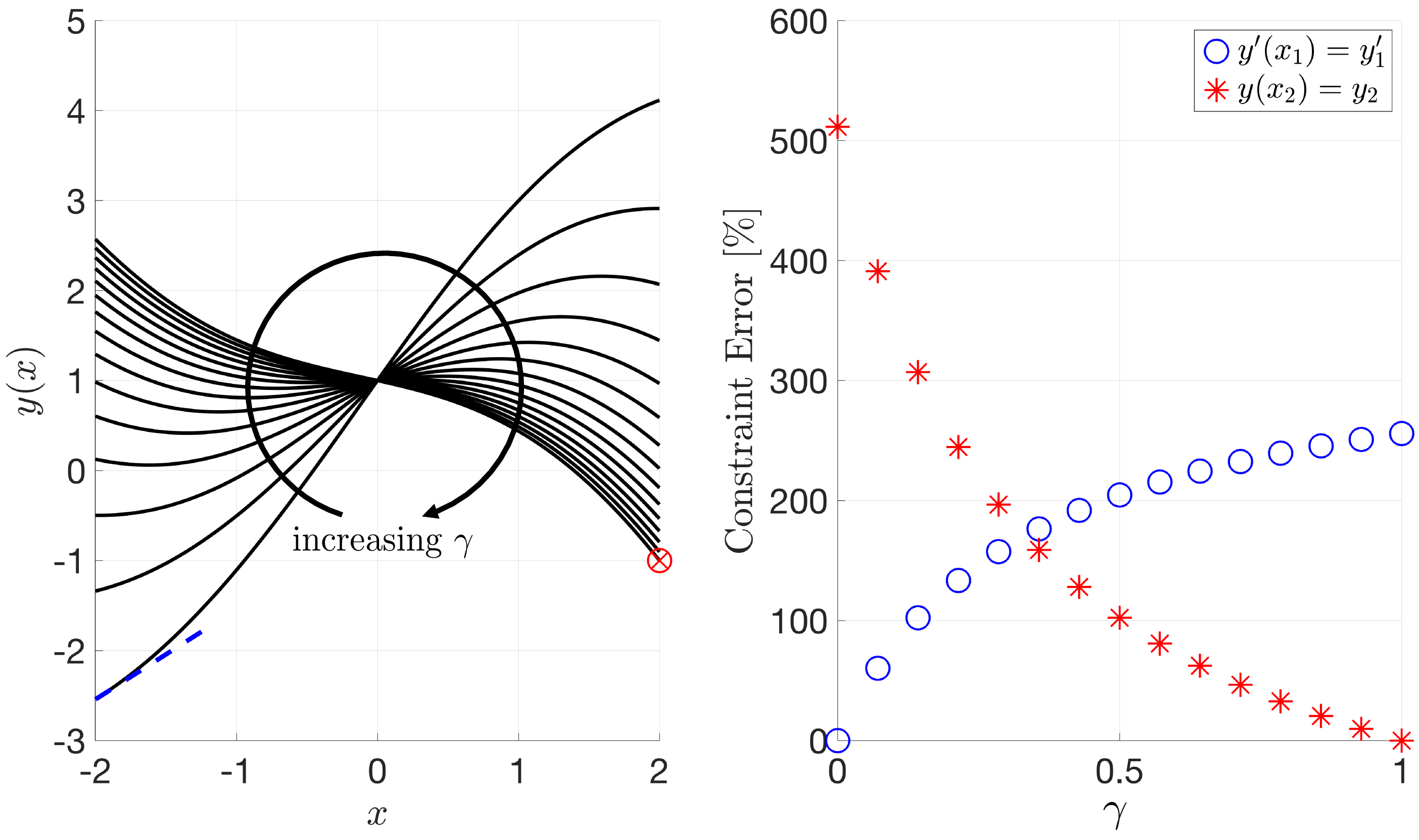}
    \caption{Equation (\ref{eq:forgot}) for varying weight values $\gamma$ using the free function $g (x) = \sin x + \cos (x/3)$.}
    \label{fig:point_and_derivative}
\end{figure}

\begin{example}{Constraints on a function and its derivative}
For further analysis, let us consider the case of constraints on a function and its derivative such that $i = 1$ and $j = 0$,
\begin{equation*}
    y_x(x_1) = y_{1_x} \quad \text{and} \quad y(x_2) = y_2,
\end{equation*}
which reduces Equation (\ref{eq:ce_2_solved}) to,
\begin{align}\label{eq:ce_p_and_d}
    y(x,g(x)) = g(x) &+ \left[\frac{s(x) s_x(x_1) w_1}{(s_x(x_1))^2 w_1 + (s(x_2))^2 w_2} \right] \Big(y_{1_x} - g_x(x_1)\Big) \nonumber\\&+ \left[\frac{s(x) s(x_2) w_2}{(s_x(x_1))^2 w_1 + (s(x_2))^2 w_2}\right] \Big(y_2 - g(x_2)\Big).
\end{align}
For this case, since $i = 1$ is the largest derivative, then $s$ must be defined such that $s_x  \ne 0$. The simplest case is to set $s = x$. Using this definition, Equation (\ref{eq:ce_p_and_d}) becomes,
\begin{equation}\label{eq:forgot}
    y  = g  + \left(\dfrac{w_1 \, x}{w_1 + w_2 \, x_2^2}\right) \Big(y_{1_x} - g_x(x_1)\Big) + \left(\dfrac{w_2 \, x_2 \, x}{w_1 + w_2 \, x_2^2}\right) \Big(y_2 - g(x_2)\Big).
\end{equation}
\end{example}
A similar test can be conducted for this case where $W$ is defined according to Equation (\ref{eq:alpha_weights}). For this particular case, let us define $g(x) = \sin (x) + \cos(x/3)$. Figure \ref{fig:point_and_derivative} shows the transformation from the initial derivative constraint to the final point constraint for various values of $\gamma$.

Additionally, Figure \ref{fig:point_and_derivative} shows the relative constraint error for each constraint as a function of the $\gamma$ parameter. In the next section, this method will be applied to three constraints with two degrees of freedom.

\subsection{Three constraints with two degrees of freedom}\label{sec:s2a_threePoints}

Now consider a constrained expression with two degrees of freedom defined as,
\begin{equation}\label{eq:ce2}
    y(x,g(x)) = g(x) + \varphi_1(x) \rho_1(x,g(x)) + \varphi_2(x) \rho_2(x,g(x)),
\end{equation}
where $s_1(x)$ and $s_2(x)$ are assigned functions and the constraints are defined such that,
\begin{equation*}
    \begin{cases} y^{(i)}(x_1) = y^{(i)}_1 \\ y^{(j)}(x_2) = y^{(j)}_2 \\ y^{(k)}(x_3) = y^{(k)}_3\end{cases} \qquad \text{where} \qquad i, j, k \in \mathbb{Z}.
\end{equation*}
Applying these constraints leads to the system of equations,
\begin{equation}\label{eq:delta}
   W \begin{bmatrix} s_1^{(i)} (x_1) & s_2^{(i)} (x_1) \\ s^{(j)}_1 (x_2) & s^{(j)}_2 (x_2) \\ s^{(k)}_1 (x_3) & s^{(k)}_2 (x_3)\end{bmatrix} \begin{bmatrix} \alpha_{11} & \alpha_{12} & \alpha_{13} \\ \alpha_{21} & \alpha_{22} & \alpha_{23} \end{bmatrix} = W \begin{bmatrix} 1 & 0 & 0 \\ 0 & 1 & 0 \\ 0 & 0 & 1\end{bmatrix}.
\end{equation}
Although Equation (\ref{eq:delta}) is expressed for three constraints, there is no upper limit on the number of constraints that can be incorporated. 

\begin{example}{Three constraints with two degrees of freedom}
Now, let us use the specific formulation, given by Equation \eqref{eq:ce2}, to derive an over-constrained expression with three point constraints.  Incorporating these constraints ($i = j = k = 0$), the system of equations in Equation \eqref{eq:delta} reduces to,
\begin{equation*} 
    W \begin{bmatrix} s_1 (x_1) & s_2 (x_1) \\ s_1 (x_2) & s_2 (x_2) \\ s_1 (x_3) & s_2 (x_3)\end{bmatrix} \begin{bmatrix} \alpha_{11} & \alpha_{12} & \alpha_{13} \\ \alpha_{21} & \alpha_{22} & \alpha_{23} \end{bmatrix} = W .
\end{equation*}
For this problem, let us define $s_1(x)  = 1$, $s_2(x) = x$, and the diagonal weight matrix as,
\begin{equation*}
     W = \begin{bmatrix} w_1 & 0 & 0 \\ 0 & w_2 & 0 \\ 0 & 0 & w_3\end{bmatrix}.
\end{equation*}
By solving the system using the weighted least-squares technique the over-constrained switching functions become,
\begin{align*}
    \varphi_1(x) &= s_1(x) \alpha_{11} + s_2(x) \alpha_{21} =  \frac{w_1}{D}\Big(w_2 (x_2-x) \Delta_{21} + w_3 (x_3-x) \Delta_{31}\Big) \\
    \varphi_2(x) &= s_1(x) \alpha_{12} + s_2(x) \alpha_{22} = \frac{w_2}{D} \Big(w_1 (x-x_1) \Delta_{21} + w_3 (x_3-x) \Delta_{32}\Big) \\
    \varphi_3(x) &= s_1(x) \alpha_{13} + s_2(x) \alpha_{23} = \frac{w_3}{D} \Big(w_1 (x-x_1) \Delta_{31} + w_2 (x-x_2) \Delta_{32}\Big)
\end{align*}
such that $D := \Delta_{21}^2 w_1 w_2 + \Delta_{31}^2 w_1 w_3 + \Delta_{32}^2 w_2 w_3$ and $\Delta_{ij} := x_i - x_j$. Therefore, the over-constrained expression becomes
\begin{equation}\label{eq:three_points}
    y(x,g(x)) = g(x) + \varphi_1(x) (y_1 - g(x_1)) + \varphi_2(x) (y_2 - g(x_2)) + \varphi_3(x) (y_3 - g(x_3)).
\end{equation}
\end{example}

First, let us analyze the simplification when the weights are prescribed as $w_1 = w_2 = 1$ and $w_3 = 0$. Using these weights, Equation (\ref{eq:three_points}) reduces to,
\begin{equation*}
    y(x,g(x)) = g(x)  + \left(\frac{x_2 - x}{x_2 - x_1}\right) (y_1 - g(x_1)) + \left(\frac{x - x_1}{x_2 - x_1}\right)(y_2 - g(x_2)) \\ + [0] (y_3 - g(x_3)),
\end{equation*}
which is the exact constrained expression obtained for the constraints $y(x_1) = y_1$ and $y(x_2) = y_2$ when using the methods developed earlier. Since the constraints are analytically embedded in Equation (\ref{eq:three_points}), the $g(x)$ function represents the solution space that satisfies the three constraints by weighted least-squares. 

While this section simply introduces the over-constrained expression concept, in Section \ref{sec:s3_overConDE}, we will look into using this framework to solve over-constrained differential equations.
%

\chapter{A GENERAL FORMULATION OF THE UNIVARIATE \emph{THEORY OF FUNCTIONAL CONNECTIONS}\label{chap:tfc_general}}

This section rigorously defines the TFC \ce\ and provides some relevant proofs. First, the definition of a functional and properties of a functional are defined.
\begin{definition}{}
A functional, e.g. $f(x,g(x))$, has independent variable(s) and function(s) as inputs, and produces a function as an output.   
\end{definition}
\noindent Note that a function as defined here coincides with the computer science definition of a functional. One can think of a functional as a map for functions. That is, the functional takes a function, $g(x)$, as its input and produces a function, $f^*(x) = f(x,g(x))$ for any specified $g(x)$, as its output. Since this body of work is focused on constraint embedding, or in other words, functional interpolation, we will not concern ourselves with the domain/range of the input and output functions. Rather, we will discuss functionals only in the context of their potential input functions, hereon referred to as the domain of the functional, and potential output functions hereon referred to as the codomain of the functional. 

\begin{figure}[ht!]
    \centering\includegraphics[width=.65\linewidth]{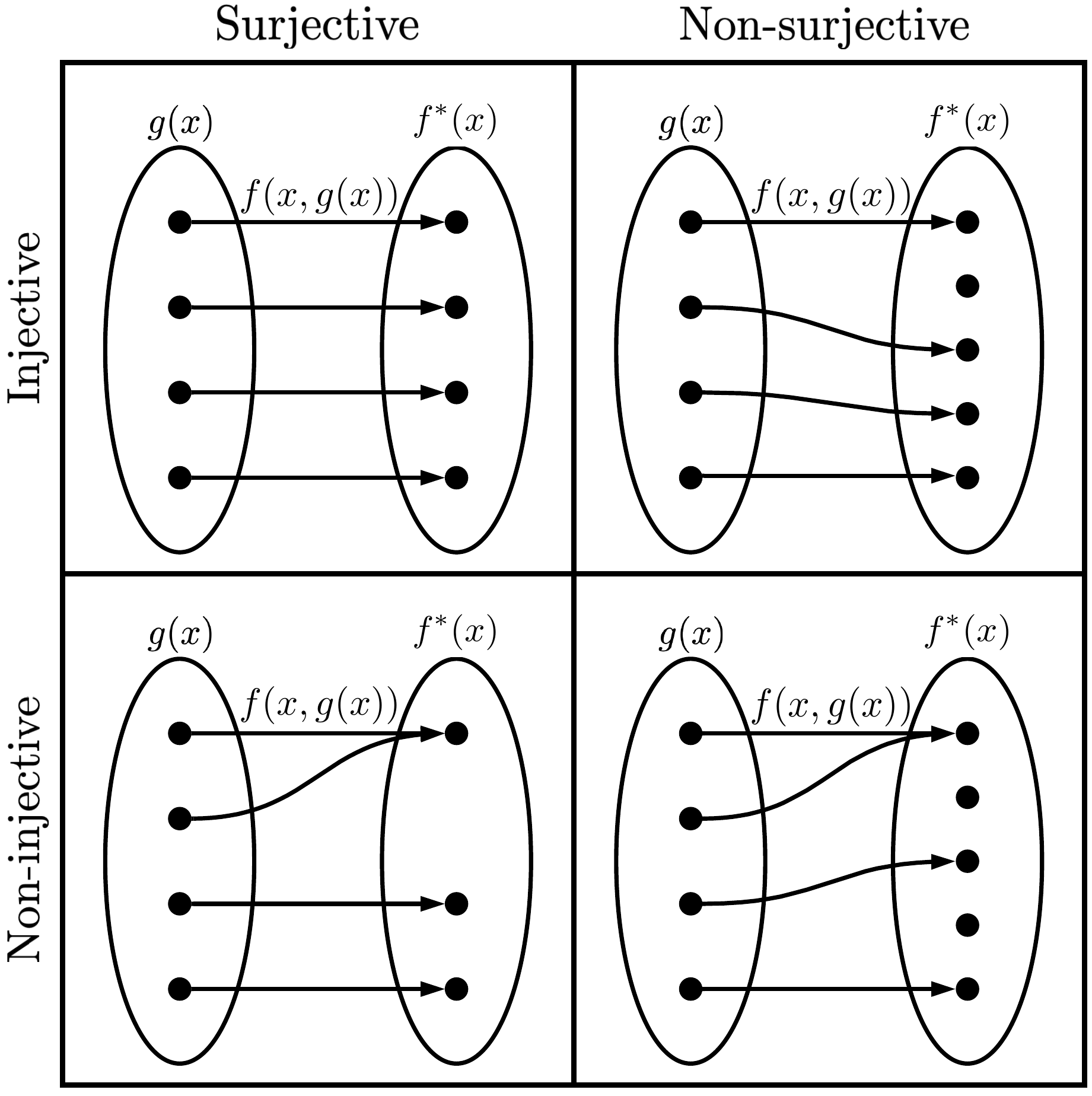}
    \caption{Graphical representation of injective and surjective functionals.}
    \label{fig:Injective_Surjective}
\end{figure}

Next, the definitions of injective, surjective, and bijective are extended from functions to functionals. 
\begin{definition}{}
A functional, $f(x,g(x))$, is said to be injective if every function in its \\codomain is the image of at most one function in its domain. 
\end{definition}
\begin{definition}{}
A functional, $f(x,g(x))$, is said to be surjective if for every function in the codomain, $f^*(x)$, there exists at least one $g(x)$ such that $f^*(x) = f(x,g(x))$.
\end{definition}
\begin{definition}{}
A functional, $f(x,g(x))$, is said to be bijective if it is both injective and surjective. 
\end{definition}
To elaborate, Figure \ref{fig:Injective_Surjective} gives a graphical representation of each of these functionals, and examples of each of these functionals follow. Note that the phrase ``smooth functions'' is used here to denote continuous, infinitely differentiable, real-valued functions. Consider the functional $f(x,g(x)) = e^{-g(x)}$ whose domain is all smooth functions and whose codomain is all smooth functions. The functional is injective because for every $f^*(x)$ in the codomain there is at most one $g(x)$ that maps $f(x,g(x))$ to $f^*(x)$.
However, the functional is not surjective, because the functional does not span the space of the codomain. For example, consider the desired output function $f^*(x) = -2$: there is no $g(x)$ that produces this output. Next, consider the functional $f(x,g(x)) = g(x)-g(0)$ whose domain is all smooth functions and whose codomain is the set of all smooth functions $f^*(x)$ such that $f^*(0)=0$. This functional is surjective because it spans the space of all smooth functions that are $0$ when $x=0$, but it is not injective. For example, the functions $g(x) = x$ and $g(x) = x+3$ produce the same result, i.e., $f(x,x) = f(x,x+3) = x$. Finally, consider the functional $f(x,g(x)) = g(x)$ whose domain is all smooth functions and whose codomain is all smooth functions. This functional is bijective, because it is both injective and surjective. 

Also, the notion of projection is extended to functionals. Consider the typical definition of a projection matrix $P^n = P$ for some $n\in\mathbb{Z}^+$. In other words, when $P$ operates on itself, it produces itself: a projection property for functionals can be defined similarly. 
\begin{definition}{}
A functional is said to be a projection functional if it produces itself when operating on itself.
\end{definition}
For example, consider a functional operating on itself, $f(x,f(x,g(x)))$.  Then, if \newline $f(x,f(x,g(x))) = f(x,g(x))$, then the functional is a projection functional. Note that proving $f(x,f(x,g(x))) = f(x,g(x))$ automatically extends to a functional operating on itself $n$ times: for example, $f(x,f(x,f(x,g(x))) = f(x,f(x,g(x))) = f(x,g(x))$, and so on.

Now that a functional and some properties of a functional have been defined, the notation used in the prior section can be leveraged to rigorously define TFC related concepts. First, it is useful to define the constraint operator, denoted by the symbol $\C{}$.
\begin{definition}{}
The constraint operator, $\C{i}$, is a linear operator that, when operating on a function, returns the function evaluated at the $i$-th specified constraint. 
\end{definition}
As an example, consider the linear constraint $3 = 2 y(2) + \pi y_{xx}(0)$, and suppose it is the first constraint in the set  ($i = 1$). For this constraint, the constraint operator operates as follows,
\begin{equation*}
    \C{1} [y(x)] =  2 y(2) + \pi y_{xx}(0).
\end{equation*}
The constraint operator is a linear operator, as it satisfies the two properties of a linear operator:
\begin{enumerate}
    \item $\C{i} [f(x) + g(x)] = \C{i}[f(x)] + \C{i}[g(x)]$
    \item $\C{i}[a g(x)] = a\C{i}[g(x)]$
\end{enumerate}
For example, again consider the linear constraint $3 = 2 y(2) + \pi y_{xx}(0)$,
\begin{align*}
    \C{1} [f(x)+g(x)] &= \C{1} [f(x)] + \C{1}[g(x)] = 2 f(2) + \pi f_{xx}(0) + 2 g(2) + \pi g_{xx}(0) \\
    \C{1}[a f(x)] &= a \C{1} [f(x)] = a \Big( 2 f(2) + \pi f_{xx}(0)\Big).
\end{align*}
Naturally, the constraint operator has specific properties when operating on the support functions, switching functions, and projection functionals.

\begin{property}{}\label{prop:co_on_s}
The constraint operator acting on the support functions $s_j (x)$ produces the support matrix \begin{equation*}
    \mathbb{S}_{ij} = \C{i}[s_j(x)].
\end{equation*}
\end{property}
For example, consider the two constraints, $y(1)=y(0)$ and $3 = 2 y(2) + \pi y_{xx}(0)$. Applying the constraint operator,
\begin{align*}
    \mathbb{S}_{ij} = \C{i}[s_j(x)] &= \begin{bmatrix}\C{1}[s_1(x)] & \C{1}[s_2(x)] \\ \C{2}[s_1(x)] & \C{2}[s_2(x)] \end{bmatrix} \\ &= \begin{bmatrix}s_1(1) - s_1(0) & s_2(1) - s_2(0) \\ 2s_1(2) + \pi s_{1_{xx}}(0) & 2s_2(2) + \pi s_{2_{xx}}(0) \end{bmatrix}.
\end{align*}
In fact, the support matrix $\mathbb{S}_{ij}$ is simply the matrix multiplying the $\alpha_{ij}$. Therefore, it follows that, $\mathbb{S}_{ij} \, \alpha_{jk} =  \alpha_{ij} \, \mathbb{S}_{jk} = \delta_{ik}$, where $\delta_{ik}$ is the Kroneker delta, and the solution of the $\alpha_{ij}$ coefficients are precisely the inverse of the constraint operator operating on the support matrix. 

\begin{property}{} 
The constraint operator acting on the switching functions $\phi_j(x)$ produces the Kronecker delta.  
\begin{equation*}
    \C{i}[\phi_j(x)] = \delta_{ij}
\end{equation*}
\end{property}
\noindent
This property is just a mathematical restatement of the linguistic definition of the switching function given earlier. One can intuit this property from the switching function definition, since they evaluate to $1$ at their specified constraint condition (i.e., $i=j$) and to $0$ at all other constraint conditions (i.e., $i \neq j$).

Using this definition of the constraint operator, one can define the projection functional in a compact and precise manner.
\begin{definition}{}\label{def:projection_function}
Let $g(x)$ be the free function where $g(x): \R \to \R$, and let $\kappa_i\in\mathbb{R}$ be the numerical portion of the $i^{th}$ constraint. Then,
\begin{equation*}
    \rho_i(x,g(x)) = \kappa_i - \C{i}[g(x)].
\end{equation*}
\end{definition}
Again, consider the linear constraint $3 = 2 y(2) + \pi y_{xx}(0)$. The projection function is,
\begin{align*}
    \rho_1(x,g(x)) &= \kappa_1 - \C{1}[g(x)] \\&= 3 - 2 g(2) - \pi g_{xx}(0).
\end{align*}
Moving forward, we look to leverage the definitions and properties of the TFC formulation to prove a few aspects of the TFC constrained expression that will be useful during numerical implementation.

\begin{theorem}{}\label{thrm:UniCe}
For any function, $f (x)$, satisfying the constraints, there exists at least one free function, $g (x)$, such that the TFC \ce\ $y(x,g(x)) = f(x)$.

\noindent\hrulefill

\noindent \textbf{Proof:}
As highlighted in Properties \ref{prop:proj1}, \ref{prop:proj2}, \ref{prop:proj3}, and \ref{prop:proj4}, the projection functionals are equal to zero whenever $g(x)$ satisfies the constraints. Thus, if $g(x)$ is a function that satisfies the constraints, then the \ce\ becomes,
\begin{align*}
    y (x, g (x)) &= g (x) + \rho_i (x, g(x)) \, \phi_i (x) \\&= g(x) + 0 \, \phi_i (x) \\&= g (x).
\end{align*}
Hence, by choosing $g (x) = f (x)$, the \ce\ becomes $y (x, f (x)) = f (x)$. Therefore, for any function satisfying the constraints, $f(x)$, there exists at least one free function $g (x) = f (x)$, such that the constrained expression is equal to the function satisfying the constraints, i.e., $y (x, f (x)) = f (x)$.
\end{theorem}

\begin{theorem}{}\label{thrm:ProjUni}
The TFC univariate \ce\ is a projection functional.

\noindent\hrulefill

\noindent \textbf{Proof:}
To prove Claim \ref{thrm:ProjUni}, one must show that $y (x, y (x, g (x))) = y (x, g (x))$. By definition, the constrained expression returns a function that satisfies the constraints. In other words, for any $g (x)$, $y (x, g (x))$ is a function that satisfies the constraints. From Claim \ref{thrm:UniCe}, if the free function used in the \ce\ satisfies the constraints, then the \ce\ returns that free function exactly. Hence, if the \ce\ functional is given itself as the free function, it will simply return itself.
\end{theorem}

\begin{theorem}{}\label{thrm:NonUniG}
For a given function, $f (x)$, satisfying the constraints, the free function, $g (x)$, in the TFC \ce\ $y(x,g(x)) = f(x)$ is not unique. In other words, the TFC \ce\ is a surjective functional. 

\noindent\hrulefill

\noindent \textbf{Proof:}
Consider the free function choice $g (x) = f (x) + \beta_j \,  s_j (x)$ where $\beta_j$ are scalar values on $\R$ and $s_j (x)$ are the support functions used to construct the switching functions $\phi_i (x)$.
\begin{equation*}
    y (x,g(x)) = g (x) + \phi_i (x) \, \rho_i (x, g (x)).
\end{equation*}
Substituting the chosen $g (x)$ yields,
\begin{equation*}
    y (x,g(x)) = f (x) + \beta_j \, s_j (x) + \phi_i (x) \, \rho_i (x, f (x) + \beta_j \, s_j (x)).
\end{equation*}
Now, according to Definition \ref{def:projection_function} of the projection functional,
\begin{equation*}
    y (x,g(x)) = f (x) + \beta_j \, s_j(x) + \phi_i (x) \Big(\kappa_i - \C{i} [f(x) + \beta_j \, s_j(x)]\Big).
\end{equation*}
Since the constraint operator $\C{i}$ is a linear operator,
\begin{equation*}
    y(x,g(x)) = f(x) + \beta_j s_j(x) +  \phi_i(x)\Big(\kappa_i - \C{i}[f(x)] -  \C{i}[s_j(x)]\beta_j\Big).
\end{equation*}
Since $f (x)$ is defined as a function satisfying the constraints, then $\C{i} [f(x)] = \kappa_i$, and,
\begin{equation*}
    y(x,g(x)) = f(x) + \beta_j s_j(x) -  \phi_i(x)\C{i}[s_j(x)]\beta_j.
\end{equation*}
Now, according to Property \ref{prop:co_on_s} of the constraint operator, and by decomposing the switching functions $\phi_i(x)$,
\begin{equation*}
    y(x,g(x)) = f(x) + \beta_j \, s_j(x) -  \alpha_{ki} \, s_k(x) \mathbb{S}_{ij} \, \beta_j.
\end{equation*}
Collecting terms results in,
\begin{equation*}
    y(x,g(x)) = f(x) + \beta_j \Big(\delta_{jk} - \alpha_{ki} \, \mathbb{S}_{ij}\Big) s_k(x).
\end{equation*}
However, $\mathbb{S}_{ki} \alpha_{ij} = \delta_{kj}$ because $\alpha_{ij}$ is the inverse of $\mathbb{S}_{ki}$. Therefore, by the definition of inverse,
$\mathbb{S}_{ki} \alpha_{ij} =  \alpha_{ki} \mathbb{S}_{ij}   = \delta_{kj}$, and thus,
\begin{equation*}
    y(x,g(x)) = f(x) + \beta_j \Big(\delta_{jk} - \delta_{jk} \Big) s_k(x).
\end{equation*}
Simplifying yields the result,
\begin{equation*}
    y(x,g(x)) = f(x),
\end{equation*}
which is independent of the $\beta_js_j(x)$ terms in the free function. Therefore, the free function is not unique.
\end{theorem}

Notice that the non-uniqueness of $g(x)$ depends on the support functions used in the \ce, which has an immediate consequence when using \ces\ in optimization. If any terms in $g (x)$ are linearly dependent to the support functions used to construct the constrained expression, their contribution is negated and thus arbitrary. For some optimization techniques, it is critical that the linearly dependent terms that do not contribute to the final solution be removed; else, the optimization technique becomes impaired. For example, prior research focused on using this method to solve ODEs \cite{LDE,NDE} through a basis expansion of $g(x)$ and least-squares, and the basis terms linearly dependent to the support functions had to be omitted from $g(x)$ to maintain full rank matrices in the least-squares.

The previous proofs coupled with the functional definitions and properties given earlier provide a more rigorous definition for the TFC \ce: the TFC \ce\ is a surjective, projection functional whose domain is the space of all real-valued functions that are defined at the constraints and whose codomain is the space of all real-valued functions that satisfy the constraints. It is surjective because it spans the space of all functions that satisfy the constraints, its codomain, based on Claim \ref{thrm:UniCe}, but is not injective because Claim \ref{thrm:NonUniG} shows that functions in the codomain are the image of more than one function in the domain: the functional is thus not bijective either because it is not injective. Moreover, the TFC \ce\ is a projection functional, as shown in Claim \ref{thrm:ProjUni}. 

This formal definition of the univariate TFC is simple yet powerful, as its claims apply to any combination of the constraints introduced previously, and it can easily be extended to $n$-dimensions; The multivariate TFC is the topic of Carl Leake's dissertation \cite{KarlDissertation} and was first introduced in Leake, Johnston, and Mortari \cite{M-TFC-new}

\begin{mybox}
\begin{center}
\begin{huge}
Part 2\\
\vspace{0.2in}
Application
\end{huge}

\vspace{0.5in}

Can you truly appreciate how special or\\
beautiful something is if you don't know\\
what it took to get it? If you never had\\
to work for it?\\
\vspace{0.1in}
--- Unravel, \emph{ColdWood Interactive}
\end{center}
\end{mybox}
\pagebreak{}
%

\chapter{APPLICATION TO THE SOLUTION OF ORDINARY DIFFERENTIAL EQUATIONS\label{chap:ode}}

In the prior sections, we developed a technique to derive functionals, called constrained expressions, which represented all possible functions satisfying a given set of constraints. One of the obvious applications of these expressions is to the solution of differential equations. In general, differential equations (DEs) are used as numerical models to describe physical phenomena throughout engineering and science. The solution of these equations is vital for design, predictive modeling, and optimization, and therefore, fast and accurate solutions are vital.

In the following section, the process to solve these equations using the TFC framework is introduced and used to solve various differential equations of varying complexity. Furthermore, while this work focuses explicitly on the solution of ordinary differential equations, the technique is easily extended to partial differential equations and was first covered in detail in Leake, Johnston, and Mortari \cite{M-TFC-new} and Schiassi et al. \cite{XTFC}. Again, for a complete development of multivariate TFC and the solution of partial differential equations, the reader is directed to the dissertation of Carl Leake \cite{KarlDissertation}.

Moving forward, we must first understand the two main approaches used to solve these types of problems. First, due to the structure of some types of problems, a differential equation can sometimes be solved analytically, and thereby, admit a closed-form solution. However, in most practical applications, the differential equations to be solved are complex, and numerical techniques become important when a solution, albeit approximated, is needed.

\section{Analytical methods to solve ODEs}
As mentioned above, some differential equations can be solved analytically to provide a closed-form solution to the equations. This solution is exact and suffers no associated error; however, these solutions are limited to a class of differential equations and do not encompass all differential equations. For example, for first-order differential equations, analytical techniques exist for the solutions of classes such as directly integrable, linear, separable, homogeneous, exact, and Bernoulli, etc. In fact, resources, including References \cite{analyical_ODE_1,analytical_ODE_2}, provide an extensive list of closed-form solutions to many classes of ordinary differential equations. However, the advancement and widespread use of computers has increased the emphasis on research towards solving these equations numerically. Additionally, since many numerical models are associated with complex differential equations, numerical solutions are sometimes the only available avenue to solve the problem.

\section{Numerical methods to solve ODEs}
The techniques to solve (or approximate) DEs are littered throughout literature, spanning almost all science, engineering, and mathematics fields. To understand how the TFC based method fits into the existing literature, let us look into the most popular numerical methods to solve ODEs, summarized in the following sections.

\subsection{Runge-Kutta family}
Some of the most widely used techniques are based on the Runge-Kutta family of integrators. Examples of these integrators include lower-order methods such as the Euler Method (first-order), Midpoint Method (second-order), and the Runge-Kutta Method (fourth-order) \cite{RK}. To highlight the general idea of these approaches, let us look at an example of solving the ODE, $y_x = f(x,y)$ subject to $y(x_0) = y_0$.

\begin{blankBox}{Low order Runge-Kutta methods}
\noindent Methods based on the Runge-Kutta method are forward-propagation schemes that, in general, rely on estimating the next value of the solution (i.e., the $k+1$ value) by an approximation involving the evaluation of the function $f(x,y)$ and some step size. The specific propagation equations for Euler, Midpoint, and Runge-Kutta methods are provided below:
\\\noindent\textbf{Euler Method}
\begin{equation*}
    y_{k+1} = y_k + h k_1 + \mathcal{O}(h^2)
\end{equation*}

\noindent\textbf{Midpoint Method}
\begin{equation*}
    y_{k+1} = y_k + h k_2 + \mathcal{O}(h^3)
\end{equation*}

\noindent\textbf{Runge-Kutta Method (RK4)}
\begin{equation*}
     y_{k+1} = y_k + \frac{1}{6} h\Big( k_1 + 2k_2 + 2k_3 + k_4\Big) + \mathcal{O}(h^5)  
\end{equation*}
where $k$ is the current time step, $k+1$ is the next time step, and $h$ is the step size. Additionally, $\mathcal{O}$ signifies the truncation order and is omitted in the numerical solution. In these equations, the values of $k_1$, $k_2$, $k_3$, and $k_4$ are intermediate calculations based on the order of the method, and are as follows,
\begin{align*}
    k_1 &= f(x_k,y_k)\\
    k_2 &= f(x_k + \frac{h}{2},y_k + \frac{h}{2}k_1)\\
    k_3 &= f(x_k + \frac{h}{2},y_k + \frac{h}{2}k_2)\\
    k_4 &= f(x_k + h,y_k + hk_3).
\end{align*}
\end{blankBox}
A typical approach to solving differential equations using the Runge-Kutta method is the RK45 technique, which combines an RK4 and RK5 method to adaptively select the step size $h$. This technique, called the Runge-Kutta-Fehlberg method, compares the difference between the value obtained from the 4th order and 5th order method to determine the optimal step size $h_{\text{best}}$. A summary of the RK45 algorithm is summarized below.
\begin{blankBox}{Runge-Kutta-Fehlberg method}
\noindent\textbf{4th order Runge-Kutta method}

\begin{equation*}
    \p{4}{y}_{k+1} = y_k + \frac{25}{216} h k_1 + \frac{1408}{2565} h k_3 + \frac{2197}{4104} h k_4 - \frac{1}{5} h k_5
\end{equation*}

\noindent\textbf{5th order Runge-Kutta method}

\begin{equation*}
    \p{5}{y}_{k+1} = y_k + \frac{16}{135} h k_1 + \frac{6656}{12825} h k_3 + \frac{28561}{56430} h k_4 - \frac{9}{50} h k_5 + \frac{2}{55} h k_6
\end{equation*}
where the pre-superscripts denote the order of the method, and the coefficients are as follows,
\begin{align*}
    k_1 &= f(x_k,y_k)\\
    k_2 &= f\Big(x_k + \frac{h}{4},y_k + \frac{h}{4}k_1\Big)\\
    k_3 &= f\Big(x_k + \frac{3h}{8},y_k + \frac{3}{32}h k_1 + \frac{9}{32} h k_2\Big)\\
    k_4 &= f\Big(x_k + \frac{12}{13}h, y_k + \frac{1932}{2197}hk_1 - \frac{7200}{2197}hk_2 + \frac{7296}{2197}hk_3\Big)\\
    k_5 &= f\Big(x_k + h, y_k + \frac{439}{216}hk_1 - 8hk_2 + \frac{3680}{513}hk_3 - \frac{845}{4104}hk_4\Big)\\
    k_6 &= f\Big(x_k + \frac{h}{2}, y_k - \frac{8}{27}hk_1 + 2hk_2 - \frac{3544}{2565}hk_3 + \frac{1859}{4104}hk_4 - \frac{11}{40}hk_5\Big).
\end{align*}
The optimal step size is then defined by
\begin{equation*}
    h_{\text{opt}} = h_{\text{last}} \Bigg(\dfrac{\varepsilon h_{\text{last}} }{2\Large|\p{4}{y}_{k+1} - \p{5}{y}_{k+1}\Large|}\Bigg)^{1/4} \approx \frac{\text{desired error}}{\text{actual error}}.
\end{equation*}
\end{blankBox}
The above technique is similar to what is implemented in algorithms such as MATLAB's \verb"ode45()" \cite{ode45} and the Python package SciPy's \verb"scipy.integrate.ode()" \cite{scipy}. In many numerical tests in this chapter, we will use the RK45 solution as the baseline to compare against the TFC method in terms of speed and accuracy.

\subsection{Gauss-Jackson}
Another technique widely used in the astrodynamics community is the Gauss-Jackson method, a multistep predictor-corrector method. First introduced in a 1924 paper by Jackson \cite{GJ1}, this technique has been further studied in References \cite{GJ2,GJ3}. In general, this method is a summed form of the Stormer-Cowell integrator \cite{SC_integration}.

In order to understand the fundamentals of this method, consider the ordinary differential equation of the form $y_{xx} = f(x,y,y_x)$. The Gauss-Jackson technique first predicts the solution value $y(x)$ for the next step and evaluates the function $f(x,y,y_x)$ at this point. Then, this predicted function value is added to the backpoints, i.e., prior calculated points. A corrector formula is utilized to revise this set of data and refine the prediction of $y(x)$. The general implementation of these algorithms can be grouped into two methods, 1) Predict-Evaluate-Correct (PEC) and 2) Predict-Evaluate-Correct-Evaluate (PECE), where the latter performs a second evaluation step to increase accuracy. Furthermore, these processes can perform additional iterations to meet some tolerance. 

The following example box provides a summary of the major equations in the Gauss-Jackson method.
\begin{blankBox}{Gauss-Jackson method}
Consider the ordinary differential equation where $y_{xx} = f(x,y)$, subject to the initial conditions $y(x_0) = y_0$ and $y_x(x_0) = y_{0_x}$. The Gauss-Jackson correction and prediction formulas are as follows where $H.O.T$ stands for higher order terms.

\noindent\textbf{Gauss-Jackson corrector formula}
\begin{equation*}
    y_k = h^2\Big[\nabla^{-2}y_{k_{xx}} + \Big(\frac{1}{12} - \frac{1}{240}\nabla^2 - \frac{1}{240}\nabla^3 - \frac{221}{60480}\nabla^4 + \hdots + H.O.T.\Big) y_{k_{xx}} \Big]
\end{equation*}

\noindent\textbf{Gauss-Jackson predictor formula}
\begin{equation*}
    y_{k+1} = h^2\Big[\nabla^{-2}y_{k_{xx}} + \Big(\frac{1}{12} + \frac{1}{12}\nabla + \frac{19}{240}\nabla^2 + \frac{3}{40}\nabla^3 + \hdots + H.O.T.\Big)y_{k_{xx}}\Big]
\end{equation*}
where $\nabla$ is the backwards difference operator such that $\nabla f_k = f_k - f_{k-1}$. The higher-order difference operators can be easily derived and are provided in Reference \cite{GJ3}. Additionally, the predication and correction of the first derivative, $y_x$, is given by the summed Adams method as,
\begin{align*}
    y_{k_{x}} &= h\Big[\nabla^{-1} - \frac{1}{2} - \frac{1}{12}\nabla - \frac{1}{24}\nabla^2 - \frac{19}{720}\nabla^3 - \frac{3}{160}\nabla^4 \hdots + H.O.T.\Big) y_{k_{xx}} \Big]\\
     y_{k+1_{x}} &= h\Big[\nabla^{-1} + \frac{1}{2} + \frac{5}{12}\nabla + \frac{3}{8}\nabla^2 + \frac{251}{720}\nabla^3 + \frac{95}{288}\nabla^4 + \hdots + H.O.T.\Big)y_{k_{xx}}\Big].
\end{align*}
\end{blankBox}
When solving differential equations using the Gauss-Jackson method (and other predictor -corrector methods), the main hurdle is initialization. Since the initial conditions are given at some epoch $x_0$, there are no backpoints, and these must be calculated before the algorithm is used. One way to initialize these backpoints is to use a single-step integrator such as the Runge-Kutta methods described in the prior section.

\subsection{Modified Chebyshev-Picard Iteration}
Modified Chebyshev-Picard Iteration \cite{MCPI1,MCPI2,MCPI3} is a path-length integral approximation that has been recently proven to be highly effective. This technique has been successfully applied to initial- and boundary-value problems in orbit propagation. The following summarizes the main parts of the method.
\begin{blankBox}{Modified Chebyshev-Picard Iteration}
Given a differential equation \begin{equation*}
    y_x = f(x,y)
\end{equation*}
where $y(x_0) = y_0$, the domain is first transformed to that of the closed interval of the Chebyshev polynomials [-1, +1],
\begin{equation*}
    x = \omega_1 + \omega z \quad \omega_1 = \frac{x_f + x_0}{2} \quad \omega_2 = \frac{x_f - x_0}{2}.
\end{equation*}
This transformation allows us to rewrite the differential equation as,
\begin{equation*}
    y_z = q(z,y) = \omega_2 f(\omega_1 + \omega_2 z, y).
\end{equation*}
The solution to this equation is provided by Picard iteration where,
\begin{equation*}
    y^i(z) = y_0 + \int_{-1}^{z} q(s,y^{i-1}(s))\dd s \quad i = 1, 2, \hdots 
\end{equation*}
Next, the state $y^i$ and the integrand are approximated by a sum of Chebyshev polynomials with unknown coefficients discretized at $(N+1)$ Chebyshev-Gauss-Lobatto (CGL)
nodes such that,
\begin{equation*}
    z_j = \cos\Big(\frac{j\pi}{N}\Big) \quad j = 0,1,2\hdots,N
\end{equation*}
The forcing function is approximated by Chebyshev polynomials through,
\begin{align*}
    q(z,y^{i-1}(z)) &\approx \sum_{k=0}^{k=N} \phantom{}' F_k^{i-1} T_k(z) \\
    &\equiv \frac{1}{2}F_0^{i-1}T_0(z) + F_1^{i-1}T_1(z) + F_2^{i-1}T_2(z) + \hdots F_N^{i-1}T_N(z).
\end{align*}
The discrete orthogonality of Chebyshev polynomials \cite{MCPI_Cheb} allows for the direct computation of $F_k$,
\begin{align*}
    F_k^{i-1} &= \frac{2}{N} \sum_{j=0}^{N} \phantom{}'' q(z_j,y^{i-1}(z_j))T_k(z_j) \\
    &= \frac{1}{N} q(z_0,y^{i-1}(z_0)) T_k(z_0) + \phantom{\hdots} \\ & \hspace{15.5pt} \frac{2}{N} q(z_1,y^{i-1}(z_1)) T_k(z_1) + \hdots + \frac{1}{N} q(z_N,y^{i-1}(z_N)) T_k(z_N)
\end{align*}
where $\sum \phantom{}'$ denotes that the first term is halved and $\sum \phantom{}''$ represents
that both the first and last terms are halved. Plugging this into the Picard iteration equation leads to,
\begin{equation*}
    y^{i} = y_0 + \sum_{r=0}^{N} \phantom{}' F_r^{i-1} \int_{-1}^{z}T_r(s)\dd s \equiv \sum_{k=0}^{N} \phantom{}' \beta^{i}_k T_k(z)
\end{equation*}
where the updated equations for the coefficients are derived in detail in Reference \cite{Bai_MCPI} and summarized below,
\begin{align*}
    \beta_k^i &= \frac{1}{2k} (F_{k-1}^{i-1} - F_{k+1}^{i-1}) \quad k = 1,2,\hdots,N-1 \\
    \beta_N^i &= \frac{F_{N-1}^{i-1}}{2N} \\
    \beta_0^i &= 2y_0 + 2 \sum_{k=1}^{k=N}(-1)^{k+1}\beta_k^i
\end{align*}
\end{blankBox}

\subsection{Collocation and Spectral Methods}
The previously mentioned methods are based on low-order Taylor expansions, which limit the step size that can be used to propagate the solution. Additionally, a common weakness of all methods based on low-order Taylor expansion is that they are not effective in enforcing algebraic constraints. Therefore, recent research has looked for other numerical schemes. 

\subsubsection{Collocation methods}
 One of these numerical schemes is the collocation method \cite{CM1,CM2,Collocation}. In this method, the solution components are approximated by piecewise polynomials on a mesh. The mesh is made up of a number of points in the domain (called collocation points), and the problem is solved by minimizing the residual of the differential equation at the collocation points. In general, this reduces to computing the unknown coefficients of the polynomial functions. The approximation to the solution must satisfy the constraint conditions and the differential equation at the collocation points in each mesh subinterval. In the collocation methods, the placement of the collocation points is not arbitrary. A modified Newton-type method, known as quasi-linearization, is then used to solve the nonlinear equations for the polynomial coefficients. The mesh is then refined by equally distributing the estimated error over the whole interval, and therefore, an initial estimation of the solution across the mesh is required. In general, this method numerically approximates the differential equation and the specified constraints.

\subsubsection{Spectral methods}
On the other hand, spectral methods \cite{Spectral} model the differential equation's solution by a sum of ``basis functions'' with unknown coefficients that are solved according to the specific differential equation. The differential equation is then approximated by 1) discretizing the domain and 2) solving the resulting algebraic equations of the differential equation and specified constrained at these nodes. In general, this method benefits from being less computationally expensive than approaches like collocation methods, but it suffers from accuracy problems when applied to complex geometries such as discontinuities. Furthermore, spectral methods are the most similar to the TFC approach since they both are an ``assumed'' solution method. In both techniques, we assume the form of the solution (i.e., Chebyshev orthogonal polynomials) and solve for unknown coefficients that minimize the residual of the differential equation. The key difference between spectral methods and the TFC method is in spectral methods, the constraints have to be introduced into the numerical scheme and therefore have associated error, whereas, in the TFC method, the constraints are satisfied analytically via the constrained expression.

\subsection{Machine Learning}
With the current boom in machine learning and artificial intelligence spurred by the increasing capabilities of computers, researchers have looked to apply these algorithms to the numerical solution of differential equations. This method is similar to the spectral method; however, the ``basis functions'' are replaced with neural networks (NNs) and paired with a multitude of optimization algorithms to solve the problem. In fact, various authors have explored the feasibility of using Neural Networks (NNs) to solve ODEs and PDEs. 

The basis of this work leverages two main ideas. First, the \textit{Universal Approximation Theorem} \cite{NN_univApproximator_cybenko1989approximation,hornik1991approximation}, which states that NNs are universal approximators, and therefore, can potentially represent the function that is the solution of a given differential equation \cite{NN_univApproximator_cybenko1989approximation,elm_UnivApproximator} as the number of neurons go to infinity. Using these ideas, in 1995, Chen and Chen \cite{chen1995universal} were able to show that NNs could approximate nonlinear operators. Furthering this work, Pinkus \cite{approxtheory} and Lu et al. \cite{lu2019deepxde} detailed a function and its partial derivatives that could simultaneously and uniformly be approximated with a single layer NN with a sufficiently large number of hidden neurons.

Of importance to the topic of this dissertation, for ODEs, multiple NN-based solutions have been proposed, including Yang et al. \cite{yang} Legendre Neural Networks (LeNNs), Sun et al. \cite{BNN} Bernstein Neural Network (BNNs), and Mall and Chakraverty \cite{CNN} Chebyshev Neural Network (CNNs). All of these techniques use single-layer NNs where the activation functions are Legendre, Bernstein, or Chebyshev polynomials, respectively. The network is trained via the Extreme Learning Machine (ELM) algorithm, proposed by Huang et al. \cite{ELM}\footnote{The author notes that the method of Legendre, Bernstein, or Chebyshev Neural Networks paired with the ELM algorithm is exactly the method defined by the spectral method by simply using Legendre, Bernstein, or Chebyshev polynomials.}. The ELM algorithm is used for single-hidden layer feed-forward networks where the hidden input weights and biases are randomly selected, and the output weights are solved via least-squares. To satisfy the problem constraint, a constraint penalty is added to the loss function minimized during the training phase.

\section{The TFC method to solve ODEs}
\begin{figure}[ht]
    \centering\includegraphics[width=.65\linewidth]{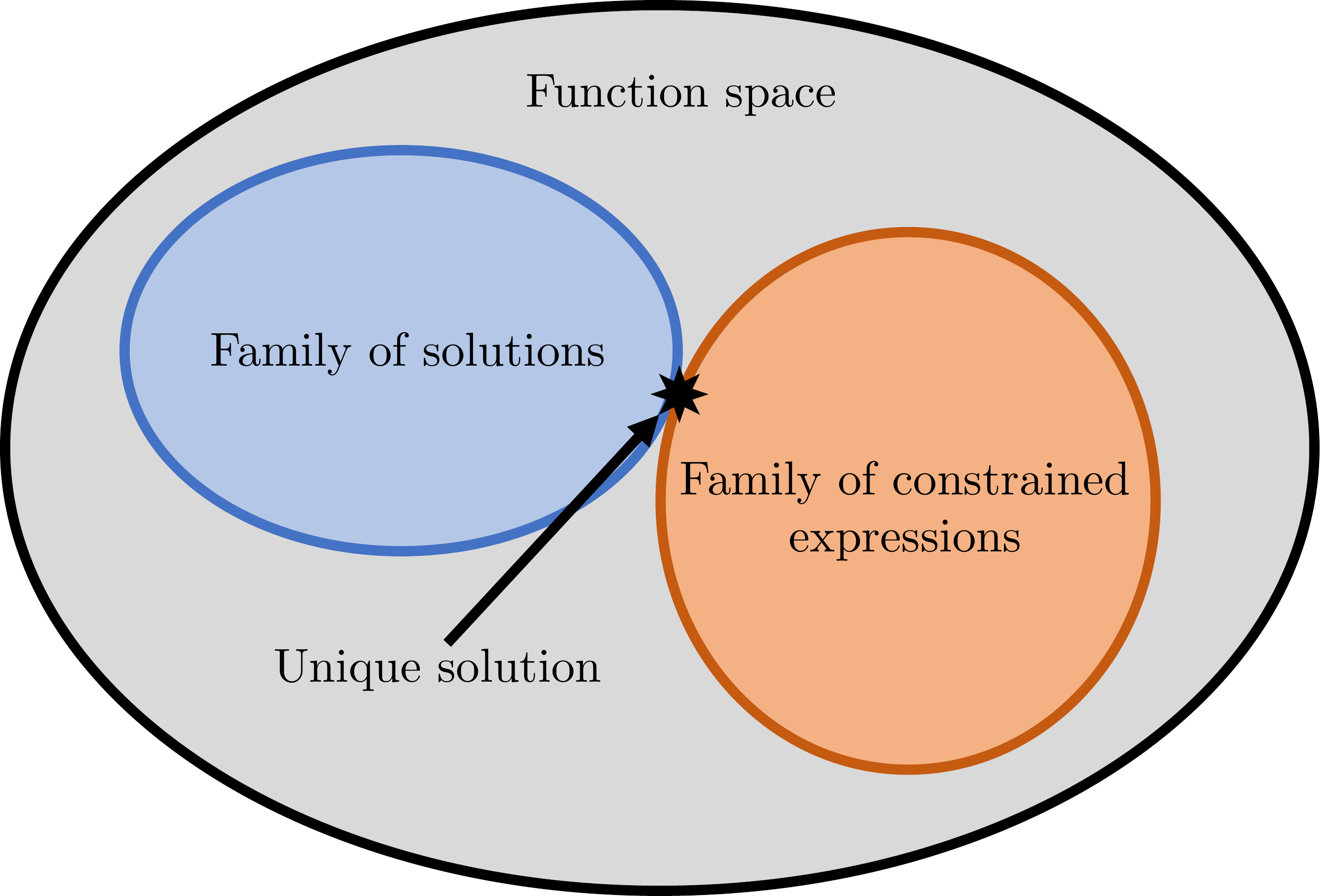}
    \caption{Diagram of function space associated with the solution of a ordinary differential equation. Note: this figure is used for conceptual purposes and is not a rigorous mathematical description. For example, in the solution of some differential equations, there could be more than one, or even infinite intersection points depending on the nature of the differential equation.}
    \label{fig:s4_venn}
\end{figure}

As we will soon see, the TFC method shares a similar approach to the collocation method, spectral method, and ELMs. However, the distinction is that the constraints are embedded analytically before the numerical approximation step. In summary, this will provide us with two unique advantages, 1) the constraints are always satisfied analytically, and 2) the loss functions only deal with the differential equation to be solved. 
In general, the TFC method is planted between the two general methods (analytical and numerical) to solve differential equations. This can be easily visualized in the diagram of Figure \ref{fig:s4_venn}.
In the prior section, we discussed the solution of DEs through analytical techniques. The analytical method is represented by the blue oval, where a family of solutions is provided. The unique solution (the black star) is then determined by applying the constraints to the differential equation. On the other hand, numerical solutions (excluding IVPs) must search the entire function space to find a unique solution. Conversely, the TFC method solves the problem in the opposite sequence of the analytical approach. First, the candidate solution is constructed by using a constrained expression. The constrained expression represents a reduction of the function space to a set only the functions satisfying the DE's constraints. Then, the codomain of the constrained expression is used to find the unique solution of the differential equation. In another sense, if we assume that our free function, $g(x)$, covers the function space of the solution, then the constrained expression is projecting this function into a reduced set of the constraints, i.e., the orange oval. It should be clear from this discussion that the solution of the differential equation is dependent on the definition of $g(x)$.

To further understand these concepts, let us consider a general differential equation,
\begin{equation}\label{eq:s3_genDE}
F\Big(x,y,\frac{\dd y}{\dd x},\frac{\dd^2 y}{\dd x^2}, \cdots, \frac{\dd^n y}{\dd x^n}\Big) = 0,
\end{equation}
subject to $n$ linear constraints. Using the TFC framework, the first step is to derive the switching and projection functions of Equation \eqref{eq:s2a_uniCeAlt}. By doing this, the constraints of Equation \eqref{eq:s3_genDE} are decoupled from the solution of the differential equation, and the differential equation is transformed into,
\begin{equation}\label{eq:s3_genDE_g}
\tilde{F}\Big(x,g,\frac{\dd g}{\dd x},\frac{\dd^2 g}{\dd x^2}, \cdots, \frac{\dd^n g}{\dd x^n}\Big) = 0,
\end{equation}
where the solution to this ``differential equation''\footnote{The use of quotations around the word differential equation is used because the resulting expression is: 1) technically not a differential equation and 2) cannot be solved using the analytical techniques to solve differential equations. To date, this type of equation has not been rigorously defined.} is obtained by finding the function $g(x)$ satisfying Equation \eqref{eq:s3_genDE_g}. In order to solve this new equation, four major steps must be taken: 1) define the free function $g(x)$ and 2) determine the derivatives of the free function $g(x)$ 3) discretize the domain, and 4) solve the resulting algebraic equation. The following sections elaborate on these steps.

\subsection{Defining the free function}\label{sec:s3_definingG}
For our definition of the free function, we will allow the domain of this function, $z$, to be different from the differential equation problem domain $x$. Ultimately, we will need to map between the domains with some function $z = z(x)$; however, allowing for different basis and problem domains is necessary in most cases since some numerical bases are defined on closed domains, e.g., Chebyshev orthogonal polynomials are defined on $z \in [-1,+1]$. This will be made clear in Section \ref{sec:s3_deriv_free}.

In selecting a free function, we are essentially looking for the best (differentiable) function approximator. A simple definition of $g(x)$ could be the monomial expansion of $m$ terms,
\begin{equation}\label{eq:s3_mono}
    g(x) = \sum_{k = 0}^{m-1} a_k \, z^k,
\end{equation}
where $a_k$ are coefficients and $z$ is simply the independent variable. According to Claim \ref{thrm:NonUniG}, the terms linearly dependent to the support functions used in the \ce\ must be removed. While this definition is valid, a linear combination of orthogonal polynomials can be leveraged for their advantageous numerical properties. 

Consider the definition of Chebyshev polynomials of the first kind,
\begin{equation}\label{eq:s3_CP}
    g(x) = \sum_{k = 0}^{m-1} a_k \, T_k(z),
\end{equation}
where again $a_k$ are coefficients and $T_k(z)$ are the Chebyshev polynomials terms. Again, Claim \ref{thrm:NonUniG} must be considered in this expansion. It has been shown that Chebyshev polynomials of the first kind produce a function that minimizes the maximum error in its application. In fact, these polynomials are part of a special class well suited for function approximation \cite{NM_OP}. Furthermore, this expansion also provides a simple way to estimate the solution's accuracy by observing the size of the coefficients of latter terms in the expansion (i.e., the coefficients of the highest-order terms), which is justified by the convergence properties of Chebyshev polynomials. An even better approximation is obtained by comparing the sets of coefficients obtained when the number of basis terms is varied \cite{AA}.

Additionally, the Legendre orthogonal polynomials, defined as,
\begin{equation}\label{eq:s3_LeP}
    g(x) = \sum_{k = 0}^{m-1} a_k \, L_k(z),
\end{equation}
where $a_k$ are coefficients and $L_k$ are polynomial terms, are another important expansion, which has been used extensively in function approximation and the solution of differential equations with beneficial error properties for the approximation of smooth functions \cite{LePDiss}. In fact, both orthogonal polynomials types mentioned have been extensively used in spectral methods \cite{Spectral}.

Moreover, our definition of $g(x)$ can even extend to machine learning where the function is defined as a neural network where we would express
\begin{equation*}
    g(x) = \mathcal{N}(z,\theta),
\end{equation*}
where the architecture is based on the independent variable $z$ and trainable parameters $\theta$, such as the weights and the biases. A complete study of the use of neural networks is out of the scope of this work, and interested readers are directed to Leake and Mortari \cite{deepToC} for a more detailed look into applying TFC in this field.

In addition to the general neural networks, one specific architecture has shown promising results which is based on the theory of the ELM \cite{ELM}. ELMs are a single-layer feed-forward NN where in the univariate definition,
\begin{equation}\label{eq:s3_ELM}
    g(x) = \sum_{k=0}^{m-1} a_k \, \sigma \left(w_k z + b_k \right).
\end{equation}
In this equation, $m$ is the number of hidden neurons, i.e., similar to the number of basis functions, and $\sigma$ is a user-defined activation function, e.g., sigmoid, tanh, swish, etc. The terms $w_k$ and $b_k$ are the associated weights and biases for the nodes and are selected randomly according to any continuous probability distribution proven in Theorems 2.1 and 2.2 in G.-B. Huang et al. \cite{ELM}. Therefore, it makes the unknown coefficients, $a_k$, linear in the form of Equation \eqref{eq:s3_ELM} similar to Equations \eqref{eq:s3_mono}, \eqref{eq:s3_CP}, and \eqref{eq:s3_LeP}.

Moving forward we will only consider the free function defined in terms of the Chebyshev polynomials Equation \eqref{eq:s3_CP}, Legendre polynomials Equation \eqref{eq:s3_LeP}, and ELMs, Equation \eqref{eq:s3_ELM}. Since all functions are linear in their unknown coefficients, $a_k$, let us write the general expansion as,
\begin{blankBox}{General Basis Expansion}
\begin{equation}\label{eq:s3_basis}
    g (x) = \B{\xi} \T\B{h}(z) \quad \text{where} \quad z = z(x)
\end{equation}
where  $\B{\xi} = \{ a_0, \cdots, a_k, \cdots, a_{m-1} \}\T$ and $\B{h}(z)$ is a vector function of the $m$ functions. 
\end{blankBox}
\subsection{Derivatives of the free function}\label{sec:s3_deriv_free}
In most cases the domain of the free function will not coincide with the domain of the problem. For example, for the orthogonal polynomials mentioned, the domain is defined for $z \in [-1,+1]$ and most of the time it is desirable to scale the input which may be different than our problem domain, $x \in [x_0,x_f]$. Therefore, these functions must be linearly mapped to the independent variable $x$. This can be done using the equations, 
\begin{equation}\label{eq:s3_x2z}
z = z_0 + \frac{z_f-z_0}{x_f-x_0}(x - x_0) \quad \longleftrightarrow \quad x = x_0 + \frac{x_f-x_0}{z_f-z_0}(z - z_0),
\end{equation}
where $x_f$ represents the upper integration limit. The subsequent derivatives of the free function are defined as,
\begin{equation*}
    \frac{\dd^{n} g}{\dd x^{n}} = \B{\xi} \T  \frac{\dd^{n} \B{h}(z)}{\dd z^{n}} \left(\frac{\dd z}{\dd x}\right)^{n},
\end{equation*}
where by defining,
\begin{equation}\label{eq:s3_linearMapping}
c := \frac{\dd z}{\dd x} = \frac{z_f - z_0}{x_f - x_0}
\end{equation}
the expression can be simplified to, 
\begin{blankBox}{Derivatives of the free function}
\begin{equation*}
    \frac{\dd^{n} g}{\dd x^{n}} = c^{n} \B{\xi} \T  \frac{\dd^{n} \B{h}(z)}{\dd z^{n}},
\end{equation*}
\end{blankBox}
which defines all mappings of the free function. By defining the free function according to the form of Equation \eqref{eq:s3_basis}, our transformed differential equation, Equation \eqref{eq:s3_genDE_g}, that was derived earlier reduces to,
\begin{equation}\label{eq:s3_genDE_xi}
    \tilde{F}(x,\B{\xi}) = 0.
\end{equation}
Next, the problem domain, $x$, must be discretized to eventually solve for the unknown coefficients and ultimately solve the differential equation. Therefore, a specific discretized scheme is needed.

\subsection{Discretization of the domain}
Since the ultimate goal is to solve Equation \eqref{eq:s3_genDE} computationally, the problem domain (and therefore the basis function domain) must be discretized. In the case of defining $g(x)$ using an ELM, the discretization can simply be selected as uniformly spaced points. However, when using Chebyshev and Legendre orthogonal polynomials, the discretization scheme is slightly more involved. For these polynomials, the optimal discretization scheme is Chebyshev-Gauss-Lobatto nodes \cite{Colloc,ChebCol}. For $N+1$ points, the discrete points are calculated as,
\begin{blankBox}{Discretization scheme for Chebyshev-Gauss-Lobatto nodes}
\begin{equation}\label{eq:collo}
    z_j = -\cos\left(\frac{j \pi}{N}\right) \quad \text{for} \quad j = 0, 1, 2, \cdots, N.
\end{equation}
\end{blankBox}
Compared with the uniform distribution, this distribution results in a much slower increase of the condition number of the matrix to be inverted in the least-squares as the number of basis functions, $m$, increases. The nodes can be realized in the problem domain through the relationship provided in Equation \eqref{eq:s3_linearMapping}.

By discretizing the domain according to the specific free function used, Equation \eqref{eq:s3_genDE_xi} becomes a system of equations that is linear if Equation \eqref{eq:s3_genDE} is linear and nonlinear if Equation \eqref{eq:s3_genDE} is nonlinear. This can be written as a loss vector at the discretized points,
\begin{equation}\label{eq:s3_genDE_dis}
    \mathbb{L}(\B{\xi}) = \begin{Bmatrix}\tilde{F}(x_0, \B{\xi}) \\ \vdots \\ \tilde{F}(x_k, \B{\xi}) \\ \vdots \\ \tilde{F}(x_f,\B{\xi}) \end{Bmatrix} = \B{0}
\end{equation}
where $x_k$, and therefore  $z_k$, are defined by Equation \eqref{eq:s3_x2z} and Equation \eqref{eq:collo}.
\subsection{Solving the resulting algebraic equation}
For a linear differential equation $F$ (and therefore a linear differential equation $\tilde{F}$), the \ce\ and its derivatives will show up linearly, and therefore, will remain linear in the unknown $\B{\xi}$ term. This leads to the form,
\begin{equation}\label{eq:s3_Axb}
    \mathbb{A} \B{\xi} + \B{b} = 0,
\end{equation}
where the matrix $\mathbb{A}$ is composed of a linear combination of the terms linear in the unknown coefficients. Written in terms of the loss function $\tilde{F}$, $\mathbb{A}$ is simply the Jacobian of the loss vector Equation \eqref{eq:s3_genDE_dis}, 
\begin{equation*}
    \mathbb{J}(\B{\xi}) = \frac{\partial \mathbb{L}}{\partial \B{\xi}} = \begin{bmatrix} \dfrac{\partial \tilde{F}(x_0,\B{\xi})}{\partial \B{\xi}} \\ \vdots \\ \dfrac{\partial \tilde{F}(x_k,\B{\xi})}{\partial \B{\xi}} \\ \vdots \\ \dfrac{\partial \tilde{F}(x_N,\B{\xi})}{\partial \B{\xi}} \end{bmatrix}.
\end{equation*}
Since the loss function is linear in $\B{\xi}$, it will be independent of $\B{\xi}$. Additionally, the vector $\B{b}$ is simply the loss vector evaluated at $\B{\xi} = \B{0}$,
\begin{equation*}
    \B{b} =  \mathbb{L}(\B{0}) = \begin{Bmatrix} \tilde{F}(x_0,\B{0}) \\ \vdots \\ \tilde{F}(x_k,\B{0}) \\ \vdots \\ \tilde{F}(x_f,\B{0}) \end{Bmatrix}.
\end{equation*}
Therefore, Equation \eqref{eq:s3_Axb} can also be realized as,
\begin{equation}\label{eq:s3_updateLinear}
     \mathbb{J}(\B{0}) \, \B{\xi} = - \mathbb{L}(\B{0}).
\end{equation}
In these linear cases, Equation \eqref{eq:s3_updateLinear} can be solved directly using any available least-squares technique. A summary of these numerical schemes are provided in Appendix \ref{chap:LS}. However, in the case of a nonlinear differential equations, Equation \eqref{eq:s3_genDE_dis} will be nonlinear in the $\B{\xi}$ coefficients. This system can be solved by an iterative least-squares method similar to Equation \eqref{eq:s3_updateLinear}; however, now a multivariate Newton's method is used to solve the nonlinear system for the change in the $\B{\xi}$ parameter denoted by $\Delta \B{\xi}$,
\begin{equation}\label{eq:s3_updateNonlinear}
     \mathbb{J}(\B{\xi}_i) \, \Delta \B{\xi}_i = - \mathbb{L}(\B{\xi}_i)
\end{equation}
\begin{blankBox}{Parameter update equations}
The parameter update of $\B{\xi}$ is provided by,
\begin{equation*}
    \B{\xi}_{i+1} = \B{\xi}_i + \Delta \B{\xi}_i
\end{equation*}
where the $\Delta \B{\xi}_i$ can be defined using classic least-squares,
\begin{equation*}
    \Delta \B{\xi}_i = -\Big(\mathbb{J}(\B{\xi}_i)\T \mathbb{J}(\B{\xi}_i) \Big)^{-1} \mathbb{J}(\B{\xi}_i)\T \mathbb{L}(\B{\xi}_i),
\end{equation*}
or any other least-squares technique provided in Appendix \ref{chap:LS}. This process is repeated until some stopping criteria are met. The original work on the solution of nonlinear differential equations by Mortari, Johnston, and Smith \cite{NDE} used the $L_2$ norm of the loss function and the $L_2$ norm of the least-squares step ($\Delta \B{\xi}$) such that it was below some tolerance, $\varepsilon$, according to the following equations,
\begin{equation*}
    L_2 [\mathbb{L}(\B{\xi}_i)] < \varepsilon \qquad \text{or} \qquad
    L_2 [\Delta \B{\xi}_{i}] < \varepsilon.
\end{equation*}
However, the work presented in this dissertation utilizes a slightly different stopping condition to reduces computational overhead such that,
\begin{equation*}
    \max [\mathbb{L}(\B{\xi}_i)] < \varepsilon \qquad \text{or} \qquad
    \max [\Delta \B{\xi}_{i}] < \varepsilon.
\end{equation*}
\end{blankBox}

In all, the solution of a linear versus a nonlinear ordinary differential equation is reduced simply to the difference between Equation \eqref{eq:s3_updateLinear} and Equation \eqref{eq:s3_updateNonlinear}, where the linear case only requires ``one'' iteration compared to the nonlinear equations. This similarity is highlighted in Section \ref{sec:s3_L}, where the problem is formulated according to both notations.

Additionally, since the constraints are embedded in the constrained expression before forming the loss vector, the numerical scheme does not change between boundary conditions. In other words, an initial-value problem is solved in the same manner as a boundary-value problem. We will soon see the power of this when applying TFC to the solution of boundary-value problems.

\subsection{The TFC roadmap}
\begin{figure}[ht]
    \centering\includegraphics[width=.9\linewidth]{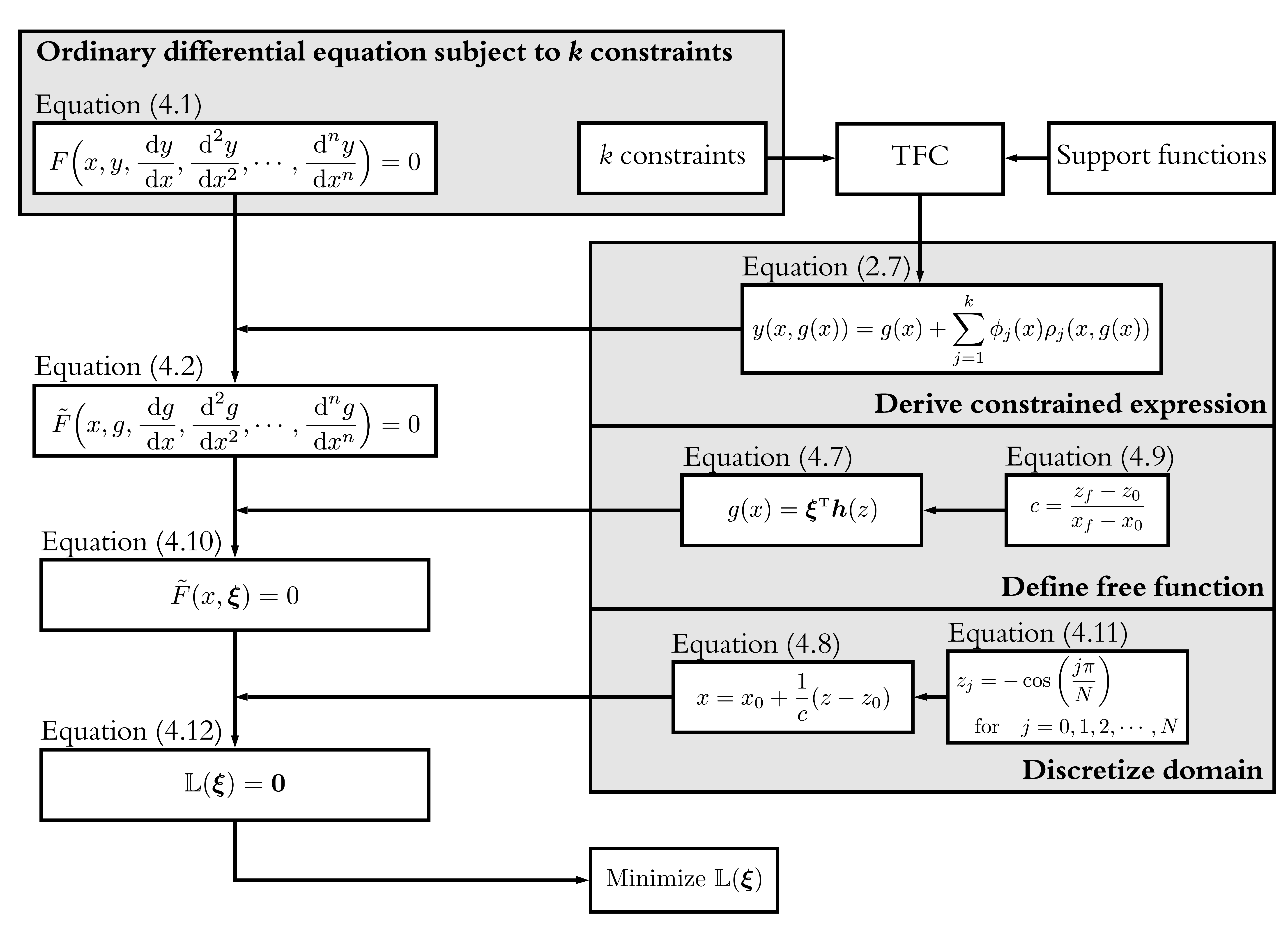}
    \caption{Flowchart of the TFC method applied to solving an ordinary differential equation in the form of Equation \eqref{eq:s3_genDE}.}
    \label{fig:s4_tfcSumaary}
\end{figure}

Before moving to our numerical examples, it is useful to summarize the entire process of solving differential equations using the TFC approach. This is provided in the flowchart in Figure \ref{fig:s4_tfcSumaary}, where the process is summarized with all major equations.
First, given the differential equation, Equation \eqref{eq:s3_genDE}, subject to $k$ constraints, we embed these constraints into the constrained expression, Equation \eqref{eq:s2a_uniCeAlt}, by selecting acceptable support functions and deriving the projection functionals and switching functions. The constrained expression and its derivative are substituted into Equation \eqref{eq:s3_genDE}, which transforms the differential equation subject to $k$ constraints to one which is unconstrained and denoted by $\tilde{F}(x,\B{\xi})$, Equation \eqref{eq:s3_genDE_g}. After this, the free function $g(x)$ is expressed by one of the many function approximation methods discussed in Section \ref{sec:s3_definingG} using Equation \eqref{eq:s3_basis}. By doing this, the differential equation is transformed into an algebraic equation with the unknown vector $\B{\xi}$. Next, we discretize the basis function domain according to Equation \eqref{eq:collo} when using Chebyshev or Legendre polynomials, and uniformly when using ELMs, and connect these to the problem domain by Equation \eqref{eq:s3_x2z}. By evaluating Equation \eqref{eq:s3_genDE_xi}, the loss function, at these discretization nodes and stacking them in a loss vector we are led to Equation \eqref{eq:s3_genDE_dis}. Finally, Equation \eqref{eq:s3_genDE_dis} is minimized using least-squares or nonlinear least-squares, depending on the linearity of the original differential equation, Equation \eqref{eq:s3_genDE}. Note, we are not limited to least-squares techniques, and in fact, any numerical minimization scheme can be used to solve the system $\mathbb{L}(\B{\xi}) = \B{0}$. With that said, the work in this dissertation focuses specifically on least-squares techniques for numerical simplicity and speed advantages. However, with the increasing complexity of problems, least-squares can become prohibitive, and the use of different optimizers is an area of future research summarized in Section \ref{sec:con_OtherOpts}. 

\section{Numerical Implementation}
To demonstrate how the TFC approach is used to solve differential equations, we will start with two simple examples covering a linear initial-value problem (Section \ref{sect:s3_laneEm}) and a nonlinear boundary-value problem (Section \ref{sect:s3_bvp}). These problems provide the full derivation and explicitly provide the Jacobian of the loss vector directly in the text for clarity. After these problems, all analytical Jacobians are not provided directly in the main text but collected in Appendix \ref{chap:app_jacob}. Following this, a brief discussion is provided on how systems of differential equations (or a subclass, vector equations) can be solved in the same manner. Lastly, Section \ref{sec:s3_extensions} discusses two adjustments to the theory to solve problems with discontinuous dynamics and unknown final times. Additionally, all numerical results were produced on a MacBook Pro (2016) macOS Version 10.15, with a 3.3 GHz Dual-Core Intel\textsuperscript{\textregistered} Core\texttrademark \, i7 and with 16 GB of RAM.

\section{Lane-Emden equation}\label{sect:s3_laneEm}
As a motivating example, let us consider the Lane-Emden equation where, 
\begin{equation}\label{eq:laneEm}
    y_{xx} + \frac{2}{x} y_{x} + y^a = 0 \quad \text{such that} \quad (x > 0, a \geq 0) \quad \text{subject to:} \quad \begin{cases} y(0) = 1 \\ y_x(0) = 0 \end{cases}
\end{equation}
For this differential equation, an exact solution exists for $a =$ 0, 1, and 5. We can see, regardless of the value of $a$, the constrained expression will be the same. Therefore, whether the equation is linear or nonlinear does not affect the derivation of the constrained expression. This should be obvious since the TFC approach decouples the problem's constraints from the solution of the differential equation. Using the theory developed earlier, the \ce\ for this problem can be solved by defining the projection functionals as,
\begin{equation*}
    \rho_1(x,g(x)) = 1 - g(0) \qquad \text{and} \qquad \rho_2(x,g(x)) = - g_x(0)
\end{equation*}
and the switching functions are determined by choosing the support functions $s_1(x) = 1$ and $s_2(x) = x$ and solving for the coefficients $\alpha_{ij}$,
\begin{align*}
    \begin{bmatrix} 1 & 0 \\ 0 & 1 \end{bmatrix} \begin{bmatrix} \alpha_{11} & \alpha_{12}\\  \alpha_{21} & \alpha_{22} \end{bmatrix} &= \begin{bmatrix} 1 & 0 \\ 0 & 1 \end{bmatrix}
\end{align*}
where it can easily be seen that $\alpha_{ij} = \delta_{ij}$. Thus, the switching functions are $\phi_1(x) = 1$ and $\phi_2(x) = x$, and the final \ce\ is,
\begin{equation}\label{eq:laneCE}
    y(x,g(x)) = g(x) + (1 - g(0)) + x(-g_x(0)).
\end{equation}
The simplicity of this expression is due in part to the second-order initial value constraints. See Appendix \ref{chap:commonSwithcingFunctions} for a summary of the associated switching functions and projection functionals for other typical constraint cases. The constrained expression, Equation \eqref{eq:laneCE}, always satisfies the constraints of Equation \eqref{eq:laneEm}.

Now, by defining $g(x)$ according to Equation \eqref{eq:s3_basis}, the constrained expression and its derivatives can be written as a linear function of the unknown coefficients,
\begin{align}
    y(x,\B{\xi}) &= \Big(\B{h} - \B{h}(z_0) - x \, c \B{h}_z(z_0)\Big)\T \B{\xi} + 1 \label{eq:s3_laneEmCEy}\\
    y_x(x,\B{\xi}) &= \Big(c \B{h}_z -  c \B{h}_z(z_0)\Big)\T \B{\xi}\label{eq:s3_laneEmCEyp} \\
    y_{xx}(x,\B{\xi}) &= \Big(c^2 \B{h}_{zz}\Big)\T \B{\xi}\label{eq:s3_laneEmCEypp}
\end{align}
In the following sections, we will use our description of the constrained expression to solve each case of the Lane-Emden equation.

\subsection{Linear differential equations}\label{sec:s3_L}

First, let us consider the solution of the linear differential equation associated with setting $a=0$ in the Lane-Emden equation,
\begin{equation*}
    y_{xx} + \frac{2}{x} y_{x} + 1 = 0 \quad \text{subject to:} \begin{cases} y(0) = 1 \\ y_x(0) = 0 \end{cases}.
\end{equation*}
This equation is singular at the initial value of $x = 0$ due to the coefficient function $\frac{2}{x}$. However, we can avoid this by multiplying both sides of the equation by the variable $x$. Hence, the differential equation becomes,
\begin{equation*}
    x \, y_{xx} + 2 \, y_x + x = 0,
\end{equation*}
which when evaluated at $x=0$ gives us the initial derivative constraint. By substituting the constrained expression into the differential equation, we are left with an algebraic equation with unknowns $\B{\xi}$,
\begin{equation}\label{eq:laneEm_xiDE}
    \Big[x \, c^2 \B{h}_{zz} + 2 \, \Big(c \B{h}_z -  c \B{h}_z(z_0)\Big)\Big]\T \B{\xi} = -x,
\end{equation}
where the coefficient $c$ comes from our mapping between the basis function domain and problem domain (recall Equation \eqref{eq:s3_linearMapping}). Now, by discretizing the domains, Equation \eqref{eq:laneEm_xiDE} can be written as a linear system of equation such that,
\begin{equation*}
    \begin{bmatrix}\Big[x_0 \, c^2 \B{h}_{zz}(z_0) + 2 \, \Big(c \B{h}_z(z_0) -  c \B{h}_z(z_0)\Big)\Big]\T \\ \vdots \\ \Big[x_k \, c^2 \B{h}_{zz}(z_k) + 2 \, \Big(c \B{h}_z(z_k) -  c \B{h}_z(z_0)\Big)\Big]\T \\ \vdots \\ \Big[x_f \, c^2 \B{h}_{zz}(z_f) + 2 \, \Big(c \B{h}_z(z_f) -  c \B{h}_z(z_0)\Big)\Big]\T \end{bmatrix} \B{\xi} = \begin{Bmatrix} -x_0 \\ \vdots \\ -x_k \\ \vdots \\ -x_f \end{Bmatrix}
\end{equation*}
which is of the form $A \B{x} = \B{b}$ and can be solved with any least-squares technique. While the construction of this linear system was straightforward, there is another formalization that will be consistent between linear and nonlinear differential equations. To realize this, consider rewriting the differential equation as the loss function, 
\begin{equation*}
    \tilde{F} =  \Big[x \, c^2 \B{h}_{zz} + 2 \, \Big(c \B{h}_z -  c \B{h}_z(z_0)\Big)\Big]\T \B{\xi} + x = 0
\end{equation*}
which can be written as a loss vector which is the discretization of $\tilde{F}$ at the collocation nodes,
\begin{equation*}
    \mathbb{L}(\B{\xi}) = \begin{Bmatrix} \tilde{F}(x_0,\B{\xi}) \\ \vdots \\ \tilde{F}(x_k,\B{\xi}) \\ \vdots \\ \tilde{F}(x_f,\B{\xi}) \end{Bmatrix} =  \begin{Bmatrix} \Big[x_0 \, c^2 \B{h}_{zz}(z_0) + 2 \, \Big(c \B{h}_z(z_0) -  c \B{h}_z(z_0)\Big)\Big]\T \B{\xi} + x_0 \\ \vdots \\ \Big[x_k \, c^2 \B{h}_{zz}(z_k) + 2 \, \Big(c \B{h}_z(z_k) -  c \B{h}_z(z_0)\Big)\Big]\T \B{\xi} + x_k \\ \vdots \\ \Big[x_f \, c^2 \B{h}_{zz}(z_f) + 2 \, \Big(c \B{h}_z(z_f) -  c \B{h}_z(z_0)\Big)\Big]\T \B{\xi} + x \end{Bmatrix} = \B{0}
\end{equation*}
with the Jacobian term of,
\begin{equation*}
    \mathbb{J}(\B{\xi}) = \begin{bmatrix} \Big[x_0 \, c^2 \B{h}_{zz}(z_0) + 2 \, \Big(c \B{h}_z(z_0) -  c \B{h}_z(z_0)\Big)\Big]\T \\ \vdots \\ \Big[x_k \, c^2 \B{h}_{zz}(z_k) + 2 \, \Big(c \B{h}_z(z_k) -  c \B{h}_z(z_0)\Big)\Big]\T \\ \vdots \\ \Big[x_f \, c^2 \B{h}_{zz}(z_f) + 2 \, \Big(c \B{h}_z(z_f) -  c \B{h}_z(z_0)\Big)\Big]\T \end{bmatrix}
\end{equation*}
where the equation,
\begin{equation*}
    \mathbb{J}(\B{\xi} =\B{0}) \B{\xi} = - \mathbb{L}(\B{\xi} =\B{0}) \quad \Longleftrightarrow \quad A \B{x} = \B{b};
\end{equation*}
however, this is the same as the first iteration of the nonlinear least-squares approach. Therefore, writing all problems (linear or nonlinear) using the loss function and Jacobian allows us to use the same process and simplify notation.

Next, for the Lane-Emden equation where $a=1$, the loss function becomes
\begin{equation*}
    \tilde{F} =  \Big[x \, c^2 \B{h}_{zz} + 2 \, \Big(c \B{h}_z -  c \B{h}_z(z_0)\Big) + x \, \Big(\B{h} - \B{h}(z_0) - x \, c \B{h}_z(z_0)\Big) \Big]\T \B{\xi} = 0
\end{equation*}
making the loss vector,
\begin{align*}
    \mathbb{L}(\B{\xi}) &= \begin{Bmatrix} \tilde{F}(x_0,\B{\xi}) \\ \vdots \\ \tilde{F}(x_k,\B{\xi}) \\ \vdots \\ \tilde{F}(x_f,\B{\xi}) \end{Bmatrix} \\ &= \begin{Bmatrix} \Big[x_0 \, c^2 \B{h}_{zz}(z_0) + 2 \, \Big(c \B{h}_z(z_0) -  c \B{h}_z(z_0)\Big) + x_0 \, \Big(\B{h}(z_0) - \B{h}(z_0) - x_0 \, c \B{h}_z(z_0)\Big) \Big]\T \B{\xi}  \\ \vdots \\ \Big[x_k \, c^2 \B{h}_{zz}(z_k) + 2 \, \Big(c \B{h}_z(z_k) -  c \B{h}_z(z_0)\Big) + x_k \, \Big(\B{h}(z_k) - \B{h}(z_0) - x_k \, c \B{h}_z(z_0)\Big) \Big]\T \B{\xi}  \\ \vdots \\ \Big[x_f \, c^2 \B{h}_{zz}(z_f) + 2 \, \Big(c \B{h}_z(z_f) -  c \B{h}_z(z_0)\Big) + x_f \, \Big(\B{h}(z_f) - \B{h}(z_0) - x_f \, c \B{h}_z(z_0)\Big) \Big]\T \B{\xi}  \end{Bmatrix}
\end{align*}
with Jacobian,
\begin{equation*}
    \mathbb{J}(\B{\xi}) = \begin{bmatrix}\Big[x_0 \, c^2 \B{h}_{zz}(z_0) + 2 \, \Big(c \B{h}_z(z_0) -  c \B{h}_z(z_0)\Big) + x_0 \, \Big(\B{h}(z_0) - \B{h}(z_0) - x_0 \, c \B{h}_z(z_0)\Big) \Big]\T  \\ \vdots \\ \Big[x_k \, c^2 \B{h}_{zz}(z_k) + 2 \, \Big(c \B{h}_z(z_k) -  c \B{h}_z(z_0)\Big) + x_k \, \Big(\B{h}(z_k) - \B{h}(z_0) - x_k \, c \B{h}_z(z_0)\Big) \Big]\T \\ \vdots \\ \Big[x_f \, c^2 \B{h}_{zz}(z_f) + 2 \, \Big(c \B{h}_z(z_f) -  c \B{h}_z(z_0)\Big) + x_f\, \Big(\B{h}(z_f) - \B{h}(z_0) - x_f \, c \B{h}_z(z_0)\Big) \Big]\T \end{bmatrix}.
\end{equation*}.

\subsection{Nonlinear ordinary differential equations}\label{sec:s3_NL}
Now, let us consider the nonlinear cases of the Lane-Emden equation where, 
\begin{equation*} 
    y_{xx} + \frac{2}{x} y_{x} + y^a = 0 \quad \text{such that} \quad (x > 0, a \geq 2) \quad \text{subject to:} \quad \begin{cases} y(0) = 1 \\ y_x(0) = 0 \end{cases}
\end{equation*}
Again, since the constraints are the same as the linear instance of the differential equation, the constrained expression is the same as in Equation \eqref{eq:s3_laneEmCEy}. Now, the approach is exactly the same as the linear cases. First, we form the loss vector such that,
\begin{equation*}
    \tilde{F} = x \, y_{xx} + 2 y_{x} + x \, y^a  = 0
\end{equation*}
where for clarity the terms $y$, $y_x$, and $y_{xx}$ are not expanded. These equations are defined by Equation \eqref{eq:s3_laneEmCEy}, Equation \eqref{eq:s3_laneEmCEyp}, and Equation \eqref{eq:s3_laneEmCEypp}, respectively.
This produces the loss vector,
\begin{equation*}
    \mathbb{L}(\B{\xi}) =  \begin{Bmatrix} \tilde{F}(x_0,\B{\xi}) \\ \vdots \\ \tilde{F}(x_k,\B{\xi}) \\ \vdots \\ \tilde{F}(x_f,\B{\xi}) \end{Bmatrix} = \begin{Bmatrix}  x_0 \, y_{xx}(x_0,\B{\xi}) + 2 y_{x}(x_0,\B{\xi}) + x_0 \, y^a(x_0,\B{\xi}) \\ \vdots \\ x_k \, y_{xx}(x_k,\B{\xi}) + 2 y_{x}(x_k,\B{\xi}) + x_k \, y^a(x_k,\B{\xi}) \\ \vdots \\ x_f \, y_{xx}(x_f,\B{\xi}) + 2 y_{x}(x_f,\B{\xi}) + x_f \, y^a(x_f,\B{\xi})  \end{Bmatrix} 
\end{equation*}
Additionally, it follows that the Jacobian is,
\small
\begin{equation*}
    \mathbb{J}(\B{\xi}) = \begin{bmatrix}\Big[x_0 \, c^2 \B{h}_{zz}(z_0) + 2 \, \Big(c \B{h}_z(z_0) -  c \B{h}_z(z_0)\Big) + x_0 \, a \, y^{a-1}(x_0,\B{\xi}) \Big(\B{h}(z_0) - \B{h}(z_0) - x_0 \, c \B{h}_z(z_0)\Big) \Big]\T  \\ \vdots \\ \Big[x_k \, c^2 \B{h}_{zz}(z_k) + 2 \, \Big(c \B{h}_z(z_k) -  c \B{h}_z(z_0)\Big) + x_k \, a \, y^{a-1}(x_k,\B{\xi}) \Big(\B{h}(z_k) - \B{h}(z_0) - x_k \, c \B{h}_z(z_0)\Big) \Big]\T \\ \vdots \\ \Big[x_f \, c^2 \B{h}_{zz}(z_f) + 2 \, \Big(c \B{h}_z(z_f) -  c \B{h}_z(z_0)\Big) + x_f\, a \, y^{a-1}(x_f,\B{\xi}) \Big(\B{h}(z_f) - \B{h}(z_0) - x_f \, c \B{h}_z(z_0)\Big) \Big]\T \end{bmatrix}
\end{equation*}.
\normalsize

Again, following the same process, the nonlinear least-squares method is used to update the $\B{\xi}$ coefficient vector and ultimately solve the differential equation. In the proceeding section, we look at the accuracy obtained for this problem.

\subsection{Numerical results of the Lane-Emden equation}
The Lane-Emden equation has an analytical solution for the following values of $a$,
\begin{equation*}
    a = \begin{cases} 
    0 \quad\longrightarrow\quad y = 1 - \dfrac{x^2}{6} \\ \\
    1 \quad\longrightarrow\quad y = \dfrac{\sin(x)}{x}\\ \\
    5 \quad\longrightarrow\quad y = \dfrac{1}{\sqrt{1 + \dfrac{x^2}{3}}}
    \end{cases}
\end{equation*}
In the following examples, we will solve this differential equation for these values of $a$ to directly compare with the analytical solution. This will allow us to analyze the accuracy of the TFC method compared to others in the literature.

\begin{whiteExample}{Lane-Emden (\texorpdfstring{$a = 0$}{a = 0})}
In this example, the Lane-Emden equation is solved for $a=0$ on the domain $x \in [0, 10]$. The results given in Figure \ref{fig:s3_landEmden_type_0_sweep} detail the TFC method's accuracy compared to the spectral method using either Chebyshev polynomials or ELMs with the sigmoid function. Figure \ref{fig:s3_landEmden_type_0_compare} compares the TFC method to spectral both expressed using Chebyshev polynomials to directly quantify the maximum accuracy of these two methods. Finally, Figure \ref{fig:s3_laneEmden_type_0_time} provides a speed versus accuracy comparison of the techniques mentioned above along with the RK45 technique using SciPy's \verb"scipy.integrate.solve_ivp" algorithm \cite{scipy}.

\begin{figure}[H]
	\centering
    \includegraphics[width=.7\linewidth]{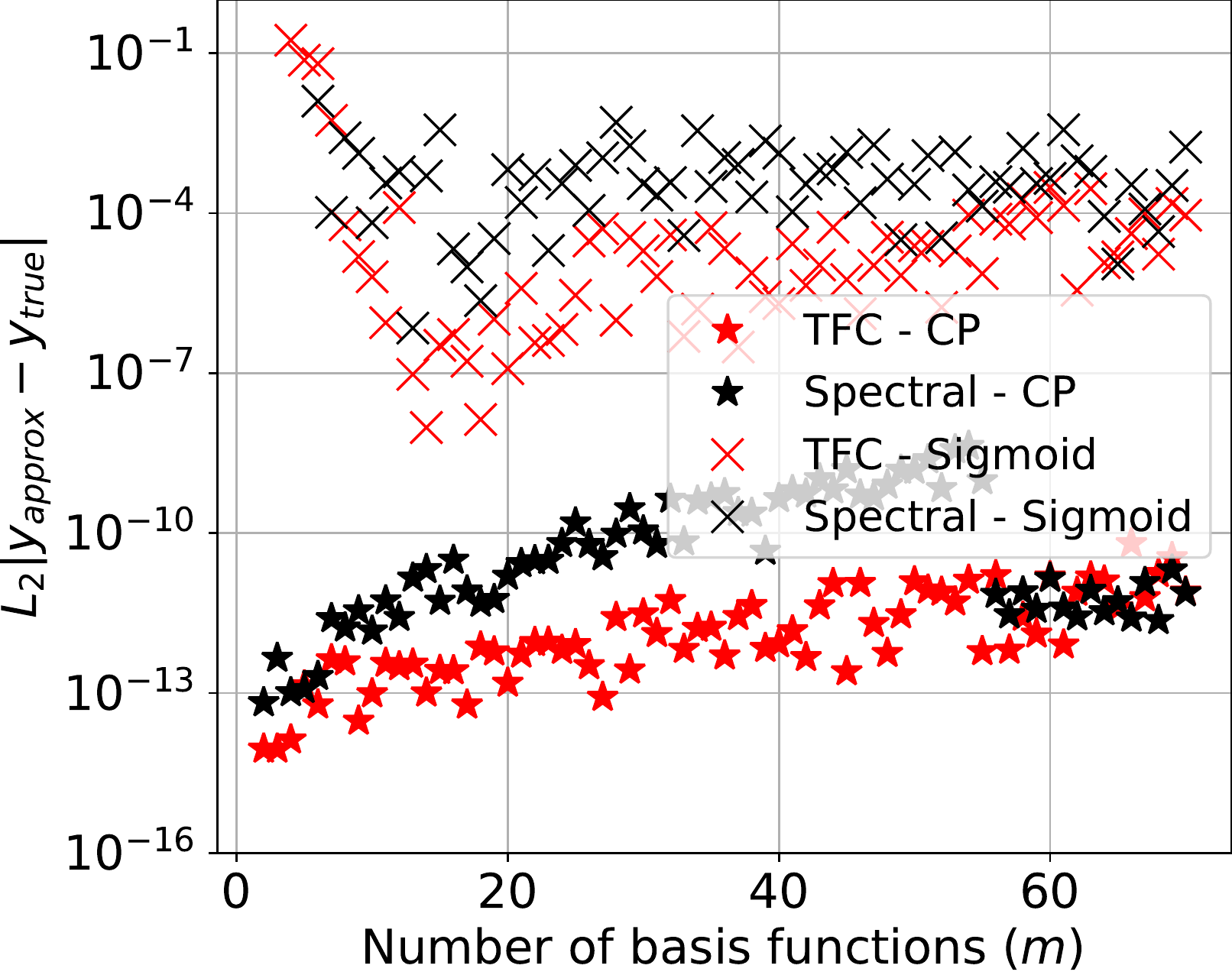}
    \caption{Accuracy of TFC and spectral method for varying number and types of basis functions for the Lane-Emdem equation ($a=0)$.}
    \label{fig:s3_landEmden_type_0_sweep}
\end{figure}

Looking at the TFC based solutions given in Figure \ref{fig:s3_landEmden_type_0_sweep}, it can be seen that the orthogonal polynomial definition of the free function provides dramatic accuracy gain at a lower number of terms. Furthermore, even by adding basis terms, the ELM based free function (sigmoid) does no match the accuracy of the Chebyshev orthogonal polynomials. In fact, for the solution of ordinary differential equations, ELMs are never more accurate than the orthogonal basis set. 

Looking at the comparison of the TFC method with the spectral method given in Figure \ref{fig:s3_landEmden_type_0_compare}, we can see a slight accuracy gain when using TFC versus a spectral method that increases as the number of basis functions increases. However, this gain of accuracy at higher basis functions is misleading because overall, both methods lose accuracy with this increase. This can be explained by looking at the analytical solution for $a=0$, which is a quadratic polynomial. This means that an expression (either with the spectral method or TFC) based on the orthogonal polynomials should have the best solution at $m=2$. Any terms past this only contribute noise to the solution of the differential equations. Regardless, at $m=2$, the TFC method is about an order of magnitude more accurate than the spectral method. 

\begin{figure}[H]
	\centering
    \includegraphics[width=.7\linewidth]{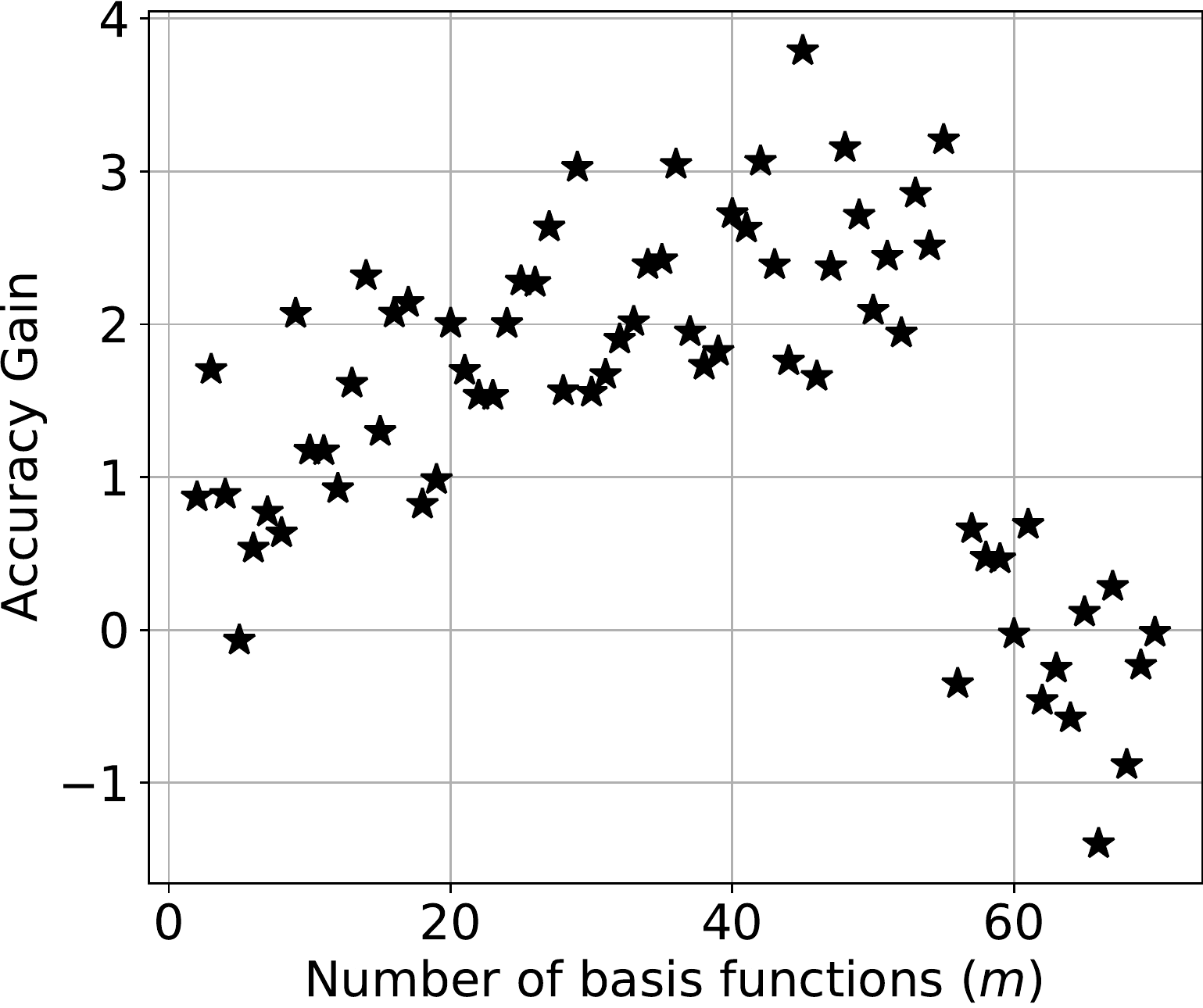}
    \caption{Accuracy gain of TFC vs. spectral for the Solution of Lane-Emdem ($a=0)$. The accuracy gain is quantified in terms of $\log_{10}(\frac{\text{spectral method error}}{\text{TFC error}})$ and therefore, the $y$-axis is by orders of magnitude. For example, when this value is greater than zero, TFC is more accurate, and vice-versa.}
    \label{fig:s3_landEmden_type_0_compare}
\end{figure}

Figure \ref{fig:s3_laneEmden_type_0_time} shows that when more solution accuracy is needed, the RK45 method requires more time to solve the problem, while the spectral and TFC method see little change in computation time. However, comparing spectral and TFC method, there seems to be little difference in accuracy versus speed, with TFC maintaining only a slight advantage.

\begin{figure}[H]
	\centering
    \includegraphics[width=.7\linewidth]{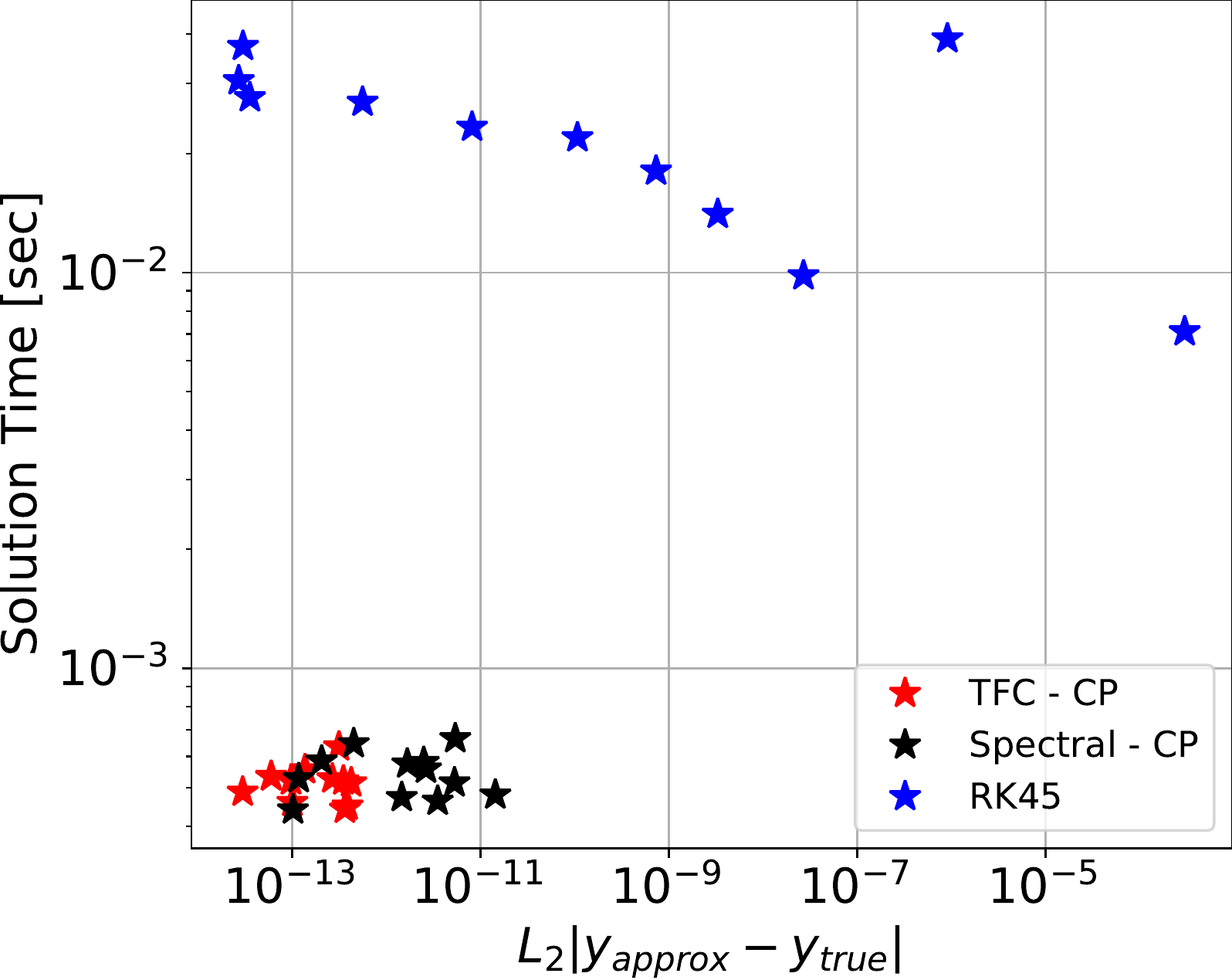}
    \caption{Timed solution of Lane-Emdem ($a=0$).}
    \label{fig:s3_laneEmden_type_0_time}
\end{figure}

\end{whiteExample}

\begin{whiteExample}{Lane-Emden (\texorpdfstring{$a = 1$}{a = 1})}
In this example, the Lane-Emden equation is solved for $a=1$ on the domain $x \in [0, 10]$. The results given in Figures \ref{fig:s3_landEmden_type_1_sweep}-\ref{fig:s3_landEmden_type_1_compare} compare the TFC method to both spectral method and ELMs based on the number of basis terms used. Additionally, Figure \ref{fig:s3_laneEmden_type_1_time} provides a speed versus accuracy comparison of the techniques mentioned above along with the RK45 technique using SciPy's \verb"scipy.integrate.solve_ivp" algorithm \cite{scipy}.

\begin{figure}[H]
	\centering
    \includegraphics[width=.7\linewidth]{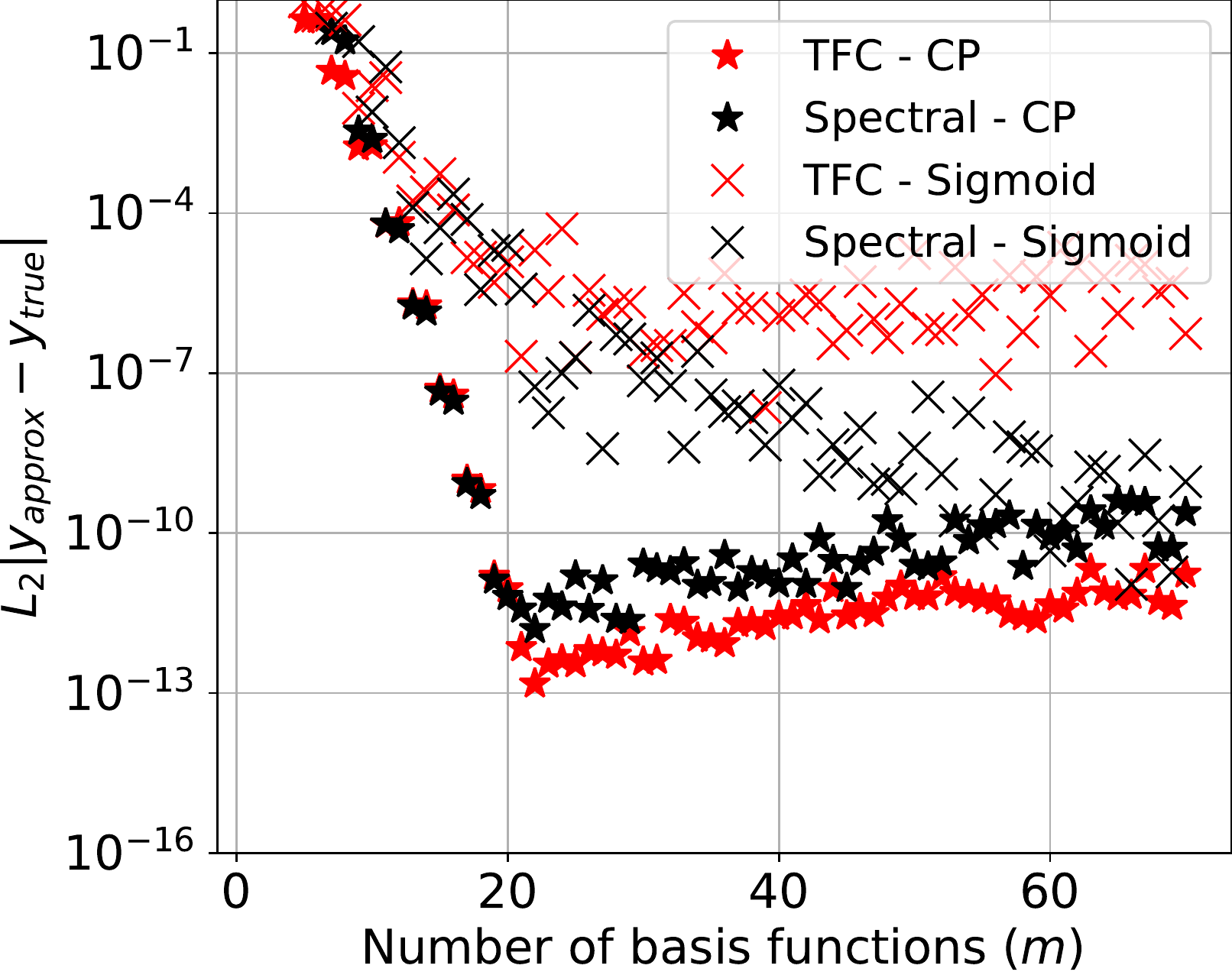}
    \caption{Accuracy of TFC and spectral method for varying number and types of basis functions for the Lane-Emdem equation ($a=1$).}
    \label{fig:s3_landEmden_type_1_sweep}
\end{figure}

Looking at the TFC based solutions given in Figure \ref{fig:s3_landEmden_type_1_sweep}, it can be seen that the orthogonal polynomial definition of the free function quickly reaches a minimum at 22 basis terms. Furthermore, even by adding basis terms, the ELM-based free functions (sigmoid) do not match the Chebyshev orthogonal polynomials' accuracy.

Looking at the comparison of the TFC method with the spectral method in Figure \ref{fig:s3_landEmden_type_1_compare}, we can see that the TFC method is always more accurate than the spectral method when more than 20 basis terms are used. However, at lower basis terms, TFC and the spectral method are comparable in terms of accuracy.

\begin{figure}[H]
	\centering
    \includegraphics[width=.7\linewidth]{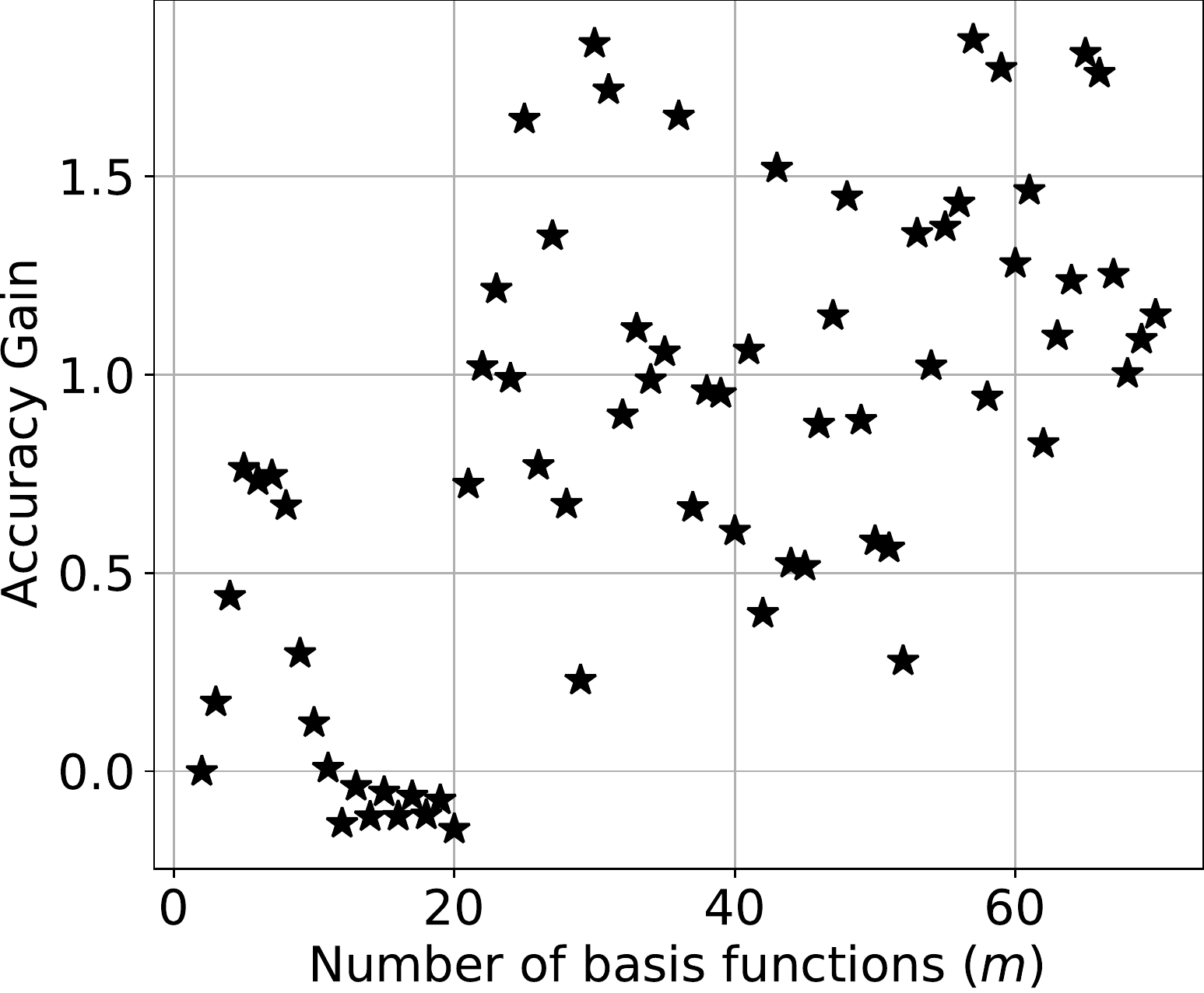}
    \caption{Accuracy gain of TFC vs. spectral for the solution of Lane-Emdem ($a=1)$. The accuracy gain is quantified in terms of $\log_{10}(\frac{\text{spectral method error}}{\text{TFC error}})$, and therefore, the $y$-axis is by orders of magnitude. For example, when this value is greater than zero, TFC is more accurate, and vice-versa.}
    \label{fig:s3_landEmden_type_1_compare}
\end{figure}

Figure \ref{fig:s3_laneEmden_type_1_time} shows that when more solution accuracy is needed, the RK45 method requires more time to solve the problem, while the spectral and TFC method see little change in computation time. In this case, TFC is slightly more accurate and faster than the spectral method.

\begin{figure}[H]
	\centering
    \includegraphics[width=.7\linewidth]{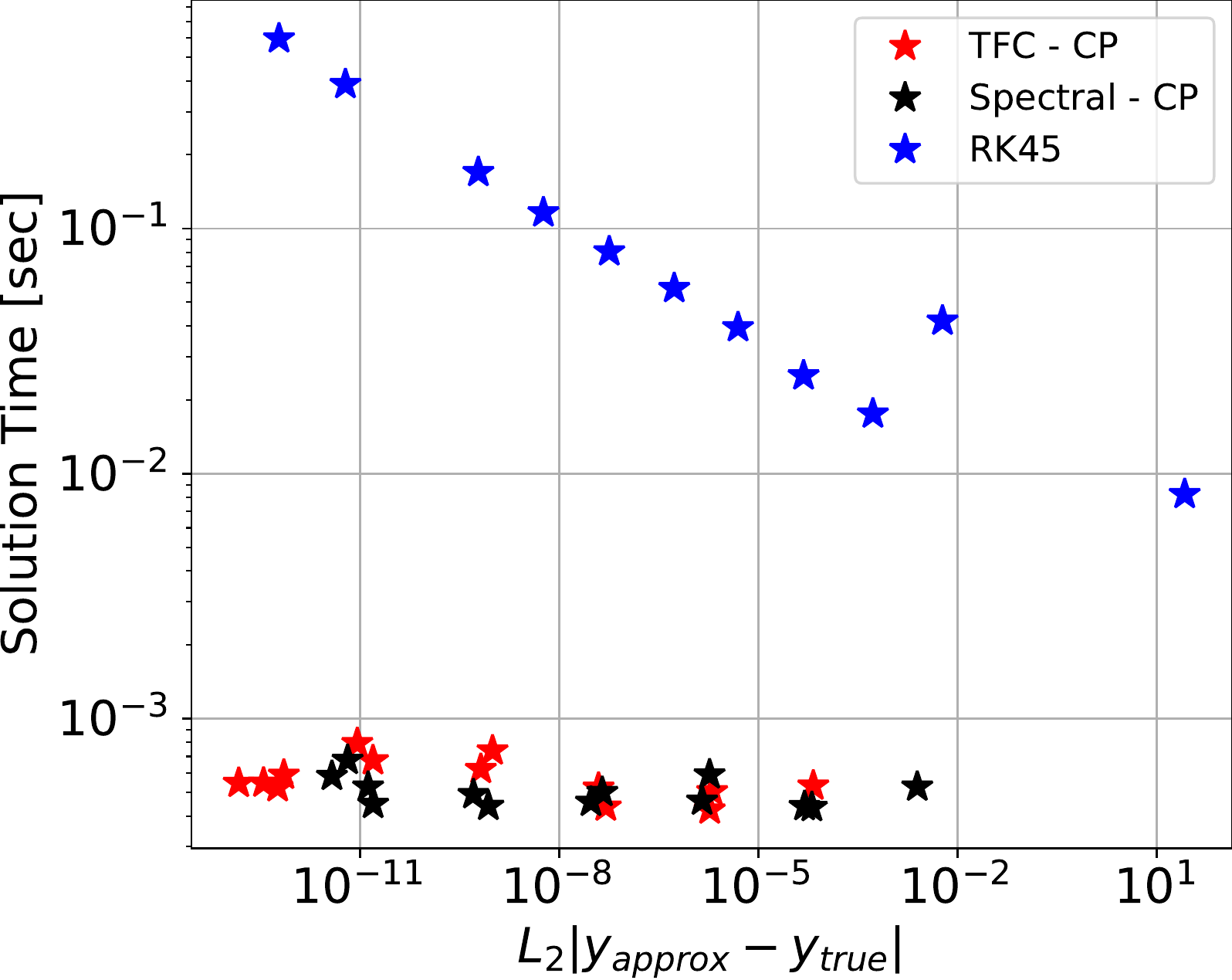}
    \caption{Timed solution of Lane-Emdem ($a=1$).}
    \label{fig:s3_laneEmden_type_1_time}
\end{figure}
\end{whiteExample}

\begin{whiteExample}{Lane-Emden (\texorpdfstring{$a = 5$}{a = 5})}
In this example, the Lane-Emden equation is solved for $a=5$ on the domain $x \in [0, 10]$. The results given in Figures \ref{fig:s3_landEmden_type_5_sweep}-\ref{fig:s3_landEmden_type_5_compare} compare the TFC method to the spectral method with a varying number of basis terms. Additionally, Figure \ref{fig:s3_laneEmden_type_5_time} provides a speed versus accuracy comparison of the techniques mentioned above, along with the RK45 technique.

\begin{figure}[H]
	\centering
    \includegraphics[width=.7\linewidth]{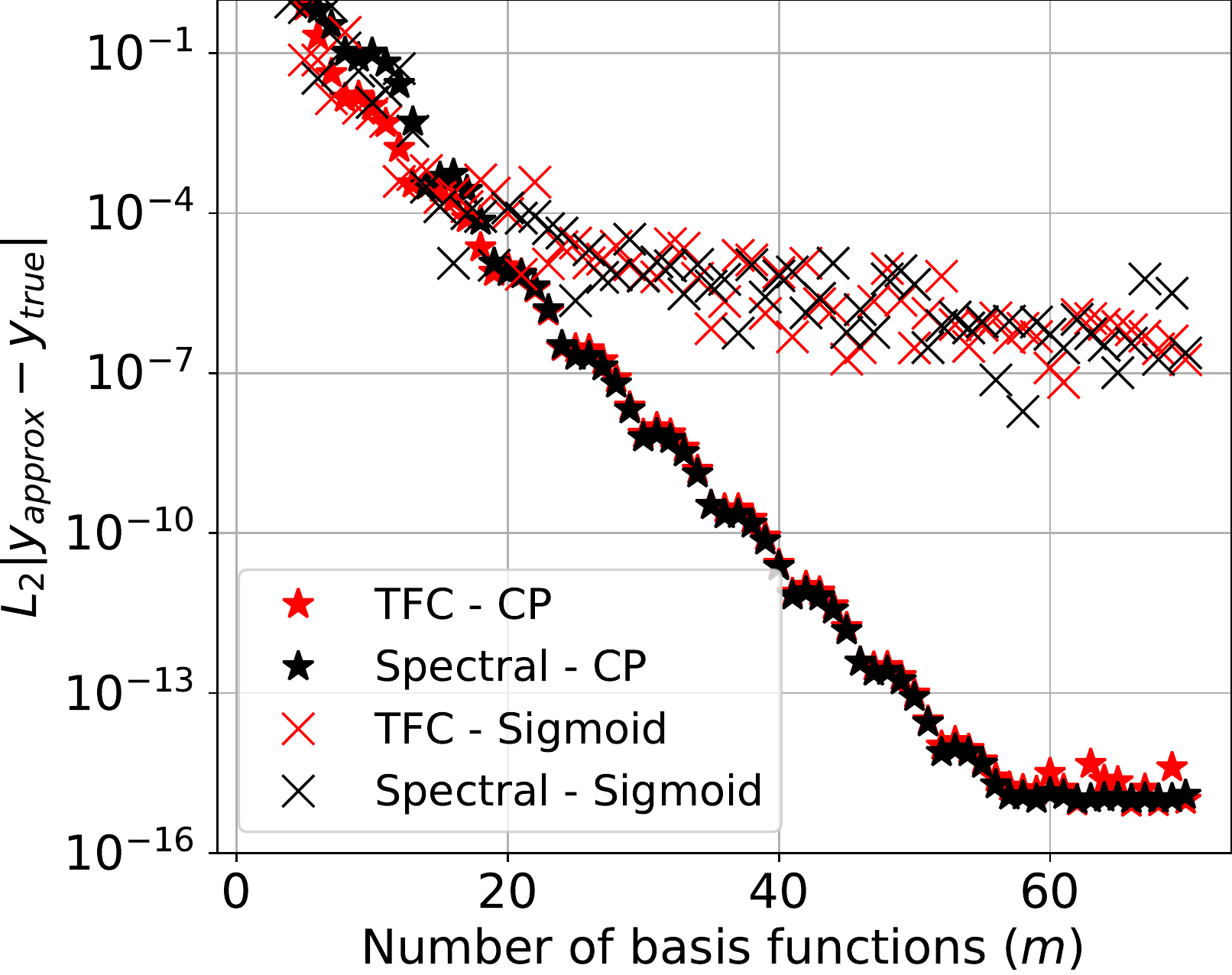}
    \caption{Accuracy of TFC and spectral method for varying number and types of basis functions for the Lane-Emdem equation ($a=5$).}
    \label{fig:s3_landEmden_type_5_sweep}
\end{figure}

Looking at the TFC based solutions given in Figure \ref{fig:s3_landEmden_type_5_sweep}, it can be seen that the orthogonal polynomial definition of the free function quickly reaches a minimum at 62 basis terms. Furthermore, even by adding basis terms, the ELM-based free functions do not match the Chebyshev orthogonal polynomials' accuracy. In fact, the solution with the sigmoid function is seven orders of magnitude less accurate.

Next, comparing the TFC method with the spectral method in Figure \ref{fig:s3_landEmden_type_5_compare}, we can see a slight accuracy gain for the TFC method until about 20 terms, where spectral and TFC method are the same in terms of accuracy. Then, around 60 terms, the spectral method has a slight accuracy gain.

\begin{figure}[H]
	\centering
    \includegraphics[width=.7\linewidth]{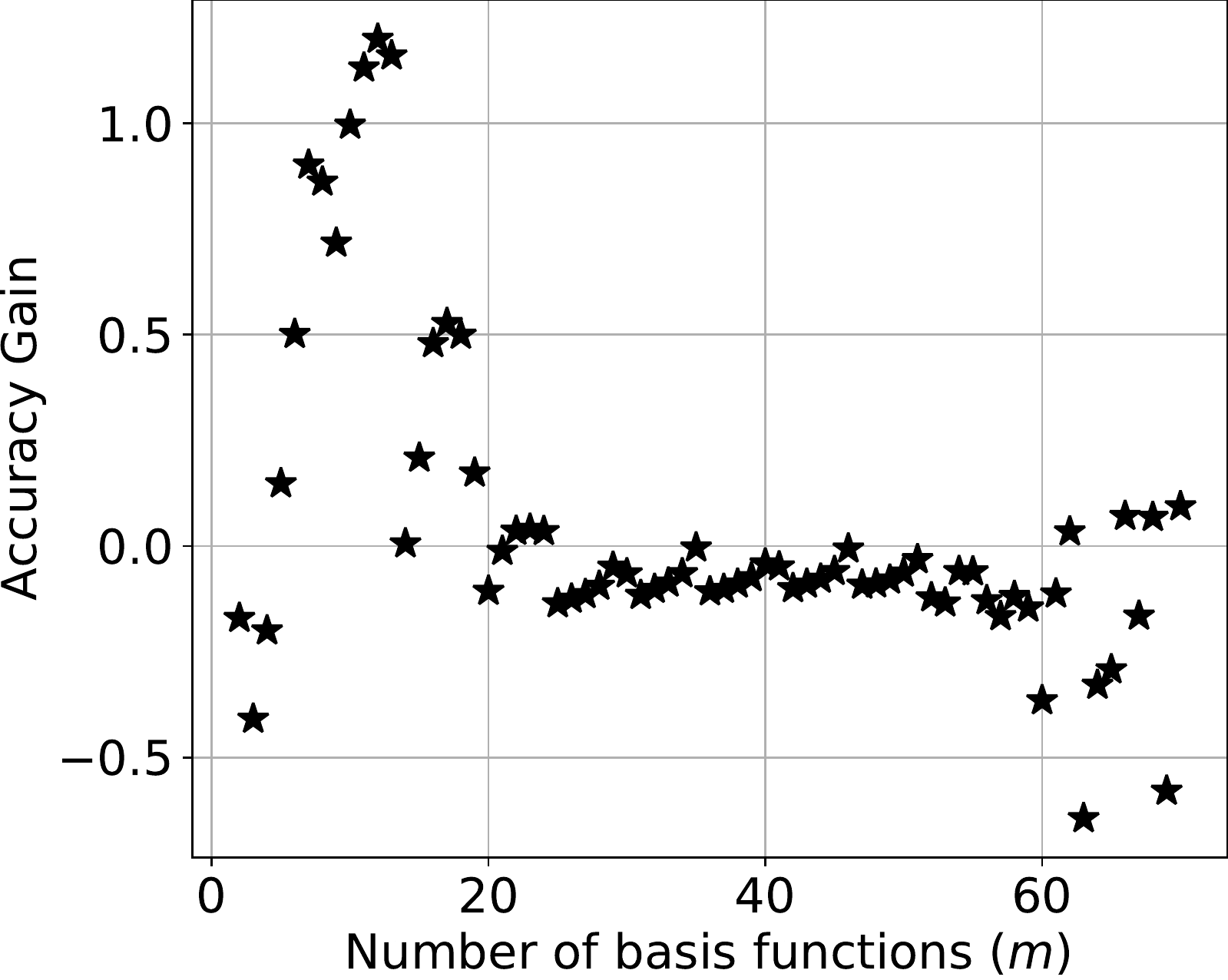}
    \caption{Accuracy gain of TFC vs. spectral for the solution of Lane-Emdem ($a=5)$. The accuracy gain is quantified in terms of $\log_{10}(\frac{\text{spectral method error}}{\text{TFC error}})$, and therefore, the $y$-axis is by orders of magnitude. For example, when this value is greater than zero, TFC is more accurate, and vice-versa.}
    \label{fig:s3_landEmden_type_5_compare}
\end{figure}

In Figure \ref{fig:s3_laneEmden_type_5_time}, we can see when more solution accuracy is needed, the RK45 method requires more time to solve the problem; however, in this case, we do see a similar trend in the spectral and TFC method where the speed is reduced for more accurate solutions. As for the comparison between spectral and TFC method, in this test, the TFC method is slightly faster to converge.
\begin{figure}[H]
	\centering
    \includegraphics[width=.7\linewidth]{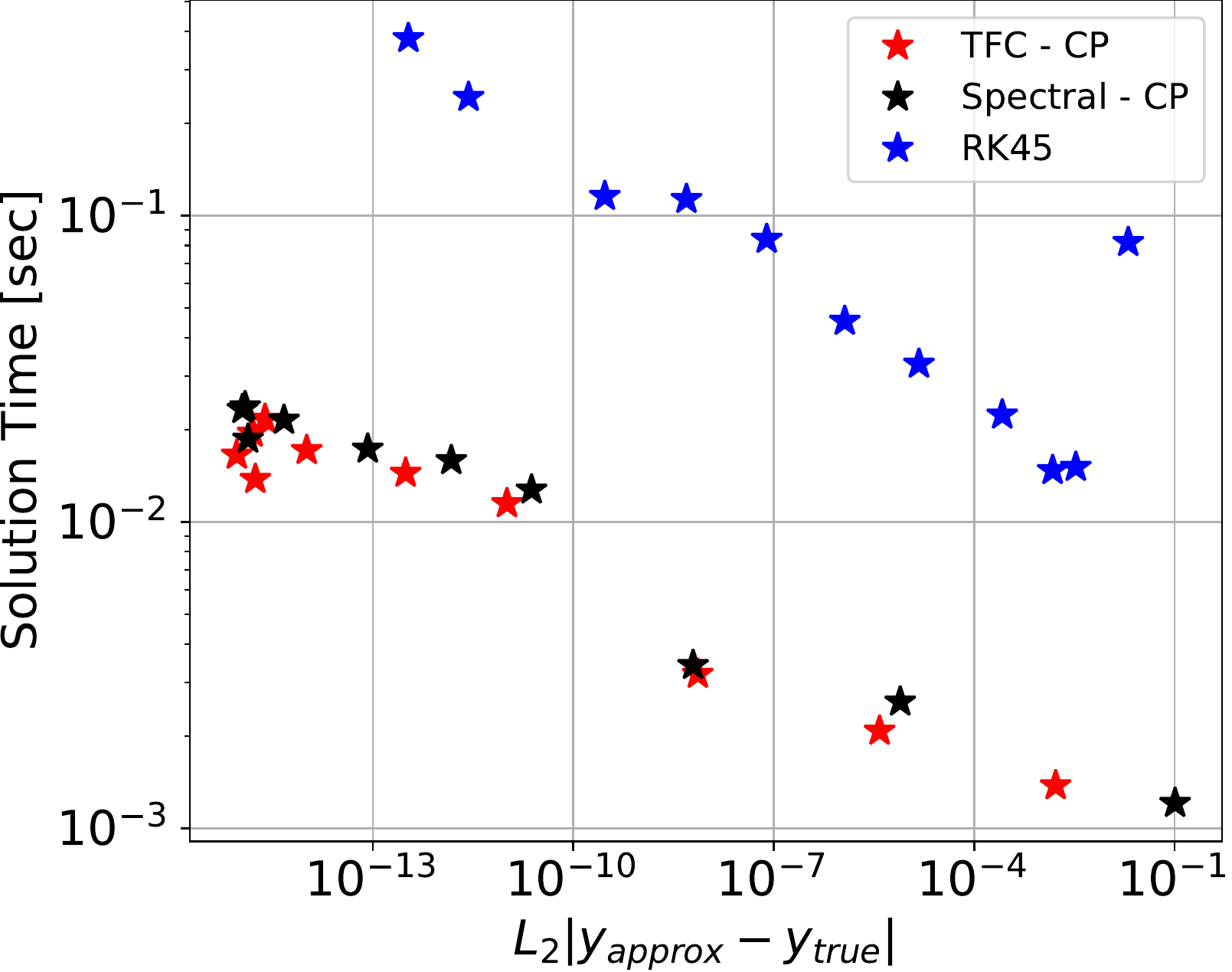}
    \caption{Timed solution of Lane-Emdem ($a=5$).}
    \label{fig:s3_laneEmden_type_5_time}
\end{figure}
\end{whiteExample}

\section{Boundary-value problem}\label{sect:s3_bvp}
From Section \ref{sec:s3_L} and Section \ref{sec:s3_NL} is was observed that solving nonlinear differential equations with TFC is the same as solving linear differential equations with one exception: the nonlinear case requires multiple iterations to solve for $\B{\xi}$. In fact, the TFC approach is a unified approach to solve differential equations, meaning that the solution method is the same regardless of the constraints. This property results from the constrained expression, which decouples the differential equation constraints from the dynamics. To highlight this, let's consider the solution of a two-point boundary value problem,
\begin{equation}\label{eq:s3_BVP}
    y_{xx} + y y_{x} = f(x) \quad \text{subject to:} \quad \begin{cases} y(0) = 0 \\ y(\pi) = 0 \end{cases}
\end{equation}
such that $f(x) = e^{-2x} \sin(x) \Big(\cos(x) - \sin(x)\Big)-2e^{-x}\cos(x)$. Using the our generalized theory, the projection functionals are,
\begin{equation*}
    \rho_1(x,g(x)) = - g(0) \qquad \text{and} \qquad \rho_2(x,g(x)) = - g(\pi).
\end{equation*}
Again, the switching functions are determined by choosing the support functions $s_1 = 1$ and $s_2 = x$ and solving for the coefficients $\alpha_{ij}$,
\begin{align*}
    \begin{bmatrix} 1 & 0 \\ 1 & \pi \end{bmatrix} \begin{bmatrix} \alpha_{11} & \alpha_{12}\\  \alpha_{21} & \alpha_{22} \end{bmatrix} &= \begin{bmatrix} 1 & 0 \\ 0 & 1 \end{bmatrix}
\end{align*}
\begin{equation*}
    \begin{bmatrix} \alpha_{11} & \alpha_{12}\\  \alpha_{21} & \alpha_{22} \end{bmatrix} = \begin{bmatrix} 1 & 0 \\ 1 & \pi \end{bmatrix}^{-1} = \frac{1}{\pi}\begin{bmatrix} \pi & 0 \\ -1 & 1 \end{bmatrix}
\end{equation*}
which leads to the switching functions,
\begin{equation*}
    \phi_1(x) = \frac{\pi - x}{\pi} \qquad \text{and} \qquad
    \phi_2(x) = \frac{x}{\pi}.
\end{equation*}
The constrained expression in terms of $\B{\xi}$ is
\begin{align*}
    y(x,\B{\xi}) &= \Big(\B{h} -\frac{\pi - x}{\pi} \B{h}(z_0) - \frac{x}{\pi} \B{h}(z_f) \Big)\T \B{\xi} \\ 
    y_x(x,\B{\xi}) &= \Big(c \, \B{h}_{z} + \frac{1}{\pi} \B{h}(z_0) - \frac{1}{\pi} \B{h}(z_f) \Big)\T \B{\xi} \\ 
    y_{xx}(x,\B{\xi}) &= \Big(c^2 \, \B{h}_{zz} \Big)\T \B{\xi}.
\end{align*}
Just like the Lane-Emden initial-value problem, the constraints are embedded, and we have a transformed differential equation subject to no constraints. Therefore, the last step is to form the loss vector and Jacobian and solve for the coefficients using our nonlinear least-squares method. Therefore, it should now be clear by this example that the process of solving the differential equations is unaffected by different constraint types. For completeness, the associated loss function, loss vector, and Jacobian are provided below.
\begin{equation*}
    \tilde{F} = y_{xx} + y \, y_{x} - f(x) = 0
\end{equation*}
\begin{equation*}
    \mathbb{L}(\B{\xi}) =  \begin{Bmatrix} \tilde{F}(x_0,\B{\xi}) \\ \vdots \\ \tilde{F}(x_k,\B{\xi}) \\ \vdots \\ \tilde{F}(x_f,\B{\xi}) \end{Bmatrix} = \begin{Bmatrix}  y_{xx}(x_0,\B{\xi}) + y(x_0,\B{\xi}) \, y_x(x_0,\B{\xi}) - f(x_0) \\ \vdots \\ y_{xx}(x_k,\B{\xi}) + y(x_k,\B{\xi}) \, y_x(x_k,\B{\xi}) - f(x_k) \\ \vdots \\ y_{xx}(x_f,\B{\xi}) + y(x_f,\B{\xi}) \, y_x(x_f,\B{\xi}) - f(x_f)  \end{Bmatrix} 
\end{equation*}

\scriptsize
\begin{equation*}
    \mathbb{J}(\B{\xi}) = \begin{bmatrix}\Big[c^2 \, \B{h}_zz(z_0) + y(x_0,\B{\xi})\Big(c \, \B{h}_{z}(z_0) + \frac{1}{\pi} \B{h}(z_0) - \frac{1}{\pi} \B{h}(z_f) \Big) + y_x(x_0,\B{\xi})\Big(\B{h}(z_0) -\frac{\pi - x_0}{\pi} \B{h}(z_0) - \frac{x_0}{\pi} \B{h}(z_f) \Big) \Big]\T  \\ \vdots \\ \Big[c^2 \, \B{h}_zz(z_k) + y(x_k,\B{\xi})\Big(c \, \B{h}_{z}(z_k) + \frac{1}{\pi} \B{h}(z_0) - \frac{1}{\pi} \B{h}(z_f) \Big) + y_x(x_k,\B{\xi})\Big(\B{h}(z_k) -\frac{\pi - x_k}{\pi} \B{h}(z_0) - \frac{x_k}{\pi} \B{h}(z_f) \Big) \Big]\T   \\ \vdots \\ \Big[c^2 \, \B{h}_zz(z_f) + y(x_f,\B{\xi})\Big(c \, \B{h}_{z}(z_f) + \frac{1}{\pi} \B{h}(z_0) - \frac{1}{\pi} \B{h}(z_f) \Big) + y_x(x_f,\B{\xi})\Big(\B{h}(z_f) -\frac{\pi - x_f}{\pi} \B{h}(z_0) - \frac{x_f}{\pi} \B{h}(z_f) \Big) \Big]\T   \end{bmatrix}
\end{equation*}.
\normalsize

\begin{whiteExample}{Solution to two-point boundary-value problem}
In this example, the two-point boundary-value problem given by Equation \eqref{eq:s3_BVP}. The results given in Figures \ref{fig:s3_BVP_TFC}-\ref{fig:s3_BVP_compare} compare the TFC method to both spectral method and ELMs based on the number of basis terms used. Additionally, Figure \ref{fig:s3_BVP_time} provides a speed versus accuracy comparison of the techniques mentioned above along a 4th order collocation algorithm with the control of residuals from SciPy's \verb"scipy.integrate.solve_bvp" algorithm \cite{scipy}.

\begin{figure}[H]
	\centering
    \includegraphics[width=.7\linewidth]{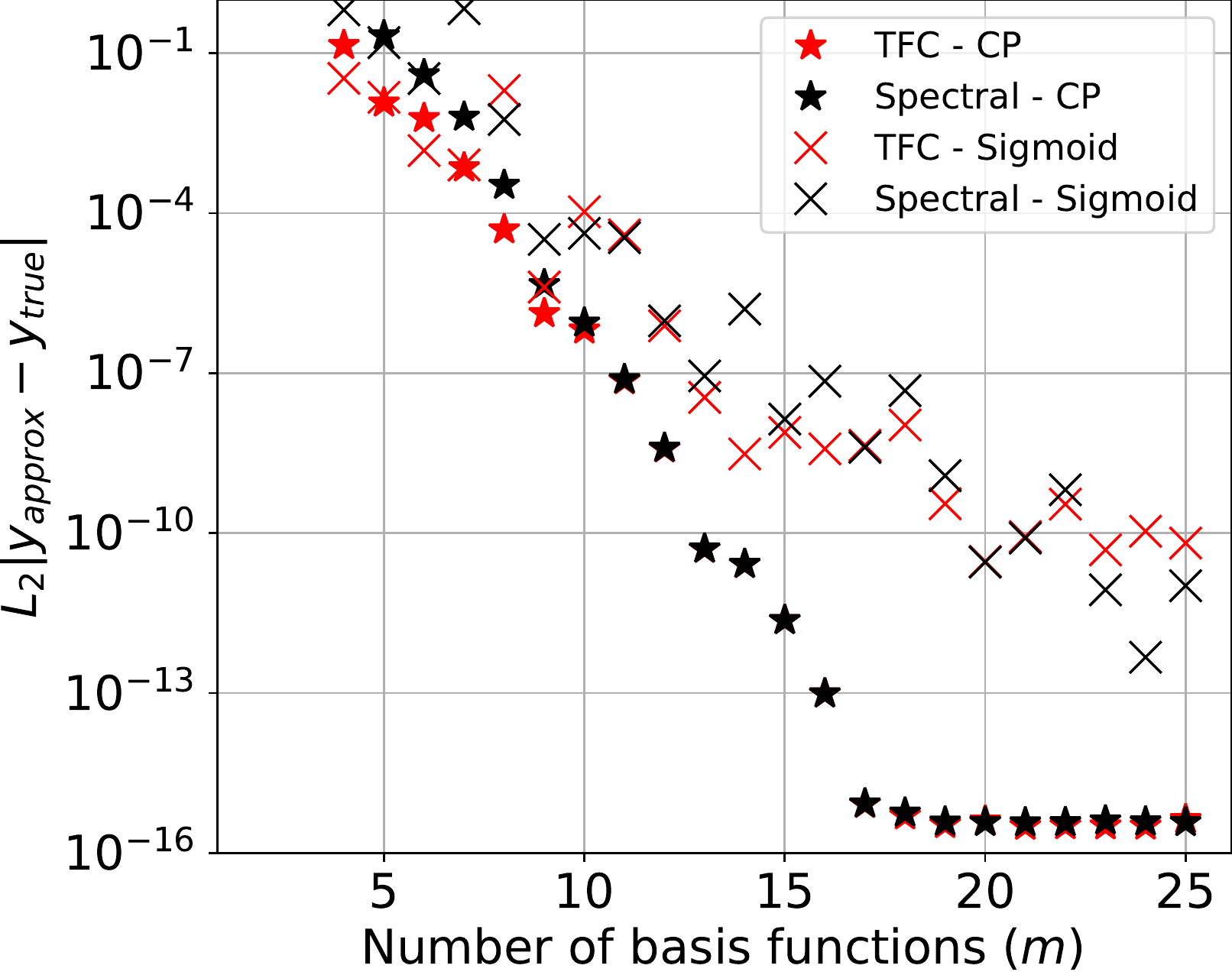}
    \caption{Accuracy of TFC and spectral method for varying number and types of basis functions for the boundary-value problem.}
    \label{fig:s3_BVP_TFC}
\end{figure}

Looking at the TFC based solutions given in Figure \ref{fig:s3_BVP_TFC}, it can be seen that similar to the solutions of the Lane-Emden differential equation, the orthogonal polynomial definition of the free function is superior. Additionally, at 22 Chebyshev basis terms, both TFC and spectral method reach a minimum with respect to solution error. Further analysis shows that that the ELM based free functions are at least 3 orders of magnitude less accurate than the orthogonal polynomials; it is clear that using TFC with orthogonal polynomials to solve ordinary differential equations is the preferred approach. Therefore, after this example, all following examples will utilize Chebyshev or Legendre polynomials as the free function. 

Next, the comparison of the TFC method with the spectral method is given in Figure \ref{fig:s3_BVP_compare}. We can see a slight accuracy gain for the TFC method until about ten terms, where the spectral and TFC method are the same in terms of accuracy. Then, around 20 terms, the TFC method has a slight accuracy gain.

\begin{figure}[H]
	\centering
    \includegraphics[width=.7\linewidth]{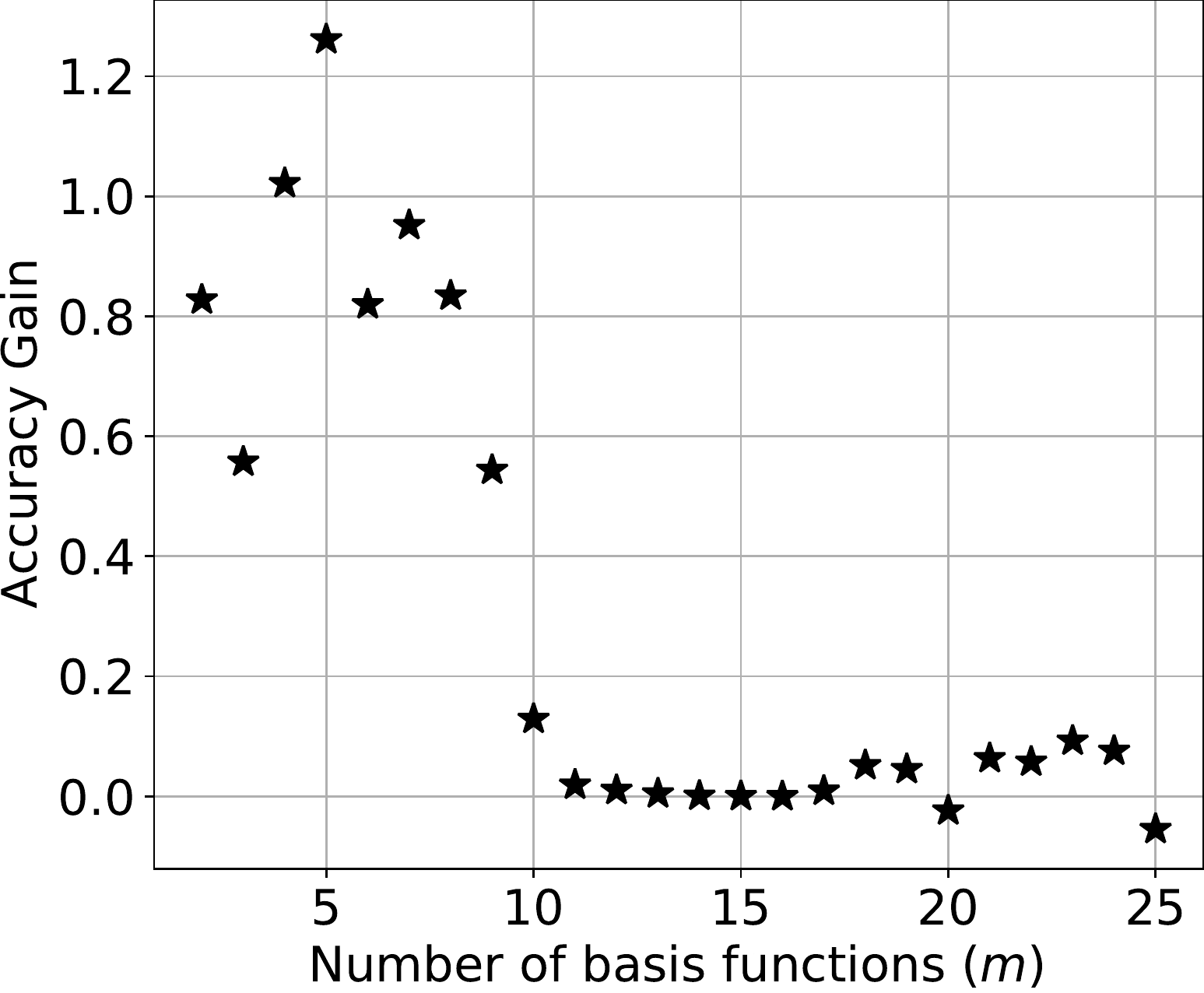}
    \caption{Accuracy gain of TFC vs. spectral method for the solution of the simple boundary-value problem. The accuracy gain is quantified in terms of $\log_{10}(\frac{\text{spectral method error}}{\text{TFC error}})$ and therefore, the $y$-axis is by orders of magnitude.}
    \label{fig:s3_BVP_compare}
\end{figure}

Finally, in Figure \ref{fig:s3_BVP_time}, a comparison of computation time is given for all of the previous techniques along with the RK45 method. Again, when more solution accuracy is needed, the RK45 method paired with a shooting method requires more time to solve the problem. Additionally, the maximum accuracy obtained from this method is on the order of $10^{-11}$. For the spectral and TFC based methods, we notice only a slight increase in computation time with increasing accuracy, and the TFC method is slightly faster.

\begin{figure}[H]
	\centering
    \includegraphics[width=.7\linewidth]{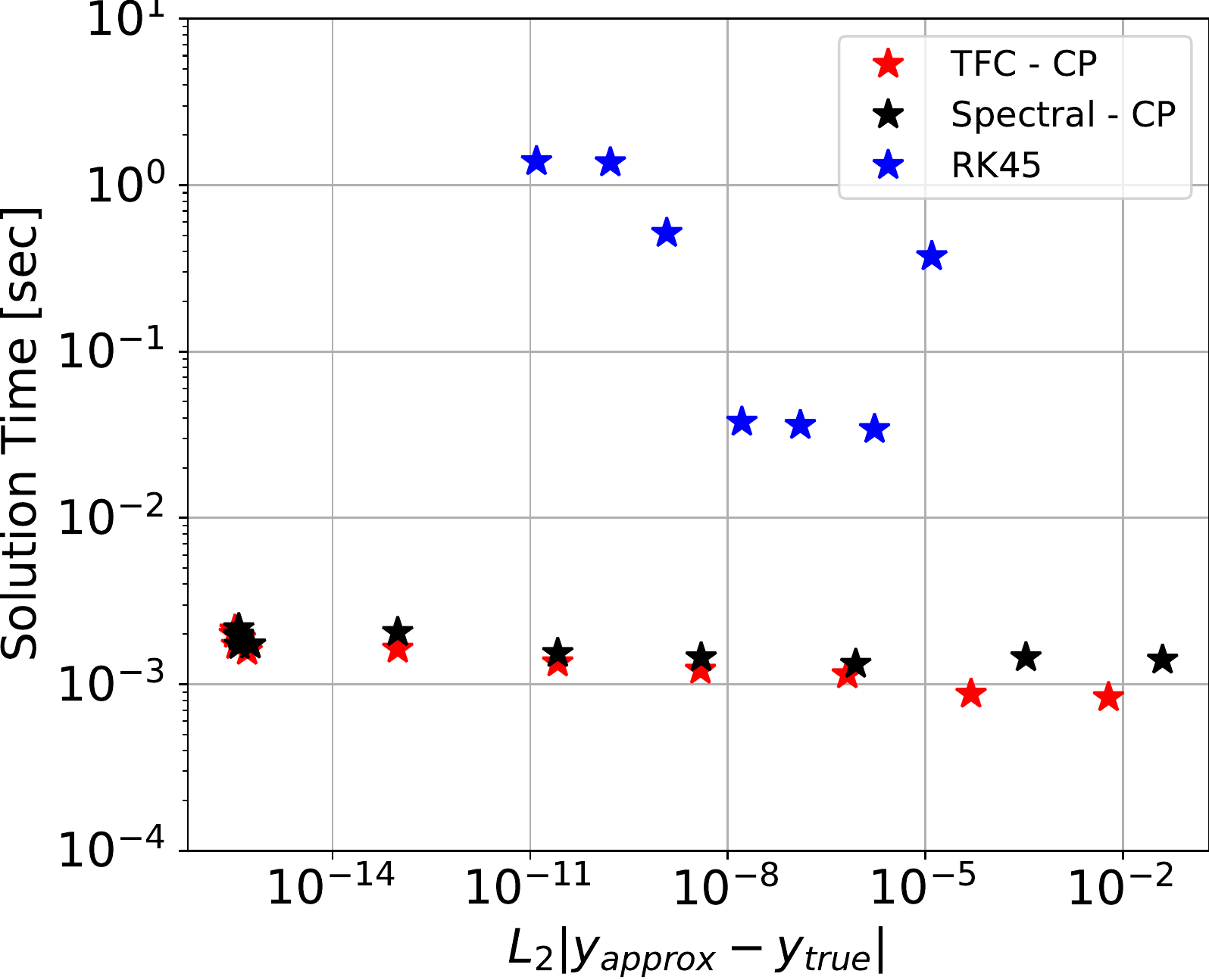}
    \caption{Timed Solution of BVP.}
    \label{fig:s3_BVP_time}
\end{figure}
\end{whiteExample}

\section{Solving systems of ordinary differential equations}

The process discussed to solve single differential equations can directly be used to solve systems of differential equations. In general, we can consider a vector function $\B{v}(t): \R \to \R^n$ where $\B{v}(t) = \{v_1(t), v_2(t), \cdots, v_n(t) \}\T$ where $v_i: \R \to \R$ or in a vector-sense, the components of the vector. This vector function is subject to some set of differential equations and constraints imposed on the $v_i$ components. Therefore, just as we have done in the single differential equation examples, a system of differential equations can be solved by deriving the constrained expressions for the $n$ component functions according to the theory provided in Chapters \ref{chap:tfc_intro} and \ref{chap:tfc_general}. In fact, if constraints are shared between components, the theory can easily incorporate these constraints (see Example \ref{sec:s2_comp}). Finally, these constrained expressions can be parameterized by defining $n$ free functions and creating a system of algebraic equations that must then be discretized and solved as usual.

\section{Two major extensions for use in optimal control problems}\label{sec:s3_extensions}
Until now, we have dealt with ordinary differential equations where 1) the free function $g(x)$ is expressed as an orthogonal polynomial set that can accurately and completely describe the solution and 2) the integration range was explicitly stated (i.e., the initial and final time of the problems were known). However, in many optimal control problems, we run into two scenarios that cause issues with the standard framework. Thus, extra theory must be developed to handle it; however, the tools and concepts developed in the earlier sections make this task an effortless step forward.

\subsection{A hybrid systems approach*}\label{sec:s3_hybridSystem}
\blfootnote{*Reprinted (along with revisions and updates unique to this dissertation) by permission from Elsevier the Journal of Computational and Applied Mathematics ``Least-squares solutions of boundary-value problems in hybrid systems,'' Johnston, H. and Mortari, D., 2021, J. Comput. Appl. Math., 393, 113524, Copyright 2021, \cite{pieceTFC}}

\begin{figure}[ht]
    \centering
    \includegraphics[width=.75\textwidth]{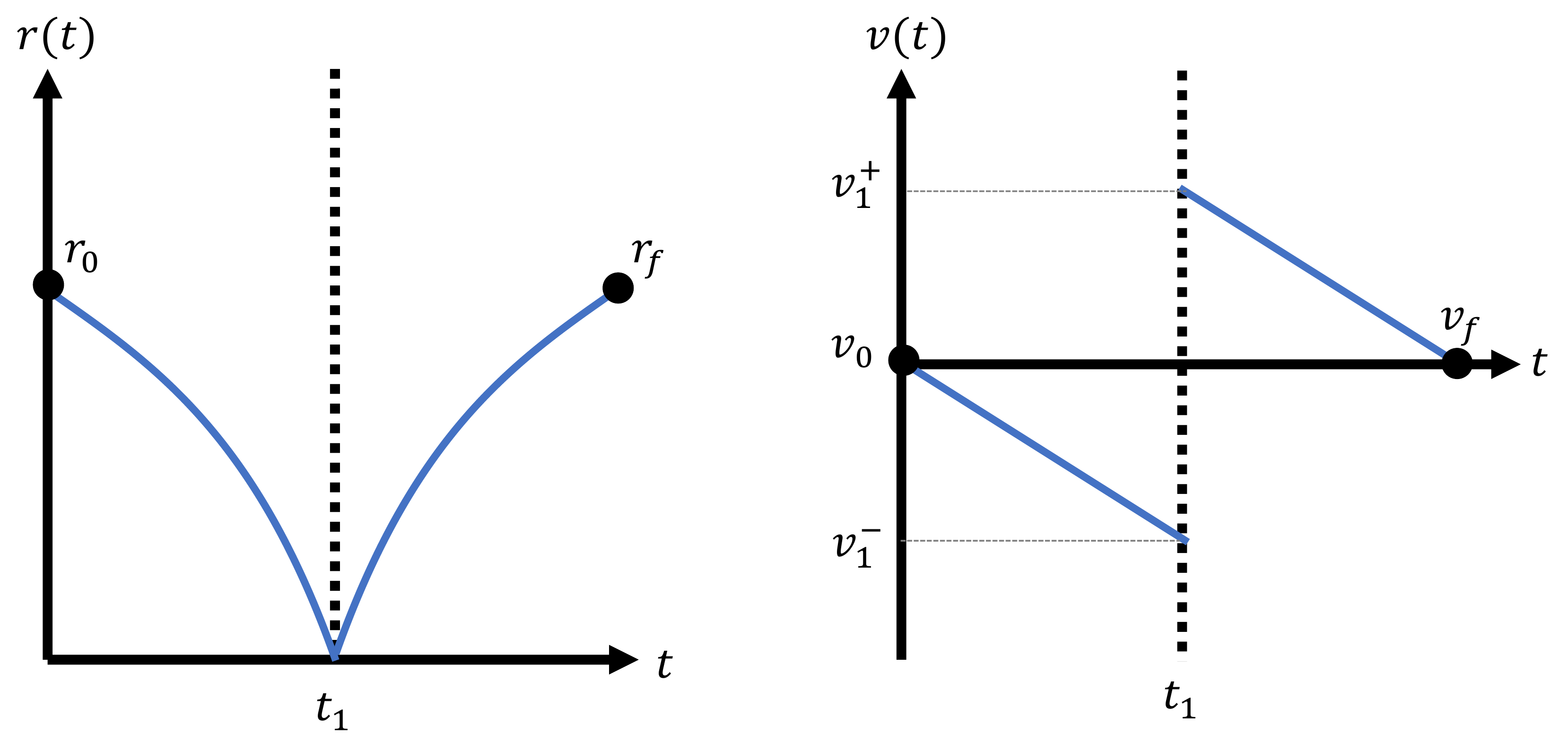}
    \caption{Graphical representation of the bouncing ball hybrid system. Reprinted with permission from \cite{pieceTFC}.}
    \label{fig:s3_bouncingBall}
\end{figure}

First, we need to adapt the constrained expression for use in hybrid systems. The original adaptation was spurred by the problem of bang-bang control structure inherent in the fuel optimal landing problem solved in Johnston et al. \cite{FOL} and explored in more detail in Johnston and Mortari \cite{pieceTFC}. By definition, hybrid systems are dynamical systems governed by a time-sequence of differential equations, either linear or nonlinear. A simple example is a bouncing ball where the motion is described by a sudden variation (or jump) in the dynamics when the ball impacts the ground, shown in Figure \ref{fig:s3_bouncingBall}. 
These systems become even more common in the study of control problems where a dynamical system is controlled by discrete controls (e.g., bang-bang control). In fact, these are considered a special case of hybrid systems called variable structure systems (VSS), and the study of the control of these systems is named variable structure control (VSC) \cite{VSS}.

\begin{figure}[ht]
    \centering
    \includegraphics[width=.5\textwidth]{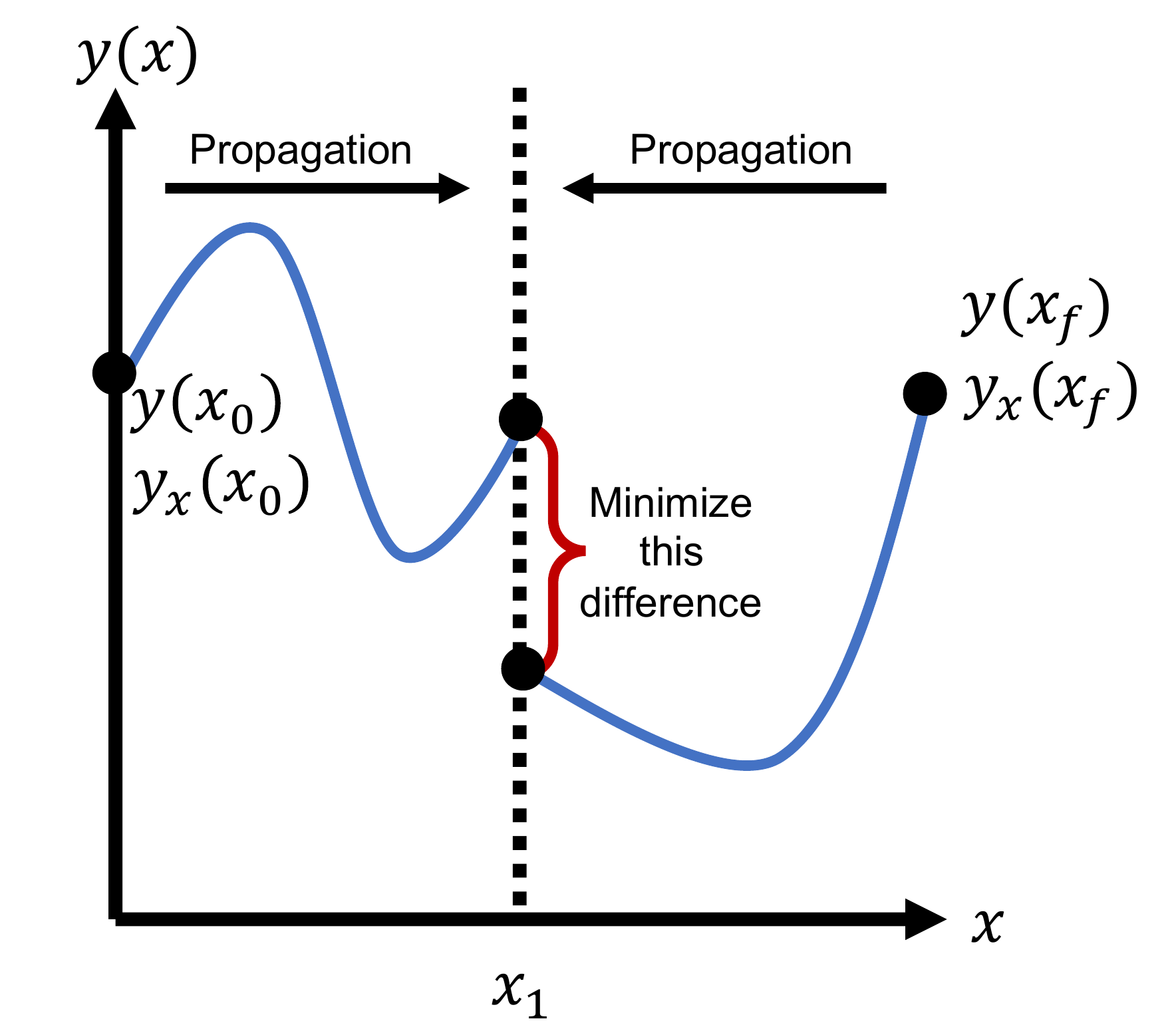}
    \caption{Graphical representation of shooting method. Reprinted with permission from \cite{pieceTFC}.}
    \label{fig:s3_shootingMethod}
\end{figure}

Initial-value problems for these systems can be easily solved by propagating the initial conditions over the domain of the first differential equation in the sequence. The final conditions can then be used as the initial conditions for the next differential equation, and the process can be repeated indefinitely (ignoring any accumulation of numerical error). However, boundary-value problems do not offer this luxury and will be the main focus of the proceeding section. The study of these problems is not new, and numerical techniques to solve these problems have existed since the 1960s, based on the shooting method \cite{shooting_TPBVP,shooting_bang,shooting_BVP,LineShoot} detailed in Figure \ref{fig:s3_shootingMethod}.
In these approaches, the interval is divided over multiple sub-intervals, and the boundary-value problem is converted to multiple initial-value problems. The unknown boundary conditions are then solved by minimizing the DE residuals and the residuals of function and derivative continuities connecting all sub-intervals. In practice, root solving techniques (bisection, Newton's method, etc.) are used to minimize all residuals. In general, even when two subsequent linear differential equations are connected, solutions based on a shooting method requires an initial guess of the unknown parameters that are used to iterate until the solution is obtained. Note that the convergence is dictated by the initial guess \cite{shooting_converg}, and it is not guaranteed. Regardless, studies have been conducted to quantify these methods' error once an approximation is obtained \cite{shooting_err1,shooting_err2}.

Other techniques for solving these problems include finite difference and finite element methods. A finite difference method where the differential equation is approximated by a difference equation that converts the problem into a system of equations that are solved using linear algebra techniques. On the other hand, in finite element methods (collocation, Galerkin, etc.) \cite{FEM}, the problem is split into smaller parts called finite elements. Simple approximated equations are used to model these elements. These elements are then assembled into a larger system of equations that model the entire problem. The finite difference and finite element method's major drawback is the number of subdivisions needed to capture large variations in the solution. 

\begin{figure}[ht]
    \centering\includegraphics[width=.55\linewidth]{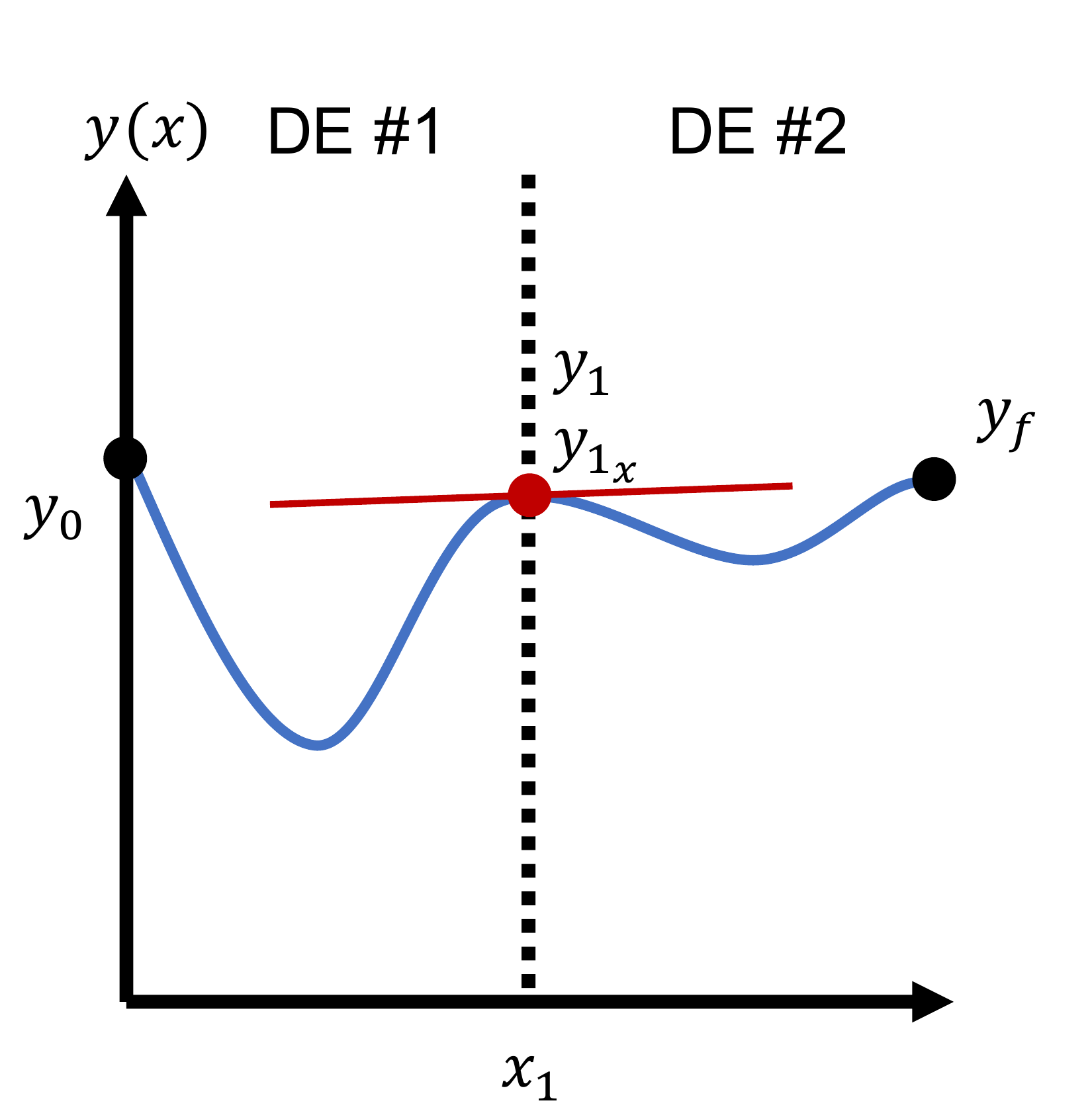}
    \caption{Illustration of piecewise TFC approach enforcing $C^1$ continuity over two segments. Reprinted with permission from \cite{pieceTFC}.}
    \label{fig:s3_seg2Template}
\end{figure}

The simplest example of a hybrid system is a differential equation with a discrete jump in the dynamic behavior at a single point along the domain. When solving a two-point BVP according to these dynamics, not only must the solution satisfy the boundary condition, but it must also preserve the $C^1$ continuity over the jump. The differential equation associated with the single switch in dynamics can be expressed in its explicit form by,
\begin{equation*}
\begin{cases} \p{1}{F}(x,y,y_x,y_{xx}) = 0 \quad \text{for } x \leq x_1\\ \p{2}{F}(x,y,y_x,y_{xx}) = 0 \quad \text{for } x > x_1 \end{cases} \quad \text{subject to:} \begin{cases} y(x_0) = y_0 \\ y(x_f) = y_f\end{cases}
\end{equation*}
where $x_1 \in (x_0, x_f)$ and $\p{1}{F}(x,y,y_x,y_{xx})$ and $\p{2}{F}(x,y,y_x,y_{xx})$ are both functions of the independent variable $x$, the function $y$, and its derivatives. For this system, a separate \ce\ for  
each segment must be derived. Additionally, at the boundary of the differential equations, in this case $x_1$, continuity must be enforced. Figure \ref{fig:s3_seg2Template} depicts the \ce\ over the two differential equation segments.
This can be done by considering each segment independently and introducing two new unknown values $y(x_1) = y_1$ and $y_x(x_1) = y_{1_x}$, which are the value and derivative of the function at the intersection. Therefore, the \ce\ over the first segment must be written for an initial value, final value, and initial derivative, while the \ce\ over the second segment must be written for an initial value, initial derivative, and final value. Using the theory already developed (and using monomial support functions), these \ces\ take the form,
\begin{align}\label{eq:s3_hybrid_ce1}
    \p{1}{y}(x,\p{1}{g}(x)) = \p{1}{g}(x) &+ \p{1}{\phi_1}(x)\Big(y_0 - \p{1}{g}(x_0)\Big) \nonumber\\&+ \p{1}{\phi_2}(x)\Big(y_1 - \p{1}{g}(x_1)\Big) + \p{1}{\phi_3}(x)\Big(y_{1_x} - \p{1}{g_x}(x_1)\Big) 
\end{align}
\begin{align}\label{eq:s3_hybrid_ce2}
    \p{2}{y}(x,\p{2}{g}(x)) = \p{2}{g}(x) &+ \p{2}{\phi_1}(x)\Big(y_1 - \p{2}{g}(x_1)\Big) \nonumber\\&+ \p{2}{\phi_2}(x)\Big(y_{1_x} - \p{2}{g_x}(x_1)\Big) + \p{2}{\phi_3}(x)\Big(y_f - \p{2}{g}(x_f)\Big)
\end{align}
where the switching functions are provided below,
\begin{align*}
    \p{1}{\phi}_1(x) &= \frac{1}{(x_1 - x_0)^2}\Big( x_1^2 - 2 x_1 x + x^2 \Big) \\ 
    \p{1}{\phi}_2(x) &= \frac{1}{(x_1 - x_0)^2}\Big(x_0(x_0-2x_1) + 2x_1 x - x^2\Big) \\ 
    \p{1}{\phi}_3(x) &= \frac{1}{x_1 - x_0}\Big( x_0 x_1 - (x_0 + x_1)x + x^2\Big)
\end{align*}
\begin{align*}
    \p{2}{\phi}_1(x) &= \frac{1}{(x_f - x_1)^2} \Big(x_f(x_f-2x_1) + 2 x_1 x - x^2 \Big)\\ 
    \p{2}{\phi}_2(x) &= \frac{1}{x_f - x_1} \Big(-x_f x_1 + (x_f + x_1) x - x^2 \Big)\\ 
    \p{2}{\phi}_3(x) &= \frac{1}{(x_f - x_1)^2} \Big(x_1^2 - 2 x_1 x + x^2 \Big).
\end{align*}

The major result of the \ces\ derived in Equations \eqref{eq:s3_hybrid_ce1} and Equation \eqref{eq:s3_hybrid_ce2} is that for all finite values of $y_1$ and $y_{1_x}$, $C^1$ continuity is satisfied. However, this formulation comes with one caveat. Since $y_1$ and $y_{1_x}$ were considered arbitrary, they are free parameters that must be solved for when solving the differential equation. Therefore, for numerical implementation, this causes the number of parameters to be solved to scale with the number of segments in the hybrid system. We will find that this is not a major issue for ordinary differential equations.

\subsubsection{Generalization for \texorpdfstring{$n$}{``n''} segments}
\begin{figure}[ht]
    \centering\includegraphics[width=.65\linewidth]{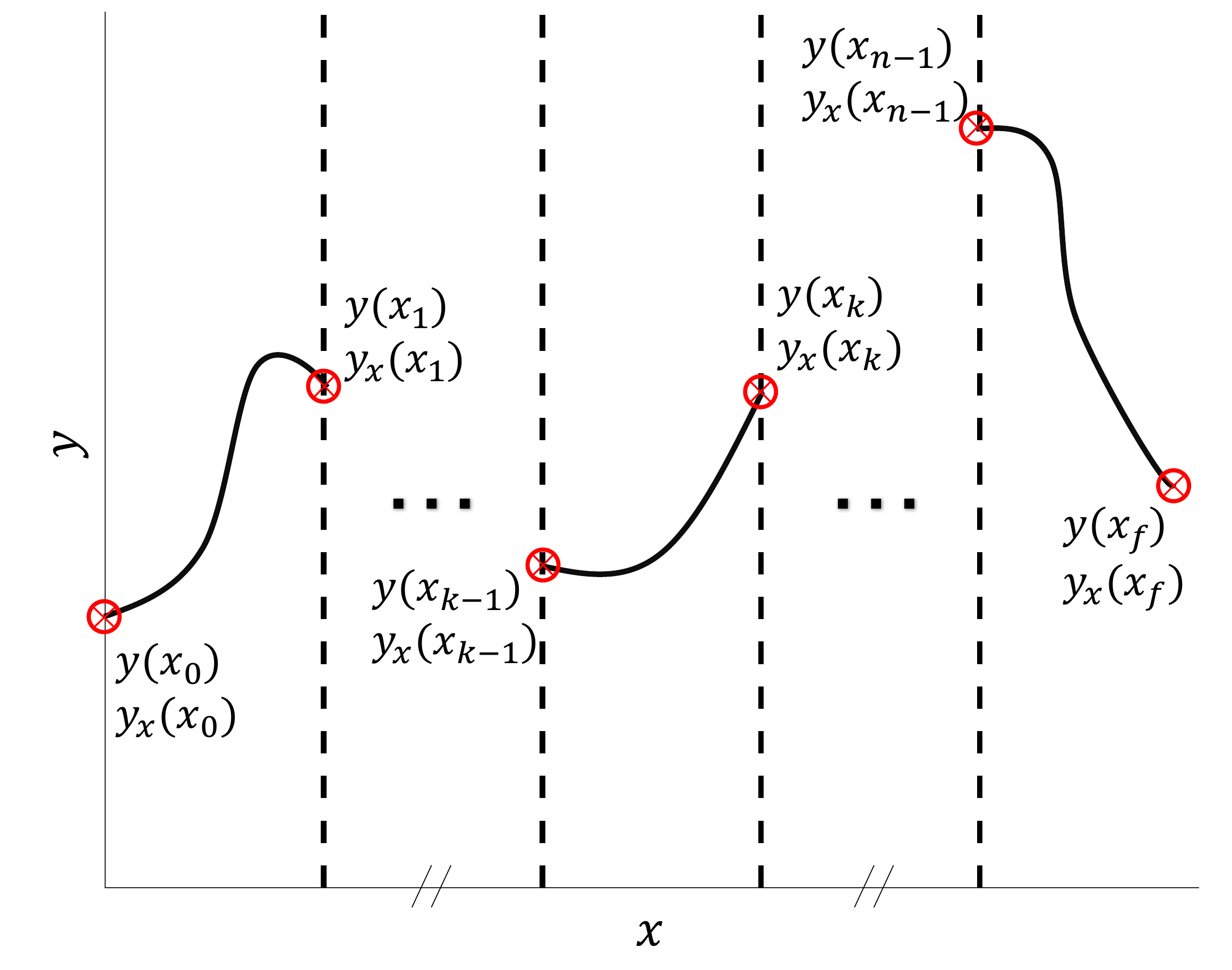}
    \caption{Illustration of segmented TFC approach to enforce $C^1$ continuity over $n$ segments. Reprinted with permission from \cite{pieceTFC}.}
    \label{fig:n_seg}
\end{figure}

Suppose the problem is subject to $n$ jumps in dynamics as detailed in Figure \ref{fig:n_seg}. This case is the generalization of the problem presented in Section \ref{sec:s3_hybridSystem}. 
Additionally, this generalization necessitates the introduction of another set of switching functions that can be derived using the TFC method. Since our future applications will focus on optimal control problems governed by second-order dynamics, we will consider each segment constrained on both sides by point and derivative constraints. The \ce\ for this constraint type produces the equation,
\begin{align*}
    \p{k}{y}(x,\p{k}{g}) = \p{k}{g}(x) &+ \p{k}{\phi}_1(x) \Big(\p{k-1}{\beta} - \p{k}{g}(x_{k-1}) \Big) + \p{k}{\phi}_2(x) \Big(\p{k}{\beta} - \p{k}{g}(x_{k}) \Big) \\ &+ \p{k}{\phi}_3(x) \Big(\p{k-1}{\beta}_x - \p{k}{g}_x(x_{k-1}) \Big) + \p{k}{\phi}_4(x) \Big(\p{k}{\beta}_x - \p{k}{g}_x(x_{k}) \Big),
\end{align*}
where $k = 0, 1, \cdots, n$, and $\p{k-1}{\beta}$, $\p{k}{\beta}$, $\p{k-1}{\beta}_x$, and $\p{k}{\beta}_x$ are the value and derivative continuity constraints when $0 < k < n$. The conditions at $k=0$ and $k=n$ are defined by the boundary constraints.
In this equation, the switching functions (when selecting the support functions as $s_1 = 1$, $s_2 = x$, $s_3 = x^2$, and $s_4 = x^3$) become,
\begin{align*}
    \p{k}{\phi}_1(x) &= \frac{1}{(x_{k} - x_{k-1})^3} \Big( -x_{k}^2 (3  x_{k-1}- x_{k}) + 6  x_{k-1}  x_{k} x  -3 ( x_{k-1}+ x_{k}) x^2 + 2x^3 \Big)\\
	\p{k}{\phi}_2(x) &= \frac{1}{(x_{k} - x_{k-1})^3} \Big( - x_{k-1}^2 ( x_{k-1}-3  x_{k}) -6  x_{k-1}  x_{k} x  + 3 ( x_{k-1}+ x_{k})x^2 -2x^3 \Big)\\
	\p{k}{\phi}_3(x) &= \frac{1}{(x_{k} - x_{k-1})^2} \Big(-x_{k-1}  x_{k}^2 +  x_{k} (2  x_{k-1}+ x_{k})x -(x_{k-1}+2  x_{k})x^2 + x^3 \Big)\\
	\p{k}{\phi}_4(x) &= \frac{1}{(x_{k} - x_{k-1})^2} \Big( - x_{k-1}^2  x_{k} + x_{k-1} ( x_{k-1}+2  x_{k}) x - (2  x_{k-1}+ x_{k}) x^2 + x^3 \Big),
\end{align*}
where $x_k$ denotes the boundaries of the segments. Lastly, by expressing the free function in the form of Equation \eqref{eq:s3_basis} and discretizing the domains, the generalization can be written in a compact block diagonal matrix of the form,
\begin{equation*}
    \B{y} = \begin{bmatrix}\mathbb{A}_1 & \vdots & \mathbb{A}_2 \end{bmatrix} \, \Xi + \mathbb{B}
\end{equation*}
where the terms of this equation are,
\begin{equation*}
    \mathbb{A}_1 = \begin{bmatrix} \p{1}{H} & \cdots & \B{0} & \cdots & \B{0}  \\ \vdots & \ddots & \vdots & \ddots & \vdots \\ \B{0} & \cdots & \p{k}{H} & \cdots & \B{0} \\ \vdots & \ddots & \vdots & \ddots & \vdots \\ \B{0} & \cdots & \B{0} & \cdots & \p{n}{H} \end{bmatrix}.
\end{equation*}
In this block diagonal matrix, the terms of $\p{k}{H}$ are matrices of the terms multiplied by $\p{k}{\B{\xi}}$. For example, the term $\p{1}{H}$ is simply,
\begin{equation*}
   \p{1}{H} = \begin{bmatrix} \Big[\B{h}(z_0) - \p{1}{\phi}_1(x_0)\B{h}(z_0) - \p{1}{\phi}_2(x_0)\B{h}(z_1) - \p{1}{\phi}_3(x_0)\, c \, \B{h}_z(z_0) - \p{1}{\phi}_4(x_0)\, c\, \B{h}_z(z_1) \Big]\T \\ \vdots \\ \Big[\B{h}(z_1) - \p{1}{\phi}_1(x_1)\B{h}(z_0) - \p{1}{\phi}_2(x_1)\B{h}(z_1) - \p{1}{\phi}_3(x_1)\, c \, \B{h}_z(z_0) - \p{1}{\phi}_4(x_1)\, c\, \B{h}_z(z_1) \Big]\T \end{bmatrix}
\end{equation*}
Next,
\begin{landscape}
\vspace*{\fill}
\begin{equation*}
    \mathbb{A}_2 = \begin{bmatrix}  \p{1}{\Phi}_2 & \p{1}{\Phi}_4 & \B{0} & \B{0} & \cdots & \B{0} & \B{0} & \B{0} & \B{0} & \cdots & \B{0} & \B{0} & \B{0} & \B{0} 
    \\
    \p{2}{\Phi}_1 & \p{2}{\Phi}_3 & \p{2}{\Phi}_2 & \p{2}{\Phi}_4 & \cdots & \B{0} & \B{0} & \B{0} & \B{0} & \cdots & \B{0} & \B{0} & \B{0} & \B{0} 
    \\ \vdots & \vdots & \vdots & \vdots & \ddots & \vdots & \vdots & \vdots & \vdots & \ddots & \vdots & \vdots & \vdots & \vdots \\
    \B{0} & \B{0} & \B{0} & \B{0}  & \cdots &  \p{k-1}{\Phi}_2 & \p{k-1}{\Phi}_4 & \B{0} & \B{0}  &\cdots & \B{0} & \B{0} & \B{0} & \B{0} 
    \\
    \B{0} & \B{0} & \B{0} & \B{0}  & \cdots &  \p{k}{\Phi}_1 & \p{k}{\Phi}_3 & \p{k}{\Phi}_2 & \p{k}{\Phi}_4 &\cdots & \B{0} & \B{0} & \B{0} & \B{0} 
    \\
    \B{0} & \B{0} & \B{0} & \B{0}  & \cdots & \B{0} & \B{0}  & \p{k+1}{\Phi}_1 & \p{k+1}{\Phi}_3 &\cdots & \B{0} & \B{0} & \B{0} & \B{0} 
    \\ \vdots & \vdots & \vdots & \vdots & \ddots & \vdots & \vdots & \vdots & \vdots & \ddots & \vdots & \vdots & \vdots & \vdots \\
    \B{0} & \B{0} & \B{0} & \B{0} & \cdots & \B{0} & \B{0} & \B{0} & \B{0} & \cdots  & \p{n-1}{\Phi}_1 & \p{n-1}{\Phi}_3 & \p{n-1}{\Phi}_2 & \p{n-1}{\Phi}_4 
    \\
    \B{0} & \B{0} & \B{0} & \B{0} & \cdots & \B{0} & \B{0} & \B{0} & \B{0} & \cdots & \B{0} & \B{0} & \p{n}{\Phi}_1 & \p{n}{\Phi}_3 \end{bmatrix}, 
\end{equation*}
\vspace{0.5in}
\vspace*{\fill}
\end{landscape}
\noindent and
\begin{equation*}
    \p{k}\Phi_j = \begin{Bmatrix} \p{k}\phi_j(x_{k-1}) & \cdots \p{k}\phi_j(x_k) \end{Bmatrix}\T,
\end{equation*}
is used for the switching functions $\p{k}\phi_j(x)$ evaluated at the discretization points. Lastly,
\begin{equation*}
    \mathbb{B} = \begin{Bmatrix}\p{1}{\Phi}_1\T & \p{1}{\Phi}_3\T & \B{0}\T & \p{n}{\Phi}_2\T & \p{n}{\Phi}_4\T \end{Bmatrix}\T,
\end{equation*}
which is a vector associated with the boundary constraints. For this system, the unknown vector is,
\begin{equation*}
\begin{split}
\Xi &= 
\{\begin{matrix} \p{1}{\B{\xi}}\T & \cdots & \p{k}{\B{\xi}}\T  & \cdots &  \p{n}{\B{\xi}}\T \end{matrix} \\
 &\qquad\qquad \begin{matrix} \p{1}{\beta} & \p{1}{\beta}_x & \cdots &  \p{k-1}{\beta} & \p{k-1}{\beta}_x  & \p{k}{\beta} & \p{k}{\beta}_x & \cdots &  \p{n-1}{\beta} & \p{n-1}{\beta}_x \end{matrix}\}\T.
\end{split}
\end{equation*}

Since this is a linear set of equations all subsequent derivatives are the derivatives of the individual components. The $d$-th order derivative of $\B{y}$ becomes, 
\begin{equation*}
    \B{y}^{(d)} = \begin{bmatrix}\mathbb{A}^{(d)}_1 & \vdots & \mathbb{A}^{(d)}_2 \end{bmatrix} \, \Xi + \mathbb{B}^{(d)},
\end{equation*}
which is also a block diagonal matrix.

Moving forward, numerical examples are provided for two cases:  1) a hybrid system governed by a linear to nonlinear differential equation sequence and 2) the one-dimensional convection-diffusion equation. The solution of the convection-diffusion highlights that this technique can also be applied outside of hybrid systems---specifically when the dynamics of two regions in the differential equation behavior drastically different.

\subsubsection{Linear-to-nonlinear differential equation sequence}\label{sect:s3_linearNonlinearDE}

Consider a second-order linear-nonlinear DE sequence such that,
\begin{equation}\label{eq:linear_nonlinear}
    y_{xx} + y (y_x)^a =  -e^{\pi - 2 x} + e^{\pi/2 - x} \quad \text{subject to: } \begin{cases} y(0) = \dfrac{9}{10} + \dfrac{1}{10} e^{\pi/2}(5 - 2 e^{\pi/2})  \\ y(\pi) = e^{-\pi/2} \end{cases}
\end{equation}
where the parameter $a$ is determined by,
\begin{equation*}
    a = \begin{cases} 0 \quad &\text{for } x \le \pi/2 \\ 1 \quad &\text{for } x > \pi/2 \end{cases}.
\end{equation*}
At the switch, $x_1 = \pi/2$, the differential equation changes from an linear differential equation to a nonlinear differential equation. This differential equation has the unique solution defined by,
\begin{equation*}
    y(x) = \begin{cases} = - \dfrac{1}{5} e^{\pi - 2 x} + \dfrac{1}{2} e^{\pi/2 - x} + \dfrac{9\cos(x) + 7\sin(x)}{10} \quad &\text{for } x \le \pi/2 \\=  e^{\pi/2 - x} \quad &\text{for } x > \pi/2 \end{cases}.
\end{equation*}
Since the sequence has a nonlinear differential equation (over the second segment), an iterative least-squares approach is necessary. For this, we define the residual of the differential equation as the loss functions such that,
\begin{equation}\label{eq:L1ex2}
    \p{1}{\tilde{F}} = \p{1}{y_{xx}} + \p{1}{y} - e^{\pi/2} + e^{\pi/2 - x}
\end{equation}
\begin{equation}\label{eq:L2ex2}
    \p{2}{\tilde{F}} = \p{2}{y_{xx}} + \p{2}{y}\p{2}{y_x} - e^{\pi/2} + e^{\pi/2 - x}
\end{equation}
where $\p{1}{y}$, $\p{1}{y_{xx}}$, $\p{2}{y}$, $\p{2}{y_{x}}$, and $\p{2}{y_{xx}}$ are defined by the \ces\ given by Equations \eqref{eq:s3_hybrid_ce1}-\eqref{eq:s3_hybrid_ce2}, which have the unknown parameters $\p{1}{\B{\xi}}$, $\p{2}{\B{\xi}}$, $y_1$, and $y_{1_x}$. By substituting these equations into Equations \eqref{eq:L1ex2} and \eqref{eq:L2ex2} and taking the partials with respect to the unknown parameters, a Jacobian can be derived and ultimately used to solve the differential equations. The analytical partials that form the Jacobian are provided in Appendix \ref{sect:app_linearNonlinearDE}.

\begin{example}{Results of linear-nonlinear differential equation sequence}
Just like all other problems, this system can be solved using an iterative least-squares approach. However, an initial guess must be provided for the iterative least-squares. In the case of BVPs using the TFC method, the initial parameters can be determined by connecting the boundary constraints using a straight line. The line initial guess is adopted here in all the hybrid numerical tests provided. Therefore, the initial estimate of $y_1$ and $y_{1_x}$ is automatically determined by this initialization. For this problem,
\begin{equation*}
    \Xi_0 = \begin{Bmatrix} \B{0}\T & \B{0}\T & \frac{y(\pi)-y(0)}{2} + y(0) & \frac{y(\pi)-y(0)}{\pi}  \end{Bmatrix}\T.
\end{equation*}
A visualization of this initial guess compared to the true solution is provided in Figure (\ref{fig:ex2_init}).
\begin{figure}[H]
    \centering\includegraphics[width=.65\linewidth]{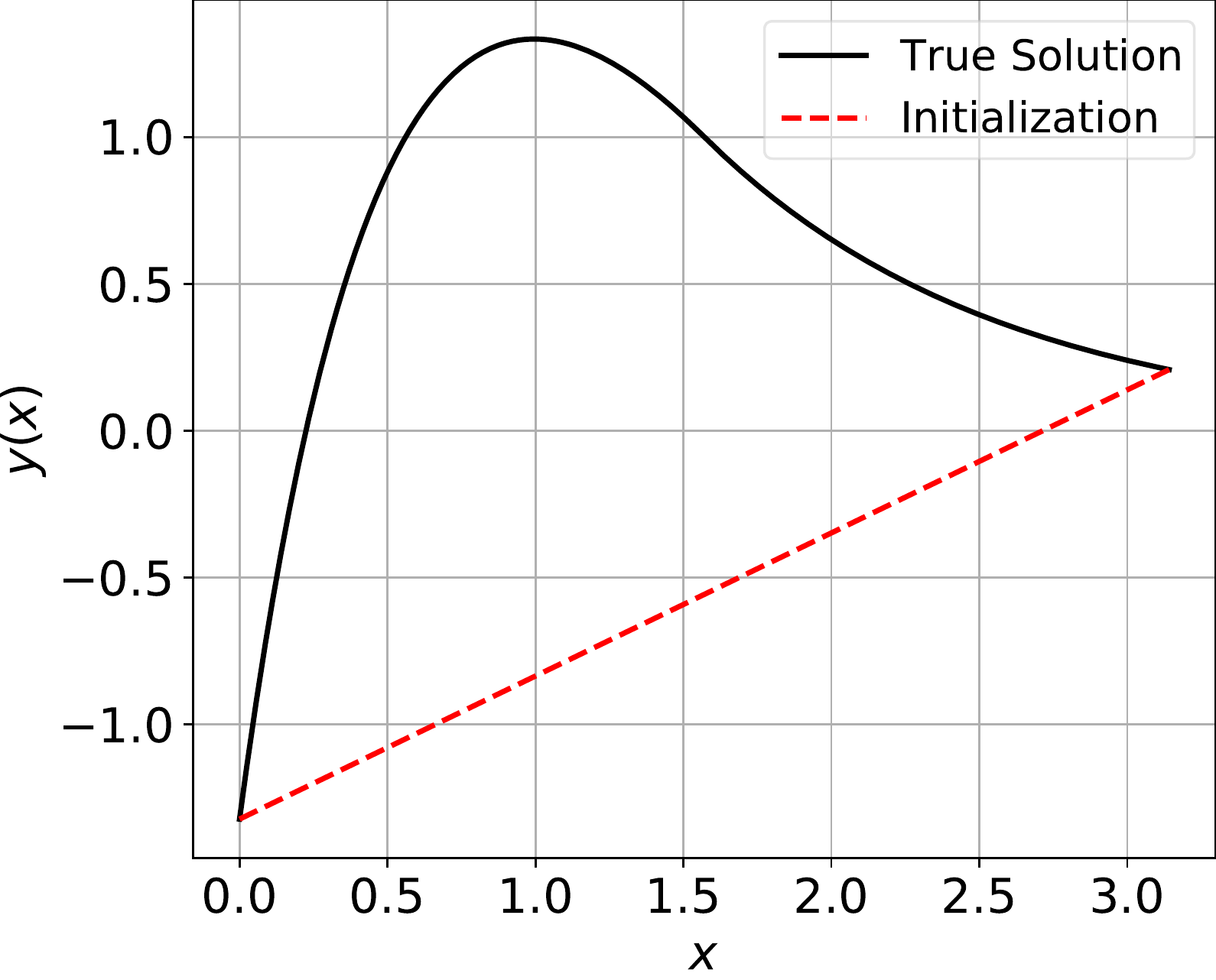}
    \caption{Initial guess and true solution for the linear-nonlinear sequence. Reprinted with permission from \cite{pieceTFC}.}
    \label{fig:ex2_init}
\end{figure} 
For the DE presented in Equation \eqref{eq:linear_nonlinear}, $N = 100$ and
$m = 16$ basis functions were used for each segment. The solution reached machine error accuracy in $15$ iterations. The results of this numerical test are shown in Figures \ref{fig:ex2:function} and \ref{fig:ex2:error}. The results show the function, its first two derivatives, and the associated absolute errors compared to the analytical solution.
\begin{figure}[H]
    \centering\includegraphics[width=0.75\linewidth]{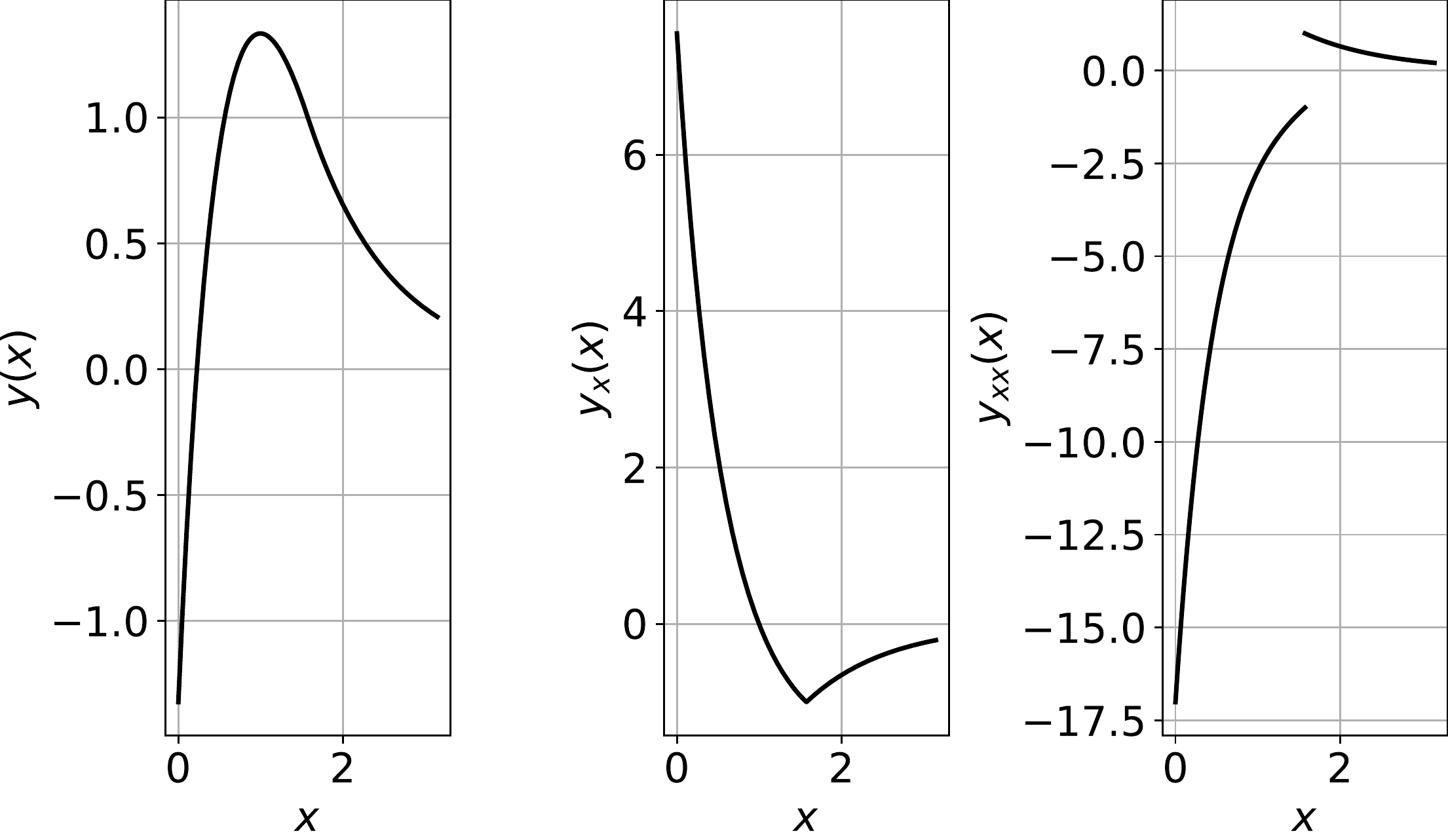}
    \caption{Solution of linear-nonlinear differential equation sequence. Reprinted with permission from \cite{pieceTFC}.}
    \label{fig:ex2:function}
\end{figure} 

\begin{figure}[H]
    \centering\includegraphics[width=0.7\linewidth]{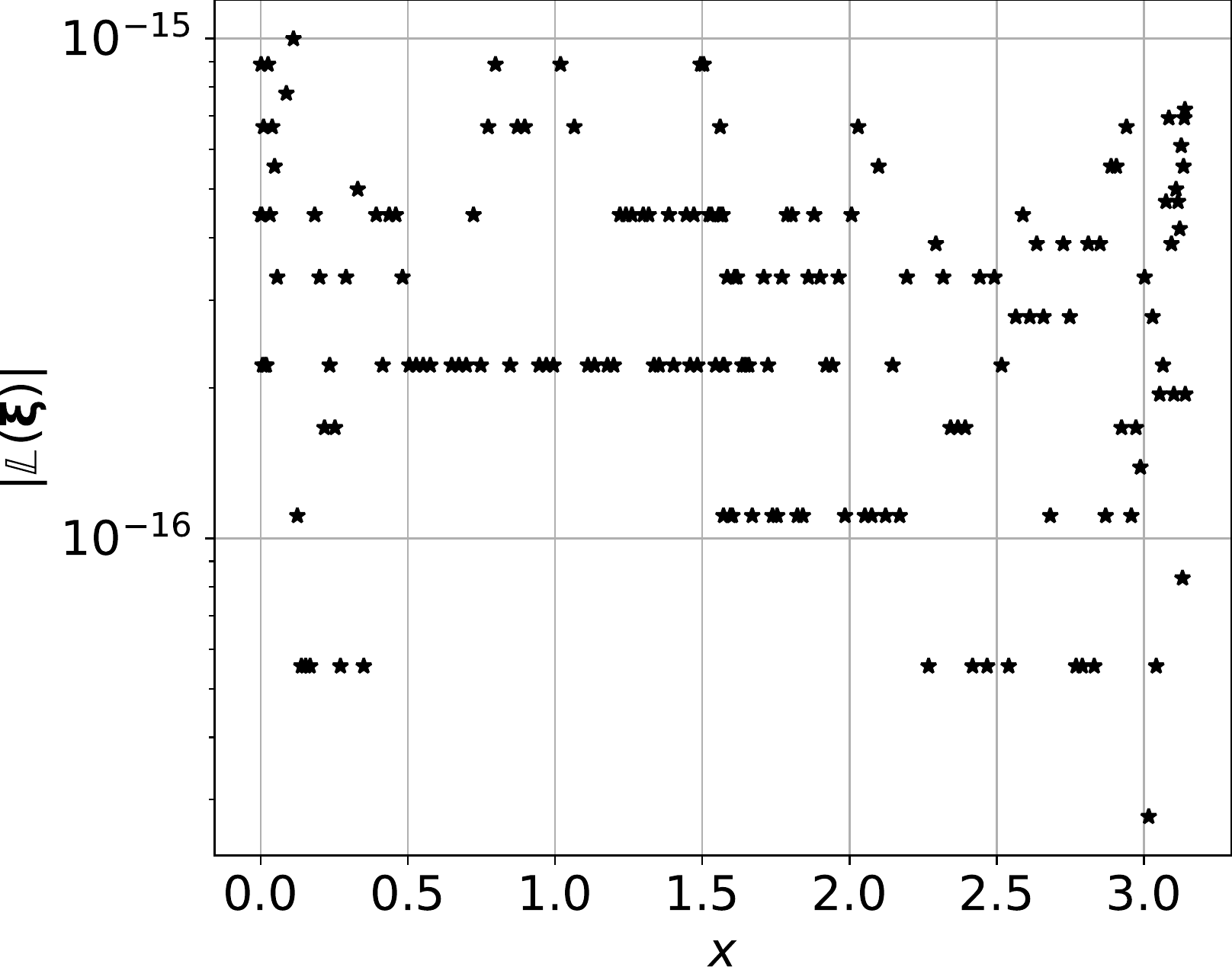}
    \caption{Absolute error of solution of linear-nonlinear differential equation sequence. Reprinted with permission from \cite{pieceTFC}.}
    \label{fig:ex2:error}
\end{figure} 

\end{example}

\subsubsection{1D convection-diffusion equation}\label{sec:s3_condif}
\begin{figure}[ht]
    \centering\includegraphics[width=0.75\linewidth]{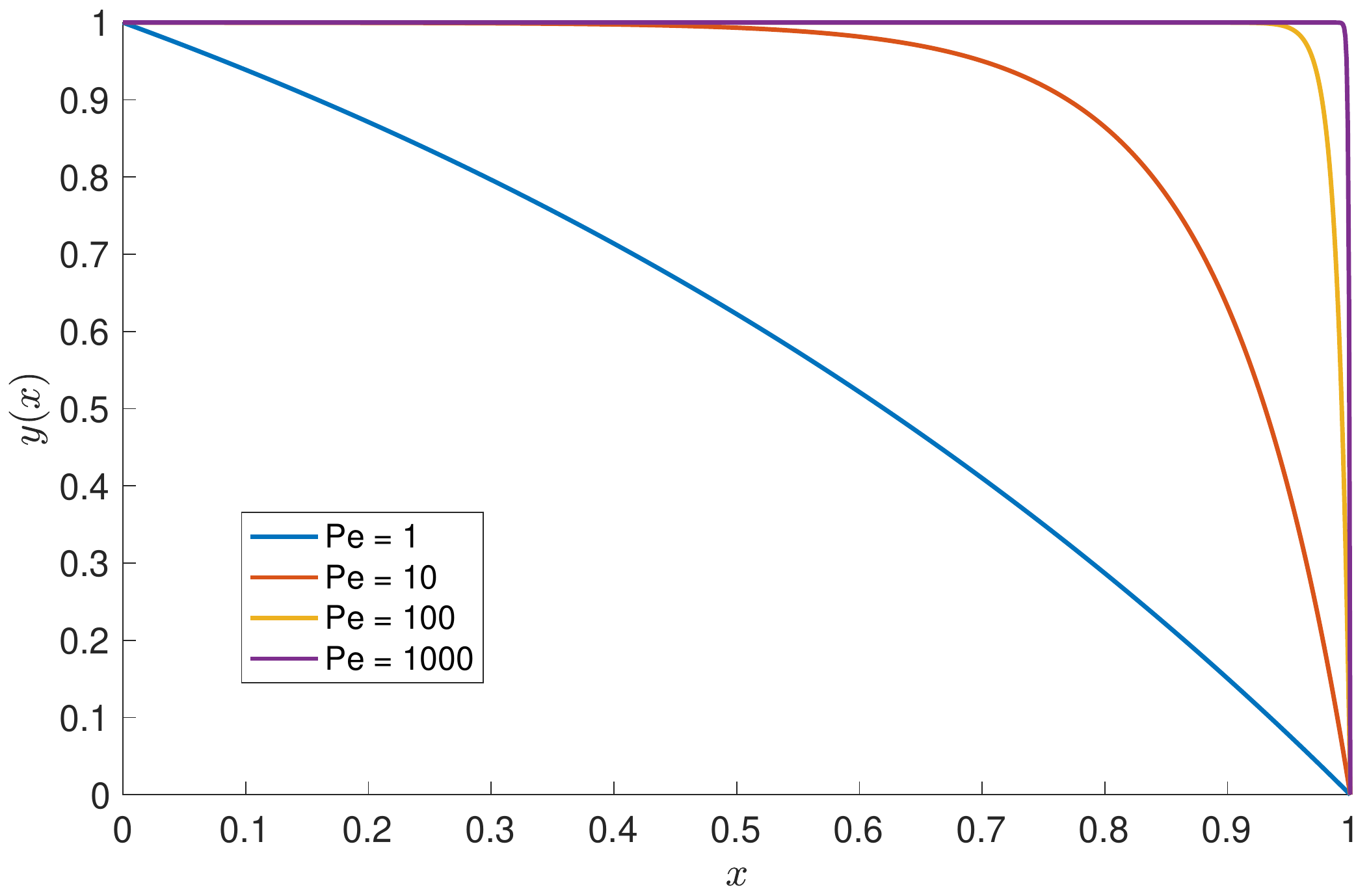}
    \caption{Solution of the 1D convection-diffusion equation for varying values of the Peclet number. As the Peclet number increases, the solution exhibits sharp transient behavior close to the endpoint. Reprinted with permission from \cite{pieceTFC}.}
    \label{fig:s3_condif}
\end{figure}

This technique doesn't just apply to hybrid systems. In fact, the concept of splitting the problem domain can be utilized when the dynamics exhibit transient behavior. To further highlight this concept, consider the example of the one-dimensional convection-diffusion equation defined by the differential equation,
\begin{equation}\label{eq:s3_condifDE}
    y_{xx} - \text{Pe} \, y_x = 0 \quad \text{subject to: } \begin{cases} y(0) = 1 \\ y(1) = 0 \end{cases}
\end{equation}
with analytical solution
\begin{equation*}
    y = \frac{1 - e^{\text{Pe}(x-1)}}{1 - e^{-\text{Pe}}}.
\end{equation*}
In these equations, Pe is the Peclet number defined by the equation,
\begin{equation*}
    \text{Pe} = \frac{uL}{k} = \text{RePr} \approx \frac{\text{heat transported}}{\text{heat conducted}},
\end{equation*}
where $u$ is the fluid velocity, $L$ is the characteristic length, $k$ is the thermal diffusivity of the fluid, Re is the Reynolds number, and Pr is the Prandtl number. We are interested in the behavior of the solution as the Peclet number increases, as shown in Figure \ref{fig:s3_condif}.
As Pe increases to around 100, the function begins to have a sharp transient behavior near the end of the domain. In order to solve this problem, let us consider a TFC solution where the domain is split into two segments such that the switch is defined at some value $x_1 \in (0, 1)$. The constrained expression follow as,
\begin{align*}
    \p{1}{y}(z,\p{1}{\B{\xi}}) = \Big(\B{h}(z) - \p{1}{\phi}_1(z) \B{h}(z_0) &- \p{1}{\phi}_2(z) \B{h}(z_f) - \p{1}{\phi}_3(z) \B{h}_z(z_f) \Big)\T \p{1}{\B{\xi}} \\&+ \p{1}{\phi}_1(z) y_0 + \p{1}{\phi}_2(z) y_1 + \p{1}{\phi}_3(z) \frac{y_{1_x}}{\p{1}{c}}
\end{align*}
\begin{align*}
    \p{2}{y}(z,\p{2}{\B{\xi}}) = \Big(\B{h}(z) - \p{2}{\phi}_1(z) \B{h}(z_0) &- \p{2}{\phi}_2(z) \B{h}_z(z_0) - \p{2}{\phi}_3(z) \B{h}(z_f) \Big)\T \p{2}{\B{\xi}} \\&+ \p{2}{\phi}_1(z) y_1 + \p{2}{\phi}_2(z) \frac{y_{1_x}}{\p{2}{c}} + \p{2}{\phi}_3(z) y_f
\end{align*}
where the segment domains are defined in the basis domain $z$ and $\p{1}{x} \in [0, x_1] \rightarrow \p{1}{z} \in [z_0, z_f]$ and $\p{2}{x} \in [x_1, x_f] \rightarrow \p{2}{z} \in [z_0, z_f]$. Since we have written the constrained expression in the basis domain, the derivative constraints must be divided by the mapping coefficient $\p{1}{c}$ and $\p{2}{c}$ to account for this. Next, the switching functions are defined as,
\begin{align*}
    \p{1}{\phi}_1(z) &= \frac{1}{(z_f - z_0)^2}\Big( z_f^2 - 2 z_f z + z^2 \Big)\\
    \p{1}{\phi}_2(z) &= \frac{1}{(z_f - z_0)^2}\Big(z_0(z_0-2z_f) + 2z_f z - z^2\Big)\\
    \p{1}{\phi}_3(z) &= \frac{1}{z_f - z_0}\Big( z_0 z_f - (z_0 + z_f)z + z^2\Big)\end{align*}
for the first segment's constrained expression, and as
\begin{align*}
    \p{2}{\phi}_1(z) &= \frac{1}{(z_f - z_0)^2} \Big(z_f(z_f-2z_0) + 2 z_0 z - z^2 \Big)\\
    \p{2}{\phi}_2(z) &= \frac{1}{z_f - z_0} \Big(-z_f z_0 + (z_f + z_0) z - z^2 \Big)\\
    \p{2}{\phi}_3(z) &= \frac{1}{(z_f - z_0)^2} \Big(z_0^2 - 2 z_0 z + z^2 \Big)
\end{align*}
for the first segment's constrained expression. Now, we can rewrite the differential equation given by Equation \eqref{eq:s3_condifDE} as,
\begin{equation*}
    \begin{cases} \text{for } x \leq x_1 : \p{1}{\tilde{F}} = \p{1}{c}^2 \p{1}{y}_{xx} - \text{Pe} \p{1}{c} \p{1}{y}_x = 0 \\ \text{for } x \geq x_1 : \p{2}{\tilde{F}} = \p{2}{c}^2 \p{2}{y}_{xx} - \text{Pe} \p{2}{c} \p{2}{y}_x = 0 \end{cases}.
\end{equation*}

To solve this differential equation, we could simply select the value of $x_1$ based on intuition and proceed with the same process as described earlier. However, it is highly likely that the selected value of $x_1$ will not be optimal and should therefore be a value that is optimized. Two methods exist to determine this value. The first method involves combining the TFC approach with an outer-loop optimizer (i.e., \verb"fsolve", a genetic algorithm, etc.) to solve for $x_1$. In this method, the TFC method supplies the estimated solution accuracy through the differential equation residuals, and the outer-loop optimizes the value of $x_1$ to minimize the residual. The second method is to include the solution of $x_1$ inside the TFC method. This can be realized by a single coefficient, since $\p{1}{c}$ and $\p{2}{c}$ are connected through the value $x_1$ by the equations
\begin{align*}
    \p{1}{c} &= \frac{z_f - z_0}{x_1 - x_0} \\ 
    \p{2}{c} &= \frac{z_f - z_0}{x_f - x_1} 
\end{align*}
and $\p{2}{c}$ can be rewritten in terms of $\p{1}{c} := \bar{c}$
\begin{equation*}
    \p{2}{c} = \frac{\bar{c}(z_f - z_0)}{\bar{c}(x_f - x_0) - z_f + z_0} = \frac{\bar{c} \Delta z}{\bar{c} - \Delta z}
\end{equation*}
and 
\begin{equation*}
    \frac{\partial \p{2}{c}}{\partial \bar{c}}= -\frac{\Delta z^2}{(\bar{c} - \Delta z)^2}.
\end{equation*}
This reduces the mapping coefficient to a single parameter that can be plugged into the constrained expressions and differential equation; however, doing so forces the system of equations to be nonlinear. The loss functions become, 
\begin{equation*}
    \mathbb{\bar{L}}(\Xi) = \begin{bmatrix} \p{1}{\mathbb{L}}(\Xi) & \p{2}{\mathbb{L}}(\Xi) \end{bmatrix}\T
\end{equation*}
with the unknown vector defined as,
\begin{equation*}
    \Xi = \begin{bmatrix} \p{1}{\B{\xi}} & \p{2}{\B{\xi}} & y_1 & y_{1_x} & \bar{c} \end{bmatrix}.
\end{equation*}
Additionally, the terms of the total loss vector and Jacobian are provided in Appendix \ref{sect:app_condif} for completeness. 

\begin{wblankBox}{Adaptation for other numerical techniques}
While the equations above show the nonlinear least-squares approach to solve the problem, the equations can easily be adapted where an external optimizer handles the estimation of the optimal $x_1$ location ($\bar{c}$ in the above equations). By removing the unknown $\bar{c}$ from the equations, we are left with a linear set of equations defined by the same loss function and the updated unknown vector,
\begin{equation*}
    \Xi = \begin{bmatrix} \p{1}{\B{\xi}} & \p{2}{\B{\xi}} & y_1 & y_{1_x} \end{bmatrix}
\end{equation*}
By defining this TFC method where the input is $x_1$ (or $\bar{c}$) and the output is some function of the loss vector (in this case we use $\max|\mathbb{L}(\Xi)|$), a suite of optimizers can be leveraged.
\end{wblankBox}

\begin{example}{Results of the 1D convection-diffusion equation}

For the numerical solution of the 1D convection-diffusion equation three methods where used: 1) a nonlinear least-squares (NLS) approach, and two approaches relying on the adaptation discuss earlier, 2) a differential evolution algorithm (DEvo) utilizing SciPy's \verb"optimize.differential_evolution()" and 3) SciPy's \verb"optimize.fsolve()" algorithm. In all cases, $N=200$ discretization points were used per segment and the basis functions were taken to the $190^{\text{th}}$ degree term ($m = 187$ basis functions). Additionally, the problem was solved for a range of Peclet numbers from $10^2$ to $10^6$ with a convergence criteria of $\varepsilon = 1 \times 10^{-13}$ and an initial guess of $x_1 = 0.75$. The results are captured in Table \ref{tab:connDiff}. For numerical stability of the algorithms, the upper bound of $x_1$ in the NLS and \verb"fsolve()" approaches was set to be 0.9990 while the DEvo was set to 0.9999990.

In general, it can be seen that the algorithms have similar maximum errors; however, the location of the estimated $x_1$ value differs considerably, and the computation time of the NLS approach is two orders of magnitude faster than the differential evolution algorithm. The difference in $x_1$ is because $x_1$ is a numerical construct based on solving the differential equation. This value does not show up naturally in the equations, and therefore there is a potential of many local minima. This is very evident in the solution of the problem for $Pe = 10^2$ and $10^3$ where the \verb"fsolve" algorithm simply chooses the initial guess as the best solution, yet has similar accuracy to the other two methods.

\vspace{0.1in}
\begin{table}[H]
\caption{Solution for convention-diffusion equation using traditional TFC with nonlinear least-squares and with a genetic algorithm to solve for $x_1$ over a span of Peclet numbers. In all test cases, the number of points was $N=200$ for each segment and the basis functions were taken to the $190^{\text{th}}$ degree term ($m = 187$ basis functions). Reprinted with permission from \cite{pieceTFC}.} 
\centering
\begin{tabular}{|c|c|c|c|c|c|}
\hline
Type & Pe & $\max|\text{Error}|$ & $\max|L(\Xi)|$ & $x_1$ & \makecell{Computation\\time [s]} \\\hline\hline
NLS &  $10^2$ & $5.13\times 10^{-15}$ & $7.28\times 10^{-12}$ & 0.91127 & 1.10\\\hline 
NLS &  $10^3$ & $5.36\times 10^{-14}$ & $4.66\times 10^{-10}$ & 0.91827 & 0.91\\\hline 
NLS &  $10^4$ & $4.97\times 10^{-13}$ & $5.96\times 10^{-8}$ & 0.99000 & 0.53\\\hline 
NLS &  $10^5$ & $4.22\times 10^{-12}$ & $7.63\times 10^{-6}$ & 0.99900 & 0.91\\\hline 
NLS &  $10^6$ & $3.10\times 10^{-11}$ & $6.10\times 10^{-4}$ & 0.99900 & 3.82\\\hline
\hline
DEvo  &  $10^2$ & $5.53\times 10^{-15}$ & $4.15\times 10^{-12}$ & 0.98374 & 9.25\\\hline 
DEvo  &  $10^3$ & $4.46\times 10^{-14}$ & $2.95\times 10^{-10}$ & 0.87589 & 10.88\\\hline 
DEvo  &  $10^4$ & $1.65\times 10^{-13}$ & $2.20\times 10^{-8}$ & 0.90771 & 11.10\\\hline 
DEvo  &  $10^5$ & $3.307\times 10^{-12}$ & $3.53\times 10^{-6}$ & 0.99838 & 10.50\\\hline 
DEvo  &  $10^6$ & $3.94\times 10^{-11}$ & $2.66\times 10^{-4}$ & 0.99945 & 9.72\\\hline
\hline
\verb"fsolve"  &  $10^2$ & $4.88\times 10^{-15}$ & $4.81\times 10^{-12}$ & 0.75000 & 2.11\\\hline 
\verb"fsolve"  &  $10^3$ & $4.71\times 10^{-14}$ & $2.60\times 10^{-10}$ & 0.75000 & 1.54\\\hline 
\verb"fsolve"  &  $10^4$ & $3.68\times 10^{-13}$ & $2.09\times 10^{-8}$ & 0.92199 & 4.51\\\hline 
\verb"fsolve"  &  $10^5$ & $4.21\times 10^{-12}$ & $4.10\times 10^{-6}$ & 0.99900 & 1.54\\\hline 
\verb"fsolve"  &  $10^6$ & $3.11\times 10^{-11}$ & $5.79\times 10^{-4}$ & 0.99900 & 1.75\\\hline
\hline 
\end{tabular}
\label{tab:connDiff}
\end{table}
\end{example}

\subsection{Dealing with unspecified time and nonlinear constraints}\label{sec:s3_freeTime}

Suppose we are faced with a problem that involves solving a differential equation subject to both linear and nonlinear boundary constraints along with an unknown final time. These conditions are typical of optimal control problems; therefore, let us consider a simple controls problem,
\begin{align*}
    \dot{x} &= \alpha x + \beta u \\
    \dot{u} &= \beta x  - \alpha u
\end{align*}
subject to $ x(0) = x_0$, $x(t_f) = x_f$, where $t_f$ is unknown\footnote{Note, this system of equation is derived from the optimal control problem $\text{min} \, J = \int_0^{t_f} \frac{1}{2}(x^2 + u^2)\dd t$ subject to the dynamics $\dot{x} = \alpha x + \beta u$ constrained such that $x(0) = x_0$ and $x(t_f) = x_f$.}. Additionally, the system must satisfy the algebraic constraint at the final time
\begin{equation*}
    \frac{1}{2} \Big(x^2(t_f) - u^2(t_f)\Big) - \frac{\alpha}{\beta} x(t_f)u(t_f) = 0.
\end{equation*}
Now, since the final time, $t_f$, is unknown, let us write the entire problem in the basis function domain $z$ and map to the problem domain $t$ using the parameter $c$ from Equation \eqref{eq:s3_linearMapping}. Therefore, the system of equations to be solved becomes,
\begin{align}
    F^x &= c \, x_z -\alpha x - \beta u = 0 \label{eq:s3_Fx}\\
    F^u &= c \, u_z - \beta x + \alpha u  = 0 \label{eq:s3_Fu}\\
    f &= \frac{1}{2} \Big(x^2(t_f) - u^2(t_f)\Big) - \frac{\alpha}{\beta} x(t_f)u(t_f) = 0 \label{eq:s3_f}
\end{align}
Now, we use the developed method, but we write all constrained expressions in the $z$ domain such that,
\begin{align*}
    x(z,g^x(z)) &= g^x(z) + \frac{z_f - z}{z_f - z_0} \Big(x_0 - g^x(z_0)\Big) + \frac{z - z_0}{z_f - z_0} \Big(x_f - g^x(z_f)\Big) \\
    u(z,g^u(z)) &= g^u(z)
\end{align*}
where the function of $u(z)$ has no linear constraints and becomes solely a function of the free function, $g_u(z)$. By discretizing Equations \eqref{eq:s3_Fx}, \eqref{eq:s3_Fu}, and \eqref{eq:s3_f}, we can construct our typical loss vectors for each function,
\begin{align*}
    \mathbb{L}^x &= \begin{Bmatrix} F^x(z_0,\B{\xi}_x, \B{\xi}_u, c) & \hdots & F^x(z_k,\B{\xi}_x, \B{\xi}_u, c) & \hdots &  F^x(z_f,\B{\xi}_x, \B{\xi}_u, c) \end{Bmatrix}\T \\
    \mathbb{L}^u &= \begin{Bmatrix} F^u(z_0,\B{\xi}_x, \B{\xi}_u, c) & \hdots & F^u(z_k,\B{\xi}_x, \B{\xi}_u, c) & \hdots &  F^u(z_f,\B{\xi}_x, \B{\xi}_u, c) \end{Bmatrix}\T \\
    \mathbb{L}^f &= f(z_f, \B{\xi}_x, \B{\xi}_u)
\end{align*}
which is collected in a total loss vector,
\begin{equation*}
    \mathbb{L} = \begin{Bmatrix}\mathbb{L}^x \phantom{}\T & \mathbb{L}^y \phantom{}\T & \mathbb{L}^f \end{Bmatrix}\T.
\end{equation*}
In this problem, not only are the coefficients $\B{\xi}_x$ and $\B{\xi}_u$ unknowns, but the final time is also unknown, which is captured in our mapping parameter $c$, such that, $\Xi = \{ \B{\xi}_x, \, \B{\xi}_u, \, c\}\T$. Therefore, our Jacobian will also be populated by partial derivatives with respect to $\Xi$. The derivation of the Jacobian is left to the reader. 

We rely again on the nonlinear least-squares approach to solve the problem since the final equation is nonlinear in the variables $x$ and $u$. Yet, note that $c$ defines a domain length and can never be negative. Therefore, let us change the definition for this variable such that $b^2 := c$. By doing this, we avoid the time domain parameter becoming negative, and the vector of unknowns becomes $\Xi = \{ \B{\xi}_x, \, \B{\xi}_u, \, b\}\T$. In summary, this simply changes Equations \eqref{eq:s3_Fx} and \eqref{eq:s3_Fu} to,
\begin{align*}
    F^x &= b^2 \, x_z -\alpha x - \beta u = 0 \\
    F^u &= b^2 \, u_z  - \beta x + \alpha u  = 0,
\end{align*}
in the development provided above. 

\begin{wblankBox}{Adaptation for other numerical techniques}
\noindent Additionally, similar to our solution of the convection-diffusion equation in Section \ref{sec:s3_condif}, we can remove the unknown value of $t_f$ (which is related to $b^2$), creating a linear system of equations to be solved in $\mathbb{L}^x$ and $\mathbb{L}^u$. After solving the system, $|\mathbb{L}^f|$ can be used as the function to be minimized.
\end{wblankBox}

\begin{example}{Solution to free-final time problem}
In this example, we have defined the coefficients as $\alpha = \beta = 1$, with the boundary conditions set as $x(0) = 1$ and $x(t_f) = 1$. Furthermore, all numerical systems were discretized with $N = 35$ points and used basis function up to the $30^{\text{th}}$ degree term (28 basis functions for $x(t)$ and 30 basis functions for $y(t)$). Lastly, the tolerance on the algorithms was set to $\varepsilon = 2.22 \times 10^{-16}$ and were initialized with $\B{\xi}_x = \B{0}$, $\B{\xi}_u = \B{0}$, and $t_f = 1$. For reference, the solution of $x(t)$ and $y(t)$ is highlighted in Figure \ref{fig:s3_unknownTimeStates}.
\begin{figure}[H]
    \centering\includegraphics[width=0.75\linewidth]{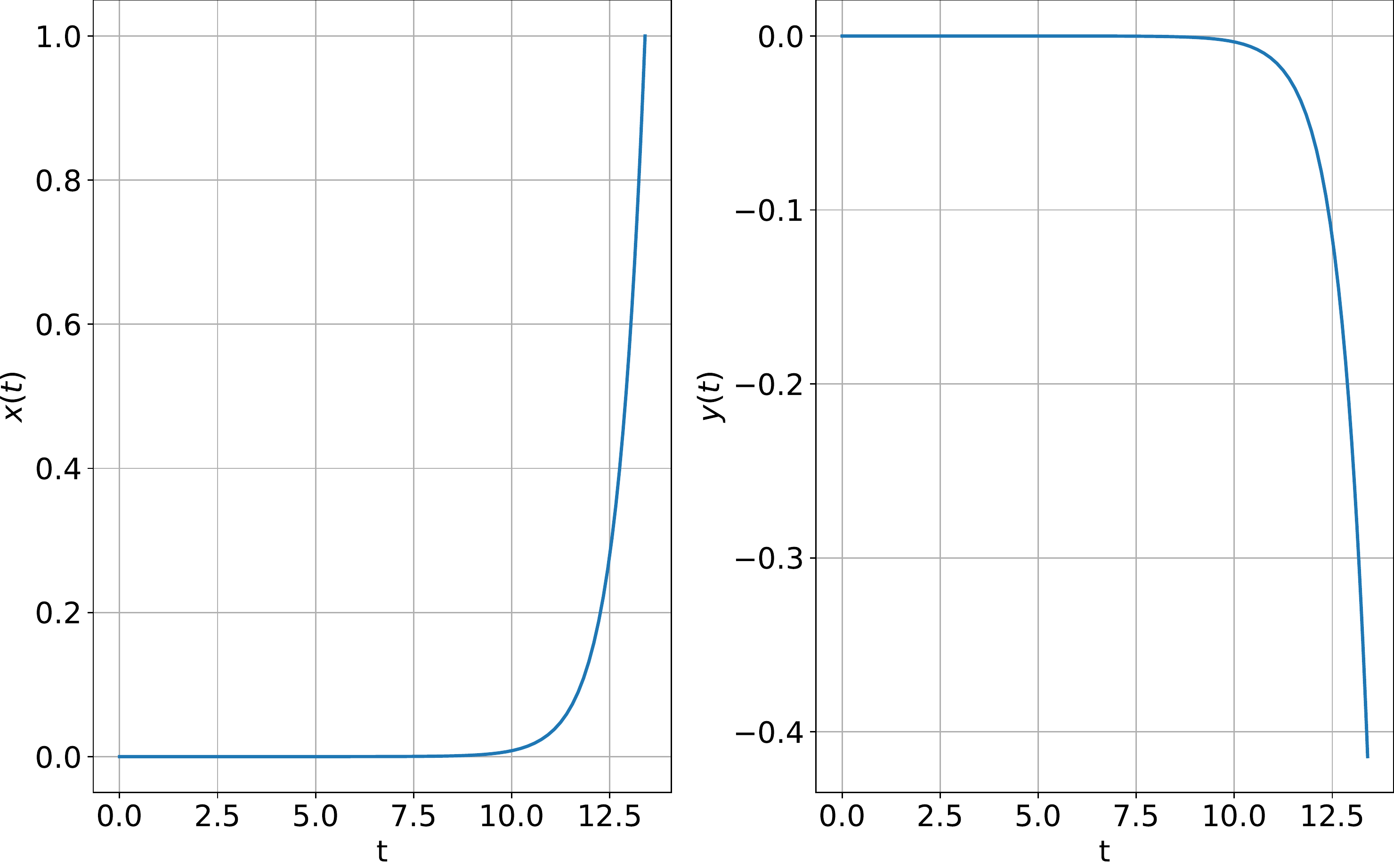}
    \caption{Time histories of the state.}
    \label{fig:s3_unknownTimeStates}
\end{figure} 
The results of this test are provided in Table \ref{tab:freeTime} where it can be seen that the $\verb"fsolve"$ method is the most accurate in terms of $\max|\mathbb{L}(\Xi)|$ and $\max|H(t)|$, which is the Hamiltonian\footnote{In this example the Hamiltonian is simply stated without definition. In Chapter \ref{chap:opt_con}, this term will be defined more rigorously.}, parameter used to derive the problem and should be zero for all times. However, we can see that all methods differ at the sixth digit of the cost, which is defined as,
\begin{equation*}
    \text{Cost} = \frac{1}{2} \int_0^{t_f} \Big(x^2(t) + u^2(t)\Big) \dd t,
\end{equation*}
and should be minimized in our case. Yet, where these solutions drastically differ is in the solution time, where the NLS approach is two orders of magnitude faster than the other approaches. This should be obvious since the NLS is the simplest approach to solving the problem. Additionally, since the cost isn't as sensitive to the final time, we see a large range of solutions for $t_f$. In all, the major consideration becomes a trade-off between solution accuracy versus computational time. In Chapters \ref{chap:eol} and \ref{chap:fol}, we will take a deeper look into this regarding optimal control problems.
\begin{table}[H]
\caption{Comparison of optimization scheme to solve the free final time problems.} 
\centering
\begin{tabular}{|c|c|c|c|c|c|c|}
\hline
Type & $\max|L(\Xi)|$ & $\max|H(t)|$ & Cost & $t_f$ & Iterations & \makecell{Comp.\\Time [s]} \\\hline\hline
NLS           & $8.36\times 10^{-14}$ & $8.36\times 10^{-14}$ & $0.206\B{9}\B{1}$ & 11.13663 & 23 & 0.0457\\\hline 
DEvo          & $9.99\times 10^{-16}$ & $5.55\times 10^{-17}$ & $0.206\B{7}\B{8}$ & 13.92129 & 77 & 2.879\\\hline 
\verb"fsolve" & $5.18\times 10^{-16}$ & $2.17\times 10^{-16}$ & $0.206\B{8}\B{2}$ & 13.24100 & 62 & 2.520\\\hline 
\hline
\end{tabular}
\label{tab:freeTime}
\end{table}
\end{example}

\section{A Solution of Lyapunov and Halo Orbits}
According to Poincar\'e, ``periodic orbits'' provide the only gateway into the otherwise impenetrable domain of nonlinear dynamics. With the advent of space exploration, periodic orbits have become an indispensable part of missions in space. The amazing fish-like Apollo orbit was the first three-body orbit used for space missions. The second three-body orbit used for space missions was the Halo orbit, discovered by Robert Farquhar in his Ph.D. thesis \cite{halo} under John Breakwell \cite{breakwell}. In 1978, Farquhar convinced NASA and led the International Sun-Earth Explorer 3 mission (ISEE3) to study the Sun from a Halo orbit around the Earth’s L1 Lagrange point. Farquhar’s original idea was to place a satellite in Halo orbit around the Lunar L2 for telecommunication support for the backside of the Moon. Today, this is indeed part of NASA’s planned return of humans to the Moon in the next few years.

Typically, the standard method for computing periodic orbits is the differential correction method (also called the shooting method), as presented by Kathleen Howell \cite{diff_cor}. One begins with an approximate solution obtained typically from normal form expansions. Using the variational equation, the guess solution is iteratively corrected for periodicity. Assuming the initial guess is in a reasonable basin of attraction to a periodic orbit, the process converges to a periodic orbit. In Hamiltonian systems, periodic orbits occur in 1-parameter families. Often, there are multiple families nearby. Hence, the convergence may not always lead to the desired orbit. Moreover, control over the specific features of the periodic orbit, such as its period or energy, requires additional work, for example, using continuation methods to reach the exact orbit desired. Using TFC, a simpler formulation and more efficient algorithm for finding periodic orbits is possible.

\subsection{System dynamics}
\begin{figure}[ht]
	\centering
	\includegraphics[width=0.7\linewidth]{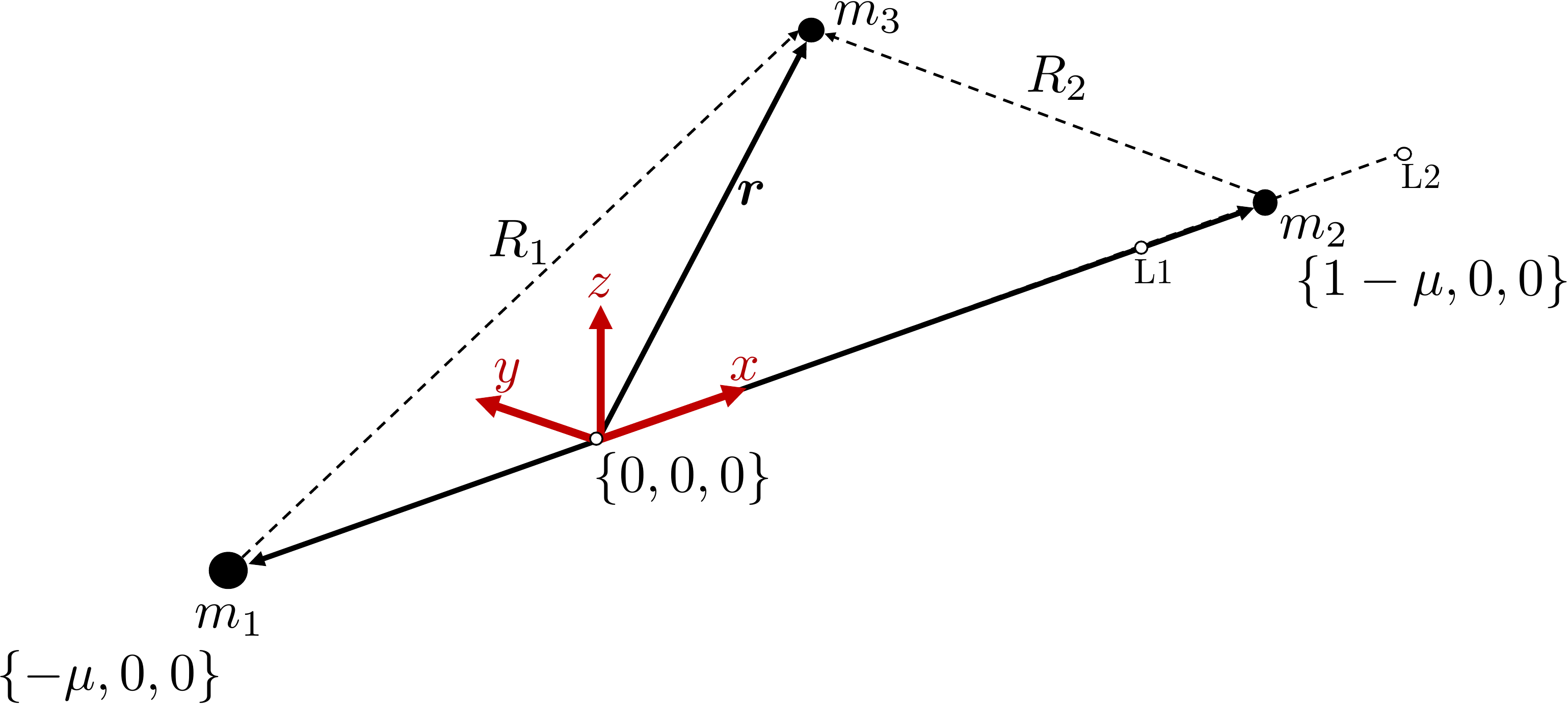}
	\caption{Schematic of the circular restricted three-body problem where the secondary body $m_2$ orbits around $m_1$ in a circular orbit. The third-body whose mass is $m_3 \ll m_2 < m_1$ is negligible and at a distance $R_1$ from $m_1$, $R_2$ from $m_2$, and $\B{r}$ from the origin, which is the system barycenter (the system's center of mass).}
	\label{fig:threeBP}
\end{figure}

The circular-restricted three-body problem is a dynamical model used to describe the motion of a particle $\B{r} = \{x, y, z\}\T$ of negligible mass under the influence of a primary body of mass $m_1$ and secondary body of mass $m_2$. Furthermore, the orbits of $m_1$ and $m_2$ are subject to circular motion about the system's barycenter and lie in the $x$-$y$ plane; the total system is depicted in Figure \ref{fig:threeBP}. Following this, the system can be non-dimensionalized by the following scaled units; unit mass is defined as $m_1 + m_2$; unit length is taken as the separation between $m_1$ and $m_2$; the unit time is chosen such that the orbits of $m_1$ and $m_2$ about the system's barycenter is $2 \pi$. By following these steps, the system can be reduced to a single parameter called the mass parameter, $\mu$, where,
\begin{equation*}
    \mu = \dfrac{m_2}{m_1+m_2}
\end{equation*}
From this, we define the terms $\mu_1$ and $\mu_2$ as
\begin{equation*}
    \mu_1 = 1 - \mu \quad \text{and} \quad \mu_2 = \mu.
\end{equation*}
Using this definition of the system, the equations of motion can be derived in the rotating frame leading to the following system of equations,
\begin{align}
    \ddot{x} - 2 \dot{y} = \dfrac{\partial \Omega}{\partial x} \nonumber\\
    \ddot{y} + 2 \dot{x} = \dfrac{\partial \Omega}{\partial y} \label{eq:stateDE}\\\nonumber 
    \ddot{z} = \dfrac{\partial \Omega}{\partial z}
\end{align}
Additionally, $\Omega$ is defined as, 
\begin{align*}
    \Omega(x,y,z):= \dfrac{1}{2}(x^2 + y^2) + \dfrac{1-\mu}{R_1} + \dfrac{\mu}{R_2} + \dfrac{1}{2}(1-\mu)\mu
\end{align*}
where $R_1 = \sqrt{(x + \mu)^2 + y^2 + z^2}$ and $R_2 = \sqrt{(x + \mu - 1)^2 + y^2 + z^2}$ are the distances to the primaries. Furthermore, the equations of motion are Hamiltonian and independent of time, and thus have an energy integral of motion $E$, where in the celestial mechanics community the Jacobi constant is used which is $C:=-2E$ and given as,
\begin{equation}\label{eq:jacobiConst}
    C = 2 \Omega - (\dot{x} + \dot{y} + \dot{z}) = (x^2 + y^2) + 2\dfrac{1-\mu}{R_1} + 2\dfrac{\mu}{R_2} + (1-\mu)\mu - (\dot{x} + \dot{y} + \dot{z})
\end{equation}
Moving forward, we will look to solve the dynamics defined by the system of equations in Equation \eqref{eq:stateDE} such that the orbit is at a fixed energy level (or rather Jacobi constant) using Equation \eqref{eq:jacobiConst}. For our implementation, it is useful to define the residuals of these equations,
\begin{align}
    0 = F_x &:=  \ddot{x} - 2 \dot{y} - \dfrac{\partial \Omega}{\partial x} \label{eq:resX}\\
    0 = F_y &:=  \ddot{y} + 2 \dot{x} - \dfrac{\partial \Omega}{\partial y} \label{eq:resY}\\
    0 = F_z &:= \ddot{z} - \dfrac{\partial \Omega}{\partial z} \label{eq:resZ}\\
    0 = F_c &:= (x^2 + y^2) + 2\dfrac{1-\mu}{R_1} + 2\dfrac{\mu}{R_2} + (1-\mu)\mu - (\dot{x} + \dot{y} + \dot{z})  - C \label{eq:resJc}. 
\end{align}
Next, we look to generate analytical expressions for the states to guarantee a periodic orbit. 

First, since in the problem the orbital period is unknown, the problem represents an unknown final time problem where we can define the problem domain as $t \in [0, T]$ where $T$ is the period of the orbit and the basis domain is $\tau \in [-1, +1]$.\footnote{Note, we have used $\tau$ here in place of $z$ because in the common notation for this problem, $z$ represents the $z$-component of the position of the body $m_3$.} The final time (or the orbital period $T$) can be parameterized in the same manner as Section \ref{sec:s3_freeTime}.

Since we are looking for periodic orbits, we can utilize the \ce\ to satisfy the following constraints,
\begin{align*}
    r_i(\tau_0) = r_i(\tau_f) = \alpha_i \quad \text{and} \quad \frac{\dd r_i}{\dd \tau}(\tau_0) = \frac{\dd r_i}{\dd \tau}(\tau_f) = \dfrac{\beta_i}{b^2} 
\end{align*}
where we define $\B{r}(\tau) := \{r_x(\tau), r_y(\tau), r_z(\tau)\}\T = \{x(\tau), y(\tau), z(\tau)\}\T$. Since the trajectory must return to the initial state at some period $T$. The constrained expressions for the three components of position are as follows,
\begin{align}\label{eq:CE_periodicOrbit}
	 r_i (\tau, g_i (\tau)) = g_i (\tau)
	&+ \phi_1(\tau) \Big(\alpha_i - g(\tau_0) \Big) 
	+ \phi_2(\tau) \Big(\alpha_i - g(\tau_f) \Big)  \nonumber\\
	&+ \phi_3(\tau) \Big(\dfrac{\beta_i}{b^2}  - g_{\tau}(\tau_0) \Big) 
	+ \phi_4(\tau) \Big(\dfrac{\beta_i}{b^2}  - g_{\tau}(\tau_f) \Big) 
\end{align}
where
\begin{align*}
    &\phi_1(\tau)   = \dfrac{1}{4} \Big(2  -3\tau + \tau^3\Big) 
	 &&\phi_2(\tau) = \dfrac{1}{4} \Big(2  +3\tau - \tau^3\Big) \\
	&\phi_3(\tau)   = \dfrac{1}{4} \Big(1  - \tau - \tau^2 + \tau^3\Big) 
	 &&\phi_4(\tau) = \dfrac{1}{4} \Big(-1 - \tau  + \tau^2 + \tau^3\Big).
\end{align*}
By their definition, the projection functionals follow as,
\begin{align*}
    &\rho_1(x,g_i(x))  = \alpha_i             - g_i(\tau_0)
	&&\rho_2(x,g_i(x)) = \alpha_i             - g_i(\tau_f)\\
	&\rho_3(x,g_i(x))  = \dfrac{\beta_i}{b^2}  - g_{\tau_i}(\tau_0)
	&&\rho_4(x,g_i(x)) = \dfrac{\beta_i}{b^2}  - g_{\tau_i}(\tau_f).
\end{align*}
Then, as usual, the constrained expressions defined by Equation \eqref{eq:CE_periodicOrbit} are used to evaluate the three differential equations and one algebraic equation given in Equations \eqref{eq:resX}, \eqref{eq:resY}, \eqref{eq:resZ}, and \eqref{eq:resJc} at the discretization points, which are ultimately used to construct a loss vector of the residuals of these equations.
\begin{equation*}
	\mathbb{L}_i(\Xi) = \begin{Bmatrix} \tilde{F}_i(\tau_0, \Xi), & \hdots, & \tilde{F}_i(\tau_k, \Xi), & \hdots, & \tilde{F}_i(\tau_f, \Xi)  \end{Bmatrix}\T = \B{0}\T_{N \times 1}
\end{equation*}
with the total loss vector of 
\begin{equation*}
	\mathbb{L}(\Xi) = \begin{Bmatrix} \mathbb{L}_x \T(\Xi), & \mathbb{L}_y \T(\Xi), & \mathbb{L}_z \T(\Xi), & \mathbb{L}_{c} \T(\Xi) \end{Bmatrix}\T  = \B{0}\T_{4N \times 1}
\end{equation*}
where the unknown vector is defined as,
\begin{equation*}
	\Xi = \begin{Bmatrix} \B{\xi}_x\T, & \B{\xi}_y\T, & \B{\xi}_z\T, & \B{\alpha}\T, & \B{\beta}\T, & b \end{Bmatrix}\T  = \B{0}\T_{(3m+7) \times 1}.
\end{equation*}

\subsection{Numerical Test}
We consider the Earth-Moon system with the parameters given in Table \ref{tab:paramLyap}.
\begin{table}[ht]
    \centering
    \caption{Earth-Moon system parameters}
    \begin{tabular}{cc}
        \toprule
        Variable & Value \\ \midrule
        Earth mass $m_1$ [kg] & $5.9724 \times 10^{24}$ \\
        Moon  mass $m_2$ [kg] & $7.346 \times 10^{22}$\\ \bottomrule
    \end{tabular}
    \label{tab:paramLyap}
\end{table}
Additionally, for the TFC implementation, the parameters used are summarized in Table \ref{tab:TFCLyap}.  
\begin{table}[ht]
    \centering
    \caption{TFC algorithm parameters}
    \begin{tabular}{cc}
        \toprule
        Variable & Value \\ \midrule
        $N$ [number of points] & $140$ \\
        $m$ [basis terms] & $130$\\
        $\varepsilon$ [tolerance] & $2.22 \times 10^{-16}$\\
        Maximum iterations & $20$\\
        \bottomrule
    \end{tabular}
    \label{tab:TFCLyap}
\end{table}

For all numerical tests, the unknown vector must be initialized. First, the terms $\B{\xi}_x$, $\B{\xi}_y$, and $\B{\xi}_z$ were all initialized by a null vector, which ultimately represents the simplest interpolating expression for the state variables. This initialization represents the worst-case scenario when there is no estimation of the trajectory. Next, the other unknown values of $\B{\alpha}$, $\B{\beta}$, and $b$ (which are associated with the position, velocity, and the period of the orbit) were initialized using Richardson's third-order analytical method for Halo-type periodic motion \cite{rich}.

This initialization was used to find the first orbit of the specified Jacobi constants. For the following orbits, the desired Jacobi constant was incrementally increased, and the converged values from the prior Jacobi constant level were used to initialize each step. 

This same process was utilized for the differential corrector method, which was implemented as a point of comparison to TFC. In the differential corrector inner-loop, the desired Jacobi constant was obtained by an iterative least-squares approach to update the initial guess.

\begin{example}{Lyapunov orbits around L1 \& L2 Lagrange points}

First, the method was used to explore the computation of Lyapunov orbits, which lie in the $x$-$y$ plane, or rather in the plane of the two primaries. For our test, the Lyapunov orbits were computed over a range of Jacobi constants, starting close to the equilibrium point's specific energy levels up to a Jacobi constant of 2.92. The associated trajectories for the orbits around L1 and L2 are provided in Figure \ref{fig:trajLyap_CP_L1} and Figure \ref{fig:trajLyap_CP_L2}.

\begin{figure}[H]
	\centering
	\includegraphics[width=0.5\linewidth]{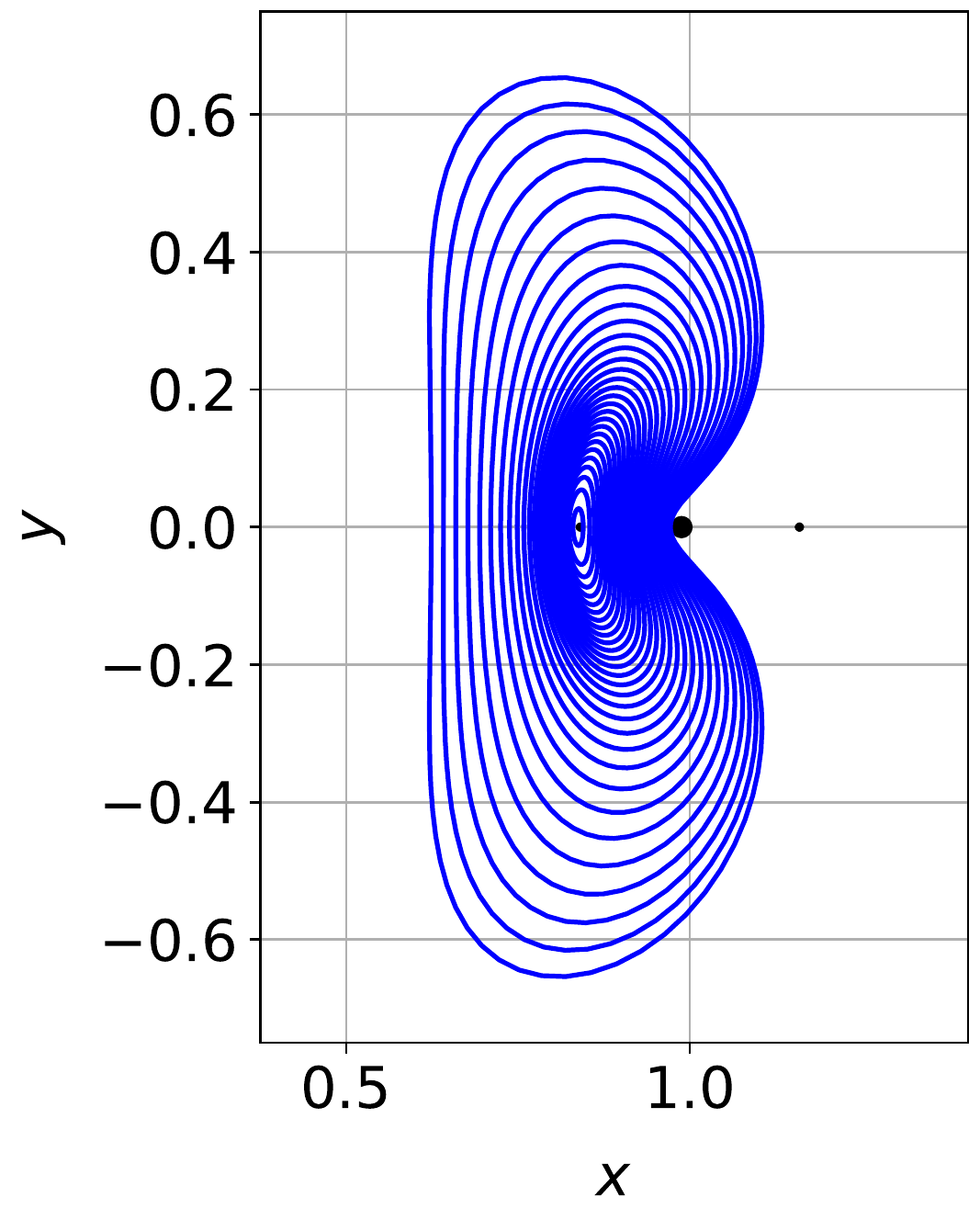}
	\caption{Lyapunov orbits for the Earth-Moon system for Jacobi constant values ranging from the energy of L1 to 2.92.}
	\label{fig:trajLyap_CP_L1}
\end{figure}

Additionally, a comparison with the differential corrector method (Reference \cite{diff_cor}) is provided in terms of speed and accuracy. Figure \ref{fig:residualsLyap_CP_L1} compares the residuals for both methods where it can be seen that the TFC approach is around 2 orders of magnitude more accurate than the differential corrector at higher Jacobi constants. Furthermore, the computation of the TFC solution is slightly faster, a little over 0.25 seconds in the extreme case, as displayed in Figure \ref{fig:compTimeLyap_CP_L1}.

\begin{figure}[H]
	\centering
	\includegraphics[width=0.7\linewidth]{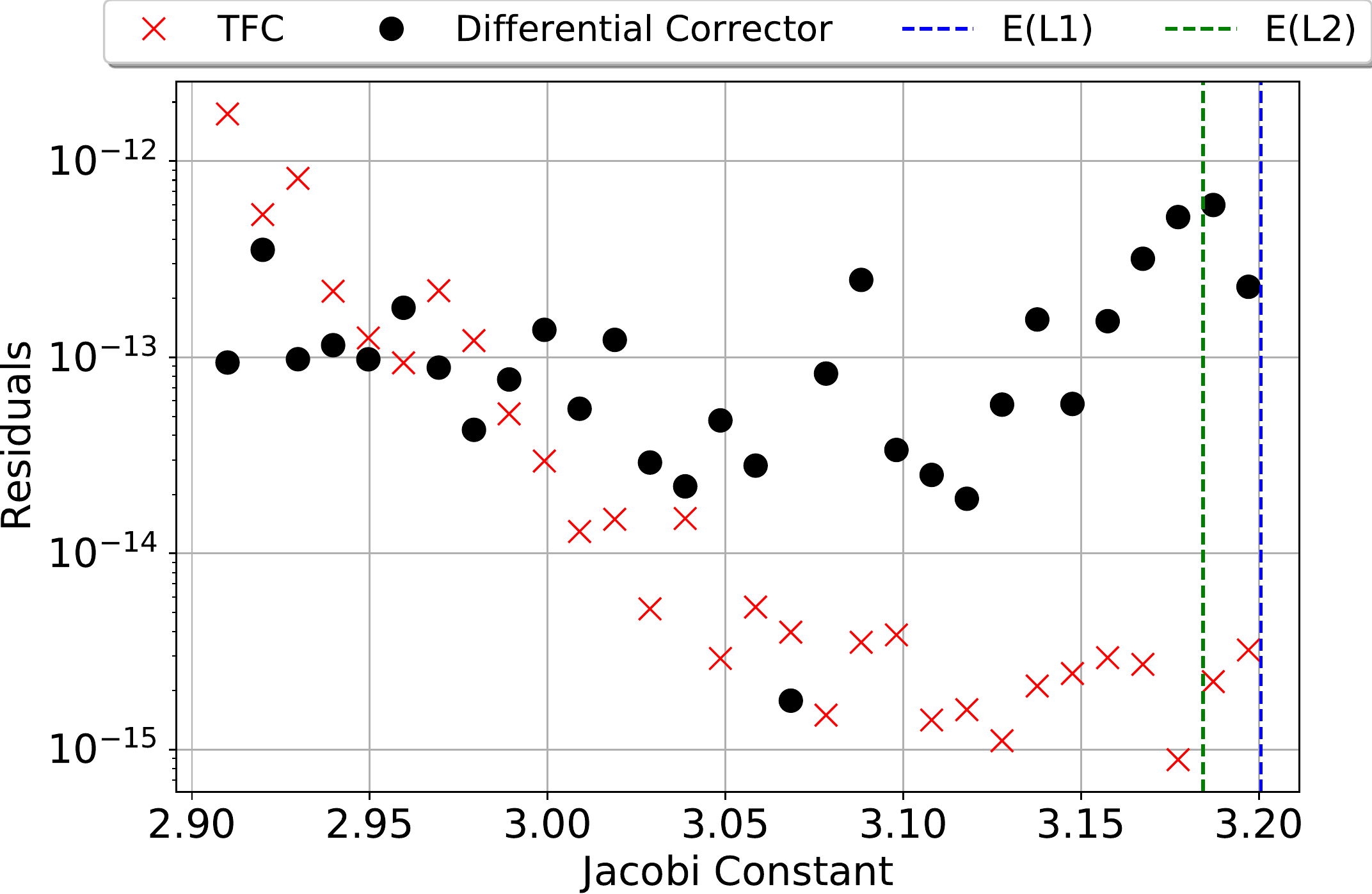}
	\caption{Maximum residuals of the loss vector for the TFC method solving for the trajectories plotted in Fig. \ref{fig:trajLyap_CP_L1} compared to that of the differential corrector. The lines of E(L1) and E(L2) represent the energy of the L1 and L2 Lagrange points respectively. The TFC approach has a slight accuracy advantage (an order-of-magnitude) as compared to the differential corrector method at higher Jacobi constants. }
	\label{fig:residualsLyap_CP_L1}
\end{figure}
\begin{figure}[H]
	\centering
	\includegraphics[width=0.7\linewidth]{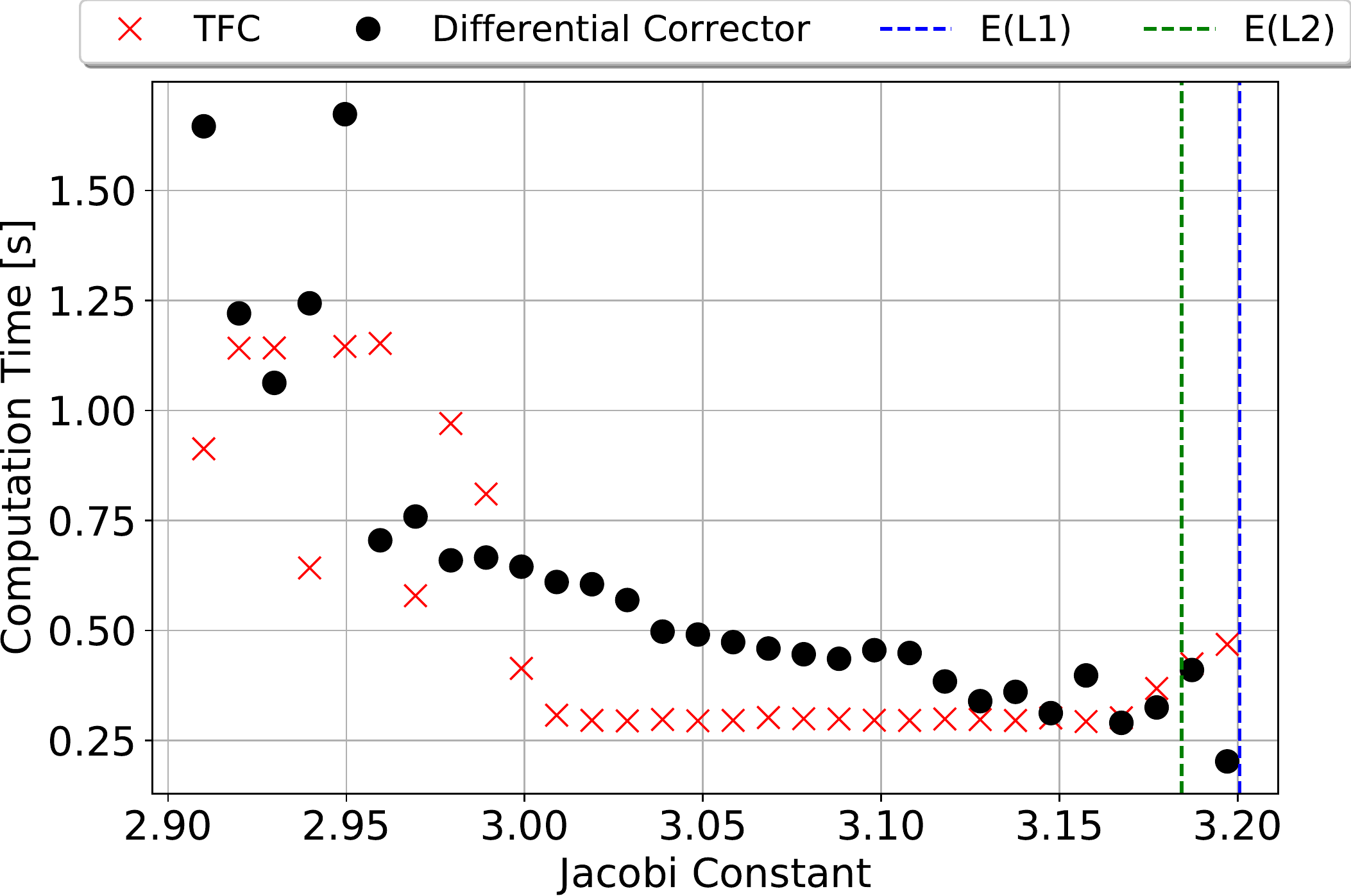}
	\caption{Computational time of the the TFC method for the trajectories plotted in Fig. \ref{fig:trajLyap_CP_L1} compared to that of the differential corrector. The TFC method holds a slight speed gain over the differential corrector.}
	\label{fig:compTimeLyap_CP_L1}
\end{figure}

Similar to the test for the L1 Lagrange point, Figure \ref{fig:trajLyap_CP_L2} displays the computed trajectories around the L2 Lagrange point. Additionally, like Figures \ref{fig:residualsLyap_CP_L1} and \ref{fig:compTimeLyap_CP_L1}, the accuracy and computation time for these tests are provided in Figures \ref{fig:residualsLyap_CP_L2} and \ref{fig:compTimeLyap_CP_L2}.
\begin{figure}[H]
	\centering
	\includegraphics[width=0.5\linewidth]{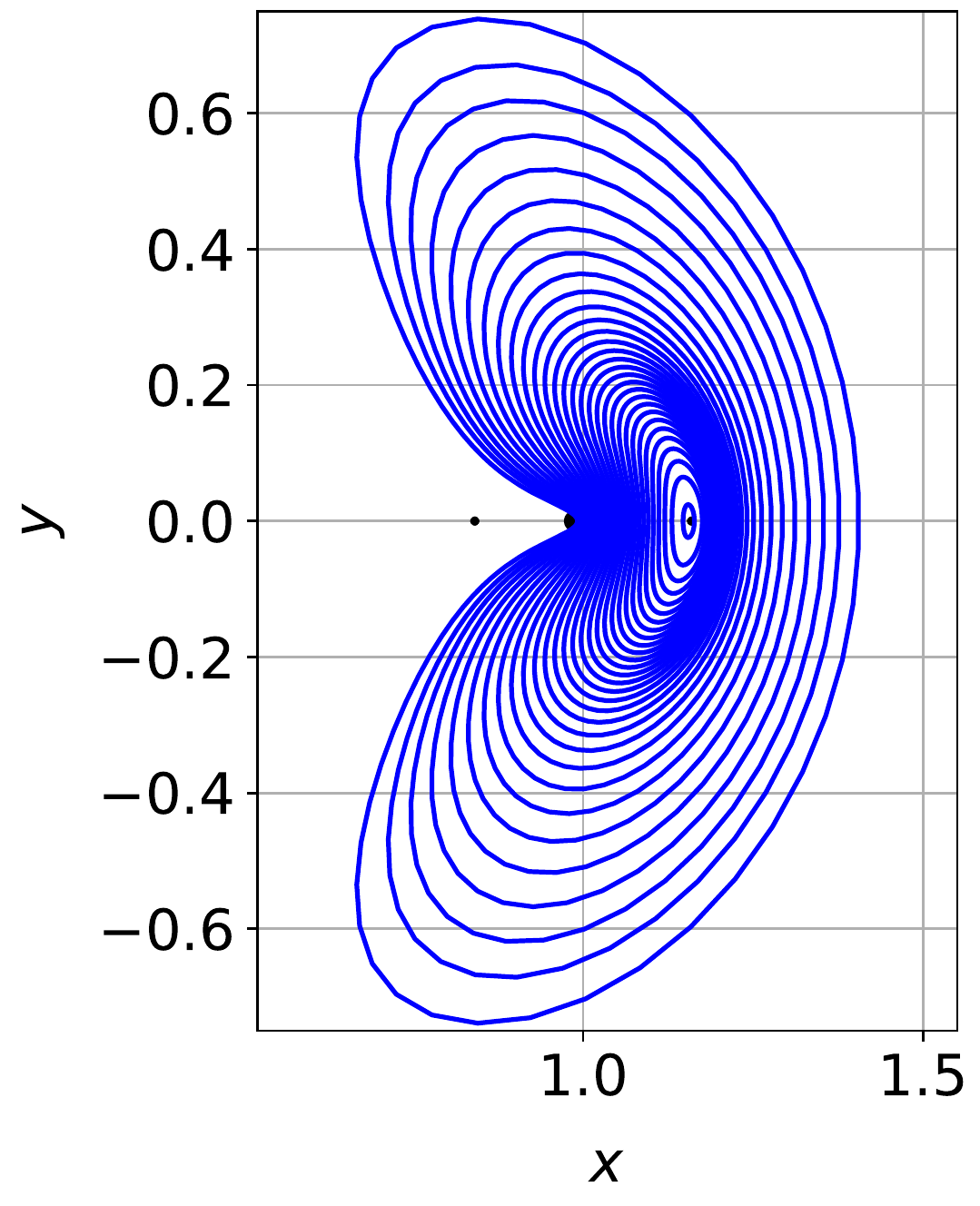}
	\caption{Lyapunov orbits for the Earth-Moon system for Jacobi constant values ranging from the energy of L2 to 2.92.}
	\label{fig:trajLyap_CP_L2}.
\end{figure}
In Figure \ref{fig:residualsLyap_CP_L2}, the TFC method is more accurate, albeit only slightly. At a Jacobi constant level of about 3.00 and above, the differential corrector method does not converge, as shown by the jump in accuracy. At this Jacobi constant level, the TFC method's accuracy starts to decrease before failing to converge at the Jacobi constant value of 2.92.
\begin{figure}[H]
	\centering
	\includegraphics[width=0.7\linewidth]{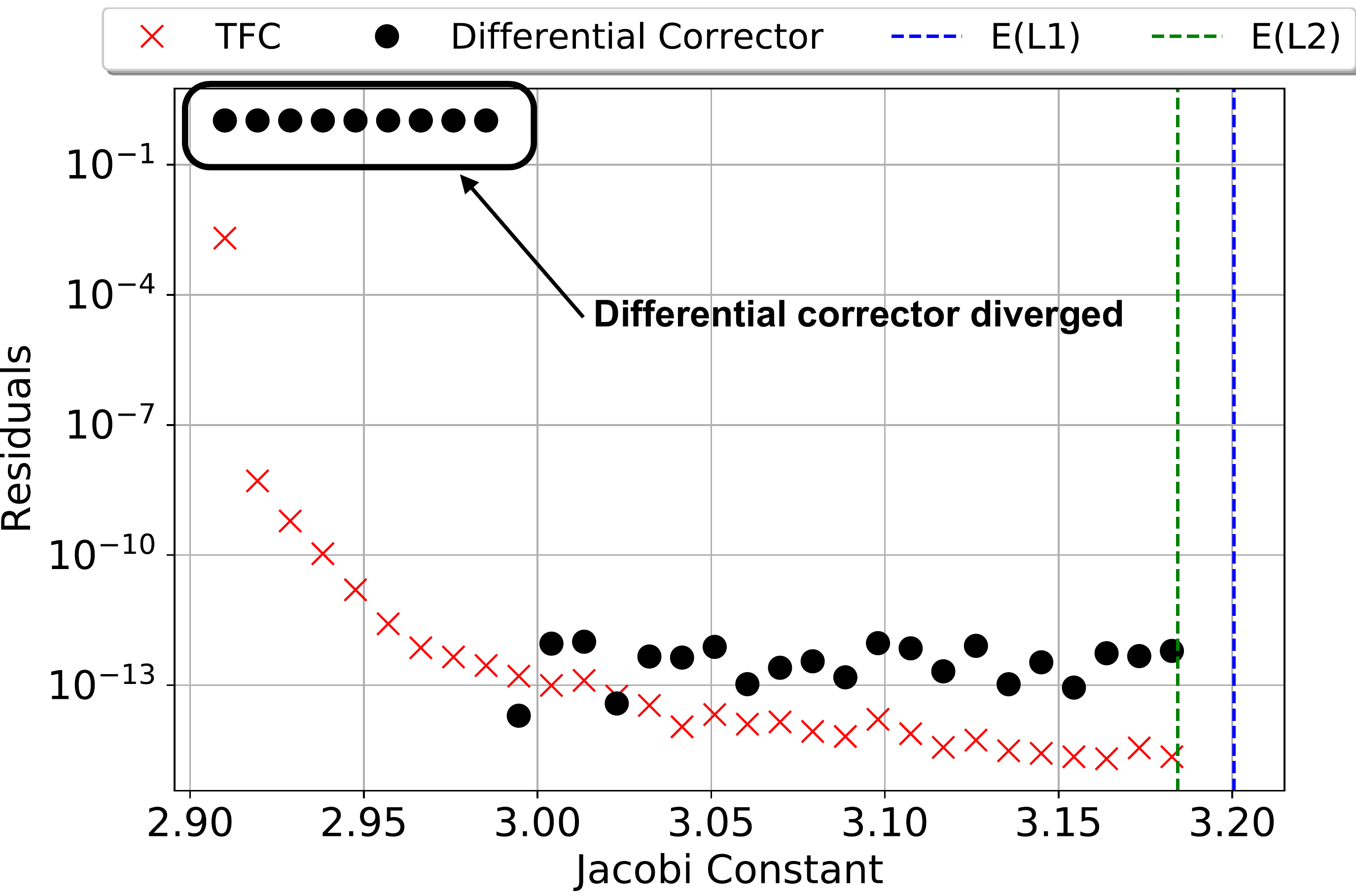}
	\caption{Maximum residuals of the loss vector for the TFC method solving for the trajectories plotted in Fig. \ref{fig:trajLyap_CP_L2} as compared to the differential corrector. For the trajectories around L2, the differential corrector diverged around a Jacobi constant level of 3.00, while the TFC method was able to solve the problem with diminishing accuracy. The black box highlights the diverged cases.}
	\label{fig:residualsLyap_CP_L2}
\end{figure}
\begin{figure}[H]
	\centering
	\includegraphics[width=0.7\linewidth]{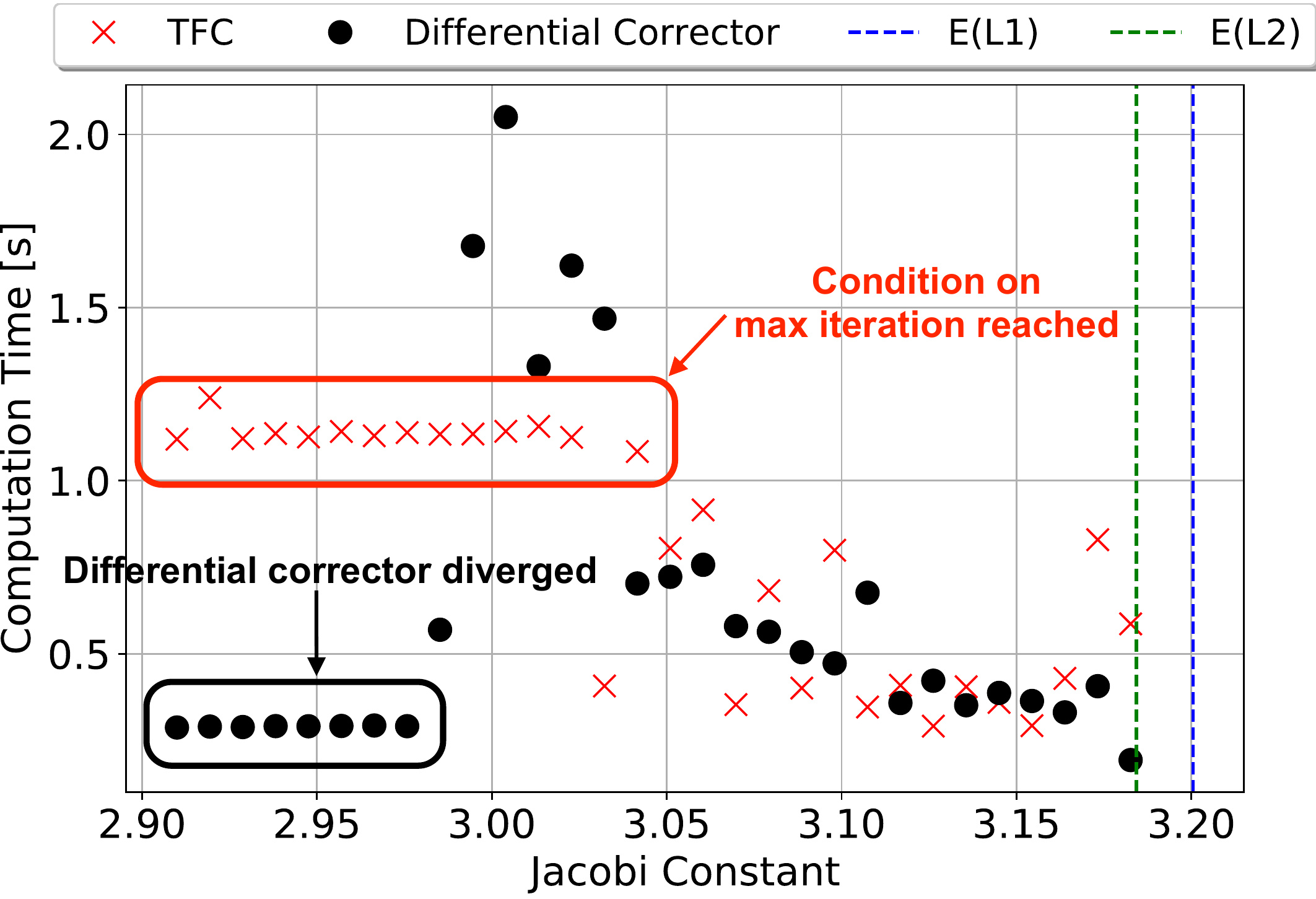}
	\caption{Computational time of the TFC method for the trajectories plotted in Fig. \ref{fig:trajLyap_CP_L2} compared to that of the differential corrector. Again, the black box highlights where the differential corrector diverged. Additionally, the red box shows where the TFC method reached its maximum allowed iterations of 20. These cases are correlated to the reduction of accuracy seen in Fig. \ref{fig:residualsLyap_CP_L2}.}
	\label{fig:compTimeLyap_CP_L2}
\end{figure}
\end{example}

\begin{example}{Halo Orbits around L1 \& L2 Lagrange points}
Next, the proposed technique was utilized to compute Halo orbits around L1 and L2. These orbits differ from Lyapunov orbits because they are \emph{not} restricted to the $x$-$y$ plane and become three-dimensional. In fact, this family of orbits is a bifurcation of the Lyapunov orbits computed in the previous section and are characterized by ``northern'' and ``southern'' bifurcations. However, when using the TFC method to compute these Halo orbits, the only thing that changes is the initialization of the $\B{\alpha}$, $\B{\beta}$, and $b$ parameters. First, we look at the computation of the ``northern'' family of Halo orbits around L1 and L2 as plotted in Figure \ref{fig:trajHalo_CP_N}. 

In these plots, we can see that around the L1 equilibrium point, the method converged to Lyapunov orbits for higher Jacobi constants. However, as the Jacobi constant decreases below 3.025, the method does not converge to a periodic orbit, as shown by the increase in residuals around the Jacobi constant value of 3.025. In fact, at the other L2 equilibrium point, the method converges for all Jacobi constant values; however, at values below 3.025, the solution jumps to a circular orbit around the Moon---therefore, these orbits were not plotted.

\begin{figure}[H]
        \begin{subfigure}[b]{0.45\textwidth}
                \centering \includegraphics[width=\linewidth]{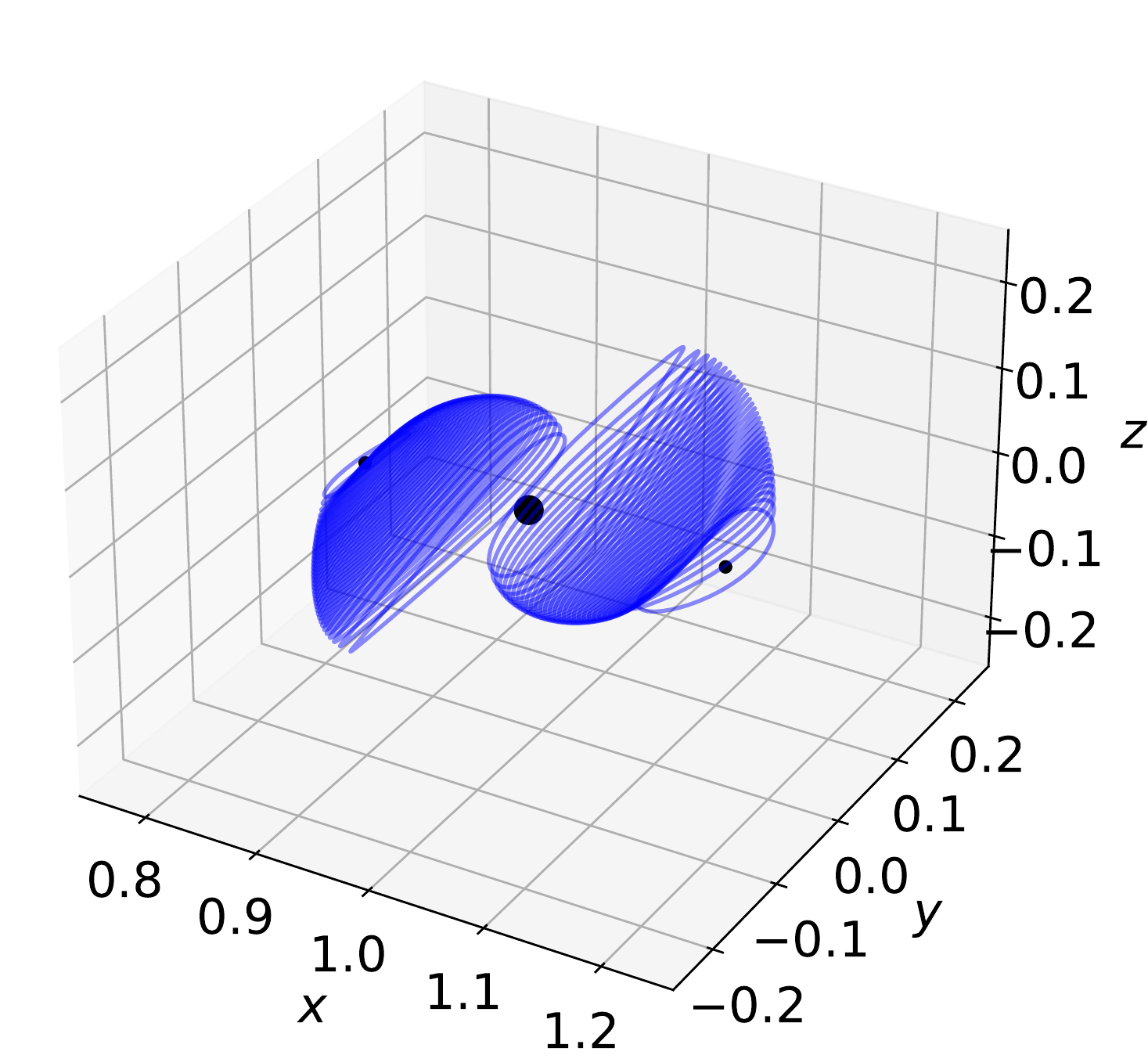}
        \end{subfigure}%
        \begin{subfigure}[b]{0.55\textwidth}
                \centering \includegraphics[width=\linewidth]{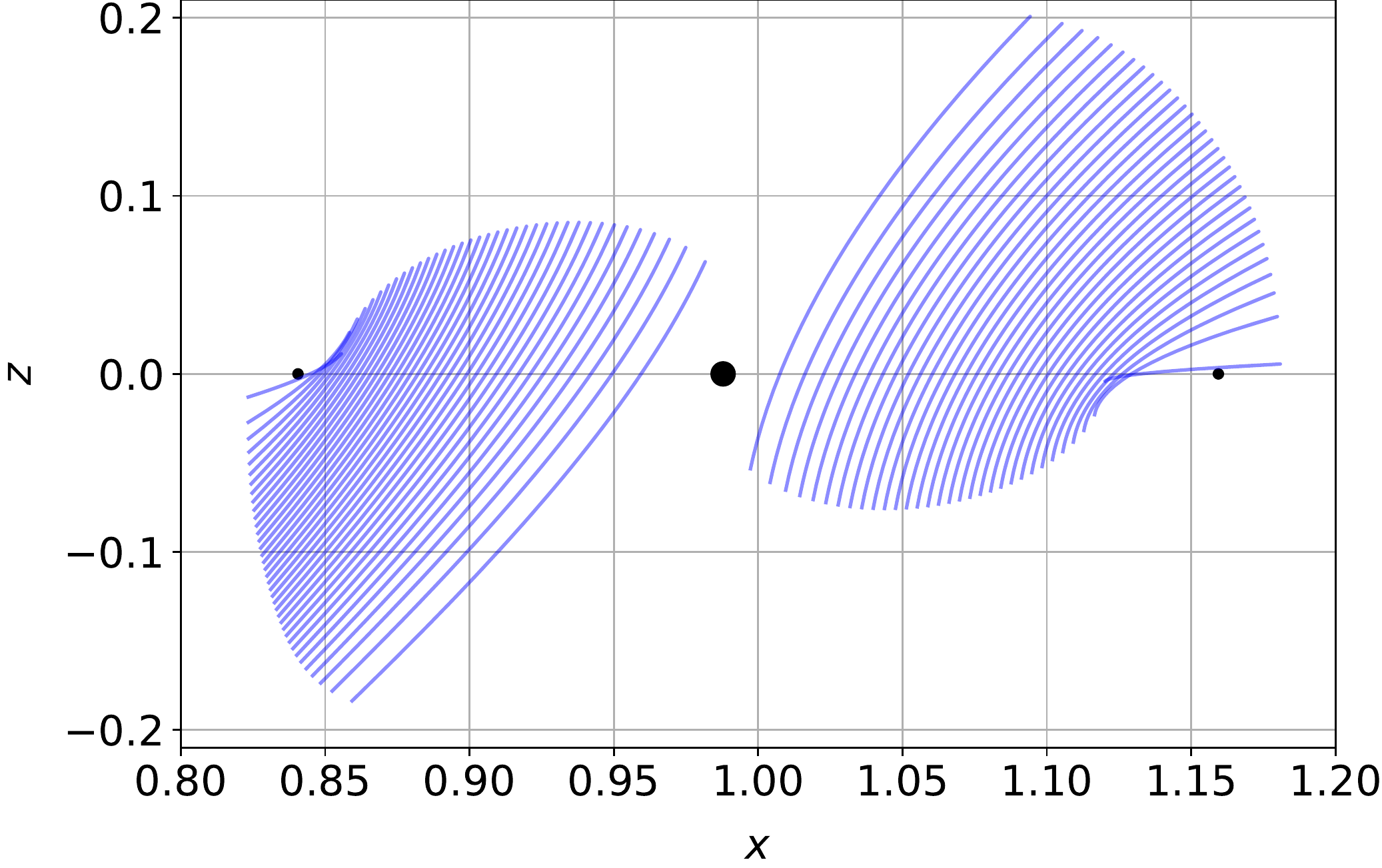}
        \end{subfigure}
        \caption{Halo orbits of the ``northern'' bifurcation around both L1 and L2 Lagrange points.}\label{fig:trajHalo_CP_N}
\end{figure}

Like the Lyapunov orbit tests, the loss vector's maximum residual was recorded and is plotted in \ref{fig:residualsHalo_CP_N}. All converged solutions were on the order of $\mathcal{O}(10^{-14})$. The convergence for Halo orbits took longer, with some cases taking 8 seconds, while on average, the solution time was around 2 seconds, as shown in Figure \ref{fig:compTimeHalo_CP_N}. 

Like the ``northern'' Halo orbits plotted in Figure \ref{fig:trajHalo_CP_N}, the Halo orbits of the ``southern'' bifurcation were computed with similar findings and, therefore, omitted from this paper for brevity.

\begin{figure}[H]
	\centering
	\includegraphics[width=0.62\linewidth]{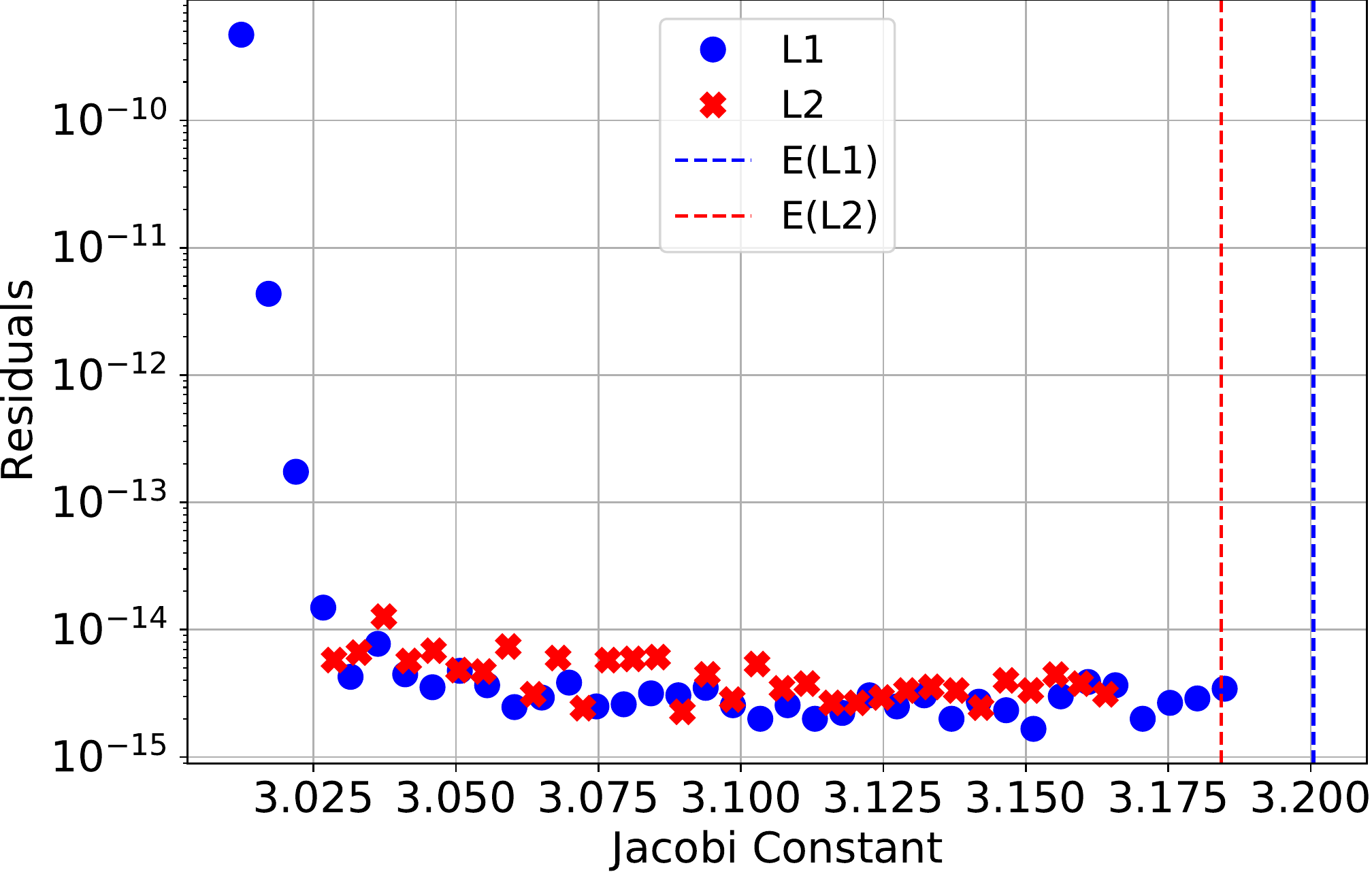}
	\caption{Maximum residuals of the loss vector for the TFC method solving for the trajectories plotted in Fig. \ref{fig:trajHalo_CP_N}. For almost all cases, the solution accuracy is on the order, $\mathcal{O}(10^{-14})$. However, around a Jacobi constant level of 3.025, the accuracy decreases for orbits around L1. The solutions for orbits around L2 lower than 3.025 are not plotted because, while they converged to a valid period orbit with high accuracy, it was not a Halo-type orbit.}
	\label{fig:residualsHalo_CP_N}
\end{figure}

\begin{figure}[H]
	\centering
	\includegraphics[width=0.59\linewidth]{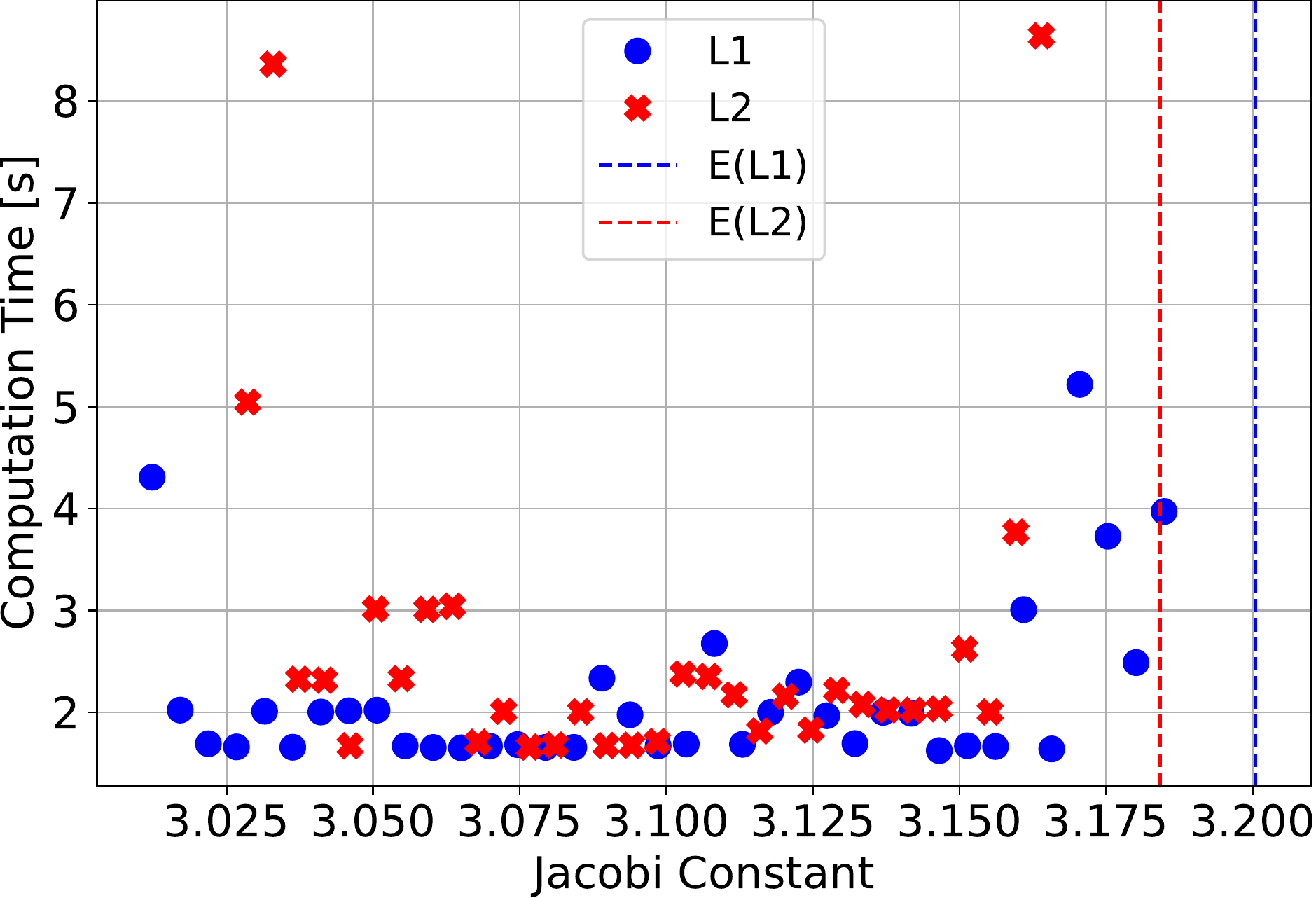}
	\caption{Computational time of the TFC method for the solution of ``northern'' Halo orbits around L1 and L2 plotted in Fig. \ref{fig:trajHalo_CP_N}. At first glance, it can easily be seen that the computation of these orbits too about twice as long to compute as the Lyapunov orbits. One cause of increased computation time is that the system of equations increased since more points and basis functions were need in the computation of these orbits.}
	\label{fig:compTimeHalo_CP_N}
\end{figure}
\end{example}

\section{Over-constrained differential equations}\label{sec:s3_overConDE}

In the following section, we revisit the theory developed in Section \ref{sec:s2a_overCON} and apply some of the over-constrained expressions to specific applications. First, a problem considering the interpolation of a trajectory based on noisy measurements augmented by a differential equation is explored. After this, a differential equation is analyzed by solving the continuous transformation from an initial-value problem to a boundary-value problem.

\subsection{Merging data with dynamics}

Consider a scenario where a trajectory is observed multiple times over its path. One question may arise about how this observational data (subject to measurement noise) can be incorporated along with the dynamical model to predict the object's actual path. The following example considers the merging of data with dynamics by using an over-constrained expression.

\begin{example}{Merging data with dynamics}
Consider a trajectory governed by the following differential equation,
\begin{equation*} 
    y_{xx} + 2y_x + y = 0,
\end{equation*}
with the analytical solution of the form $y(x) = e^{-x}(c_1 x + c_2)$, which was used to check the final answer and create the true trajectory. Additionally, assume that an object is ``observed'' under the influence of this dynamical system  at three points $x =[-1, -0.5, +1]$, and these measurements are subject to normally distributed noise such that,
\begin{equation*}
    y(-1) = \mathcal{N}(y_1^{true},\sigma_1^2), \quad y(-0.5) = \mathcal{N}(y_2^{true},\sigma_2^2), \quad \text{and} \quad \mathcal{N}(y_3^{true},\sigma_3^2).
\end{equation*}
To solve this problem, we can utilize Equation \eqref{eq:three_points} from Section \ref{sec:s2a_threePoints}. Additionally, since the measurement data has associated accuracy in terms of $\sigma_1$, $\sigma_2$, and $\sigma_3$, the weight matrix is defined by the variances such that,
\begin{equation*}
    W = \begin{bmatrix} \sigma_1^{-2} & 0 & 0 \\ 0 & \sigma_2^{-2} & 0 \\ 0 & 0 & \sigma_3^{-2} \end{bmatrix}.
\end{equation*}
For the given problem, we assume $\sigma_1 = \sigma_2 = \sigma_3 = 1$ and $y_1^{true} = 5$, $y_2^{true} = 4.515$, and $y_3^{true} = 2$. Using the development in Section \ref{sec:s2a_threePoints} this differential equation can incorporate information from all three observations even though the differential equation is only second order. Furthermore, after this step, the process to solve the differential equation is exactly the same as all prior examples and has been omitted for brevity. 

For this specific test, a Monte Carlo simulation of $10,000$ trials was conducted to determine the space that the function $y(x)$ could occupy given the ``observation'' uncertainty and subject to the governing dynamics of the differential equation. Figure \ref{fig:s3_three_points} shows the solution.
\begin{figure}[H]
	\centering
    \includegraphics[width=0.95\linewidth]{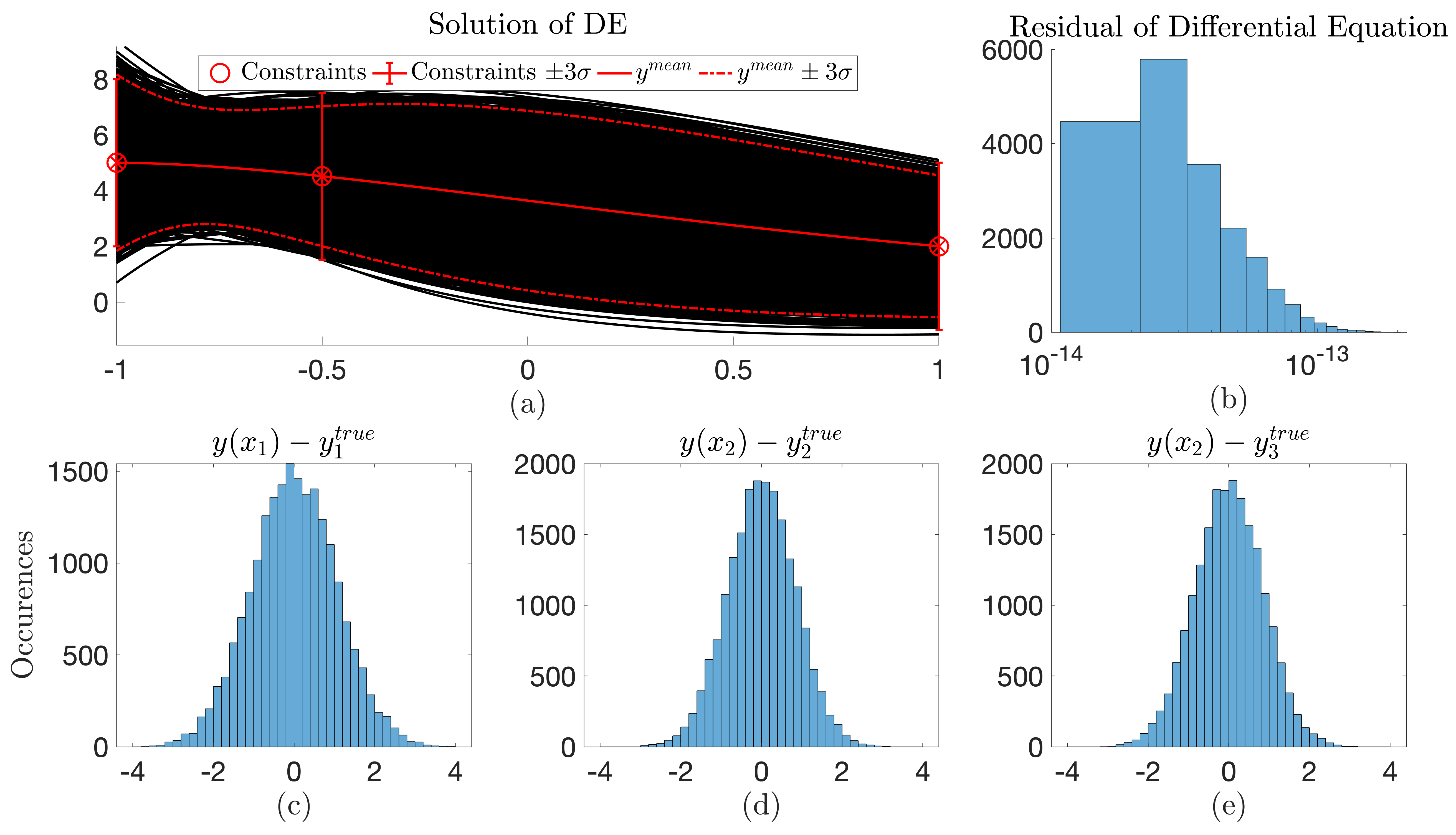}
    \caption{Monte Carlo test for $10,000$ trials. Plot (a) shows the differential equation solution space given the observation uncertainty. Plot (b) highlights the residuals of the differential equation over the entire simulation. It can be seen that the residuals of all solutions are between $10^{-13}$ to $10^{-14}$. Plots (c), (d), and (e) display the distribution of the constraint points around the true value. Note, these values are sampled from the solutions of the differential equation and not the constraints specified in the constrained expression.}
    \label{fig:s3_three_points}
\end{figure}
A probability bound for the differential equation can be produced through this test, along with an estimated mean. For all solutions, the differential equation residuals remained less than $10^{-13}$ verifying the accuracy of the method. Additionally, an interesting result of this test is in the final estimated solutions. This is most evident when observing the solution trajectories near the constraint points, $x =-0.5$ and $x =+1$; it can be seen that the $3 \sigma$ of the differential equation is less than that of $3\sigma$ associated with the constraints. This happens because the loss function in the TFC method minimizes the residuals of the differential equation; in the over-constrained TFC method, the residuals are minimized simultaneously with the weighted least-squares of the observations (or constraints).
\end{example}

\subsection{Initial to boundary value problem transformation}
The development of the over-constrained expression led to this question: if a differential equation can be solved with more constraints than its order, what is the connection between an initial- and boundary-value problem? 

\begin{example}{Initial to boundary value problem transformation}
Consider the second-order, linear differential equation given by,
\begin{equation*}
    y_{xx} + \left[\cos (3 x^2) - 3 x + 1\right] y_x + \left[6 \sin(4 x^2) - e^{\cos(3 x)}\right] y = 2 \dfrac{[1 - \sin(3 x)] (3 x - \pi)}{4 - x}
\end{equation*}
subject to the three constraints
\begin{equation*}
    y (-1) =-2, \qquad y_x (-1) =-2, \qquad \text{and} \qquad y (+1) = 2.
\end{equation*}
By using Equation \eqref{eq:delta} from Section \ref{sec:s2a_threePoints} and defining $s_1(x) = 1$ and $s_2(x) = x$, we can write,
\begin{equation*}
   W \begin{bmatrix} 1 & x_0 \\ 0 & 1 \\ 1 & x_f \end{bmatrix} \begin{bmatrix} \alpha_{11} & \alpha_{12} & \alpha_{13} \\ \alpha_{21} & \alpha_{22} & \alpha_{23} \end{bmatrix} = W.
\end{equation*}
Next, before solving this system, let us also define the weight matrix as,
\begin{equation*}
    W = \begin{bmatrix} 1 & 0 & 0 \\ 0 & 1 - \gamma & 0 \\ 0 & 0 & \gamma\end{bmatrix}
\end{equation*}
where $\gamma$ is a weight parameter transforming the problem from IVP to BVP as $\gamma\in[0, 1]$. Now, solving for the system we get the pseudo-switching functions,
\begin{align*}
    \varphi _1 = s_1 \alpha_{11} + s_2 \alpha_{21} &= \frac{1}{1 + 4\gamma - \gamma^2}\Big( (1 + \gamma) - 2 \gamma x \Big)\\
    \varphi _2 = s_1 \alpha_{12} + s_2 \alpha_{22} &= \frac{1}{1 + 4\gamma - \gamma^2}\Big( (1 - \gamma)^2 + (1-\gamma^2)x \Big)\\
    \varphi _3 = s_1 \alpha_{13} + s_2 \alpha_{23} &= \frac{1}{1 + 4\gamma - \gamma^2} \Big( -\gamma(\gamma-3) + 2 \gamma x \Big)
\end{align*}
so the total over-constrained expression takes the form,
\begin{equation*}
    y(x,g(x)) = g(x) + \varphi _1(x) \Big(-2 - g(-1)\Big) + \varphi_2(x) \Big(-2 - g_x(-1)\Big) +  \varphi_3(x) \Big(2 - g(1)\Big).
\end{equation*}
Again, by utilizing the numerical techniques discussed earlier, we can define the free function, plug the resulting expression into the differential equation to create our loss function, discretize the domain at the collocation nodes, and solve the system via least-squares. Figure \ref{fig:s3_IVP2BVP_soln} shows this transformation ``surface'' along with the residuals of the differential equation for validation of the method. Figure \ref{fig:s3_IVP2BVP_res} shows that the mean residual over all of the $\gamma$ values are on the order of $10^{-14}$ with a standard deviation on the same order.
\begin{figure}[H]
	\centering \includegraphics[width=0.7\linewidth]{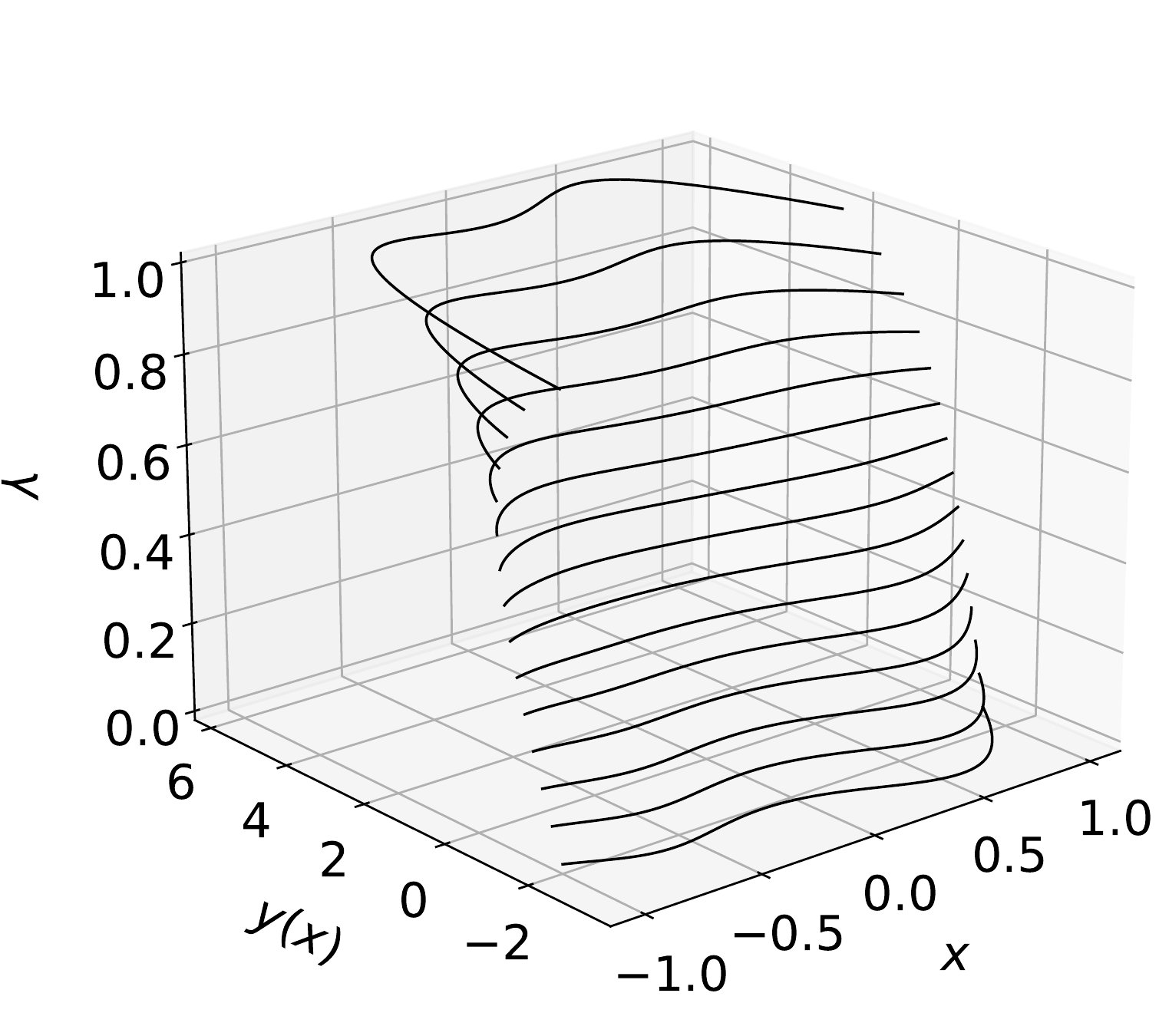}
    \caption{IVP to BVP differential equation parametric transformation. These plots shows the solution of the differential equation, $y (x)$, continuously morphing from IVP constraints to BVP constraints.}
    \label{fig:s3_IVP2BVP_soln}
\end{figure}

\begin{figure}[H]
	\centering \includegraphics[width=0.7\linewidth]{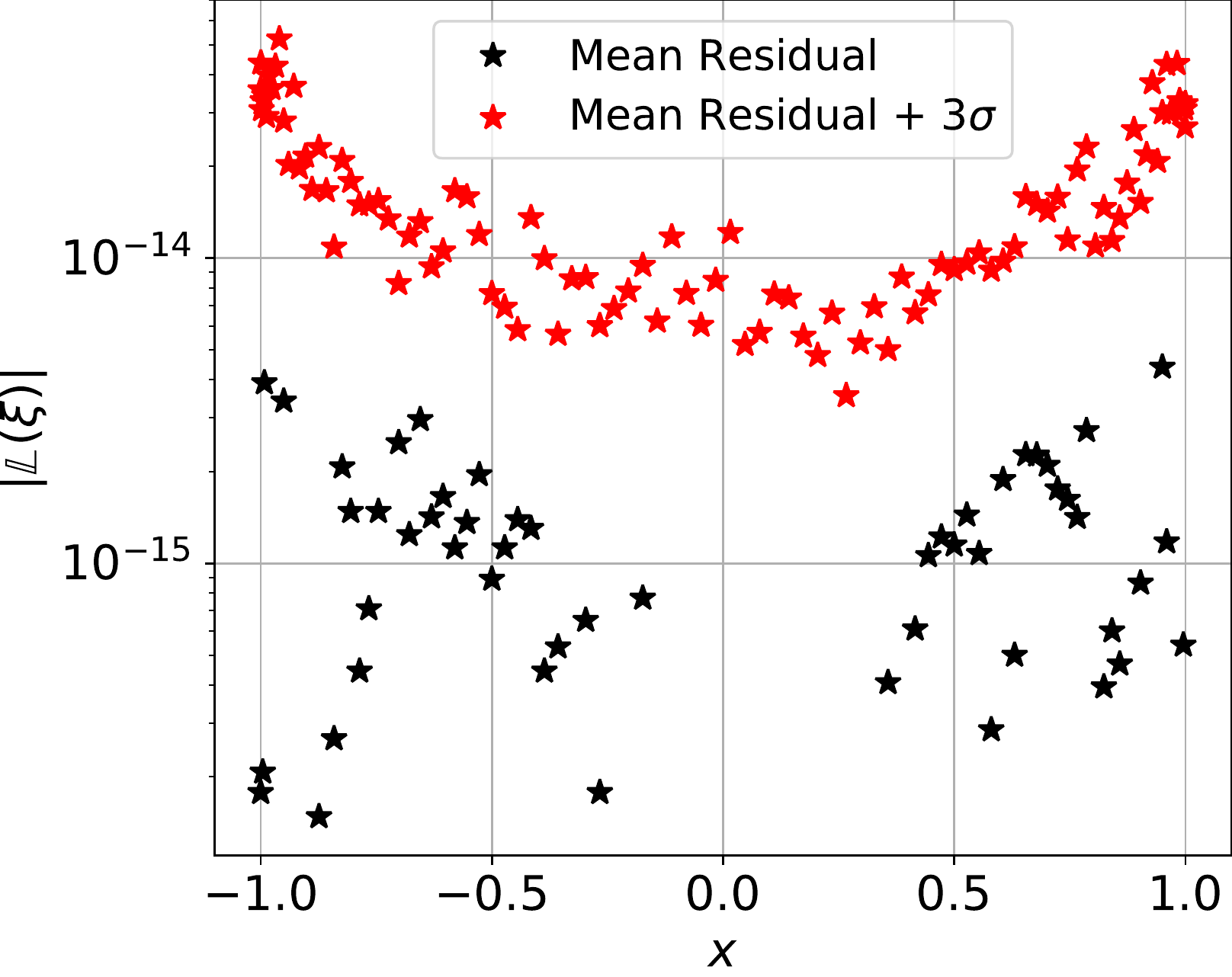}
    \caption{Residuals of loss vectors for IVP to BVP differential equation parametric transformation. In all cases, the residual of the differential equation is on the order of $10^{-14}$.}
    \label{fig:s3_IVP2BVP_res}
\end{figure}
\end{example}

\begin{mybox}
\begin{center}
\begin{huge}
Part 3\\
\vspace{0.2in}
Optimal Control
\end{huge}

\vspace{0.5in}

Some days feel warm no matter how cold\\
they are, and some things are fun no\\
matter how old you are, and sometimes\\
you wish a visit could just last forever...\\
\vspace{0.1in}
--- Unravel, \emph{ColdWood Interactive}
\end{center}
\end{mybox}
\pagebreak{}
%

\chapter{USE FOR REAL-TIME OPTIMAL CONTROLLERS IN AEROSPACE SYSTEMS\label{chap:opt_con}}

Over the previous sections, we have explored the Theory of Functional Connections to build the constrained expression and solve differential equations subject to constraints. In this section, we will take everything we have learned thus far and explore its application to the field of optimal control, and specifically, real-time optimal control, which is an active field of research. It should be clear from the examples given in Chapter \ref{chap:ode} that TFC is an effective method to solve differential equations.

Transitioning from theoretic equations to the physical world, many problems arise affecting the accuracy and robustness of controllers, including unmodelled dynamics and sensor measurement noise, which can result in a deviation from the desired optimal trajectory. Classically, this problem is overcome by deriving a closed-loop controller that tracks the optimal reference trajectory (e.g., Mars Science Laboratory guidance \cite{singh2007guidance}). While the closed-loop controller may be optimal in following the reference trajectory, it will be sub-optimal in the global problem since a disturbance in the state should redefine the full optimal trajectory.  Solving for the new optimal solution would involve computing a single-open loop trajectory consisting of the optimal state and optimal control program history. However, as mentioned above, disturbances and measurement noise will cause a deviation from this solution. Therefore, this computation would have to be done during each guidance cycle of the computer allowing for an updated solution based on the state. 

\begin{figure}[ht]
    \centering\includegraphics[width=0.95\linewidth]{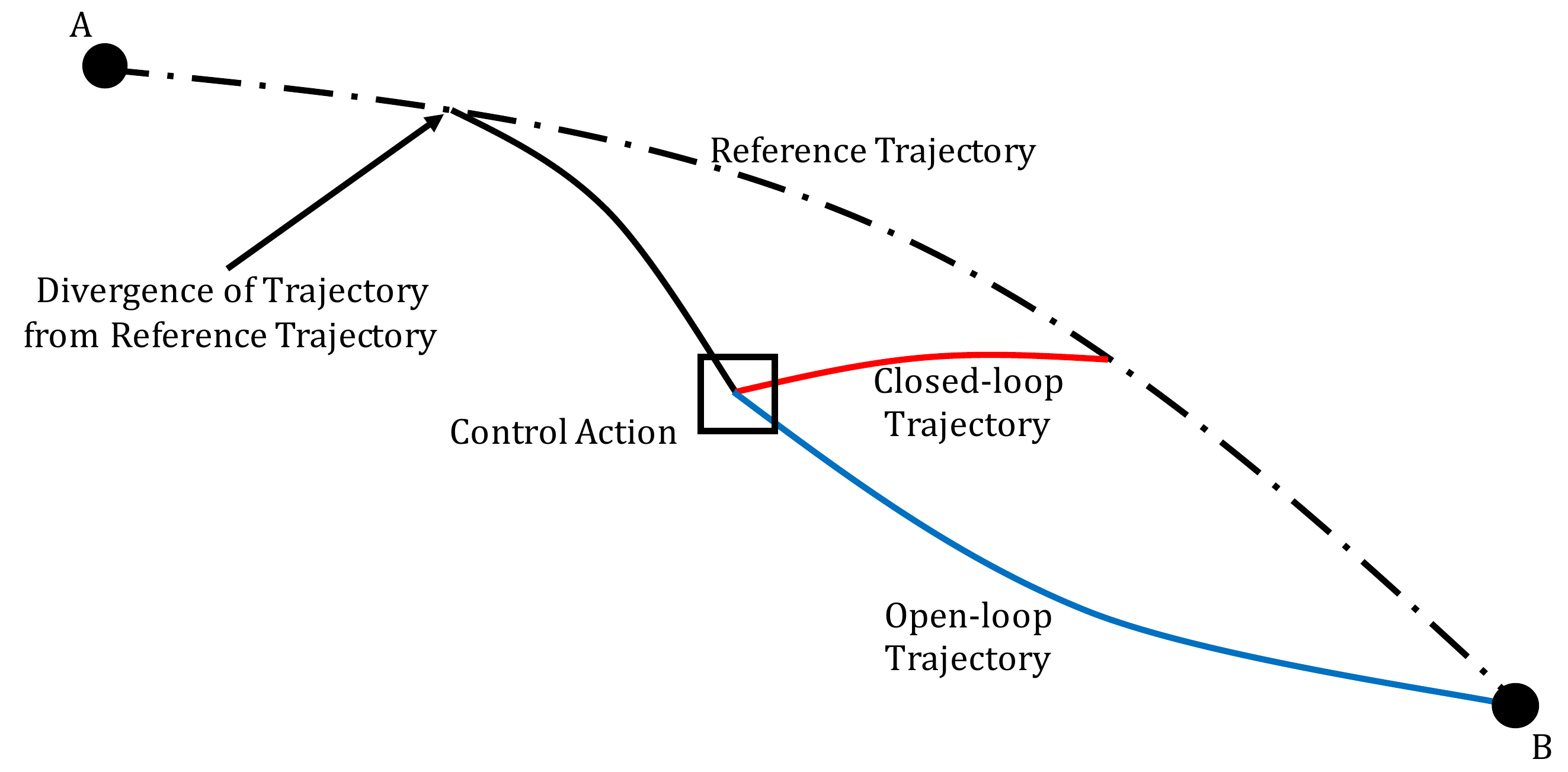}
    \caption{Trajectory going from Point A to Point B. The dashed line represents the reference trajectory. In this situation, the true trajectory deviates from the reference trajectory. At the guidance computer cycle, the closed-loop controller acts optimally to return the trajectory (red line) to the reference trajectory. On the other hand, the open-loop solution provides the optimal path from the new point and the resulting trajectory follows this path (blue line).}
    \label{fig:s4_controllers}
\end{figure}

The difference between the two methods mentioned above is easily visualized with a simple example provided in Figure \ref{fig:s4_controllers}.  
Consider some optimal control problem where it is desired that an object's trajectory goes from point A to point B subject to some cost function. Over the course of the trajectory, the true path can deviate from the reference trajectory due to such things as unmodelled dynamics, disturbances, etc. In practice, the control for this reference trajectory is followed until the next guidance cycle, signified by the black box in Figure \ref{fig:s4_controllers}. At this point, sensors provide some information on the state, i.e., position, velocity, etc., and a control action is determined. In the case of a closed-loop control law, the computed control action will be the one that optimally returns the object to the reference trajectory. Conversely, an open-loop control law will recompute a new optimal trajectory from the current state, producing a trajectory that could be drastically different reference trajectory. 

Contrary to the example in Figure \ref{fig:s4_controllers}, in actual implementation, the frequency of the guidance cycle is drastically higher, and therefore, the control is updated more often. For example, it is reported all guidance functions on the Mars Science Laboratory \cite{singh2007guidance} are within 60 to 70 Hz ($\sim$ 14 to 17 ms). While other applications, this can exceed 100 Hz. Additionally, in the case of the open-loop solution, this implies that a new solution must be computed at this frequency. 

Clearly, to enable such technology, real-time solutions must be obtained as quickly as possible to implement the recomputed optimal trajectory and control. With the exponential increase in computational power, this computation has become more feasible for onboard implementation, and researchers have started to explore the possibility of rapid and real-time trajectory generation for guidance application \cite{rao2018mesh,dennis2019computational,ross2006issues} through open-loop solutions. Additionally, the issues associated with real-time optimal control have also been recently explored in Reference \cite{ross2006issues}. Overall, the idea is to generate an optimal feedback control that can be constructed by continuously generating computational open-loop optimal trajectories quickly and efficiently \cite{ross2008optimal,rao2018mesh,dennis2019computational}. 

Therefore, with this being the ultimate goal, the following sections will focus on studying the solution of the single open-loop optimal control problems using the TFC approach, where we are interested in determining the limits of the method's speed and accuracy. By developing a fast, accurate, and robust solver, this smaller algorithm can be eventually incorporated into the larger problem, as mentioned above. In the following section, we will discuss the current techniques to solve the open-loop optimal control problem.

\section{Techniques to solve optimal control problems: direct vs. indirect method}
Usually, two methods are available to solve optimal control problems, direct and indirect methods. Direct methods are based on discretizing the continuous states and controls to transform the continuous problem into a nonlinear programming (NLP) problem \cite{darby2011hp,fahroo2002direct,ross2004pseudospectral}. The latter can be cast as a finite constrained optimization problem that can be solved via any of the available numerical algorithms that have the potential to find a local minimum, e.g., trust-region method \cite{byrd2000trust}. Whereas direct methods have been applied to solve a large variety of optimal control problems \cite{josselyn2003rapid,graham2015minimum,miller2017rapid,jiang2019integrated}, the general NLP problem is considered NP-hard, i.e., non-deterministic polynomial-time hard. NP-hard problems imply that the required amount of computational time needed to find the optimal solution does not have a predetermined bound, i.e., a bound cannot be determined a priori. NP-hard problems are such that the computational time necessary to converge to the solution is not known. As a consequence, the lack of assured convergence may result in questioning the reliability of the proposed approach. Since for optimal, closed-loop space guidance, most problems require computing numerical solutions onboard and in real-time; general algorithms that solve NLP problems cannot be reliably implemented. More recently, researchers have been experimenting with transforming optimal control problems from a general non-convex formulation into a convex optimization problem \cite{acikmese2007convex,blackmore2010minimum}. Here, the goal is to take advantage of the assured convex convergence properties. Indeed, convex optimization problems are shown to be computationally tractable as their related numerical algorithms guarantee convergence to a globally optimal solution in polynomial time. The general convex methodology requires that the optimal guidance problem is formulated as convex optimization whenever appropriate or convexification techniques are applied to transform the problem from a non-convex problem into a convex one. Such methodologies have been proposed and applied to solve optimal guidance and control via the direct method in a large variety of problems including, planetary landing \cite{acikmese2007convex,blackmore2010minimum}, entry atmospheric guidance \cite{wang2016constrained,wang2018autonomous}, rocket ascent guidance \cite{zhang2019rapid}, and low thrust \cite{wang2018minimum}.

Alternatively, a second approach to solve optimal control and guidance problems, called the indirect method, has been generally applied to various optimal control problems. This approach applies optimal control theory (i.e., Pontryagin Minimum Principle, PMP) to formally derive the first-order necessary conditions that must be satisfied by the optimal solution (state and control). The problem is cast as a two-point boundary value problem (TPBVP) that must be solved to determine the time evolution of state and costate from which the control generally depends. For general nonlinear problems, the necessary conditions result in a complicated set of equations and conditions. Additionally, the resulting TPBVP tends to be highly sensitive to the initial guess on the costates making the problem very hard to solve. Although indirect methods are known to yield more accurate optimal solutions, they are tough to implement and tend to be less used in practice with respect to direct methods. For this problem, we attempt to alleviate the sensitivity of initialization by TFC constrained expressions.

In the next section, we will look at the derivation of the TPBVP from the indirect method, starting with first principles. Additionally, we will explore how the TFC constrained expression reduces the number of algebraic equations to be solved.

\section{Summary of the indirect method}
To thoroughly understand the application of the TFC method to solve optimal control problems, a basic understanding of optimal control theory, and more specifically, the indirect method based on the calculus of variation, is needed. For the reader's convenience, the mathematical foundation for a general optimal control problem is provided in this section. For an extensive look into a plethora of optimal control problem types solved using the indirect method, the reader is directed to ``Applied Optimal Control'' by Bryson and Ho \cite{Bryson_Ho}.

\begin{figure}[ht]
    \centering\includegraphics[width=0.75\linewidth]{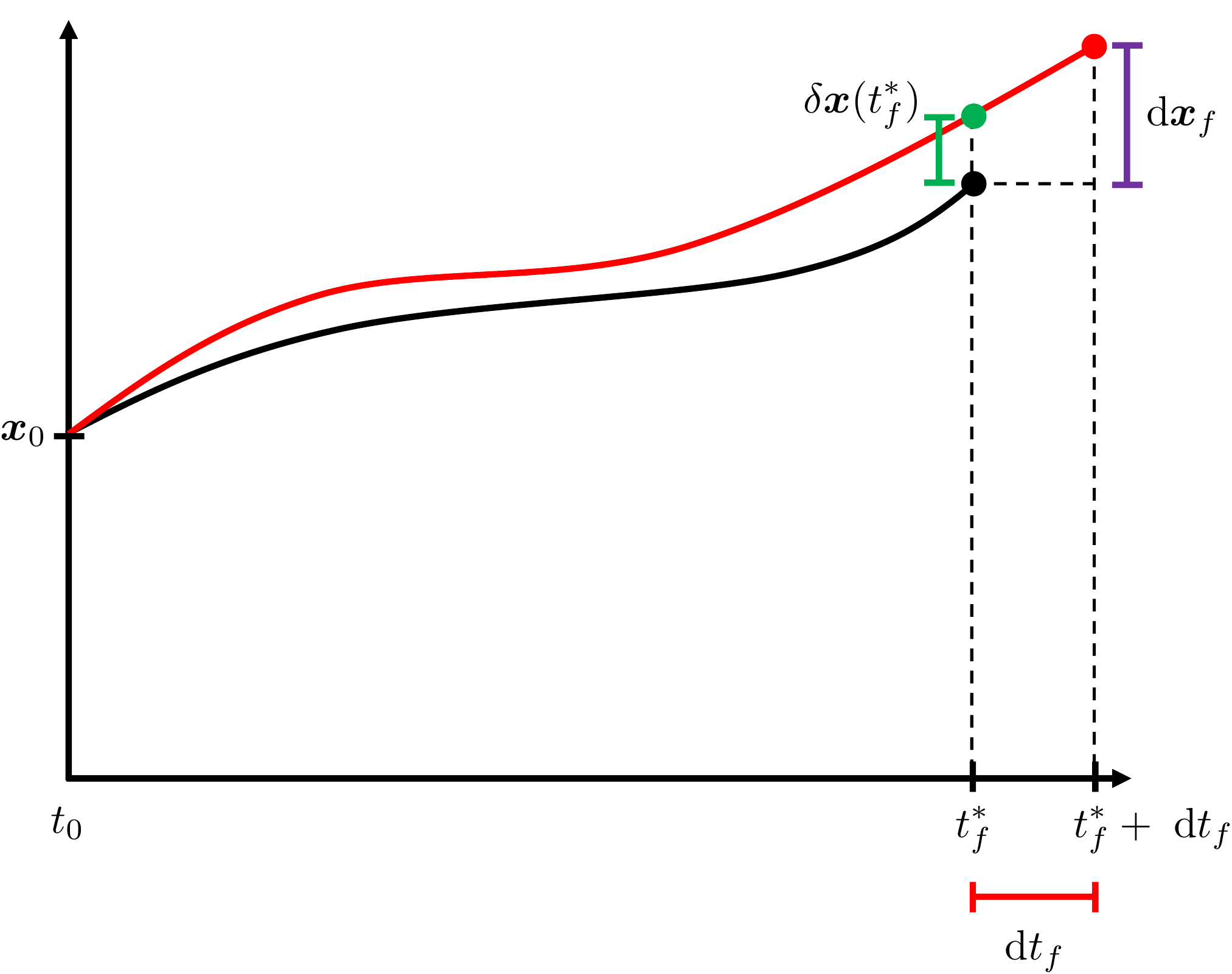}
    \caption{Graphical representation of the admissible variation, $\delta \B{x}(t^*_f)$, which is the state's variation with respect to the optimal trajectory's (black line) final condition, $\B{x}^*_f$.}
    \label{fig:s4_skewVariation}
\end{figure}

In general, a continuous-time dynamical optimization problem can be posed as a minimization of the cost functional (known as the Bolza Problem),
\begin{equation}\label{eq:s4_genCost}
    J = \Phi\left(\B{x}(t_f), t_f\right) + \int_{t_0}^{t_f} \mathcal{L}\left(\B{x}(t),\B{u}(t),t \right) \dd t
\end{equation}
where $\B{x}(t)$ is the state vector and $\B{u}(t)$ is the control vector, both a function of the independent variable of time, $t$. In this formulation, $\Phi$ is a function is the cost associated with the terminal state values and $\mathcal{L}$ is cost over the trajectory. In addition to Equation \eqref{eq:s4_genCost}, the states' dynamics are governed by a general nonlinear equation,
\begin{equation}\label{eq:s4_DE}
    \bdot{x} = \B{f}\left(\B{x}(t),\B{u}(t),t\right)
\end{equation}
with the boundary constraints 
\begin{equation}\label{eq:s4_boundary1}
    \B{x}(t_0) = \B{x}_0
\end{equation}
\begin{equation*} 
    \Psi \left(\B{x}(t_f),t_f\right) = \B{0}.
\end{equation*}
By adjoining the system of differential equations given by Equation \eqref{eq:s4_DE} with the Lagrange multiplier functions $\B{\lambda}(t)$, called the costate functions, the augmented cost function becomes,
\begin{align}\label{eq:s4_Jaug1}
    J_a = \Phi(\B{x}(t_f), t_f) &+ \B{\nu}\T \Psi(\B{x}(t_f), t_f) \nonumber\\ &+ \int_{t_0}^{t_f} \Big( \mathcal{L}(\B{x}(t),\B{u}(t),t) + \B{\lambda}\T(t) \B{f}(\B{x}(t),\B{u}(t),t) - \B{\lambda}\T(t)\bdot{x} \Big)\dd t
\end{align}
In optimal control theory, the first two terms in the integral are defined as the scalar function $H$ called the Hamiltonian,
\begin{equation}\label{eq:s4_Ham}
    H(\B{x}(t),\B{u}(t),\B{\lambda}(t), t) = \mathcal{L}\left(\B{x}(t),\B{u}(t),t \right) + \B{\lambda}\T(t) \B{f}(\B{x}(t),\B{u}(t),t)
\end{equation}
Substituting Equation \eqref{eq:s4_Ham} into Equation \eqref{eq:s4_Jaug1} and dropping the function arguments for clarity yields,
\begin{equation*}
    J_a = \Phi + \B{\nu}\T \Psi + \int_{t_0}^{t_f}\Big(H - \B{\lambda}\T \bdot{x}\Big) \dd t.
\end{equation*}
Consider the variation of the augmented cost function $J_a$ about the optimal solution of $J(\B{x}^*, \B{u}^*, t_f^*)$ where the ($^*$) signifies the optimal solution,
\begin{align*}
    \delta J_a = \frac{\partial \Phi}{\partial \B{x}^*(t_f^*)}\T \dd \B{x}^*_f &+ \frac{\partial \Phi}{\partial t_f^*}\dd t_f + \B{\nu}\T \frac{\partial \Psi}{\partial \B{x}^*(t_f^*)} \dd \B{x}_f + \B{\nu}\T \frac{\partial \Psi}{\partial t_f^*} \dd t_f + \dd \B{\nu}\T \Psi(\B{x}^*(t_f^*),t_f^*) \\
    &+ \int_{t_0}^{t_f^*}\Big[\frac{\partial H}{\partial \B{x}^*(t)}\T \delta \B{x} + \frac{\partial H}{\partial \B{u}^*(t)}\T \delta \B{u} + \frac{\partial H}{\partial \B{\lambda}^*}\T \delta \B{\lambda} - \delta\B{\lambda}\T \bdot{x}^* - \B{\lambda}^{*\mbox{\tiny T}} \delta \bdot{x}^* \Big]\dd t \\ &+ \Big[H(\B{x}^*, \B{u}^*, \B{\lambda}^*, t_f^*) - \B{\lambda}^{*\mbox{\tiny T}}(t_f^*) \bdot{x}^*(t_f^*) \Big] \dd t_f.
\end{align*}
Collecting terms and rewriting $\B{\lambda}^{*\mbox{\tiny T}} \delta \bdot{x}^*$ using integration by parts leads to,
\begin{align*}
    \delta J_a =& \Big[\frac{\partial \Phi}{\partial \B{x}^*(t_f^*)} + \frac{\partial \Psi}{\partial \B{x}^*(t_f^*)}\T \B{\nu} \Big]\T \dd \B{x}_f + \Big[\frac{\partial \Phi}{\partial t_f^*} + \B{\nu}\T \frac{\partial \Psi}{\partial t_f^*} \Big] \dd t_f + \dd \B{\nu}\T \Psi(\B{x}^*(t_f^*),t_f^*) \\ &- \B{\lambda}^{*\mbox{\tiny T}}\delta \B{x}\vert_{t_0}^{t_f^*} + \int_{t_0}^{t_f^*}\Big\{ \Big[\frac{\partial H}{\partial \B{x}^*} + \bdot{\lambda}^* \Big]\T \delta \B{x} + \frac{\partial H}{\partial \B{u}^*}\T \delta \B{u} + \Big[\frac{\partial H}{\partial \B{\lambda}^*} - \bdot{x}^* \Big]\T \delta \B{\lambda} \Big\}\dd t \\ &+ H(\B{x}^*, \B{u}^*, \B{\lambda}^*, t_f^*) \dd t_f - \B{\lambda}^{*\mbox{\tiny T}}(t_f^*) \bdot{x}^*(t_f^*) \dd t_f
\end{align*}
Now, to simplify the problem further, we need to consider the admissible variation of the state vector, $\delta \B{x}(t)$, shown in Figure \ref{fig:s4_skewVariation}
Here we can define the skew variation, $\dd \B{x}_f$, as,
\begin{equation*}
    \dd \B{x}_f = \delta \B{x}(t^*_f) + \Big(\bdot{x}^*(t_f) + \hspace{-15pt} \overbrace{\delta \bdot{x}(t_f) \Big)\dd t_f}^{\text{neglect second order term}}
\end{equation*}
which we simplify to,
\begin{equation*}
    \dd \B{x}_f = \delta \B{x}(t^*_f) + \bdot{x}^*(t_f) \dd t_f.
\end{equation*}
Using this relationship along with the fact that for $\B{x}(t_0) = \B{x}_0$ the variation of the state at the initial condition is equal to zero and thus $\B{\lambda}^{*\mbox{\tiny T}}\delta \B{x}\vert_{t_0}^{t_f^*} = \B{\lambda}^{*\mbox{\tiny T}}\delta \B{x}(t_f^*)$, we can simplify the expression of $\delta J_a$ to,
\begin{equation*} 
\begin{aligned}
    \delta J_a =& \Big[\frac{\partial \Phi}{\partial \B{x}^*(t_f^*)} + \frac{\partial \Psi}{\partial \B{x}^*(t_f^*)}\T \B{\nu} - \B{\lambda}^*(t_f^*)\Big]\T \dd \B{x}_f + \Big[\frac{\partial \Phi}{\partial t_f^*} + \B{\nu}\T \frac{\partial \Psi}{\partial t_f^*} + H(\B{x}^*, \B{u}^*, \B{\lambda}^*, t_f^*)\Big] \dd t_f \\ &+ \dd \B{\nu}\T \Psi(\B{x}^*(t_f^*),t_f^*) + \int_{t_0}^{t_f^*}\Big\{ \Big[\frac{\partial H}{\partial \B{x}^*} + \bdot{\lambda}^* \Big]\T \delta \B{x} + \frac{\partial H}{\partial \B{u}^*}\T \delta \B{u} + \Big[\frac{\partial H}{\partial \B{\lambda}^*} - \bdot{x}^* \Big]\T \delta \B{\lambda} \Big\}\dd t.
\end{aligned}
\end{equation*}
The extrema of this equation can be found by finding the conditions such that $\delta J_a$ is equal to zero. In order for this to occur, the square bracketed terms must go to zero. From this, we are lead to a set of equations that must be satisfied simultaneously that are referred to as the first-order necessary conditions for optimality (since they are based on the first variation of the augmented cost function). Dropping the (*) notation, the optimal solution is defined by the following set of differential and algebraic equations,
\begin{align}
    &\bdot{x} = \frac{\partial H}{\partial \B{\lambda}} \label{eq:s4_state}\\
    &\bdot{\lambda} = -\frac{\partial H}{\partial \B{x}} \label{eq:s4_costate}\\
    &\frac{\partial H}{\partial \B{u}} = \B{0} \label{eq:s4_control}\\
    &\Psi(\B{x}(t_f),t_f) = \B{0} \label{eq:s4_terminal}\\
    &\B{\lambda}(t_f) = \frac{\partial \Phi}{\partial \B{x}(t_f)} + \frac{\partial \Psi}{\partial \B{x}(t_f)}\T \B{\nu} \label{eq:s4_free_xf}\\
    &H(t_f) + \frac{\partial \Phi}{\partial t_f} + \B{\nu}\T \frac{\partial \Psi}{\partial t_f} = 0 \label{eq:s4_free_tf}
\end{align}
By looking at our definition of the Hamiltonian, Equation \eqref{eq:s4_Ham}, the first necessary condition simply reiterates the dynamics of the system, $\bdot{x} = f(\B{x},\B{u},t)$. Furthermore,
Equation \eqref{eq:s4_boundary1} constrains the initial values, Equations \eqref{eq:s4_costate} is a differential equation governing the costate values, and Equation \eqref{eq:s4_control} is the necessary condition for the control vector. Finally, Equations \eqref{eq:s4_terminal}, \eqref{eq:s4_free_xf}, and \eqref{eq:s4_free_tf} are necessary for the following cases and are sometimes referred to as transversality conditions,
\begin{itemize}
    \item For constraints on the final state and/or time, Equation \eqref{eq:s4_terminal} must be satisfied.
    \item For the components of $\B{x}(t_f)$ that are unconstrained (or free), Equation \eqref{eq:s4_free_xf} is used to determine the final value of the associated costate, i.e., $\B{\lambda}(t_f)$.
    \item For unconstrained (or free) final time, Equation \eqref{eq:s4_free_tf} must also be satisfied.
\end{itemize}
In all problems, Equations \eqref{eq:s4_state}, \eqref{eq:s4_costate}, and \eqref{eq:s4_control} will always be applicable, while Equations \eqref{eq:s4_terminal}, \eqref{eq:s4_free_xf}, and \eqref{eq:s4_free_tf} are dictated by the constraints of the final state and time according to the bullet points above.

\section{Addition of control inequality constraint}\label{sect:s4_ineq}
It can be seen from the prior section that the first-order necessary conditions derived from the indirect method rely on the formulation of the Hamiltonian, $H$, based on the cost function (a functional of $\Phi$ and $\mathcal{L}$), along with any terminal constraints ($\Psi$). In many problems, as is the case with the fuel-optimal landing problem explored in Chapter \ref{chap:fol}, it is necessary to constrain the control by some function. Therefore, consider the constraint,
\begin{equation*}
    \B{C}(\B{u}(t),t) \leq \B{0},
\end{equation*}
where $\B{C}$ is a vector function. The method to apply this constraint is to adjoin the constraint to Equation \eqref{eq:s4_Ham},
\begin{equation*}
    H = \mathcal{L} + \B{\lambda}\T \B{f} + \B{\mu}\T \B{C}
\end{equation*}
where $\B{\mu}$ are Lagrange multipliers that have the requirement,
\begin{equation*}
    \mu_i \quad \begin{cases} \leq 0, \quad C_i = 0, \\ = 0, \quad C_i < 0 \end{cases}
\end{equation*}
where $i$ denotes the specific constraint. By doing this, the only equation that changes in our prior derivation is Equation \eqref{eq:s4_control} since $\B{C}$ is a function of the control variable. It follows that,
\begin{equation}\label{eq:s4_controlCon}
    \B{0} = \frac{\partial H}{\partial \B{u}} = \frac{\partial \mathcal{L}}{\partial \B{u}} + \B{\lambda}\T \frac{\partial \B{f}}{\partial \B{u}} + \B{\mu}\T\frac{\partial \B{C}}{\partial \B{u}}.
\end{equation}
In general, Equation \eqref{eq:s4_controlCon} defines the set of conditions for the control based on the inequality constraints, $\B{C}$, and the state and costate values. We will revisit the application of inequality constraints in Chapter \ref{chap:fol}.

\section{Adjustment using the TFC approach and constrained expressions}

In general, through the indirect method, the optimal control problem is converted into a two-point boundary-value problem, Equations \eqref{eq:s4_state} and \eqref{eq:s4_costate}, with additional linear and nonlinear constraints, Equations \eqref{eq:s4_boundary1}, \eqref{eq:s4_control}, \eqref{eq:s4_terminal}, \eqref{eq:s4_free_xf}, and \eqref{eq:s4_free_tf}. In the case of control constraints described in Section \ref{sect:s4_ineq}, Equation \eqref{eq:s4_control} is replaced by Equation \eqref{eq:s4_controlCon}. In all, these equations represent the first-order necessary conditions that must be satisfied simultaneously.

As it should be clear from the development of the TFC approach in Sections \ref{chap:tfc_intro} and \ref{chap:tfc_general}, the benefit of this method is the ability to analytically embed linear constraints. Of the necessary conditions, the initial value constraint, Equation \eqref{eq:s4_boundary1}, and any linear terminal constraints, Equation \eqref{eq:s4_terminal}, can be easily embedded into a constrained expression for the state. To distinguish between the linear and nonlinear components of $\Psi$, let $\Psi$ be the composition of the linear and nonlinear portions,
\begin{align*}
    \Psi = \begin{Bmatrix} \Psi_{\ell} \\ \Psi_{n\ell} \end{Bmatrix} = \B{0} 
\end{align*}
where the linear terms $\Psi_{\ell}$ are embedded into the state constrained expressions and $\Psi_{n \ell}$ replaces the $\Psi$ term in Equations \eqref{eq:s4_terminal}, \eqref{eq:s4_free_xf}, and \eqref{eq:s4_free_tf}. Doing this reduces the length of the $\B{\nu}$ coefficient vector, and therefore, reduces the search space of the numerical optimization algorithm. However, in most cases, and both landing problems presented in this work, the terminal constraints are all linear, and thus, the $\Psi$ term can be eliminated. The result of the application of the TFC constrained expression is summarized in the following equations. 
\begin{align*}
    \bdot{x} = \frac{\partial H}{\partial \B{\lambda}}  \qquad &\longrightarrow& \bdot{x} = \frac{\partial H}{\partial \B{\lambda}}\\
    \bdot{\lambda} = -\frac{\partial H}{\partial \B{x}} \qquad &\longrightarrow& \bdot{\lambda} = -\frac{\partial H}{\partial \B{x}}\\
    \frac{\partial H}{\partial \B{u}} = \B{0}  \qquad &\longrightarrow& \frac{\partial H}{\partial \B{u}} = \B{0} \\
    \B{x}(t_0) = \B{x}_0 \qquad &\longrightarrow& -- \\
    \Psi(\B{x}(t_f),t_f) = \B{0}  \qquad &\longrightarrow& --\\
    \B{\lambda}(t_f) = \frac{\partial \Phi}{\partial \B{x}(t_f)} + \frac{\partial \Psi}{\partial \B{x}(t_f)}\T \B{\nu}  \qquad &\longrightarrow& \B{\lambda}(t_f) = \frac{\partial \Phi}{\partial \B{x}(t_f)} \\
    H(t_f) + \frac{\partial \Phi}{\partial t_f} + \B{\nu}\T \frac{\partial \Psi}{\partial t_f} = 0 \qquad &\longrightarrow& H(t_f) + \frac{\partial \Phi}{\partial t_f} = 0
\end{align*}

\section{Connection with the existing literature and difference between local and global collocation methods}
Over the past few decades, optimal control and trajectory optimization have been very active and interconnected fields of research. Solving optimal control problems is becoming increasingly important in developing G\&C algorithms that can effectively enable system autonomy and autonomous operations. Indeed, the recently coined term \textit{computational guidance and control} \cite{lu2017introducing} refers to a paradigm shift in which computation has a central role in defining and executing G\&C functions for aerospace systems. Newly defined algorithms tend to rely extensively on onboard computation, where numerical algorithms replace closed-loop G\&C and closed-loop predefined laws. Indeed, the vast majority of optimal control problems of interest for space systems do not have a closed-form solution and must rely on numerical methods. The latter are generally divided into two classes, i.e., direct and indirect methods.

Direct methods, sometimes called to as \textit{direct transcription methods} \cite{betts1998mesh}, refer to a class of numerical optimal control methodologies where the continuous optimal control problem is transcribed into an NLP optimization problem via proper approximation of the state and/or control. The most fundamental direct method is the single or multiple shooting method (e.g., \cite{diehl2006fast}), where the control is parametrized using a specified functional form, and the equations of motions are satisfied by direct integration. The resulting NLP minimizes the discretized cost function subject to path and/or interior-point constraints. 

In contrast, the alternative and more popular class of direct methods is the \textit{direct collocation} method. Here, both state and control are approximated using a defined functional form (e.g., a linear combination of Chebyshev polynomials). Such methods are generally divided into \textit{local and global collocation}. Local collocation divides the interval into many subintervals and enforces continuity across the interfaces. The resulting problem is further discretized using Runge-Kutta (implicit) methods (e.g., References \cite{schwartz1996consistent,hager2000runge}) or orthogonal collocation methods, where the collocation points are selected as roots of a family of orthogonal polynomials (e.g., References \cite{reddien1979collocation,herman1996direct}). Conversely, global collocation methods employ \textit{global} polynomials to approximate state and control with collocation executed at specified points across the desired time interval. 

The most popular set of global collocation methods for optimal control are named \textit{pseudospectral} methods. Indeed, there are different ways to approximate state and control. Historically, the first class of pseudospectral methods were developed by expanding state and control in a set of Chebyshev polynomials of degree $N$ \cite{fahroo2002direct,vlassenbroeck1988chebyshev}. Eventually, this approach was abandoned in favor of a linear combination of Lagrange polynomials using alternative collocation points such as Gauss-Lobatto \cite{elnagar1995pseudospectral} and Gauss-Lobatto-Radau \cite{fahroo2008pseudospectral,garg2011direct}. Such formulations were preferred mainly because the isolation condition was automatically satisfied \cite{rao2009survey} and yielded simpler conditions for collocation.

Many advancements have been made to develop both theory and practical implementation of pseudospectral methods for direct transcription of optimal control problems. Theoretical understanding in the convergence properties and connection with indirect methods \cite{gong2008connections,gong2008spectral,kang2007convergence,kang2008pseudospectral,hager2016convergence,hager2019convergence} coupled with pseudospectral algorithmic advancements to deal with a large class of smooth and non-smooth problems \cite{ross2004pseudospectral,fahroo2008pseudospectral,ross2003legendre,huntington2008comparison,agamawi2020mesh} has been paving the way to the potential application of such approaches for real-time implementation \cite{ross2008optimal,rao2018mesh,dennis2019computational}. Importantly, a new class of adaptive pseudospectral methods capable of automatically determining the number of segments and order of polynomial expansion has been recently developed \cite{darby2011hp}. Such an approach eventually led to the development of the GPOPS-II numerical platform \cite{patterson2014gpops}, which has been widely employed in trajectory optimization and control in a few applications such as low-thrust \cite{graham2015minimum}, solar sail \cite{peloni2018automated}, and rocket ascent \cite{miller2017rapid}. An in-depth review of pseudospectral methods applied to optimal control can be found in \cite{rao2009survey,ross2012review,rao2014trajectory}.

On the other end, indirect methods rely on developing the first-order necessary conditions by directly applying PMP or by the calculus of variations. The necessary conditions result in a TPBVP that must be generally resolved by application of numerical techniques such as single and multiple shooting methods \cite{keller1976numerical,stoer2013introduction}, orthogonal collocation \cite{oh1977use}, or pseudospectral methods \cite{fahroo2000trajectory}. The proposed method falls under this category, as the optimal guidance problem is cast as TPBVP that is solved via TFC.

At first glance, the proposed technique might seem similar to some of the above mentioned numerical schemes, namely, collocation methods \cite{ChebCol} and indirect pseudospectral methods \cite{fahroo2000trajectory}. This similarity is because the free function $g(t)$ is approximated using orthogonal polynomials discretized over the local or global domain, depending on the selected technique. However, there is a fundamental difference and a numerical benefit that the TFC approach adds, which is absent in previously developed techniques. For example, in indirect orthogonal collocation methods, the state and costates are parameterized using piecewise polynomial functions, transforming the problem into a nonlinear system of equations that must be solved.

Similarly, in indirect pseudospectral methods, the global spectral approach mandates that the state and costate are expanded via some basis functions. While it is true that the function $g(t)$ may be defined in the same fashion, the fundamental difference lies in how the TFC approach handles the problem's constraints: by analytically embedding them through the use of constrained expressions. In both local and global spectral methods, such constraints become part of the optimization scheme. In contrast, the TFC approach analytically reduces the search space of the solution to those that only satisfy the constraints. As a result, a simpler optimization scheme can be employed to find the solution.

To further highlight the differences, consider the differential equation to be solved in Equation \eqref{eq:s4_DEexample},
\begin{equation}\label{eq:s4_DEexample}
    F\left(t,y,\dot{y},\ddot{y}\right) = 0 \qquad \text{subject to:} \quad  \begin{cases} y(t_0) = y_{0}\\ \dot{y}(t_0) = \dot{y}_{0}\\ y(t_f) = y_{f}\\ \dot{y}(t_f) = \dot{y}_{f} \end{cases}
\end{equation}
using the spectral method. Let the function $y(t)$ be defined in the same way as the $g(t)$ function in the TFC formulation such that,
\begin{equation*}
    y(t) = \B{\zeta} \T \B{h}(z).
\end{equation*}
The key difference is that this description does not satisfy the constraints which must be enforced by the following equations,
\begin{align*}
    y(t_0) &= y_{0} = \B{\zeta} \T \B{h}(z_0)\\
    y(t_f) &= y_{f} = \B{\zeta} \T \B{h}(z_f)\\
    \dot{y}(t_0) &= \dot{y}_{0} = \B{\zeta} \T c \B{h}_{z}(z_0)\\
    \dot{y}(t_f) &= \dot{y}_{f} = \B{\zeta} \T c \B{h}_{z}(z_f).
\end{align*}
Then, to solve the problem, these equations must be appended to the residual of the differential equation,
\begin{equation*}
    F(t,\B{\zeta}) = 0
\end{equation*}
Notice that to solve for the unknown $\B{\zeta}$ coefficient vector, all five equations must be solved simultaneously. In other words, the solution of the constraints are now coupled to the solution of the dynamics, and the coefficients of $\B{\zeta}$ contribute to the constraint satisfaction, which will have numerical approximation error. Therefore, $\B{\xi}$ from the TFC development is not the same as $\B{\zeta}$ defined through the spectral method.

It should now be clear that the major novelty when solving optimal control problems is the analytical constraint satisfaction that reduces the system of equations. Since this technique is applied before numerically approximating the solution using orthogonal polynomials, there is no numerical error associated with enforcing the boundary conditions. Importantly, the constraints and dynamics are decoupled. Additionally, the constraint satisfaction is independent of how $g(t)$ is expressed, and therefore, the proposed formulation allows for a wide range of potential approximation of the free function. It is worth noting that in pseudospectral optimal control, the selection of the weighted interpolating functions is essential for convergence, and such functions are intimately connected with the problem's boundary conditions \cite{fahroo2008pseudospectral,fahroo2008advances,ross2012review}. The TFC approach decouples the two problems and only relies on the convergence properties of the selected family of functions which approximate $g(t)$.
%

\chapter{ENERGY-OPTIMAL LANDING\label{chap:eol}}

The energy-optimal landing problem is an important step in our study of the TFC method for real-time optimal control. While mathematically simpler than the fuel-optimal landing problem, it provides a real problem for testing the algorithms. In the simplest formulation, the acceleration due to gravity, $\B{a}_g$, is considered constant (which is the cause for the terminal descent phase of landing). For this case, a feedback solution can be derived based on the calculated time-to-go function and can be solved for a problem formulated in state-space (as is the following example) \cite{EOL_chris}. The feedback law is defined as,
\begin{equation*}
    \B{u} = -\frac{6}{t^2_{go}} \B{r} - \frac{4}{t_{go}}\B{v} - \B{a}_g
\end{equation*}
where $\B{u}$ is the control acceleration, $\B{r}$ and $\B{v}$ are the position and velocity states respectively, and time-to-go is $t_{go}$.

Conversely, another feedback solution exists (although not used in this work) for this problem called Zero-Effort-Miss/Zero-Effort-Velocity (ZEM/ZEV) \cite{z1,z2,ZEM_ZEV}. In this approach, ZEM is the associated error in the final distance to the landing site if no control action is taken, and ZEV is the error on final velocity again under zero control effort. This formulation collapses to the expression,
\begin{equation*}
    \B{u} = \frac{6}{t^2_{go}}\textbf{ZEM} - \frac{2}{t_{go}}\textbf{ZEV}
\end{equation*}

Moving forward, it is important to know that TFC is by nature an open-loop optimal controller since, in practice, the problem would be solved at every computer cycle to update the trajectory. The feedback solution is only valid for a constant gravity vector, $\B{a}_g$; however, since the TFC development is general, it can be easily adjusted to solve for any gravitational model.

Although not presented here, the interested reader is directed to the application of this technique to both small and large planetary bodies presented in Reference \cite{EOL_COMP}.

\section{Dynamical model}
For the problem of energy-optimal pinpoint landing on large bodies (e.g., the Moon or Mars) the governing system dynamics during the powered descent phase can be modeled as follows,
\begin{align*}
  \bdot{r} &= \,\B{v} \\ 
  \bdot{v} &= \,\B{a}_g + \B{u},
\end{align*}
where $\B{r}$ and $\B{v}$ are position and velocity vectors of the lander with respect to the landing site fixed frame. Additionally, $\B{u} = \frac{T}{m}$ is associated with the thrust acceleration of the lander and is used to determine the thrust control $T$ for the current spacecraft mass $m$. The dynamics of the mass state are governed by the equation,
\begin{equation*}
    \dot{m} = - \alpha \, T
\end{equation*}
where $\alpha = 1 / v_{ex}$, with $v_{ex}$ being the effective exhaust velocity of the rocket engine. However, since the mass dynamics are independent of the spacecraft position and velocity, and the spacecraft acceleration is the control variable, the mass state and, in turn, the thrust value can are determined after the optimal trajectory is computed. Furthermore, acceleration due to gravity, $\B{a}_g$, is considered constant since this problem deals with the terminal descent phase. For this problem, the initial and final position and velocity, and initial mass are given:
\begin{equation*}
    \begin{cases}\B{r}(0) = \B{r}_0 \\ \B{v}(0) = \B{v}_0\end{cases}, \qquad \begin{cases} \B{r}(t_f) = \B{r}_f \\ \B{v}(t_f) = \B{v}_f\end{cases}.
\end{equation*}
The objective is to minimize the energy, which can be realized by minimizing the control used while satisfying the problem's dynamics constraints. Therefore, the problem can be posed as,
\begin{wblankBox}{Optimization problem statement}
\begin{equation*}
    \underset{t_f,\B{u}}{\text{minimize}} \quad \Gamma t_f + \frac{1}{2} \int_{t_0}^{t_f} \B{u}\T \B{u} \dd \tau
\end{equation*}
\begin{equation*}
\begin{aligned}
\text{subject to} \quad &\dot{\B{r}} = \,\B{v}, \quad \dot{\B{v}} = \,\B{a}_g + \B{u}, \\
& \B{r}(0) = \B{r}_0, \quad \B{v}(0) = \B{v}_0, \\
& \B{r}(t_f) = \B{r}_f, \quad \B{v}(t_f) = \B{v}_f
\end{aligned}
\end{equation*}
\end{wblankBox}
\noindent where $\Phi(t_f) = \Gamma \, t_f$ is the terminal cost parameter for the final time. $\Gamma$ is a scalar weight parameter on the final time and represents a trade-off between the minimum-time and minimum-energy problem. For example, if $\Gamma = 0$, we recover the minimum energy cost function.

\section{First-order necessary conditions}
Applying the PMP, the Hamiltonian takes the following form,
\begin{equation*}
    H = \mathcal{L} + \B{\lambda}\T \B{f}
\end{equation*}
which can be expanded as,
\begin{equation*} 
    H = \frac{1}{2} \B{u}\T \B{u} + \B{\lambda}_r\T \B{v} + \B{\lambda}_v\T (\B{a}_g + \B{u}).
\end{equation*}
Applying the first-order necessary conditions, the optimal control action is realized by,
\begin{equation*}
    \frac{\partial{ H}}{\partial \B{u}} = \B{u} + \B{\lambda}_v = 0 \quad \longrightarrow \B{u} = -\B{\lambda}_v.
\end{equation*}
It can be seen that the vector $\B{u}$ is opposite of the costate $\B{\lambda}_v$, and therefore, we can replace this costate term directly with the control in all following equations. The additional first-order conditions lead to
\begin{equation*}
    \begin{aligned}
        \dot{\B{r}} &= \frac{\partial{ H}}{\partial \B{\lambda}_r} = \B{v}  \qquad &&\dot{\B{\lambda}}_r = -\frac{\partial{ H}}{\partial \B{r}} = \B{0}\\
        \dot{\B{v}} &= \frac{\partial{ H}}{\partial \B{\lambda}_v} = \B{a}_g + \B{u} \qquad &&\dot{\B{\lambda}}_v = -\frac{\partial{ H}}{\partial \B{v}} = -\B{\lambda}_r\\
    \end{aligned}
\end{equation*}
and as mentioned, the differential equation associated with $\dot{\B{\lambda}}_v$ can be written as,
\begin{equation*}
    \bdot{u} = \B{\lambda}_r.
\end{equation*}
Lastly, since the problem is posed as a free final time problem, the transversality condition is given by,
\begin{equation*}
    H(t_f) + \frac{\partial \Phi}{\partial t_f} = 0
\end{equation*}
which reduces to
\begin{equation*}
    H(t_f) = - \Gamma.
\end{equation*}
Collecting all equations, a constrained, differential systems of equations is formed which must be satisfied simultaneous to obtain an optimal solution,
\begin{blankBox}{First-order necessary conditions}
\begin{eqnarray}
    \dot{\B{r}} &=& \B{v}\label{eq:s5_rdot}\\
    \dot{\B{v}} &=& \B{a}_g + \B{u}\label{eq:s5_vdot}\\
    \dot{\B{\lambda}}_r &=& \B{0}\label{eq:s5_lamrdot}\\
    \dot{\B{u}} &=& \B{\lambda}_r\label{eq:s5_lamvdot}\\
    H(t_f) + \Gamma &=& 0 \label{eq:s5_htf}
\end{eqnarray}
subject to the constraints $\B{r}(t_0) = \B{r}_0$, $\B{v}(t_0) = \B{v}_0$, $\B{r}(t_f) = \B{r}_f$, and $\B{v}(t_f) = \B{v}_f$. 
\end{blankBox}

The following section reformulates the system of equations defined by Equations \eqref{eq:s5_rdot}, \eqref{eq:s5_vdot}, \eqref{eq:s5_lamrdot}, \eqref{eq:s5_lamvdot}, \eqref{eq:s5_htf} using the techniques developed in the prior sections.

\section{Solving the problem via the TFC}
Through the use of TFC, Equations \eqref{eq:s5_rdot} through \eqref{eq:s5_htf} can be reduced. First, using the TFC approach, Equation \eqref{eq:s5_rdot} is redundant, since the constrained expression will always satisfy this condition. Furthermore, Equations \eqref{eq:s5_lamrdot} and \eqref{eq:s5_lamvdot} can be combined since Equation \eqref{eq:s5_lamrdot} shows $\B{\lambda}_r$ must be constant. Therefore, these two equations can be replaced by the equation,
\begin{equation}\label{eq:s5_controlCE}
   u_i(t,\B{\xi}_{u_i}) = \B{h}_{u}\T \B{\xi}_{u_i}, \quad \text{for} \quad i = 1,2,3
\end{equation}
where $\B{h}_{u}$ consists of the constant and linear terms of the selected basis set. Lastly, the boundary constraints are fully handled by the TFC constrained expressions of the following form,
\begin{align}\label{eq:s5_stateCE}
    r_i(t,\B{\xi}_i) = \Big(\B{h}(z) - \phi_1(t) \B{h}(z_0) - \phi_2(t) \B{h}(z_f) &- \phi_3(t) c \B{h}_z(z_0) - \phi_4(t) c \B{h}_z(z_f) \Big)\T \B{\xi}_i \nonumber\\&+ \phi_1(t) r_{0_i} + \phi_2(t) r_{f_i} + \phi_3(t) v_{0_i} + \phi_4(t) v_{f_i}.
\end{align}
Therefore, the first-order necessary conditions reduce to,
\begin{eqnarray}
    \dot{v}_i &=& a_{g_i} + u_i \label{eq:s5_vdot_new}\\
    0 &=& -\frac{1}{2} u_i^2(t_f) + u_i(t_f) a_{g_i} + \Gamma, \label{eq:s5_htf_new}
\end{eqnarray}
where the state and control are written in terms of the TFC constrained expressions. In general, the unknowns of this system are the coefficients related to the state $\B{\xi}_i$ and control $\B{\xi}_{u}$ along with the final time $t_f$. Both state and control unknowns appear linearly in the system of equations; however, the final time appears nonlinearly through the transversality equation, Equation \eqref{eq:s5_htf_new}, and can be handled in two different ways. The first method uses an Outer-loop optimizer that solves for the mapping parameter, i.e., optimizes the final time with the transversality condition. In contrast, the inner TFC loop solves the least-squares problem of Equation \eqref{eq:s5_vdot_new}. The second method leverages the theory developed in Section \ref{sec:s3_freeTime}, where the mapping parameter (which is a function of $t_f$) is solved alongside the other unknowns in a single loop. This method, however, requires an implementation of a nonlinear least-squares approach. While this section has merely summarized the relevant equations, Sections \ref{sec:s5_outLoop} and \ref{sec:s5_singLoop} discuss in detail how each method can be applied to the energy-optimal landing problem. Lastly, various tests are conducted to determine the accuracy, speed, and robustness of both techniques. The findings of these tests will help us in our study of the more complex problem of fuel-optimal landing in the following chapter. 

\subsection{Outer-loop optimizer}\label{sec:s5_outLoop}

Using the constrained expression given by Equation \eqref{eq:s5_stateCE}, for the Outer-loop method, the constrained expression is written in the problem domain (i.e., in terms of time), and thus, the switching functions are,
\begin{align*}
	\phi_1(t) &= \frac{1}{\Delta t^3} \Big( -t_f^2 (3  t_0- t_f) + 6  t_0  t_f t  -3 ( t_0+ t_f) t^2 + 2t^3 \Big) \\
	\phi_2(t) &= \frac{1}{\Delta t^3} \Big( - t_0^2 ( t_0-3  t_f) -6  t_0  t_f t  + 3 ( t_0+ t_f)t^2 -2t^3 \Big) \\
	\phi_3(t) &= \frac{1}{\Delta t^2} \Big(-t_0  t_f^2 +  t_f (2  t_0+ t_f)t -(t_0+2  t_f)t^2 + t^3 \Big) \\
	\phi_4(t) &= \frac{1}{\Delta t^2} \Big( - t_0^2  t_f + t_0 ( t_0+2  t_f) t - (2  t_0+ t_f) t^2 + t^3 \Big) 
\end{align*}
Equation \eqref{eq:s5_stateCE} and its derivatives for this method are,
\begin{align}
    r_i(t,\B{\xi}_i) = \Big(\B{h}(z) - \phi_1 \B{h}(z_0) - \phi_2 \B{h}(z_f) &- \phi_3 c\B{h}_z(z_0) - \phi_4 c\B{h}_z(z_f) \Big)\T \B{\xi}_i \nonumber\\&+ \phi_1 r_{0_i} + \phi_2 r_{f_i} + \phi_3 v_{0_i} + \phi_4 v_{f_i}.\label{eq:s5_r_timeCE} 
\end{align}
Substituting Equation \eqref{eq:s5_r_timeCE} and its second derivative, i.e. acceleration, and the definition of the control, Equation \eqref{eq:s5_controlCE}, into Equation \eqref{eq:s5_vdot_new}, the loss functions becomes,
\begin{align}
    F_i(t,\Xi) = \Big(c^2 \B{h}_{zz}(z) &- \ddot{\phi}_1(t) \B{h}(z_0) - \ddot{\phi}_2(t) \B{h}(z_f) - \ddot{\phi}_3(t) \B{h}_z(z_0) - \ddot{\phi}_4(t) \B{h}_z(z_f) \Big)\T \B{\xi}_i \nonumber \\ &+ \ddot{\phi}_1(t) r_{0_i} + \ddot{\phi}_2(t) r_{f_i} + \ddot{\phi}_3(t) v_{0_i} + \ddot{\phi}_4(t) v_{f_i} - a_{g_i} - \B{h}_{u}\T \B{\xi}_{u_i} = 0 \label{eq:s5_singleLoop_lossFreeTime}
\end{align}
where the loss vector becomes,
\begin{align*}
    \mathbb{L} &= \begin{Bmatrix} \mathbb{L}_1\T & \mathbb{L}_2\T &\mathbb{L}_3\T\end{Bmatrix}\T_{3N \times 1},
\end{align*}
where
\begin{equation*}
    \mathbb{L}_i = \begin{Bmatrix}
    F_i(t_0,\Xi) & \hdots & F_i(t_k,\Xi) & \hdots & 
    F_i(t_f,\Xi)
    \end{Bmatrix}\T,
\end{equation*}
and the unknown vector is then,
\begin{align*}
    \Xi &= \begin{Bmatrix} \B{\xi}_1\T & \B{\xi}_2\T & \B{\xi}_3\T & \B{\xi}_{u_1}\T & \B{\xi}_{u_2}\T & \B{\xi}_{u_3}\T \end{Bmatrix}\T _{(3m + 6) \times 1}.
\end{align*}

It should be seen that the loss function, given by Equation \eqref{eq:s5_singleLoop_lossFreeTime}, is linear, and therefore the loss vector is a linear system of equations. The terms of this linear system are provided in Appendix \ref{sect:app_sec:s5_outLoop}. Additionally, given this linear system, any available least-squares technique can be used to solve for the unknown coefficients. Next, once these coefficients are solved, Equation \eqref{eq:s5_htf_new} is enforced using any available root solving technique (the numerical results used \verb"NumPy's" $fsolve()$ algorithm). This process is repeated until the tolerance on the inner and outer residuals are met. 

\subsection{Single-loop approach}\label{sec:s5_singLoop}
For the single-loop approach using TFC, we take advantage of the fact that the mapping coefficient, of Equation \eqref{eq:s3_linearMapping}, is a function of the final time $t_f$. Next, the parameter is redefined such that it cannot be negative: $b^2 := c$. Then, by converting the dynamics and constraints into the basis function domain, i.e., $z \in [-1,1]$ for Chebyshev and Legendre polynomials, this parameter can simply be included in the optimization loop and solved simultaneously with the $\B{\xi}_i$ and $\B{\xi}_{u_i}$ coefficients. However, in all cases, an unknown final time will appear nonlinearly, and therefore, a nonlinear least-squares will be required regardless of whether or not the original system is linear.

The first step in this method is to write the whole problem in the basis function domain. This in turn will introduce new switching functions (and for clarity will be labeled as $\phiz{z}{\phi}$), which are,
\begin{align*}
	\phiz{z}{\phi}_1(z) &= \frac{1}{\Delta z^3} \Big( -z_f^2 (3  z_0- z_f) + 6  z_0  z_f z  -3 ( z_0+ z_f) z^2 + 2z^3 \Big)\\
	\phiz{z}{\phi}_2(z) &= \frac{1}{\Delta z^3} \Big( - z_0^2 ( z_0-3  z_f) -6  z_0  z_f z  + 3 ( z_0+ z_f)z^2 -2z^3 \Big)\\
	\phiz{z}{\phi}_3(z) &= \frac{1}{\Delta z^2} \Big(-z_0  z_f^2 +  z_f (2  z_0+ z_f)z -(z_0+2  z_f)z^2 + z^3 \Big)\\
	\phiz{z}{\phi}_4(z) &= \frac{1}{\Delta z^2} \Big( - z_0^2  z_f + z_0 ( z_0+2  z_f) z - (2  z_0+ z_f) z^2 + z^3 \Big),
\end{align*}
such that $\Delta z := z_f - z_0$. This change is also reflected in the constrained expression for the state,
\begin{align}
    r_i(z,\B{\xi}) = \Big(\B{h}(z) - \phiz{z}{\phi}_1 \B{h}(z_0) - \phiz{z}{\phi}_2 \B{h}(z_f) &- \phiz{z}{\phi}_3 \B{h}_z(z_0) - \phiz{z}{\phi}_4 \B{h}_z(z_f) \Big)\T \B{\xi}_i \nonumber \\ &+ \phiz{z}{\phi}_1 r_{0_i} + \phiz{z}{\phi}_2 r_{f_i} + \phiz{z}{\phi}_3 \frac{v_{0_i}}{b^2} + \phiz{z}{\phi}_4 \frac{v_{f_i}}{b^2}. \label{eq:s5_r_zCE}
\end{align}
Hence, the need to divide the velocity constraints by the modified mapping parameter, $b$ in Equations \eqref{eq:s5_r_zCE}. Next, our definition of $\B{u}$ remains unchanged and is defined by Equation \eqref{eq:s5_controlCE}. Following the current definition of the state and costate, the differential equation of Equation \eqref{eq:s5_vdot_new} becomes,
\begin{align*}
    F_i(z,\Xi) = b^4\Big[ \Big(\B{h}_{zz}(z) &- \phiz{z}{\phi}_{1_{zz}} \B{h}(z_0) - \phiz{z}{\phi}_{2_{zz}} \B{h}(z_f) - \phiz{z}{\phi}_{3_{zz}} \B{h}_z(z_0) - \phiz{z}{\phi}_{4_{zz}} \B{h}_z(z_f) \Big)\T \B{\xi}_i \nonumber \\ &+ \phiz{z}{\phi}_{1_{zz}} r_{0_i} + \phiz{z}{\phi}_{2_{zz}} r_{f_i} + \phiz{z}{\phi}_{3_{zz}} \frac{v_{0_i}}{b^2} + \phiz{z}{\phi}_{4_{zz}} \frac{v_{f_i}}{{b^2}} \Big] - a_{g_i} - \B{h}_{u}\T \B{\xi}_{u_i} = 0
\end{align*}
and the loss function associated with Equation \eqref{eq:s5_htf_new} can be written in terms of the unknowns as,
\begin{equation*} 
    \mathbb{L}_H(\Xi) = -\frac{1}{2} \sum_{j=1}^3\Big(u^2_j(z_f) \Big)  +  \sum_{j=1}^3\Big( u_j(z_f) a_{g_j} \Big) + \Gamma = 0
\end{equation*}
with the augmented loss function 
\begin{align*}
    \mathbb{L} &= \begin{Bmatrix} \mathbb{L}_1\T & \mathbb{L}_2\T &\mathbb{L}_3\T & \mathbb{L}_H\end{Bmatrix}\T_{(3N + 1) \times 1}
\end{align*}
where
\begin{equation*}
    \mathbb{L}_i = \begin{Bmatrix}
    F_i(z_0,\Xi) & \hdots & F_i(z_k,\Xi) & \hdots & 
    F_i(z_f,\Xi)
    \end{Bmatrix}\T.
\end{equation*}
The unknown vector is then,
\begin{align*}
    \Xi &= \begin{Bmatrix} \B{\xi}_1\T & \B{\xi}_2\T & \B{\xi}_3\T & \B{\xi}_{u_1}\T & \B{\xi}_{u_2}\T & \B{\xi}_{u_3}\T &  b \end{Bmatrix}\T _{(3m + 7) \times 1}\\ 
\end{align*}
The partial derivatives of the loss functions are provided in Appendix \ref{sect:app_sec:s5_singLoop}, and nonlinear least-squares is used to update the unknowns.

\section{Parameter initialization}
Finally, the last consideration before solving the problem using either method is to initialize the unknown parameters. In the Outer-loop method detailed in Section \ref{sec:s5_outLoop}, the inner-loop is a linear system, and therefore, $\B{\xi}_i$ and $\B{\xi}_{u_i}$ do not need to be initialized. However, an estimate of the final time $t_f$ is needed: for all numerical tests, this value was chosen to be one in the scaled time.

Next, for the single-loop method, all variables must be initialized since the system is nonlinear. As observed in the earlier section, the simplest initialization of the unknowns associated with the state constrained expression is to set them equal to zero. This is equivalent to connecting the boundary value problem with the simplest interpolating expression (i.e,. $g(x) = 0$). 

Although the initialization scheme for the unknowns associated with the control expression could follow the same process, more is known about their potential solution, which can be leveraged. Following this thought, we can initialize the parameters assuming the initial control is opposite the spacecraft velocity,
\begin{equation*}
    \B{u}_0 = -\dfrac{\B{v}_0}{||\B{v}_0||}
\end{equation*}
and the final control value is assumed to be in the direction opposite of the initial position vector,
\begin{equation*}
    \B{u}_f = -\dfrac{\B{r}_0}{||\B{r}_0||}
\end{equation*}
Using these two equations, the values of $\B{\xi}_{u_i} = \begin{Bmatrix} a_{0_i} & a_{1_i} \end{Bmatrix}\T $ become,
\begin{eqnarray*}
    a_{0_i} - a_{1_i} &= -\dfrac{v_{0_i}}{ \Big(\sum_{j=1}^3 v_{0_i}^2\Big)^{1/2} } = -\mathbb{V}_{0_i} \\
    a_{0_i} + a_{1_i} &= -\dfrac{r_{0_i}}{\Big(\sum_{j=1}^3 r_{0_i}^2\Big)^{1/2} }  = -\mathbb{R}_{0_i}.
\end{eqnarray*}
Solving this linear system yields,
\begin{eqnarray*}
    a_{0_i} &= -\dfrac{1}{2} (\mathbb{R}_{0_i} + \mathbb{V}_{0_i}) \\
    a_{1_i} &= -\dfrac{1}{2} (\mathbb{R}_{0_i} - \mathbb{V}_{0_i}).
\end{eqnarray*}

\section{Results}
First, the two proposed methods are compared to the known feedback solution presented in Reference \cite{EOL_chris} to validate the TFC method's accuracy. After this, a Monte Carlo simulation is constructed to test the Single-loop and Outer-loop method over a range of initial conditions to determine expected speed and robustness. 

For the numerical test presented in this section, the problem was scaled by the initial conditions. The unit length, $\ell^*$, and unit time, $t^*$, where calculated by the following equations, 
\begin{align*}
    \ell^* &= \max\left(|\B{r}_0|\right) \\
    t^* &= \frac{\ell^*}{\max\left(|\B{v}_0|\right)}.
\end{align*}

\begin{example}{Comparison to known feedback solution}

\begin{table}[H]
    \caption{Problem parameters for numerical test.}
    \begin{subtable}{.5\linewidth}
      \centering
        \caption{Problem specific values.}
        \begin{tabular}{|c||c|}
    \hline
    Variable & Value \\
    \hline\hline
    $\B{r}_0$ & $\{500000, 100000, 50000\}\T$ [ft] \\\hline
    $\B{v}_0$ & $\{-3000, 0, 0\}\T$ [ft/s] \\\hline
    $\B{a}_g$ & $\{0, 0, -5.31\}\T$ [ft/s$^2$]\\\hline
    $\Gamma$ & $0$ and $100$ \\\hline
    \hline
    \end{tabular}
    \end{subtable}%
    \begin{subtable}{.5\linewidth}
      \centering
        \caption{TFC parameters.}
        \begin{tabular}{|c||c|}
    \hline
    Variable & Value \\
    \hline\hline
    $N$ [\# points] & $100$ \\\hline
    $m$ [\# basis functions] & $60$ \\\hline
    $\varepsilon$ [tolerance] & $2.22 \times 10^{-16} $\\\hline
    \hline
    \end{tabular}
    \end{subtable} 
\end{table}

This comparison test shows that both TFC based methods solve the problem with almost identical results to the feedback solution. Furthermore, the TFC based method produces identical results regardless of whether the Outer-loop or Single-loop method is used. The tingle-loop method is an order of magnitude faster. Compared to the spectral method, for these specific test cases, TFC is slightly slower with regard to computation time. Looking at Tables \ref{tab:EOL_compare_Gamma0} and \ref{tab:EOL_compare_Gamma100}, the difference is close to 10 milliseconds.

\begin{table}[H]
\caption{Single case energy-optimal landing for $\Gamma = 0$.}
\centering
\begin{tabular}{|c||c|c|c|c|c|}
\hline
Parameter & \multicolumn{2}{c|}{TFC} & \multicolumn{2}{c|}{Spectral} & Feedback\\
\hline\hline
Method & Outer-loop & Single-loop & Outer-loop & Single-loop & $-$ \\\hline 
$t_f$ [sec] & 406.03 & 406.03 & 406.03 & 406.04 & 406.03\\\hline 
Cost  & 19004.12 & 19004.12 & 19004.12 & 19004.12 & 19004.19\\\hline
Comp. Time [s]  & 2.63 & 0.097 & 2.00 & 0.087 & $-$\\\hline
Iterations & 19 & 15 & 19 & 11 & $-$ \\\hline
$\max|\mathbb{L}|$ & $3.3 \times 10^{-16}$ & $2.2 \times 10^{-16}$ & $4.2 \times 10^{-16}$ & $4.3 \times 10^{-16}$ & $-$ \\\hline
\hline 
\end{tabular}
\label{tab:EOL_compare_Gamma0}
\end{table}

However, this test did not provide the full picture and was used as a first test to compare all of the results with those published in Reference \cite{EOL_chris}. Regardless, this test shows that the Single-loop approach should be the focus of further testing where TFC and the spectral method are used to solve the problem over a wide range of initial conditions.

\begin{table}[H]
\caption{Single case energy-optimal landing for $\Gamma = 100$.}
\centering
\begin{tabular}{|c||c|c|c|c|c|}
\hline
Parameter & \multicolumn{2}{c|}{TFC} & \multicolumn{2}{c|}{Spectral} & Feedback\\
\hline\hline
Method & Outer-loop & Single-loop & Outer-loop & Single-loop & $-$ \\\hline 
$t_f$ [sec] & 301.05 & 301.05 & 301.05 & 301.05 & 301.05\\\hline 
Cost  & 52569.32 & 52569.32 & 52569.32 & 52569.32 & 52568.53\\\hline
Comp. Time [s]  & 2.64 & 0.124 & 1.92 & 0.111 & $-$\\\hline
Iterations & 16 & 18 & 16 & 14 & $-$ \\\hline
$\max|\mathbb{L}|$ & $4.4 \times 10^{-16}$ & $2.8 \times 10^{-16}$ & $4.4 \times 10^{-16}$ & $3.3 \times 10^{-16}$ & $-$ \\\hline
\hline 
\end{tabular}
\label{tab:EOL_compare_Gamma100}
\end{table}
\end{example}

\begin{example}{Monte Carlos simulation for varying initial conditions}

\noindent \textbf{Test Setup:}
For this test, the TFC parameters and acceleration due to gravity remained the same as the example above. Furthermore, for the following test, we have only considered the pure energy-optimal problem such that $\Gamma=0$. 

Next, to span a large variety of initial conditions, the following process was used to define an initial position ellipse and associated velocity. Recall, the equation of the radius of an ellipse is defined as,
\begin{equation*}
    R_{\text{ellipse}} = \frac{ab}{\sqrt{a^2 \sin^2(\alpha) + b^2 \cos^2(\alpha)}}
\end{equation*}
In our case, we define these parameters,
\begin{equation*}
    \alpha = \mathcal{U}(0,2\pi), \quad
    a = 1000 \, \text{[m]}, \quad
    b = 500 \, \text{[m]}.
\end{equation*}
Using this, the sample ellipse was centered $2,000$ meters up-range and with an elevation of $1,500$ meters. This is simply the point $(-2000, 0, 1500)$. Finally, our initial conditions can be written using the following equations, where $\text{SF} = \mathcal{U}(0,1)$ is a scale factor used to span the whole area of the ellipse.
\begin{equation*}
    \B{r}_0 = \begin{cases} -2,000 + \text{SF} \cdot R_{\text{ellipse}} \cos(\alpha) \\ \text{SF} \cdot  R_{\text{ellipse}} \sin(\alpha) \\
    1,500 + \mathcal{U}(-100,100) \end{cases} \text{[m]} \quad
    \B{v}_0 = \begin{cases} 100 \cos(\beta) \\ 100 \sin(\beta) \\ -75 + \mathcal{U}(-10,10)  \end{cases} \text{[m/s]} 
\end{equation*}
where $\beta = \mathcal{U}\left(-\frac{\pi}{2},\frac{\pi}{2}\right)$. The following results compare the accuracy, speed, and robustness of the Single-loop approach for both TFC and spectral method. Recall, the Single-loop method solves all first-order necessary conditions simultaneously, albeit forcing the method to become nonlinear.
\begin{figure}[H]
    \centering\includegraphics[width=0.7\linewidth]{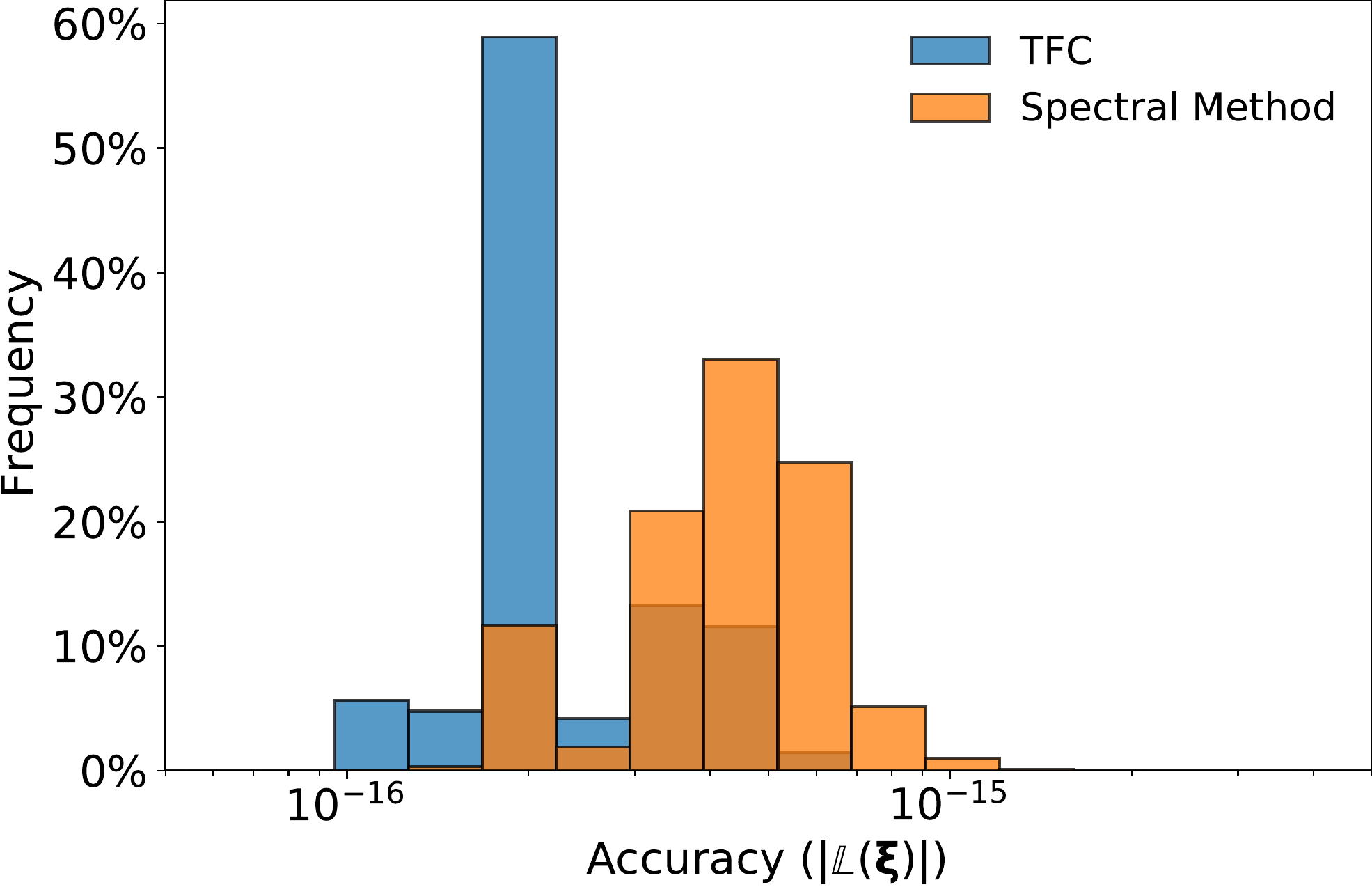}
    \caption{Histogram of the maximum residual of the loss vector.}
    \label{fig:s3_EOL_hist_L2}
\end{figure}
Over the 10,000 trial Monte Carlo simulation, the TFC method failed on three accounts or $0.03$\% of the time, and the spectral method failed 279 times or $2.79$\% of the time. Figure \ref{fig:s3_EOL_hist_L2} is a histogram of the methods' error, which shows that TFC is consistently more accurate.
\begin{figure}[H]
    \centering\includegraphics[width=0.7\linewidth]{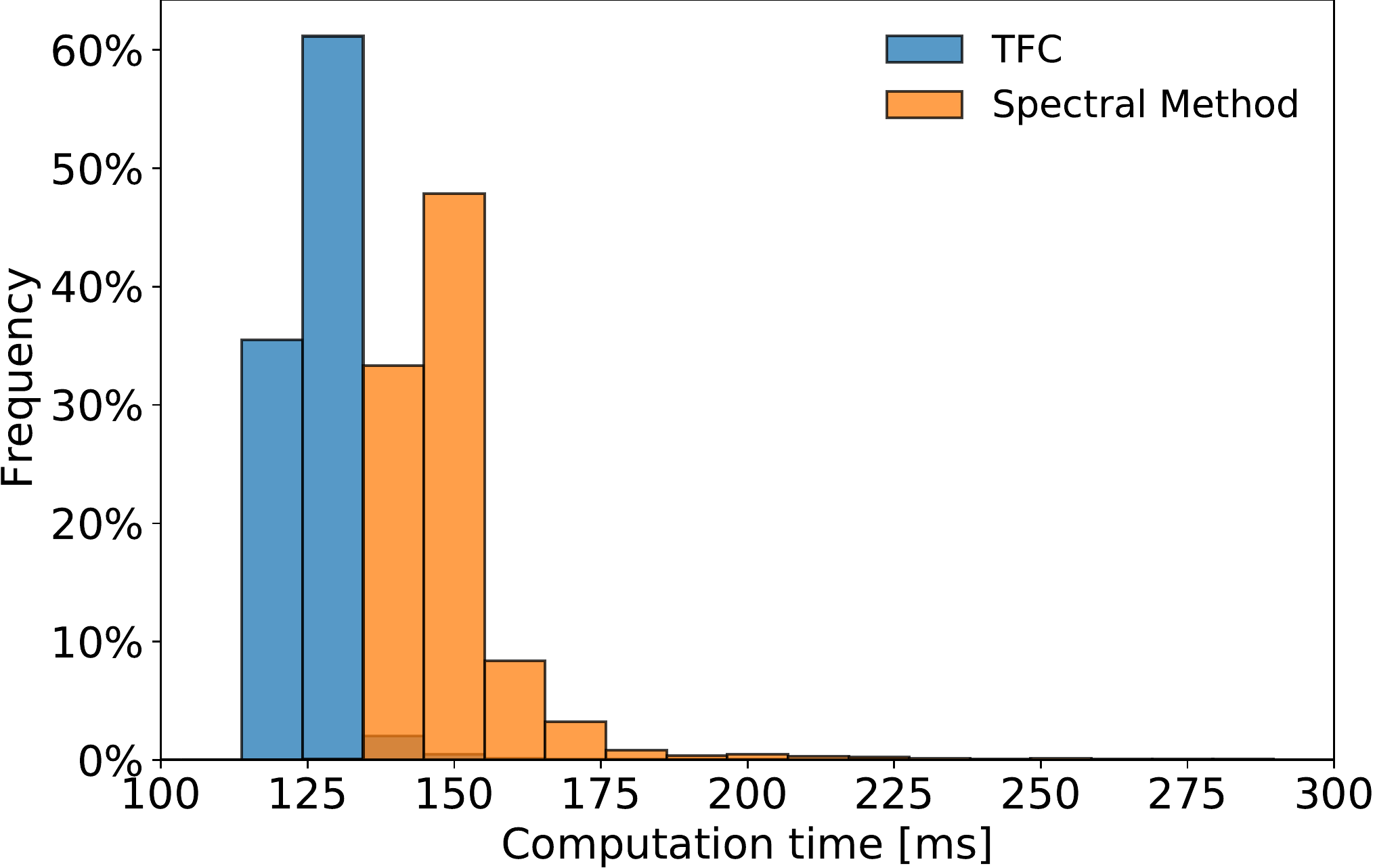}
    \caption{Histogram of the computation time of both methods.}
    \label{fig:s3_EOL_hist_time}
\end{figure}
Next, Figure \ref{fig:s3_EOL_hist_time} quantifies the computation time associated with both approaches. It can be seen that TFC produces a solution about 25-50 ms faster than the spectral method. This observation makes sense when analyzing Figure \ref{fig:s3_EOL_hist_iter}, where it can be seen that TFC usually takes between 18-20 iterations, while the spectral method can take up to 30 iterations.
\begin{figure}[H]
    \centering\includegraphics[width=0.7\linewidth]{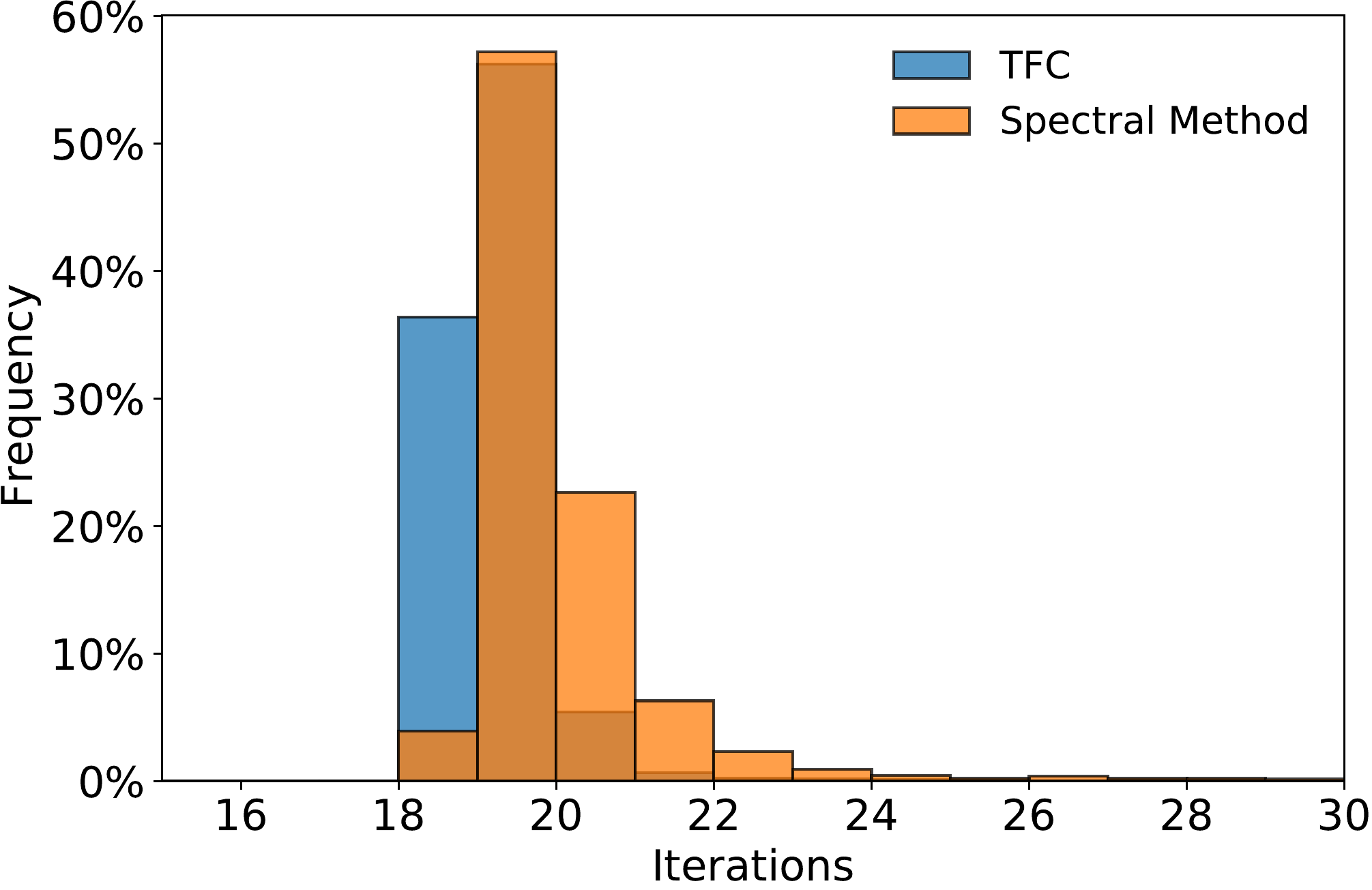}
    \caption{Histogram of the number of iterations.}
    \label{fig:s3_EOL_hist_iter}
\end{figure}
\end{example}

\section{Conclusions}
In this section, we applied TFC to the 3D energy-optimal landing problem, which has a known feedback solution for constant acceleration due to gravity $\B{a}_g$. While the TFC algorithm's implementation is relatively straightforward due to the simplicity of the optimal control problem, it gives us a major stepping stone forward in quantifying the accuracy, robustness, and speed of the TFC technique to solve realistic optimization problems. Moving forward, we will leverage what was learned from this example to make decisions in the fuel-optimal landing problem in the following section. The major takeaways from this problem are:

\begin{blankBox}{Major takeaways from energy-optimal landing tests}
\begin{enumerate}
    \item The Single-loop TFC approach requires the mapping parameter to show up nonlinearly in the dynamics. This can cause sensitivity due to initialization, which reduces the algorithm's robustness.
    \item The Outer-loop approach allows for any numerical scheme to be paired with TFC, which increases the applicability and, as seen in the prior example, can lead to increased robustness.
\end{enumerate}

\end{blankBox}
%

\chapter{FUEL-OPTIMAL LANDING*\label{chap:fol}}
\blfootnote{*Reprinted (along with revisions and updates unique to this dissertation) by permission from Springer Nature Customer Service Centre GmbH: Springer Nature The Journal of the Astronautical Sciences ``Fuel-Efficient Powered Descent Guidance on Large Planetary Bodies via Theory of Functional Connections,'' Johnston, H., Schiassi, E., Furfaro, R. et al., 2020, J Astronaut Sci 67, 1521–1552, Copyright 2020, \cite{FOL}}

The fuel-optimal (or propellant-efficient) landing is the natural extension from our solution of the energy-optimal landing problem presented in Chapter \ref{chap:eol}. This problem now introduces the mass state as another dynamic equation and inequality constraints on the spacecraft's thrust. While ultimately, we are interested in the full six-degree-of-freedom (6-DOF) solution, this 3-DOF is the natural next step where the attitude dynamics are not considered. This problem's solution is the subject of many studies, as mentioned in the literature review presented at the beginning of Chapter \ref{chap:opt_con}. Of the techniques discussed, Lu \cite{upg} has looked to solve this problem using the indirect method, which reduces the problem to a shooting method, and Acikmese and Ploen \cite{acikmese2007convex} and  Blackmore et al. \cite{blackmore2010minimum} have reformulated the problem via convex optimization to derive a solution.

\section{Dynamical model}
\begin{figure}[ht]
    \centering\includegraphics[width=\linewidth]{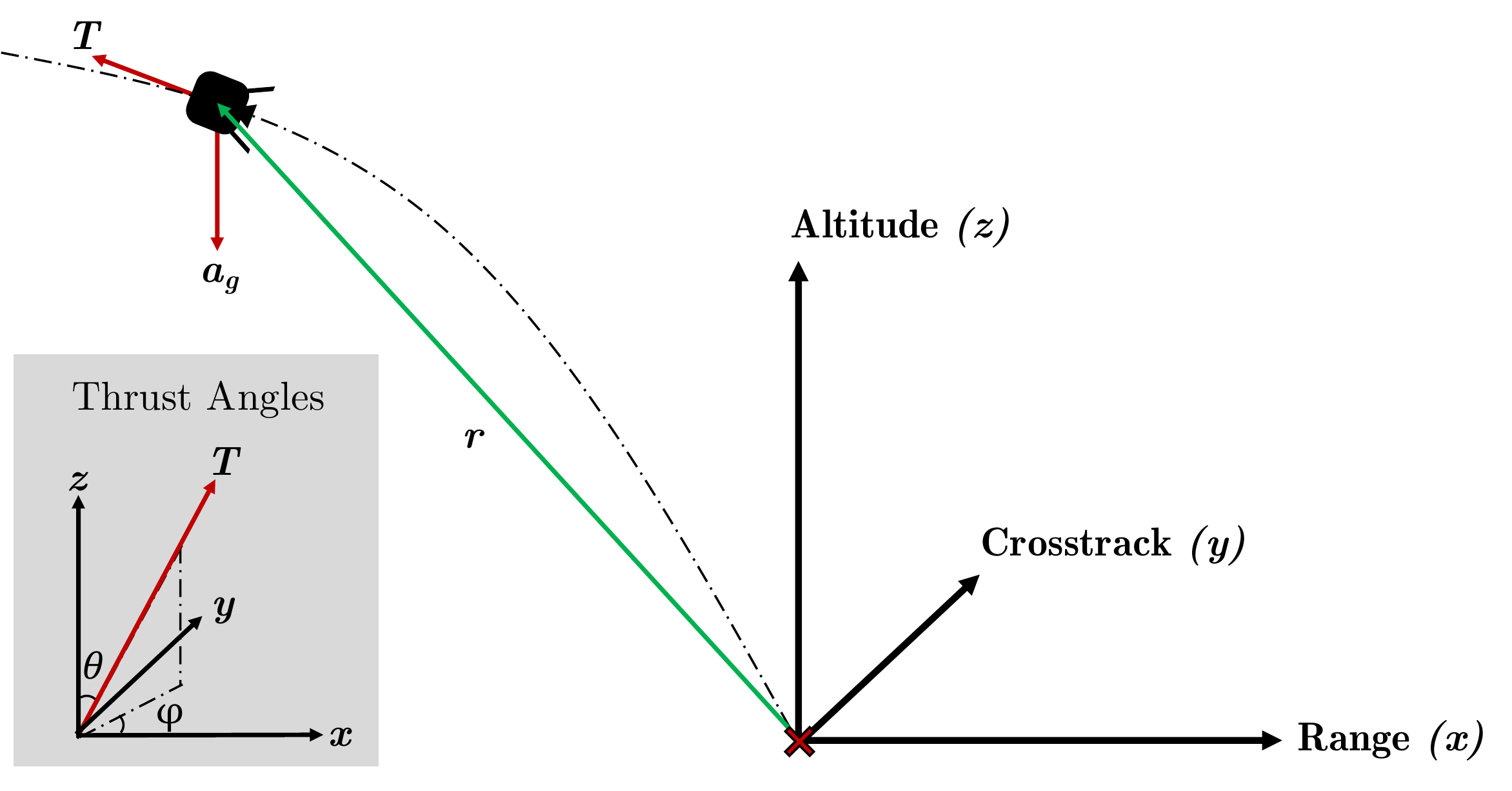}
    \caption{Coordinate frame definition for optimal powered descent pinpoint landing problem. Reprinted with permission from \cite{FOL}.}
    \label{fig:coord}
\end{figure}

For the problem of powered descent pinpoint landing guidance on large bodies (e.g., the Moon or Mars) the governing system dynamics during the powered descent phase can be modeled as follows,
\begin{align}
  \bdot{r} &= \,\B{v} \nonumber \\ 
  \bdot{v} &= \,\B{a}_g + { \dfrac{\B{T}}{m} } \nonumber\\ 
  \dot{m} &= \,-\alpha \, T \label{eq:s6_mdot}
\end{align}
where the spacecraft's state is defined by the position $\B{r}$, velocity $\B{v}$, and mass $m$. Additionally, $\alpha = 1 / v_{ex}$, where $v_{ex}$ is the effective exhaust velocity of the rocket engine that is considered constant \cite{acikmese2007convex,upg}, $T = ||\B{T}||$, and $\B{T} = T \, \hat{\B{t}}$ is the thrust and it is constrained as follows:
\begin{equation*}
\begin{aligned}
    0 \leq T_{min} &\leq T \leq T_{max} \\
    || \hat{\B{t}} ||&=1.
\end{aligned}
\end{equation*}
Furthermore, $\B{a}_g$ is the gravity acceleration, which is also considered constant. As stated in Reference \cite{upg}, this assumption is justified for short flights, as is the case for the landing's powered descent phase. A summary of the reference frame for this problem is given in Figure \ref{fig:coord}.
For the landing problem, the boundary consists of initial and final constraints on the position and velocity state and an initial constraint on the mass state,
\begin{equation*}
    \begin{cases} \B{r}(0) = \B{r}_0 \\ \B{v}(0) = \B{v}_0 \end{cases}, \qquad \begin{cases} \B{r}(t_f) = \B{r}_f \\ \B{v}(t_f) = \B{v}_f\end{cases}, \qquad \text{and} \qquad m(0) = m_0.
\end{equation*}
In all, the objective is to minimize the mass of the propellant used while satisfying the dynamics constraints of the problem. Therefore, the problem can be posed as,
\begin{wblankBox}{Optimization problem statement}
\begin{equation*}
    \underset{t_f,T}{\text{minimize}} \quad \alpha\int_{0}^{t_f}T \dd \tau
\end{equation*}
\begin{equation*}
\begin{aligned}
\text{subject to} \quad \dot{\B{r}} = \,\B{v}, \quad &\dot{\B{v}} = \,\B{a}_g + { \dfrac{\B{T}}{m} }, \quad \dot{m} = \,-\alpha \, T \\
& 0 \leq T_{min} \leq T \leq T_{max}, \quad || \hat{\B{t}} || = 1 \\
& \B{r}(0) = \B{r}_0, \quad \B{v}(0) = \B{v}_0, \quad m(0) = m_0 \\
& \B{r}(t_f) = \B{r}_f, \quad \B{v}(t_f) = \B{v}_f
\end{aligned}
\end{equation*}
\end{wblankBox}

\section{First-order necessary conditions}
From our definition of the optimization problem, we next apply the indirect method by applying the PMP dictates that the Hamiltonian takes the following form \cite{Bryson_Ho},
\begin{equation*}
    H = \mathcal{L} + \B{\lambda}\T \B{f} + \B{\mu}\T \B{C}
\end{equation*}
which can be expanded to,
\begin{equation}\label{eq:s6_ham}
    H=\alpha T + \B{\lambda}_r\T \B{v}+ \B{\lambda}_v\T  \left( \B{a}_g + \frac{T}{m} \hat{\B{t}} \right) - \lambda_m \alpha T + \mu_1 (T-T_{\max}) + \mu_2 (T_{\min}-T)
\end{equation}
where $T-T_{\max} \leq 0$ and $T_{\min}-T \leq 0$ and $\mu_1 > 0, \mu_2 > 0$. According to PMP, the optimal thrust solution is one that minimizes the Hamiltonian. Because both the thrust $T$ and mass $m$ are both non-negative, $\hat{\B{t}}$ should be in the opposite direction of of the velocity costate, i.e., $\hat{\B{t}} = -\frac{\B{\lambda}_v}{||\B{\lambda}_v||}$. This is what in Lawden's theory \cite{lawden1963optimal} is called \textit{primer's vector}. Thus Equation (\ref{eq:s6_ham}) can be rewritten as,
\begin{equation*}
H = \alpha T + \B{\lambda}_r\T \B{v}+ \B{\lambda}_v\T \B{a}_g - \frac{T}{m} ||\B{\lambda}_v|| - \lambda_m \alpha T + \mu_1 (T-T_{\max}) + \mu_2 (T_{\min}-T)
\end{equation*}
Now, to determine optimal thrust magnitude, we impose that the partial derivative of the Hamiltonian with respect to the thrust (i.e., the control) is equal to zero, which is of the form of Equation \eqref{eq:s4_control},
\begin{equation*}
    S := \frac{\partial H}{\partial T} = \underbrace{\alpha - \frac{1}{m}||\B{\lambda}_v|| - \alpha \lambda_m}_{\textstyle \text{$\sigma$} \mathstrut} + \mu_1 - \mu_2 = 0 
\end{equation*}
where there are three conditions that result in $S = 0$:
\begin{enumerate}
    \item if $\mu_1 = \mu_2 = 0 \quad (T_{\min} < T < T_{\max})$ \quad \text{then} \quad $\sigma = 0$
    \vspace{1mm}
    \item if $\mu_1 = 0$, $\mu_2 > 0 \quad (T=T_{\min})$  \ \ \quad \hspace{3pt} \quad \text{then} \quad $\sigma - \mu_2 = 0$ \quad $\rightarrow$ \quad $\sigma = \mu_2 > 0$
    \vspace{1mm}
    \item if $\mu_1 > 0$, $\mu_2 = 0 \quad (T = T_{\max})$ \: \: \quad \quad \text{then} \quad $\sigma + \mu_1 = 0$ \quad $\rightarrow$ \quad $\sigma = -\mu_1 < 0$
\end{enumerate}
Finally, one can conclude that the thrust magnitude has the following program:
\begin{equation*}
    T = \begin{cases} = T_{\max} \qquad \text{if} \qquad \sigma < 0 \\ = T_{\min} \qquad \text{if} \qquad \sigma > 0 \end{cases}
\end{equation*}
It has been demonstrated in Reference \cite{upg} that the singular case $\sigma = 0$ corresponds to a constant thrust perpendicular to the gravity vector, which is generally not possible for a powered descent problem. Therefore, a singular arc is not part of the sought optimal solution. Furthermore, it is straightforward to show that $\sigma$ changes signs at most twice and is derived in detailed in Reference \cite{upg}. Consequently, the thrust magnitude can switch between min-max twice at the most. That is, in the most general case, the thrust magnitude has a max-min-max profile. Hence, we can write the thrust magnitude as a function of time with $t_1$ and $t_2$ as parameters, where $t_1$ and $t_2$ are the times where the switches happen, i.e., $T=T(t;t_1,t_2)$. This result implies that thrust is constant between switches, and therefore, the solution of Equation \eqref{eq:s6_mdot} is a piecewise linear function in terms of  $t_1$ and $t_2$ detailed by the following equation,
\begin{equation*}
    m(t;t_1,t_2) = \begin{cases} 
    \text{if } t \leq t_1 \qquad \ : \, m_0 - \alpha \Big[T_{\max}(t-t_0)\Big]\\ 
    \text{if } t_1 \leq t \leq t_2 \ : \, m_0 - \alpha \Big[T_{\max}(t_1-t_0) - T_{\min}(t-t_1)\Big]\\ 
    \text{else} \qquad \qquad \ : \, m_0 - \alpha \Big[T_{\max}(t_1-t_0) - T_{\min}(t_2-t_1) - T_{\max}(t-t_2)\Big]. \end{cases}
\end{equation*}
In addition to these conditions, we are left with the first-order necessary conditions for the costates as given by Equation \eqref{eq:s4_costate},
\begin{eqnarray*}
    \bdot{\lambda}_r &= -\dfrac{\partial H }{\partial \B{r}} &= \B{0} \\
    \bdot{\lambda}_v &= -\dfrac{\partial H }{\partial \B{v}} &= -\B{\lambda_r} \\
    \dot{\lambda}_m &= -\dfrac{\partial H }{\partial m} &= -\frac{T}{m^2}||\B{\lambda}_v||.
\end{eqnarray*}
Finally, since the final mass state is unconstrained, Equation \eqref{eq:s4_free_xf} implies that, 
\begin{equation*}
    \lambda_m(t_f) = 0,
\end{equation*}
and likewise, since the final time of the problem is unknown, Equations \eqref{eq:s4_free_tf} leads to the condition on the final value of the Hamiltonian.
\begin{equation*}
    H(t_f) = 0.
\end{equation*}
In fact, since the Hamiltonian is not an explicit function of time, the partial derivative with respect to time is zero (i.e., $\frac{\partial H}{\partial t} = 0$), which implies a stronger condition, that the Hamiltonian should be zero for all time,
\begin{equation*}
    H(t) = 0.
\end{equation*}
We will take these conditions and look to apply the TFC method to solve all of the equations simultaneously. 


\section{Solving the problem via the TFC}\label{sect:s6_fol_jac}

With the simplifications introduced in the previous section, the following nonlinear set of equations must be solved to find the optimal state and thrust program,
\begin{wblankBox}{First-order necessary conditions}
\begin{align}
  \bdot{r} &= \B{v} \label{eq:s6_r_dot}\\
  \bdot{v} &= \B{a}_g - \beta(t)
  \dfrac{\B{\lambda_v}}{||\B{\lambda_v}||}\label{eq:s6_v_dot}\\
  \bdot{\lambda}_r &= \B{0} \label{eq:s6_lam_r_dot}\\ 
  \bdot{\lambda}_v &=  -\B{\lambda_r}\label{eq:s6_lam_v_dot} \\
  \dot{\lambda}_m &= \,- \frac{T (t;t_1,t_2)}{m^2}||\B{\lambda}_v|| \label{eq:s6_lam_m_dot}\\
 H(t_f) = 0 &= \alpha T(t_f;t_1,t_2) + \B{\lambda}_v\T(t_f) \left(\B{a}_g - \beta(t_f)\frac{\B{\lambda}_v(t_f)}{||\B{\lambda}_v(t_f)||} \right) \label{eq:s6_trans}
\end{align}
where we define
\begin{equation*}
    \beta(t) := \frac{T(t;t_1,t_2)}{m(t)}
\end{equation*}
and Equations \eqref{eq:s6_r_dot}, \eqref{eq:s6_v_dot}, and \eqref{eq:s6_lam_m_dot} are subject to
\begin{align*}
    \B{r(0)} = \B{r_0}, \quad \B{v(0)} = \B{v_0}, \quad \B{r(t_f)} = \B{r_f}, \quad \B{v(t_f)} = \B{v_f}, \quad \lambda_m(t_f) = 0.
\end{align*}
\end{wblankBox}

It must be noted that $\lambda_m$ only shows up in Equation \eqref{eq:s6_lam_m_dot}, and can therefore be solved independently. Since the transversality condition gives $\lambda_m (t_f) = 0$, Equation \eqref{eq:s6_lam_m_dot} can be solved by back propagation or by simply using the TFC method.

Since this problem's solution exhibits a bang-bang profile for thrust, the original formulation of the TFC method (i.e., as used in the Outer-loop method of the energy optimal landing problem in Section \ref{sec:s5_outLoop}) must be adjusted to accommodate switching behavior in the control. In general, this can be labeled as a hybrid system because the dynamical behavior is governed by both continuous dynamics (when the thruster is firing) and discrete dynamics (when the thrust jumps). The general theory for this extension to hybrid systems has been developed in Section \ref{sec:s3_hybridSystem} but is also fully developed in the following equations. Additionally, a few equations are redundant and can be removed completely via the TFC \ce\ to further simplify the solution of this nonlinear system of equations. As done in the last section, the differential equation expressed by Equation \eqref{eq:s6_r_dot} is unnecessary and can be disregarded. Similarly, the equations for $\bdot{\lambda}_r$ and $\bdot{\lambda}_v$ can be simplified. First, let us express the vector equations as three scalar equations, each where the index $i$ represents the individual components. Using this notation, we can expand $\B{\lambda}_v$ such that,
\begin{equation*} 
    \lambda_{v_i}  = a_{0_i} + a_{1_i} z = \B{h}_{\lambda}\T \B{\xi}_{\lambda_i}, \quad \text{for} \quad i = 1,2,3
\end{equation*}
which satisfies Equations (\ref{eq:s6_lam_r_dot}-\ref{eq:s6_lam_v_dot}) through
\begin{align*}
    \dot{\lambda}_{v_i} = c_{\lambda} \lambda_{v_i}' &= c_{\lambda} a_{1_i} \\
   -\dot{\lambda}_{v_i} = \lambda_{r_i} &= - c_{\lambda} a_{1_i}.
\end{align*}
This process reduces the problem to the solution of a single differential equation expressed by Equation \eqref{eq:s6_v_dot} and an algebraic equation for the Hamiltonian at the final time given by Equation \eqref{eq:s6_trans}. Rewriting the differential equation in indicial notation and collecting all terms on one side, a loss function based on the residuals of the differential equation can be defined,
\begin{equation}\label{eq:s6_loss}
\mathbb{L}_i = a_i - a_{g_i} + \beta(t) \, \lambda_{v_i} \left(\displaystyle\sum_{j=1}^3 \lambda^2_{v_j}\right)^{-1/2} \qquad \text{for} \qquad i = 1,2,3
\end{equation}
where $a_i := \dot{v}_i$ (or simply the acceleration of the spacecraft). Now, the only step left is to construct a \ce\ for the state variables. In the above derivation of the thrust structure, we have shown that the thrust switches at most twice, leading to a max-min-max profile. Therefore, the function $\beta(t)$ in Equation \eqref{eq:s6_loss} jumps twice along the solution trajectory. This switching causes three distinct differential equations that cannot be solved with a single polynomial expansion over the entire domain, as was done for the energy-optimal guidance in Chapter \ref{chap:eol}. Therefore, a new formulation for the TFC approach has been developed to handle these hybrid systems \cite{pieceTFC}. This process allows for the continuity between each segment of the domain.
\begin{figure}[ht]
    \centering\includegraphics[width=.75\linewidth]{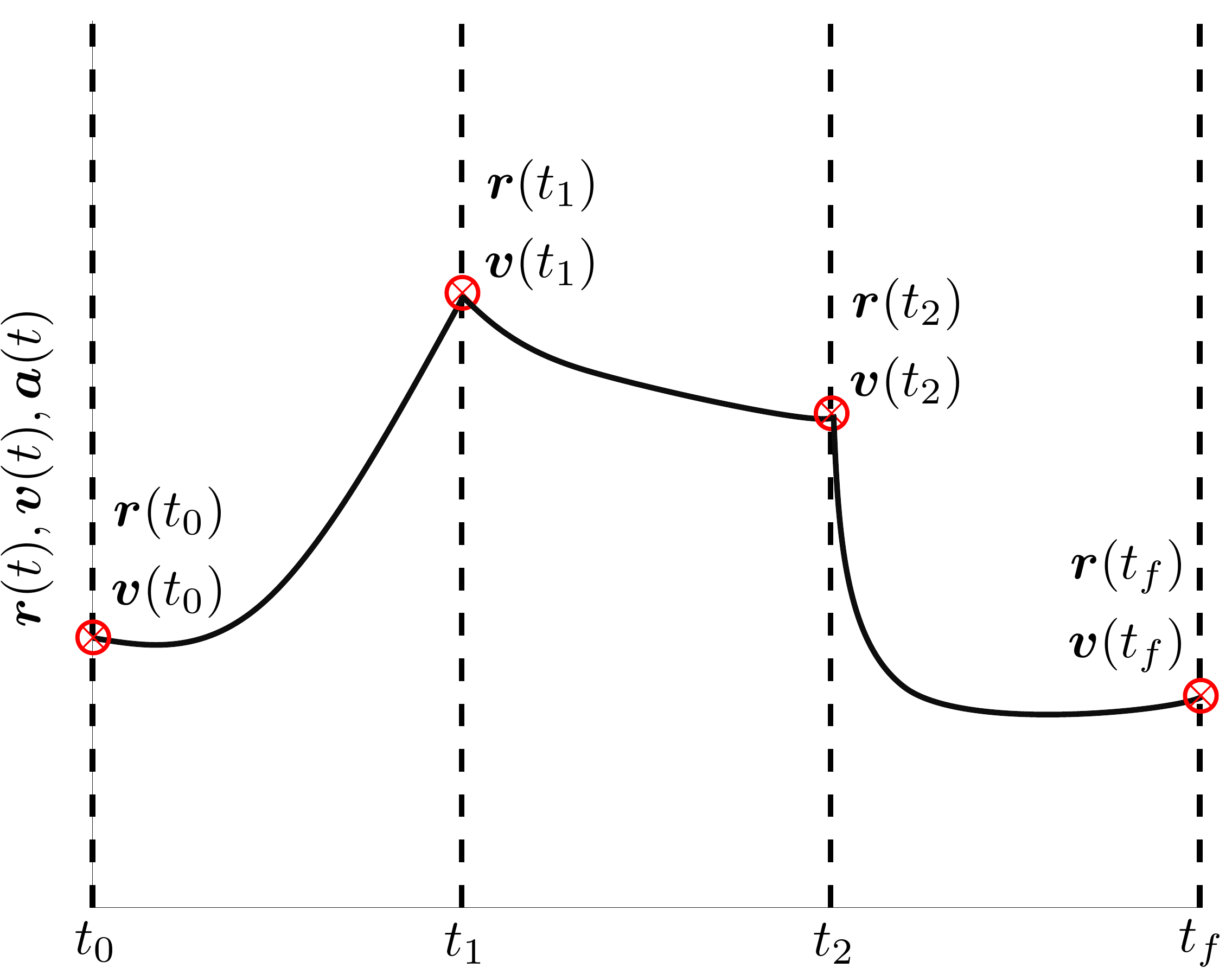}
    \caption{Visual representation of piece-wise approach using the TFC method. In this derivation, the \ces\ maintain continuity of position and velocity through embedded relative constraints. Reprinted with permission from \cite{FOL}.}
    \label{fig:piece}
\end{figure}
As shown in Figure \ref{fig:piece}, it is apparent that all sub-domains share the same constraint conditions (i.e, the initial and final position and velocity are constrained). Therefore, a single constraint expression can be derived for the case of arbitrary constraint locations and then incorporated into the sub-domains. The \ce\ for this specific case was derived in Section \ref{sec:s5_outLoop} and it is captured by Equation \eqref{eq:s5_r_timeCE}. Consequently, the position, velocity, and acceleration \ce\ can be expressed as,
\begin{align}
r_i(t,g_i(t)) = g_i(t) + \phi_1(t)(r_{0_i} - g_i(t_0)) &+ \phi_2(t) (r_{f_i} - g_i(t_f)) \nonumber\\&+ \phi_3(t) (v_{0_i} - \dot{g}_i(t_0)) + \phi_4(t)(v_{f_i} - \dot{g}_i(t_f))  \label{eq:s6_ce_r} \\ v_i(t,g_i(t)) = \dot{g}_i(t) + \dot{\phi}_1(t)(r_{0_i} - g_i(t_0)) &+ \dot{\phi}_2(t) (r_{f_i} - g_i(t_f)) \nonumber\\&+ \dot{\phi}_3(t) (v_{0_i} - \dot{g}_i(t_0)) + \dot{\phi}_4(t)(v_{f_i} -\dot{g}_i(t_f)) \nonumber \\ a_i(t,g_i(t)) = \ddot{g}_i(t) + \ddot{\phi}_1(t)(r_{0_i} - g_i(t_0)) &+ \ddot{\phi}_2(t) (r_{f_i} - g_i(t_f)) \nonumber\\&+ \ddot{\phi}_3(t) (v_{0_i} - \dot{g}_i(t_0)) + \ddot{\phi}_4(t)(v_{f_i} -\dot{g}_i(t_f)) \label{eq:s6_ce_a} 
\end{align}
The switching functions are the same as those used in the Outer-loop method for solving the energy-optimal landing problem (this is because they share the same constraint conditions) and are defined by switching functions of Section \ref{sec:s5_outLoop}. In these switching functions, $t_0$ and $t_f$ must be replaced with the respective segment's initial and final time, e.g., for the first segment $ t \in [t_0, t_1]$.

The \ce\ detailed by Equations (\ref{eq:s6_ce_r}-\ref{eq:s6_ce_a}) can be used as a template to write the \ces\ for each segment of the solution trajectory. In order to explicitly identify the segment, the pre-superscript notation will be used. For example, $\p{1}{r}_i$ describes the position \ce\ for the first segment defined on $t\in[t_0,t_1]$. For this problem, $s=1$ (where $s$ is used to denote the segment) is defined on $t\in[t_0,t_1]$, $s=2$ is defined on $t\in[t_1,t_2]$, and $s=3$ is defined on $t\in[t_2,t_f]$. Using this formulation, the \ces\ of position for each segment are,
\small
\begin{align*}
\p{1}{r}_i(t,g_i(t)) = \p{1}{g}_i(t) &+  \p{1}{\phi}_1(t)\left(r_{0_i} - \p{1}{g}_i(t_0)\right) + \p{1}{\phi}_2(t) \left(r_{1_i} - \p{1}{g}_i(t_f)\right) \\& + \p{1}{\phi}_3(t) \left(v_{0_i} - \p{1}{\dot{g}}_i(t_0)\right) + \p{1}{\phi}_4(t)\left(v_{1_i} -\p{1}{\dot{g}}_i(t_f)\right)
\end{align*}
\begin{align*}
\p{2}{r}_i(t,g_i(t)) = \p{2}{g}_i(t) &+  \p{2}{\phi}_1(t)\left(r_{1_i} - \p{2}{g}_i(t_0)\right) + \p{2}{\phi}_2(t) \left(r_{2_i} - \p{2}{g}_i(t_f)\right) \\&+ \p{2}{\phi}_3(t) \left(v_{1_i} - \p{2}{\dot{g}}_i(t_0)\right) + \p{2}{\phi}_4(t)\left(v_{2_i} -\p{2}{\dot{g}}_i(t_f)\right)
\end{align*}
\begin{align*}
\p{3}{r}_i(t,g_i(t)) = \p{3}{g}_i(t) &+ \p{3}{\phi}_1(t)\left(r_{2_i} - \p{3}{g}_i(t_0)\right) + \p{3}{\phi}_2(t) \left(r_{f_i} - \p{3}{g}_i(t_f)\right)  \\&+ \p{3}{\phi}_3(t) \left(v_{2_i} - \p{3}{\dot{g}}_i(t_0)\right) + \p{3}{\phi}_4(t)\left(v_{f_i} -\p{3}{\dot{g}}_i(t_f)\right)
\end{align*}
\normalsize
where the derivative of these functions follow the form of Equations (\ref{eq:s6_ce_r}-\ref{eq:s6_ce_a}). This allows us to collect the unknown $\B{\xi}_i$ vectors and write the \ce\ in the form,
\begin{align*}
\p{1}{r}_i(t,\p{1}{\B{\xi}}_i) = \p{1}{} \Big(\B{h}(z) - \phi_1(t)\B{h}(z_0) &- \phi_2(t)\B{h}(z_f) - \phi_3(t)c \B{h}_z(z_0)  - \phi_4(t)c \B{h}_z(z_f)  \Big)\T \p{1}{\B{\xi}}_i \\ &+ \p{1}{\phi}_1(t) r_{0_i} + \p{1}{\phi}_2(t) r_{1_i} + \p{1}{\phi}_3(t) v_{0_i} + \p{1}{\phi}_4(t) v_{1_i}
\end{align*}
\begin{align*}
\p{2}{r}_i(t,\p{2}{\B{\xi}}_i) = \p{2}{} \Big(\B{h}(z) - \phi_1(t)\B{h}(z_0) &- \phi_2(t)\B{h}(z_f) - \phi_3(t)c \B{h}_z(z_0)  - \phi_4(t)c \B{h}_z(z_f)  \Big)\T \p{2}{\B{\xi}}_i \\ &+ \p{2}{\phi}_1(t) r_{1_i} + \p{2}{\phi}_2(t) r_{2_i} + \p{2}{\phi}_3(t) v_{1_i} + \p{2}{\phi}_4(t) v_{2_i}
\end{align*}
\begin{align*}
\p{3}{r}_i(t,\p{3}{\B{\xi}}_i) = \p{3}{} \Big(\B{h}(z) - \phi_1(t)\B{h}(z_0) &- \phi_2(t)\B{h}(z_f) - \phi_3(t)c \B{h}_z(z_0)  - \phi_4(t)c \B{h}_z(z_f)  \Big)\T \p{3}{\B{\xi}}_i \\ &+ \p{3}{\phi}_1(t) r_{2_i} + \p{3}{\phi}_2(t) r_{f_i} + \p{3}{\phi}_3(t) v_{2_i} + \p{3}{\phi}_4(t) v_{f_i}
\end{align*}
Along with the linear unknowns in $\p{s}{\B{\xi}}_i$, the equations share linear unknowns in $r_{1_i},v_{1_i},r_{2_i},v_{2_i}$ which serve as the embedded relative constraints between adjacent segments. With this new formulation, we now have three separate loss functions based on the residual of the differential equation over each segment ($s$) which are as follows,
\begin{equation*} 
\p{s}{F}_i\left(t,\Xi \right) = \p{s}{a}_i - a_{g_i} + \beta(t) \, \lambda_{v_i} \, \left(\displaystyle\sum_{j=1}^3 \lambda^2_{v_j}\right)^{-1/2}.
\end{equation*}
Note that although the costate \ces\ do not need to be split into separate domains, special attention must be paid to discretizing the equations according to the segment time ranges. Again, to solve for the unknown $\B{\xi}_i$ parameters, a nonlinear least-squares technique was used, which requires computing the partials of the loss function with respect to all of the unknowns. All partial derivatives for each segment and each unknown are provided in Appendix \ref{sect:app_FOL}.

In addition to the loss functions for the problem dynamics given by Equation \eqref{eq:s6_v_dot}, a loss function associated with the transversality conditions for the Hamiltonian is defined as,
\begin{equation*}
    \mathbb{L}_H\left(t_f,\Xi\right) = \alpha T_\text{max} + \sum_{i=1}^{3}\lambda_{v_i}(t_f) a_{g_i} - \beta(t_f) \left(\sum_{i=1}^{3}\lambda^2_{v_i}(t_f)\right)^{\frac{1}{2}}.
\end{equation*}
The partial derivatives of this function are also provided in Appendix \ref{sect:app_FOL}. Next, by discretizing the domain over $N$ points, these loss functions can be organized into the loss vector, 
\begin{equation*}
    \mathbb{L} = \begin{Bmatrix} \p{1}{\mathbb{L}}_1\T & \p{1}{\mathbb{L}}_2\T & \p{1}{\mathbb{L}}_3\T & \p{2}{\mathbb{L}}_1\T & \p{2}{\mathbb{L}}_2\T & \p{2}{\mathbb{L}}_3\T & \p{3}{\mathbb{L}}_1\T & \p{3}{\mathbb{L}}_2\T & \p{3}{\mathbb{L}}_3\T & \mathbb{L}_H\end{Bmatrix}_{(\{9N+1\}\times 1)}\T
\end{equation*}
where
\begin{equation*}
 \p{s}{\mathbb{L}_i} = \begin{Bmatrix}\p{s}{F}_i(t_0,\Xi) & \hdots & \p{s}{F}_i(t_k,\Xi) & \hdots & \p{s}{F}_i(t_f,\Xi) \end{Bmatrix}\T.
\end{equation*}
Additionally, the vector of unknowns takes the form,
\begin{align*}
    \Xi = \Big\{ &\p{1}{\B{\xi}}_1\T \quad \p{1}{\B{\xi}}_2\T \quad \p{1}{\B{\xi}}_3\T \quad \p{2}{\B{\xi}}_1\T \quad \p{2}{\B{\xi}}_2\T \quad \p{2}{\B{\xi}}_3\T \quad \p{3}{\B{\xi}}_1\T \quad \p{3}{\B{\xi}}_2\T \quad \p{3}{\B{\xi}}_3\T \\ &\B{\xi}_{\lambda_1}\T \quad \B{\xi}_{\lambda_2}\T \quad \B{\xi}_{\lambda_3}\T  \quad \B{r}_1\T \quad \B{v}_1\T \quad  \B{r}_2\T \quad \B{v}_2\T  \Big\}_{(9m + 18)}\T.
\end{align*}
In general, the structure of the Jacobian is,
\begin{equation}\label{eq:s6_aug_J}
    \mathbb{J} = \begin{bmatrix} 
    \p{1}{J}_{\B{\xi}} & \B{0}_{(3N\times 3m)} & \B{0}_{(3N\times 3m)} &  \p{1}{J}_{\B{\xi}_\lambda} & \p{1}{J}_{r_1, v_1} & \B{0}_{(3N\times 6)} \\ 
    \B{0}_{(3N\times 3m)}  &   \p{2}{J}_{\B{\xi}}  & \B{0}_{(3N\times 3m)} &  \p{2}{J}_{\B{\xi}_\lambda} & \p{2}{J}_{r_1, v_1} & \p{2}{J}_{r_2, v_2}\\ 
    \B{0}_{(3N\times 3m)} &  \B{0}_{(3N\times 3m)}  &  \p{3}{J}_{\B{\xi}}  & \p{3}{J}_{\B{\xi}_\lambda} & \B{0}_{(3N\times 6)} & \p{3}{J}_{r_2, v_2} \\ \B{0}_{(1\times 3m)}  & \B{0}_{(1\times 3m)}  & \B{0}_{(1\times 3m)} & J_H & \B{0}_{(1\times 6)} & \B{0}_{(1\times 6)} \end{bmatrix}_{(\{9N+1\}\times \{9m + 18\})}.
\end{equation}
Finally, using Equation \eqref{eq:s6_aug_J} along with the augmented loss functions and unknown vector, an iterative least-squares is used to find $\Xi$.

\subsection{Jacobian properties and sparsity}
\begin{figure}[H]
    \centering\includegraphics[width=.9\linewidth]{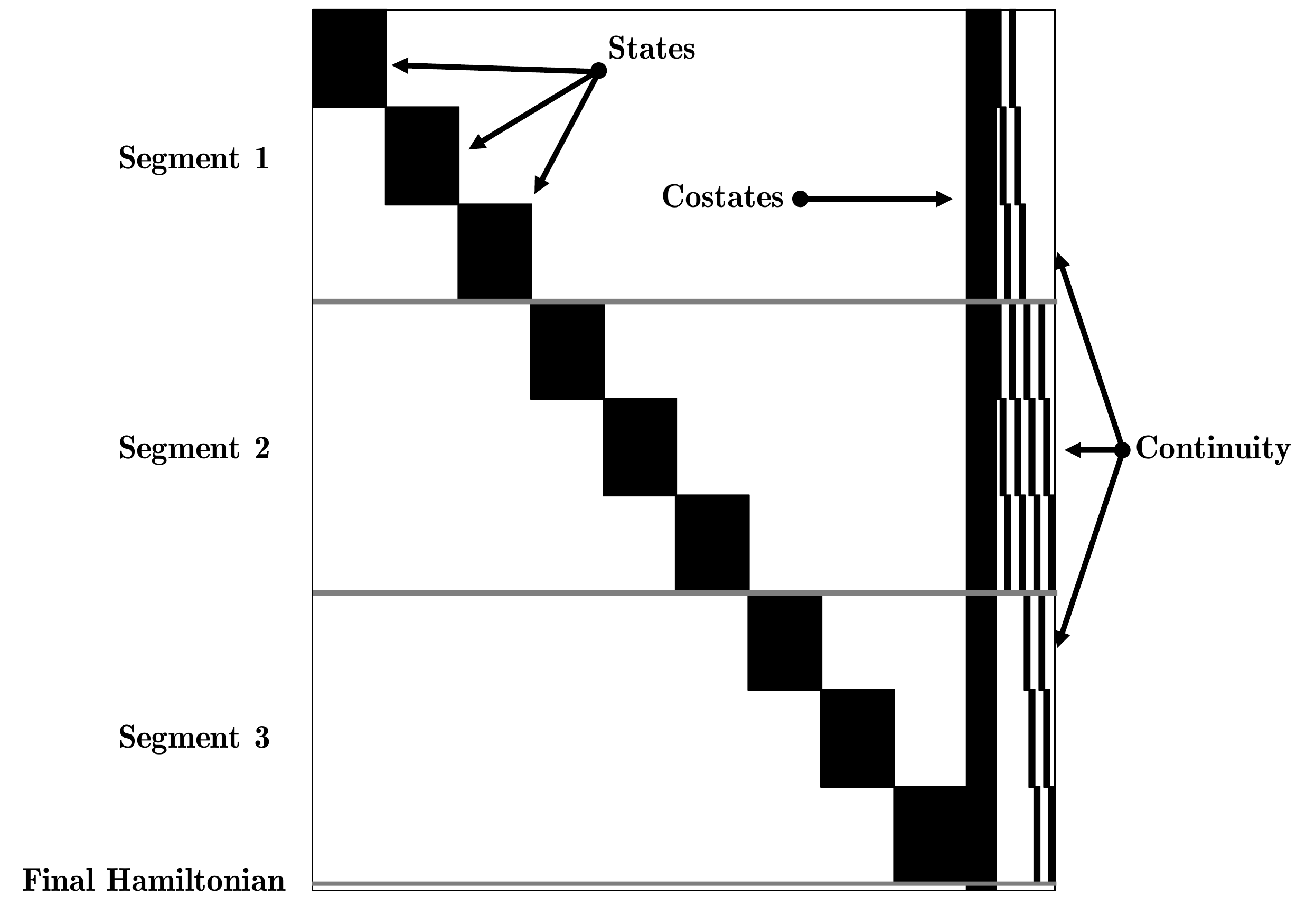}
    \caption{Visual representation of the Jacobian matrix to be inverted where the black elements represent the nonzero entries. Reprinted with permission from \cite{FOL}.}
    \label{fig:theJacob}
\end{figure}

From the prior equations, it should be evident that the Jacobian defined by Equation \eqref{eq:s6_aug_J} will need to be inverted. Therefore, Figure \ref{fig:theJacob} is provided as a visual aid to highlight the sparsity structure of this Jacobian. In addition to this structure, another property of this matrix is that the elements dealing with continuity, Jacobian terms $\p{1}{J}_{r_1,v_1}$, $\p{2}{J}_{r_1,v_1}$, $\p{2}{J}_{r_2,v_2}$, and $\p{3}{J}_{r_2,v_2}$, highlighted in the right side of Figure \ref{fig:theJacob}, are parameter independent (i.e., they are only a function of the $\phi(t)$ terms, or rather time) and therefore are constant and need only to be computed once per TFC loop.

\subsection{Initialization of parameters}
An initial estimate of the parameters is needed to initialize the iterative least-squares process. Since the problem is a boundary-value problem, the first guess for $\p{s}{\B{\xi}}_i$, $\B{r}_1$, $\B{r}_2$, $\B{v}_1$, and $\B{v}_2$ can be determined by simply connecting the initial and final position with a straight line and using this trajectory for a least-squares fitting of the \ces\ describing the $\p{s}{r}_i$ terms. Next, since $\lambda_{v_i}$ is related to the thrust direction, it can be assumed,
\begin{equation*}
    \B{\lambda}_{v_0} = \frac{\B{v}_0}{||\B{v}_0||},
\end{equation*}
similar to that presented in Reference \cite{upg} (Equation (51) in the text) However, the initialization of $\B{\lambda}_r = \B{0}$ will cases issues in the TFC method because this involves setting $\B{\xi}_\lambda$ coefficients to zeros. Therefore, in this dissertation, the coefficients are initialized using,
\begin{equation*}
    \B{\lambda}_{v_f} = -\frac{\B{r}_0}{||\B{r}_0||}.
\end{equation*}

\section{Summary of Algorithm}
\begin{figure}[H]
    \centering\includegraphics[width=\linewidth]{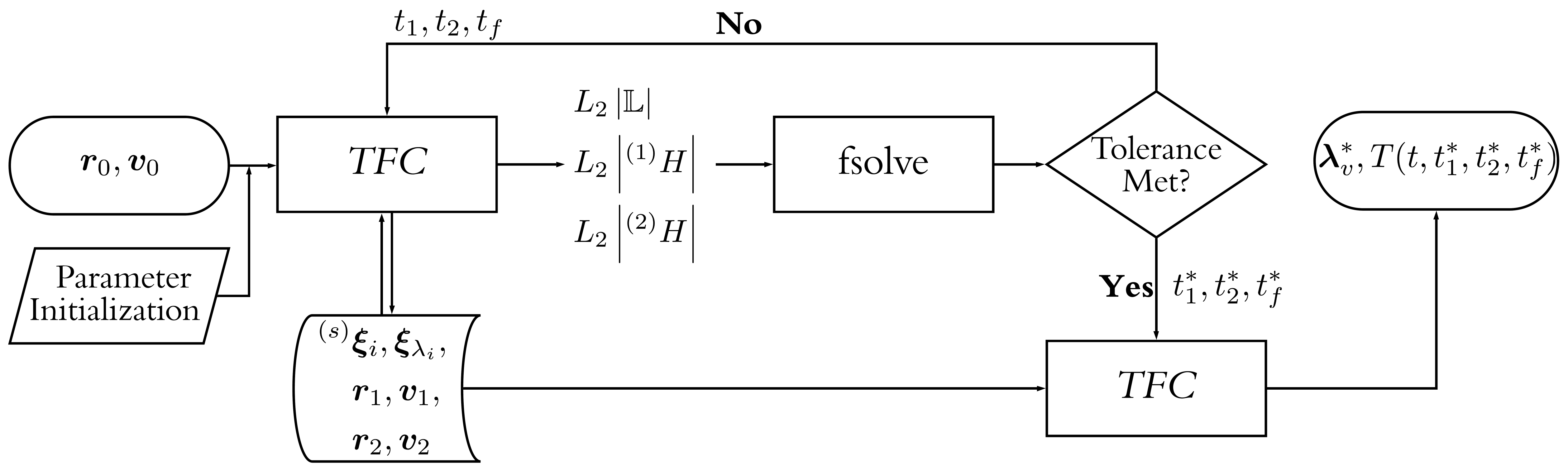}
    \caption{Summary of the full algorithm used with the TFC approach. Reprinted with permission from \cite{FOL}.}
    \label{fig:algorithm}
\end{figure}

Overall, the TFC method was used as ab inner-loop function to minimize the residuals of the first-order necessary conditions subject to a prescribed thrust profile $T(t;t_1,t_2)$, i.e., the switching times,  $t_1$ and $t_2$, and the final time, $t_f$, are assumed to be known by the TFC-based inner-loop routine. Consequently, an outer-loop routine has been developed to optimize the three time parameters $t_1, t_2, t_f$ given the $L_2$-norms of the residual of the first-order conditions, and the Hamiltonian over the first two segments (here, MATLAB's \cite{MATLAB} \verb"fsolve" was used). In other words, the following minimization problem needs to be solved for $t_1, t_2, \text{and} t_f$,
\begin{equation}\label{eq:s6_fsolve}
    \underset{t_1,t_2,t_f}\min \B{F}(t_1,t_2,t_f) = \begin{bmatrix} \max|\mathbb{L}|, &  \quad \max|\p{1}{H}(t)|, & \quad \max|\p{2}{H}(t)| \end{bmatrix}\T,
\end{equation}
where $\mathbb{L}$ is the loss function of the inner TFC loop, and $\p{1}{H}(t)$ and $\p{2}{H}(t)$ are the Hamiltonian values over the first and second segment, respectively, evaluated using the inner loop converged parameters. A flow chart of the relevant inputs and outputs is provided in Figure \ref{fig:algorithm}. 
Following the process given in Figure \ref{fig:algorithm}, the initial conditions $\B{r}_0$ and $\B{v}_0$ along with initial guesses for $t_1$, $t_2$, and $t_f$ are fed into the TFC method to minimize $\mathbb{L}$. The converged parameters are used to evaluate the Hamiltonian over the first and second segments. Using the norm of these quantities, \verb"fsolve" is used to solve the minimization problem given in Equation (\ref{eq:s6_fsolve}). If the tolerance of the outer loop is met (in all tests, the step and function tolerance of \verb"fsolve" were set to $4.4 \times 10^{-16}$), the $t_1, t_2,$ and $t_f$ are considered optimal, and the TFC loop is ran one more time to compute the optimal trajectory.\footnote{It must be noted that these low tolerances were used to quantify the baseline for speed and accuracy of the method. For implementation, the accuracy needed can be used to tune the tolerance and increase the algorithm's computational speed.}

\section{Results}
The proposed method was validated using two specific test cases based on selected initial conditions defining a powered descent guidance scenario for landing on Mars. In Example \ref{sec:s6_test1}, the algorithm is tested on initial conditions where the optimal trajectory is characterized by a min-max thrust profile. Furthermore, in Example \ref{sec:s6_test2}, the case where the optimal thrust profile is max-min-max is studied. In both cases, the results were compared with GPOPS-II solutions. The algorithm was fully implemented in MATLAB R$2019a$, and therefore not optimized for speed,

Similar to the energy-optimal landing problem in Chapter \ref{chap:eol}, the problem was scaled by the initial conditions for the numerical implementation. The unit length, $\ell^*$, and unit time, $t^*$, where calculated by the following equations, 
\begin{align*}
    \ell^* &= \max\left(|\B{r}_0|\right) \\
    t^* &= \frac{\ell^*}{\max\left(|\B{v}_0|\right)}.
\end{align*}

\subsection{Constant Test Parameters} 
We consider the trajectory optimization problem for a spacecraft performing powered descent for a pinpoint landing on Mars. The gravitational field is assumed constant, as generally, the powered descent starts below $1.5$ km. For the numerical test,  the lander parameters have been assumed to be similar to the ones presented in Reference \cite{acikmese2007convex} and reported in Table \ref{tab:test_parameters}.
\begin{table}[ht]
\caption{Constant parameters used in test cases. Reprinted with permission from \cite{FOL}.}
\centering
\begin{tabular}{|c||c|}
\hline
Variable & Value \\
\hline\hline
$\B{a}_g$ [m/s$^2$] & $\begin{Bmatrix} 0, & 0, & -3.7114\end{Bmatrix}\T$ \\\hline
$I_{sp}$ [s] & $225$ \\\hline
$g_0$ [m/s$^2$] & $9.807$\\\hline 
$\overline{T}$ [N] & $3,100$\\\hline 
$N_T$ [-] & $6$\\\hline 
$\phi_T$ [deg] & $27$ \\\hline 
\hline
\end{tabular}
\label{tab:test_parameters}
\end{table}
Thrust magnitude bounds and the $\alpha$ parameter are defined as follows:
\begin{eqnarray*}
    T_\text{min} &= 0.3 \overline{T}N_T\cos\phi_T &\approx 4,971.81 \: \: \: \text{[N]}\\
    T_\text{max} &= 0.8 \overline{T}N_T\cos\phi_T &\approx 13,258.18 \: \text{[N]}
\end{eqnarray*}
where $ \overline{T}$ is the maximum thrust for a single engine, $N_T$ is the number of thrusters in the lander, and $\phi_T$ is the cant angle of the thrusters with respect to the lander, and
\begin{equation*} 
    \alpha = \frac{1}{I_{\rm sp} \, g_0 \, \cos\phi_T} \approx 5.0863 \cdot 10^{-4} \: \text{[s/m]},
\end{equation*}
where $I_{\rm sp}$ is the engines' specific impulse and $g_0$ is Earth's gravitational constant.

\begin{example}{Test 1: Min-Max Trajectory}\label{sec:s6_test1}

For Test 1, initial conditions were selected such that the optimal thrust profile would be min-max, i.e., switch between minimum thrust to maximum thrust. Table \ref{tab:min_max} defines the boundary conditions for this test case, and Figure \ref{fig:min_max_traj} provides the converged trajectory using the TFC approach.
\begin{table}[H]
\caption{Boundary conditions for min-max trajectory profile test case. Reprinted with permission from \cite{FOL}.}
\centering
\begin{tabular}{|c||c|c|}
\hline
Variable & Initial & Final \\
\hline\hline
$\B{r}$ [m] & $\begin{Bmatrix} -900, & 10, & 1500\end{Bmatrix}\T$ & $\begin{Bmatrix} 0, & 0, & 0\end{Bmatrix}\T$ \\\hline
$\B{v}$ [m/s] & $\begin{Bmatrix} 30, & -10, & -70\end{Bmatrix}\T$ & $\begin{Bmatrix} 0, & 0, & 0\end{Bmatrix}\T$ \\\hline
$m$ [kg] & $1905$ & - \\\hline
\hline
\end{tabular}
\label{tab:min_max}
\end{table}

\begin{figure}[H]
    \centering\includegraphics[width=.75\linewidth]{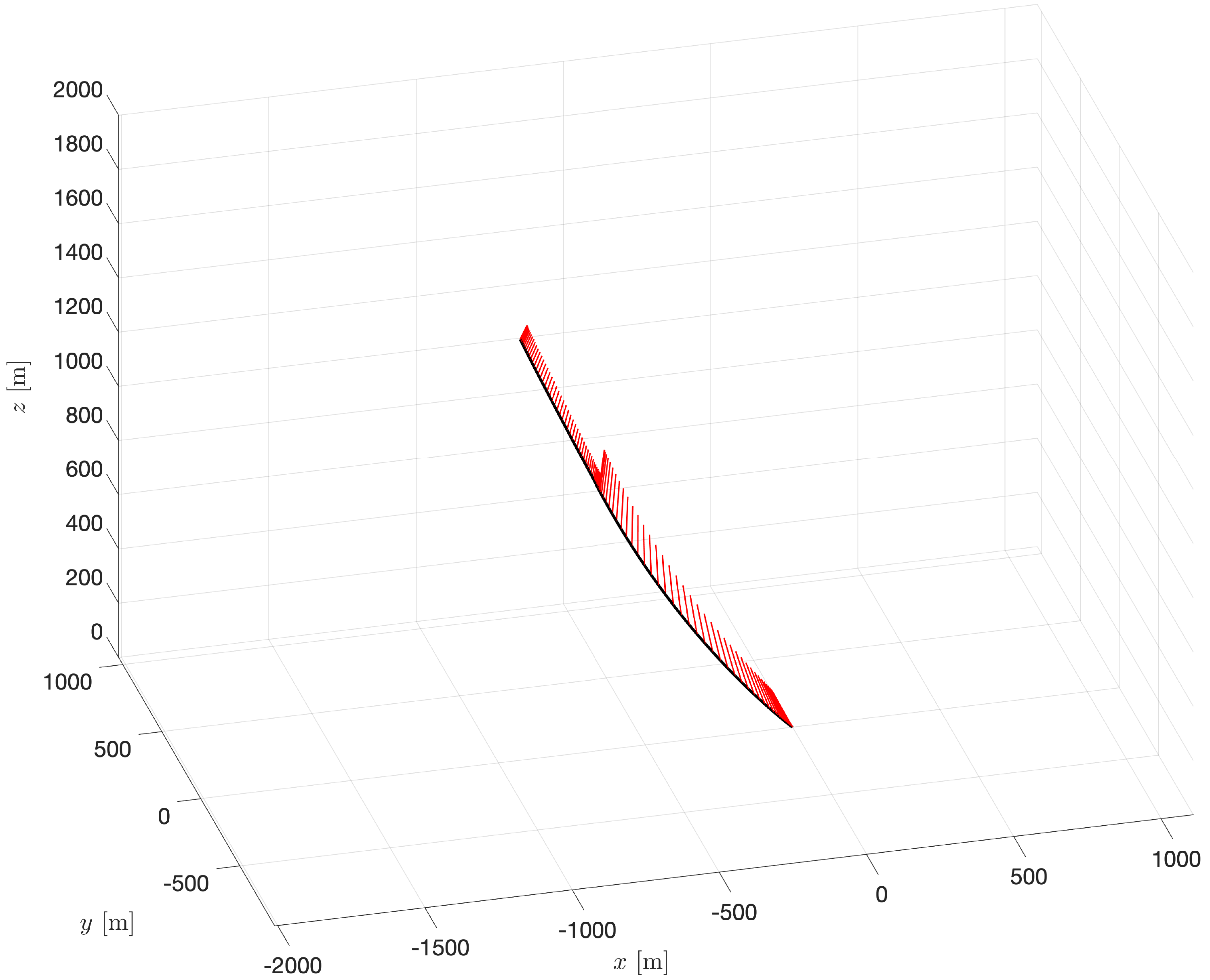}
    \caption{Landing trajectory for min-max thrust profile based on initial conditions, $\B{r}_0 = \protect\begin{Bmatrix} -900, & 10, & 1500\protect\end{Bmatrix}\T $[m], $\B{v}_0 = \protect\begin{Bmatrix} 30, & -10, & -70\protect\end{Bmatrix}\T$ [m/s], $m_0 = 1905$ [kg]. Reprinted with permission from \cite{FOL}.}
    \label{fig:min_max_traj}
\end{figure}
In addition to the trajectory, component plots of the position, velocity, and acceleration are provided in Figure \ref{fig:min_max_tfc}. Furthermore, this figure also plots the residual of the governing differential equations for mass and acceleration to quantify the method's accuracy. It can be seen that the TFC residual is about $\mathcal{O}(10^{-11})$ or less for the whole solution domain.
\begin{figure}[H]
    \centering\includegraphics[width=.95\linewidth]{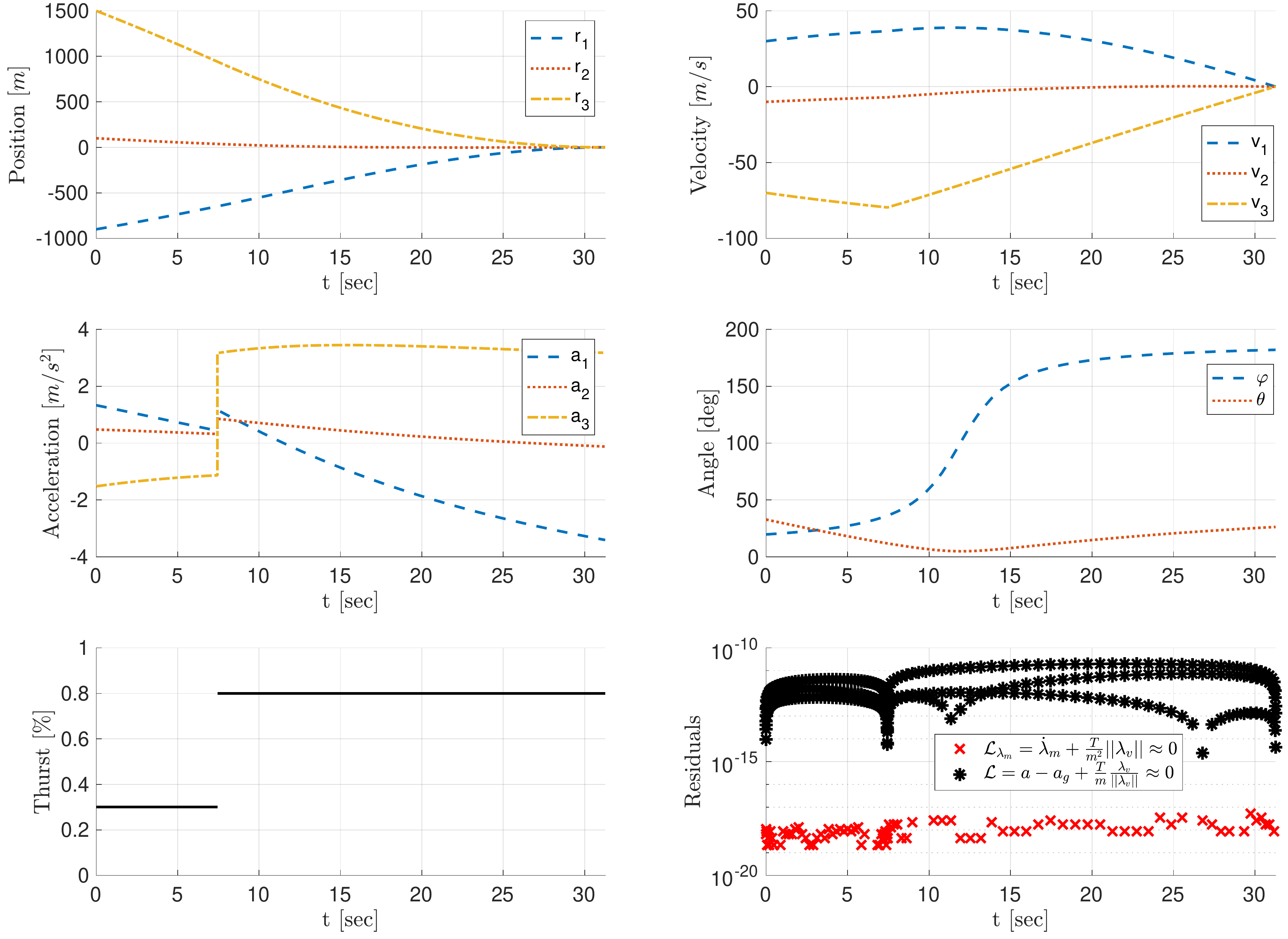}
    \caption{TFC solution of the min-max thrust profile case. The solution is presented in terms of the position, velocity, acceleration, and residuals of the differential equations. Reprinted with permission from \cite{FOL}.}
    \label{fig:min_max_tfc}
\end{figure}
The accuracy of this approach was also compared to results obtained using GPOPS-II \cite{GPOPS} and is quantified in terms of the converged parameters, the $L_2$-norms of the Hamiltonian, and propellant mass used. Moreover, to further justify the accuracy of the solution, the converged parameters of initial costate values and switching times for each method were propagated using MATLAB's \verb"ode45" with a tolerance of $2.2 \times 10^{-14}$ to check the final position and velocity error and also the final error of the $\lambda_m$ term. The tabulated values of this test are provided in Table \ref{tab:min_max_results}. In this test, \verb"fsolve" iterated 27 times with each TFC inner-loop averaging 76 ms, resulting in a total execution time of 2.1 seconds within the MATLAB implementation. Further, during this test, the TFC method converged in about 6 iterations every function call. Additionally, as a last point of comparison, the time histories of the Hamiltonian for both methods are plotted in Figure \ref{fig:min_max_H}. 
\begin{table}[H]
\caption{Converged parameters for the TFC and GPOPS-II solution for the min-max trajectory test case. The values $||\B{r}(t_f)||$, $||\B{v}(t_f)||$, and $\lambda_m(t_f)$ were determined by propagating both TFC and GPOPS-II converged solutions in order to have a one-to-one comparison on the accuracy of the converged solutions. Reprinted with permission from \cite{FOL}.}
\centering
\begin{tabular}{|c||c|c|c|}
\hline
Variable & TFC & GPOPS-II \cite{GPOPS} \\
\hline\hline
$L_2[\mathbb{L}]$ & $1.036 \cdot 10^{-10}$ & $-$  \\\hline
$L_2[H]$ & $5.488 \cdot 10^{-11}$ & $1.064 \cdot 10^{-3}$  \\\hline
$m_{\text{used}}$ [kg] & $179.447$  & $179.447$ \\\hline
$t_1$ [s] & $7.4430$ & $7.4430$ \\\hline 
$t_f$ [s] & $31.2623$ & $31.2623$  \\\hline 
$||\B{r}(t_f)||$ [m] & $2.886 \cdot 10^{-9}$ & $1.535\cdot 10^{-2}$  \\\hline 
$||\B{v}(t_f)||$ [m] & $3.166 \cdot 10^{-10}$ & $7.649\cdot 10^{-4}$  \\\hline 
$\lambda_m(t_f)$ [s] & $4.496 \cdot 10^{-14}$ & $-4.193\cdot 10^{-7}$  \\
\hline
\hline
\end{tabular}
\label{tab:min_max_results}
\end{table}

\begin{figure}[H]
    \centering\includegraphics[width=.75\linewidth]{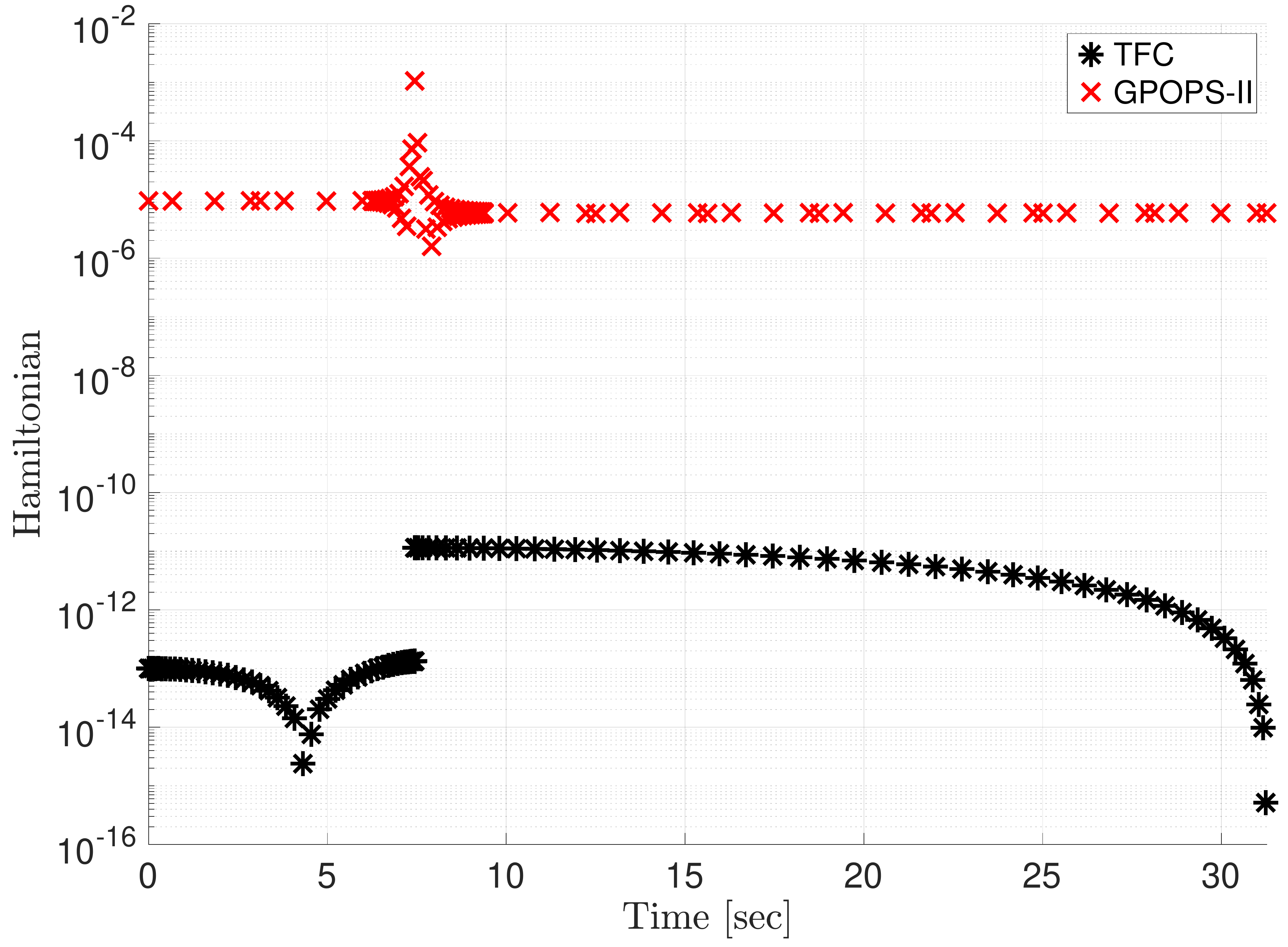}
    \caption{Comparison of Hamiltonian for TFC and GPOPS-II converged solutions for the min-max trajectory. Reprinted with permission from \cite{FOL}.}
    \label{fig:min_max_H}
\end{figure}

\end{example}

\begin{example}{Test 2: Max-Min-Max Trajectory}\label{sec:s6_test2}
In test case 2, the initial conditions were specified such that the optimal solution exhibited a max-min-max profile, i.e., the thrust switches twice, max-min and min-max. The boundary conditions for this case are provided in Table \ref{tab:max_min_max}, whereas Figure \ref{fig:max_min_max_traj} reports the shape of the trajectory computed using the TFC-based algorithm.
\begin{table}[H]
\caption{Boundary conditions for max-min-max trajectory profile test case. Reprinted with permission from \cite{FOL}.}
\centering
\begin{tabular}{|c||c|c|}
\hline
Variable & Initial & Final \\
\hline\hline
$\B{r}$ [m] & $\begin{Bmatrix} -200, & 100, & 1500\end{Bmatrix}\T$ & $\begin{Bmatrix} 0, & 0, & 0\end{Bmatrix}\T$ \\ \hline
$\B{v}$ [m/s] & $\begin{Bmatrix}85, & 50, & -65\end{Bmatrix}\T$ & $\begin{Bmatrix} 0, & 0, & 0\end{Bmatrix}\T$ \\ \hline
$m$ [kg] & $1905$ & - \\ \hline
\hline
\end{tabular}
\label{tab:max_min_max}
\end{table}
\begin{figure}[H]
    \centering\includegraphics[width=.75\linewidth]{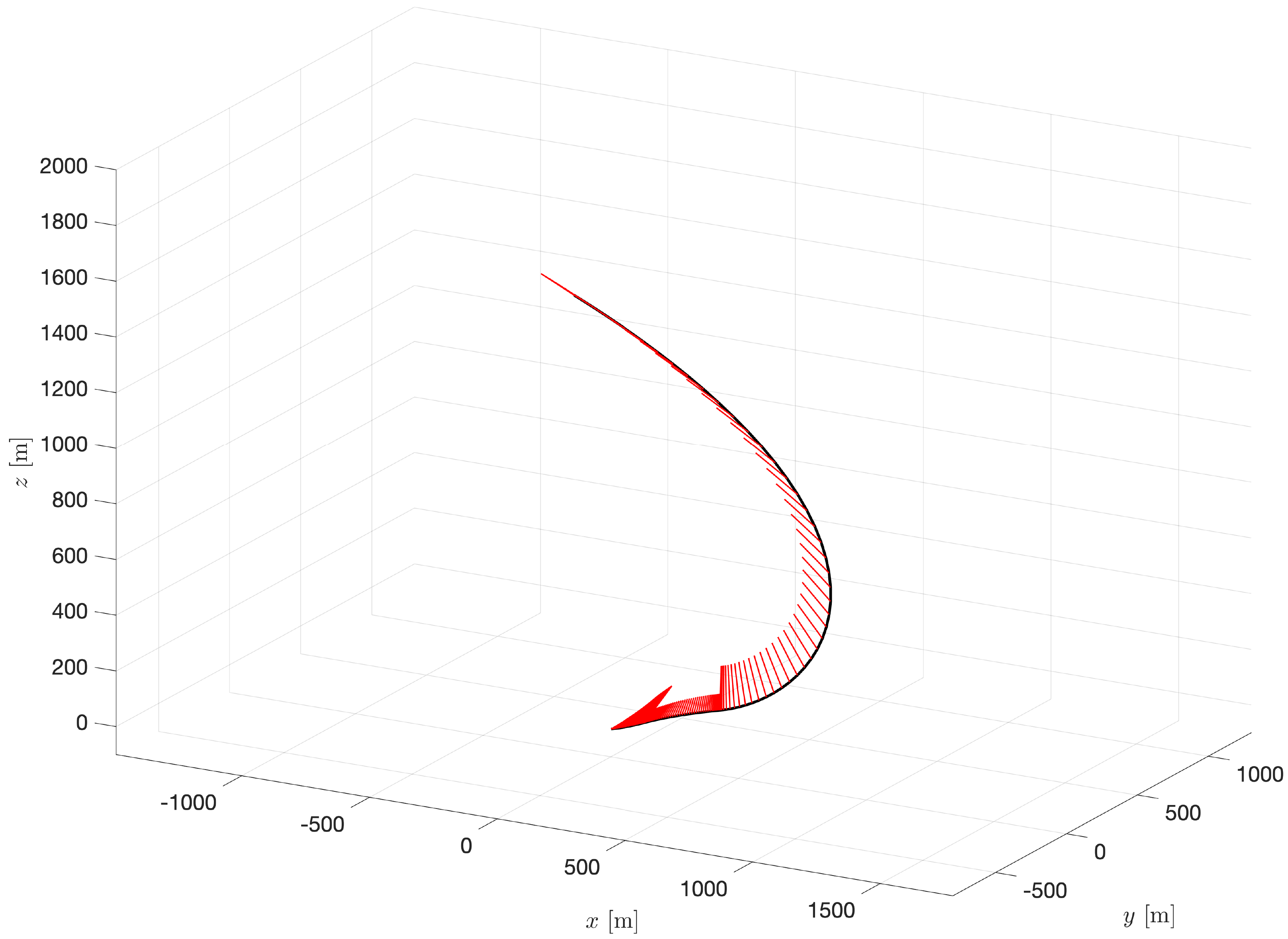}
    \caption{Landing trajectory for max-min-max thrust profile based on initial conditions, $\B{r}_0 = \protect\begin{Bmatrix} -200, & 100, & 1500\protect\end{Bmatrix}\T $[m], $\B{v}_0 = \protect\begin{Bmatrix} 85, & -50, & -65\protect\end{Bmatrix}\T$ [m/s], $m_0 = 1905$ [kg]. Reprinted with permission from \cite{FOL}.}
    \label{fig:max_min_max_traj}
\end{figure}
Again, the TFC solution history is reported for each component of position, velocity, and acceleration in Figure \ref{fig:max_min_max_tfc}. The error is quantified by the residual of the governing equation of motion and the mass costate equation. It can be seen that the TFC residual is $\mathcal{O}(10^{-12})$ or less for the whole solution domain.
\begin{figure}[H]
    \centering\includegraphics[width=.95\linewidth]{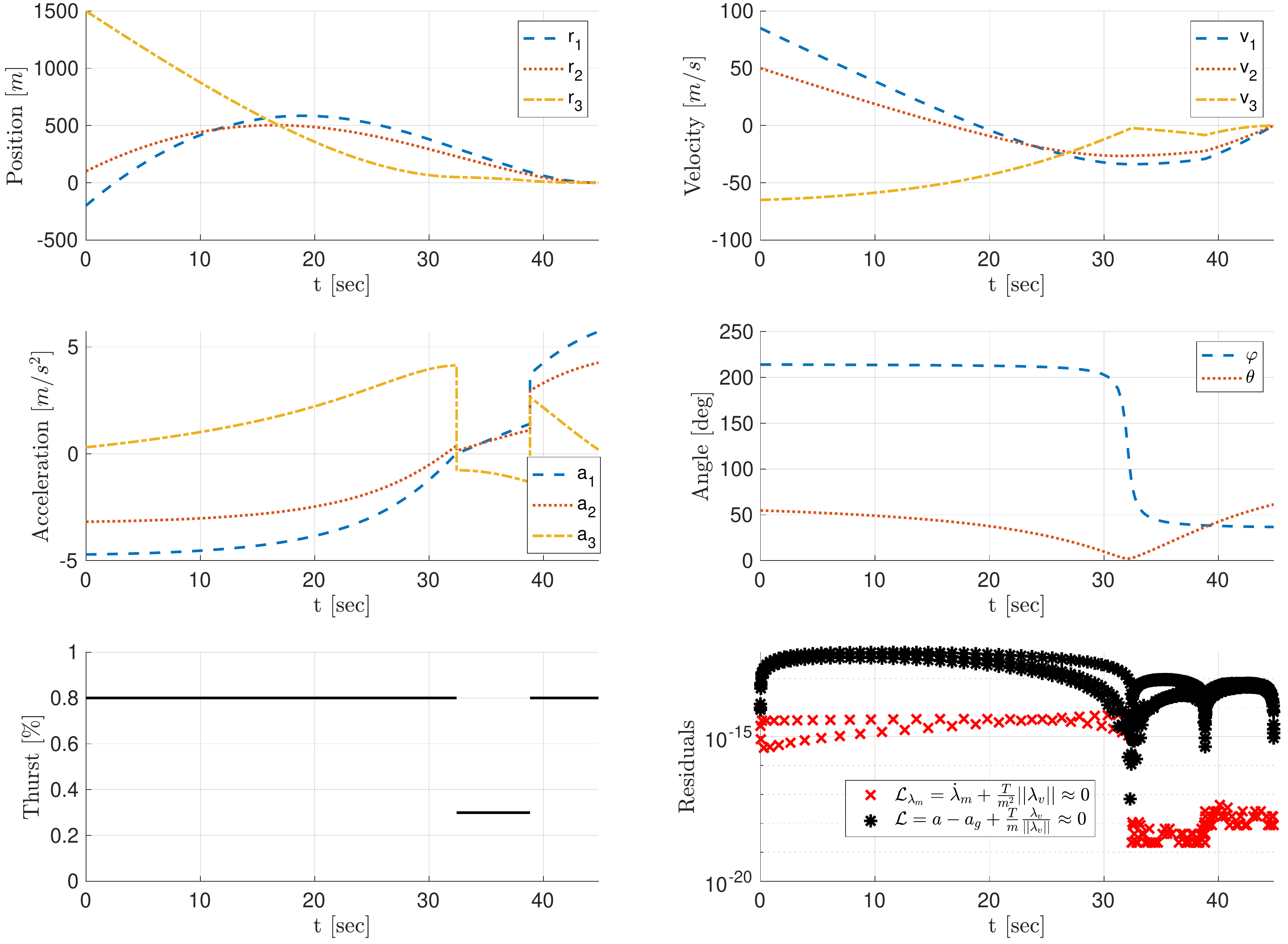}
    \caption{TFC solution of the max-min-max thrust profile case. The solution is presented in terms of the position, velocity, acceleration, and residuals of the differential equations. Reprinted with permission from \cite{FOL}.}
    \label{fig:max_min_max_tfc}
\end{figure}
Similar to test case 1, the solution is compared with the one obtained via GPOPS-II \cite{GPOPS} for all converged parameters, which now includes another switching time, $t_2$. It can be seen that the magnitude of associated errors is similar to those presented in Section \ref{sec:s6_test1}. In this test, \verb"fsolve" iterated 32 times with each TFC inner-loop averaging 81 ms, resulting in a total execution time of 2.6 seconds within the MATLAB implementation. Further, during this test, the TFC method converged in about 3 iterations every \verb"fsolve" function call. Lastly, the propagated comparison to GPOPS is provided in Table \ref{tab:max_min_max_results}, and the Hamiltonian of the two methods is plotted as a function of time in Figure \ref{fig:max_min_max_H} to highlight the optimality of both solutions.
\begin{table}[H]
\caption{Converged parameters for the TFC and GPOPS-II solution for the max-min-max trajectory test case. The values $||\B{r}(t_f)||$, $||\B{v}(t_f)||$, and $\lambda_m(t_f)$ were determined by propagating both TFC and GPOPS-II converged solutions in order to have a one-to-one comparison on the accuracy of the converged solutions. Reprinted with permission from \cite{FOL}.}
\centering
\begin{tabular}{|c||c|c|c|}
\hline
Variable & TFC & GPOPS-II \cite{GPOPS} \\
\hline\hline
$L_2[\mathbb{L}]$ & $5.654\cdot 10^{-12}$ & $-$  \\\hline
$L_2[H]$ & $8.686\cdot 10^{-8}$ & $6.418 \cdot 10^{-3}$  \\\hline
$m_{\text{used}}$ [kg]& $275.205$ & $275.206$ \\\hline
$t_1$ [s] &$32.418$& $32.417$ \\\hline 
$t_2$ [s]&$38.838$ & $38.833$  \\\hline 
$t_f$ [s]&$44.823$ & $44.823$  \\\hline
$||\B{r}(t_f)||$ [m] & $8.330\cdot 10^{-10}$ & $1.350\cdot 10^{-1}$  \\\hline 
$||\B{v}(t_f)||$ [m] & $2.812\cdot 10^{-11}$ & $2.077\cdot 10^{-2}$  \\\hline 
$\lambda_m(t_f)$ [s] & $-8.815\cdot 10^{-15}$ & $-7.354\cdot 10^{-6}$  \\
\hline
\hline
\end{tabular}
\label{tab:max_min_max_results}
\end{table}

\begin{figure}[H]
    \centering\includegraphics[width=.75\linewidth]{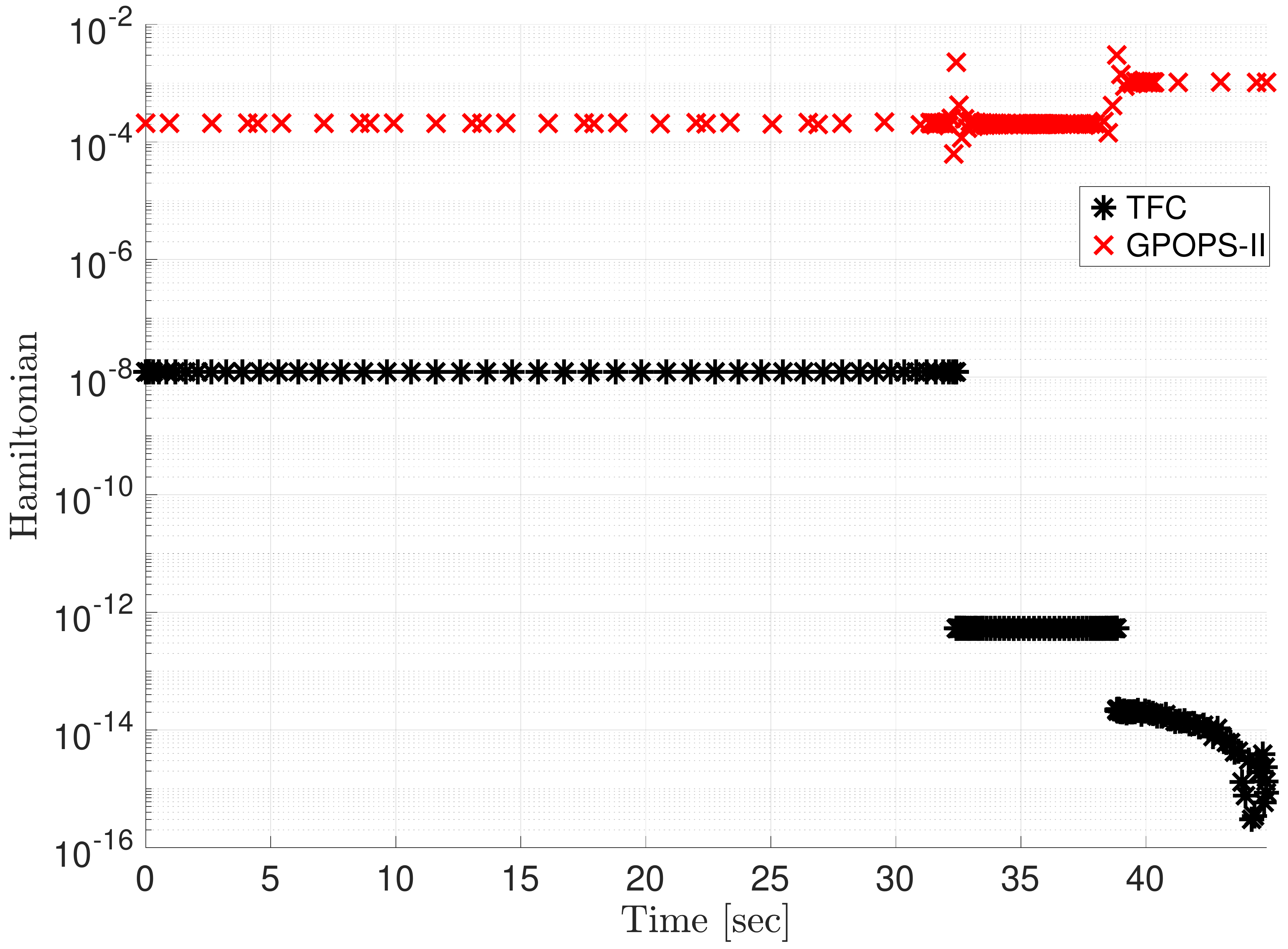}
    \caption{Comparison of Hamiltonian for TFC and GPOPS-II converged solutions for the max-min-max trajectory. Reprinted with permission from \cite{FOL}.}
    \label{fig:max_min_max_H}
\end{figure}
\end{example}

\section{Major findings and conclusions of results}

In all, the current implementation of TFC to the fuel-optimal landing problem \emph{cannot} be used in real-time applications. While the accuracy and speed, once written to a compiled language, are acceptable, the algorithm's robustness is the limiting factor. For example, the Monte Carlo test conducted in Chapter \ref{chap:eol} could not be run for this algorithm. This and other conclusions are summarized below:

\begin{blankBox}{Major takeaways from fuel-optimal landing tests}
\begin{enumerate}
    \item As illustrated in Figure \ref{fig:algorithm}, the proposed TFC-based algorithm requires an efficient implementation of the iterative least-square together with a root-finding algorithm (e.g., Trust-Region-Dogleg algorithm \cite{conn2000trust} as implemented in the \verb"fsolve" routine in MATLAB).
    \item As reported in the numerical tests presented in Examples \ref{sec:s6_test1} and \ref{sec:s6_test2}, the \verb"fsolve" routine iterates for up to 32 times with an upper bound on the execution time of about 2.6 seconds to generate one optimal trajectory, using MATLAB.
    \begin{itemize}
        \item It is known that the MATLAB programming language is about 10 times slower than a C++ executable, which is usually employed to run algorithms on the spacecraft onboard microprocessor.
        \item Therefore, a computational time gain of at least one order of magnitude is expected, thus making the algorithm attractive for real-time implementations with regards to speed.
    \end{itemize}
    \item While the problem was solved with acceptable speed and accuracy, the robustness to the initialization of the times $t_1$, $t_2$, and $t_f$ caused convergence issues that are not acceptable for real-time implementation.
    \begin{itemize}
        \item In this dissertation, two specific cases were solved for the fuel-optimal landing problem but ``hand-tuning'' was necessary for reliable convergence. 
        \item Future work could look remove the necessity of the outer-loop; however, from other studies on free final time problems, this problem may be sufficiently complex such that a single-loop least-squares, like that of Chapter \ref{chap:eol}, will not work. 
    \end{itemize}
    \item A major concern of this technique may be the trade-off in the amount of work in formulating the problem (and especially computing the constrained expression) compared to other optimization packages. While these terms are formulated analytically, the \href{https://github.com/leakec/tfc}{\textcolor{blue}{\underline{TFC GitHub}}} \cite{TfcGithub} provides a framework such after forming the loss vector, the Jacobian terms are computed through automatic differentiation and do not require analytical formulation.
\end{enumerate}
\end{blankBox}

%

\chapter{SUMMARY AND CONCLUSIONS \label{chap:summary}}
 The work presented is entitled ``A Journey from Theory to Application,'' because it represents a single route through the dense landscape of the Theory of Functional Connections. I have surely not observed, recorded, and studied all aspects along the way. However, this section is my way of creating a map for future work. Through examples presented, the reader should be familiar with the theory and how it is currently applied. To further aid the reader, the code for most of the problems and examples in this dissertation can be found for free on the \href{https://github.com/leakec/tfc}{\textcolor{blue}{\underline{TFC GitHub}}} \cite{TfcGithub}. Moving forward with this section, I look to summarize the major results of this journey along with many potential ideas I have explored. 
 
 The main route of this dissertation began with discussing the fundamentals of TFC and the process to derive \emph{constrained expressions}, which are the heart of the method. For a given set of linear constraints, the constrained expression is a functional that represents \emph{all} functions analytically satisfying the constraints, parameterized by the free function $g(x)$. While a method to derive these constrained expressions was provided in the original work on TFC (Reference \cite{U-TFC}), this dissertation presents a new formulation that exploits the main structure shared by all constrained expressions. This structure, named the switching-projection form, 1) gives a more intuitive approach to derive constrained expression, 2) provides a straightforward and general framework for the derivation of linear type constraints, 3) allows for a plethora of mathematical insights and associated claims on existence and non-uniqueness, and 4) provides a simple and elegant extension to $n$-dimensional constrained expressions. In fact, readers interested in the latter point and their application to partial differential equations are directed to Carl Leake's dissertation: ``The Multivariate Theory of Functional Connections: An $n$-Dimensional Constraint Embedding Technique Applied to Partial Differential Equations'' \cite{KarlDissertation}. In addition to many detailed examples that derive constrained expressions, the first part also provides preliminary insight for an ad-hoc method allowing inequality constraints and some discussion on over-constrained expressions. While the former has been implemented in multiple numerical solutions, the latter topic was an academic exploration that spurred from the realization that constrained expressions could also be derived using a weight least-squares approach and allow for more constraints than the number of support functions used in the derivation. In all, this topic was marginally studied, and the usefulness and potential applications are not well understood.
 
 Following the derivation of constrained expression, the second part of this dissertation focused on applying these functionals to the solution of ODEs. Compared to other numerical techniques, the one based on TFC splits the problem into two separate parts: 1) the constraints and 2) the dynamics. As should be clear from the prior sections, the TFC approach allows for the differential equation constraints to be \emph{analytically} embedded in the constrained expressions. In general, this process transforms the differential equation from a constrained optimization problem into an \emph{unconstrained} optimization problem. Next, by using the constrained expression associated with the differential equation constraints, and by 1) defining the free function, $g(x)$ as some know basis with unknown coefficients and 2) discretizing the domain, the problem is again transformed into an algebraic equation that can be solved with any optimization technique, where $\mathbb{L}(\B{\xi}) =\B{0}$. While in this dissertation, almost all problems were solved with a linear or nonlinear least-squares, except for free final time problems where \verb"fsolve" or differential evolution algorithms were also used, much fruitful research remains in the study of this technique paired with other numerical schemes. In fact, TFC \emph{is not} by itself a numerical scheme, but rather an analytical technique to reduce the computational overhead of numerically approximating the constraints. 
 
 In this part, the approach to solve differential equations was highlighted by numerous examples, starting with a simple initial-value problem and ending with complex cases such as systems of differential equations with terminal algebraic constraints and an unknown domain length. In fact, the latter examples of part two of this dissertation focused on unique corner cases that are relevant in ODEs, including 1) a technique for split domain problems and its application to 2) hybrid systems (differential equations with jumps in dynamics), 3) unknown domain length, (free final time) problems relevant in optimal control, and 4) the computation of periodic orbits, which constrained expressions provide a simple and elegant approach to tackle. Finally, some examples of the application of over-constrained constrained expressions were provided.
 
The final part of this dissertation leveraged the prior sections to solve specific aerospace engineering problems, namely terminal descent spacecraft landing on large planetary bodies. These problems were formulated using the indirect method, where the optimal control problem is transformed into a set of differential and algebraic equations that much be solved simultaneously. While this approach is known to produce more optimal solutions than the direct method, the indirect method has a few major drawbacks: 1) the size of the system is doubled with the incorporation of the costates (Lagrange multipliers), and 2) that these costates are highly sensitive to initialization. Therefore, in practice, the indirect method is used less often. Furthermore, while many other numerical approaches exist to solve these types of problems, the motivation to use TFC was that the constrained expressions would provide 1) added robustness to initialization and 2) faster solution speeds. The benefits are not as drastic as first hypothesized for the two problems studies, energy-optimal and fuel-optimal landing. While the TFC solution to the energy-optimal landing did show increased robustness, speed, and accuracy over the spectral method, the solution to the fuel-optimal landing problem lacked robustness and could only be solved for particular cases. In its current state, the TFC algorithm is not quite robust enough. Future improvements could still lead to a technique that could be leveraged to solve trajectories on-board and in real-time by recomputing the optimal trajectory at every computer guidance cycle.
 
First, the energy-optimal landing problem was analyzed for constant gravity cases. This problem has an analytical feedback solution and was used to evaluate the accuracy of the TFC method versus the spectral method and highlight the benefits of TFC. These results showed that the method built with TFC was more accurate, faster, and more robust to poor initialization. Moving forward, the lessons learned from the energy-optimal problem were translated to the fuel-optimal landing problem with one distinct difference: in the energy optimal problem, the final time was solved using a single-loop approach where all the TFC parameters were solved for simultaneously. However, it was found that this method only works for a selection of problems\footnote{The author has found that this approach also does not work for many problems in trajectory optimization, e.g., minimum-time orbit transfer with a solar sail \cite{XTFC-Traj}.}, including problems where the domain has more than one segment due to the dynamics' switching behavior, as seen in the fuel-optimal landing problem. For this reason, the fuel-optimal landing problem was solved using an inner- and outer-loop approach where the TFC method solved the problem for the fixed time cases, i.e., where the switching times and the final time where specified ($t_1, t_2, t_f$), and an outer-loop was used to determine the values of these times. The drawback of this is that the algorithm relies on an external optimizer and increases computation time; in this problem, MATLAB's \verb"fsolve" algorithm was unitized. While two specific cases were solved, showing that a solution can be obtained using TFC, the algorithm is not fit for implementation as a real-time controller in its current state. The current issues with this algorithm include 1) the lack of robustness to the initialization of $t_1$, $t_2$, and $t_f$, 2) the inability to solve the problem with a priori knowledge of the control structures, i.e., max, min-max, or max-min-max thrust arcs, and 3) no guarantees on the convergence of the algorithm. 
 
 To remedy these concerns, more research needs to be done to identify other optimization techniques that could be used in both the inner- and outer-loops of the algorithms. Additionally, the entirety of this work focuses on solving the problems derived using the indirect method. This leaves the area of direct optimization completely untouched and ripe for exploration. 
 
\section{Future research}
Based on the discussion above, I have chosen to include this section to discuss the current and most fruitful paths in the study of TFC related to the topics covered in this dissertation. In this section, I look to provide key insight into topics most likely to yield widespread improvements to the technique and its applications.

\subsection{In search of a free function} 

At the heart of TFC is the constrained expression, which can describe all functions satisfying a set of constraints. The reader should recall that the constrained expression has a free function, $g(x)$, which does not affect the constraints. In numerical applications such as solving differential equations or optimal control problems, the free function must be numerically approximated. Therefore, the representation of the free function is vital in the overall ability to solve problems; however, an in-depth study of this topic is lacking in this dissertation---along with the entire body of research of TFC.

While in this dissertation I mainly focused on the Legendre and Chebyshev orthogonal polynomials, other papers on TFC have looked into using Extreme Learning Machines \cite{XTFC} (mentioned briefly in Chapter \ref{chap:ode}) and Neural Networks (Deep-TFC) \cite{deepToC} to approximate the free function. However, the work on Deep-TFC has only used fully connected NNs up to this point, and the study of different NN architectures is an active area of research.

According to all of the research conducted to date, orthogonal polynomials for most problems are highly effective and produce solutions near machine-level precision. However, when dealing with complex problems, e.g., Naiver-Stokes equations or PDEs with sharp gradients, the Neural Network approach is more accurate. In general, the only benefit of using ELMs is in the low memory case for the solution of PDEs where the number of basis functions is reduced.

Regardless, there is major promise with the study of particular definitions of $g(x)$ leveraging some a priori knowledge of the problem dynamics. To explain this concept and shed light on a potential area of further research, consider a boundary-value problem in trajectory design that includes many revolutions (or orbits) and dynamics that are not purely Keplerian (there are perturbations due to third-body effects, the sun, etc.). In this case, to accurately determine a solution, the function of $g(x)$ must capture both the periodicity of the orbit and the orbit changes due to perturbation. One idea to solve this problem would be to use a hybrid basis composed of terms to individually capture the periodic and non-periodic portions individually.

\subsection{Other optimization schemes}\label{sec:con_OtherOpts}

Next, as mentioned in the previous section, TFC \emph{is not a numerical optimization technique}, but rather an analytical method than can be coupled with any optimization scheme that can solve $\mathbb{L}(\B{\xi}) = \B{0}$. In this dissertation, along with every paper other than Deep-TFC \cite{deepToC}, the optimization scheme used to determine the $\B{\xi}$ coefficients of the free function $g(x) = \B{\xi}\T \B{h}(x)$ were based on a simple linear or nonlinear least-squares. This was done for two reasons: 1) the simplicity of the method and the fact that 2) most problems did not require a more complex method. Outside of this dissertation, along with least-squares, Leake \cite{KarlDissertation} studied the use of three other optimizers, including Limited-memory  Broyden-Fletcher-Goldfarb-Shanno algorithm, Adam (a first-order gradient-based optimization of stochastic objective functions), and constrained support vector machines for the solution of differential equations.

However, just as I have discussed with the definition of the free function, an exploration of a wide range of numerical optimization techniques should be the focus of future work in the application of TFC. For the increasing complexity of problems, this will also be a necessity. Above all, TFC can reduce the set of admissible functions and has the potential to speed up many optimization techniques. 

One of the potential areas of research is pairing TFC with other optimization schemes within direct optimization. In this dissertation and all other work utilizing TFC, optimal control problems were solved using the indirect method. Similar to how convex optimization is used to convert nonconvex problems into convex problems to assure convergence with NLP solvers, there is potential that the TFC constrained expressions can be used to complement current NLP solvers.

\section{Additional Literature on TFC}\label{sec:outro}

In this section, I look to provide the reader with the most up-to-date capabilities of the theory and the many areas not covered in this dissertation. In all, I hope that the text is a springboard for interested researchers that provides references to all prior work and gives a clear path to more fruitful studies in this area. The list below provides a short description of each paper's contribution along with the links (the PDF file provides clickable links).

\subsection{Functional Interpolation}
\begin{itemize}
    \item Mortari, D. The Theory of Connections: Connecting Points. \emph{Mathematics} \textbf{2017}, 5(4), 57; \href{https://doi.org/10.3390/math5040057}{[\textbf{Link}]}
    \begin{adjustwidth*}{1cm}{1cm}
    This is the seminal paper on the Theory of Functional Connections. The work presented explores the fundamental idea of functional interpolation using an additive formulation. Constraint interpolation is introduced for points, derivatives, and linear combinations of them. The additive form of functional interpolation is the basis for all subsequent works.
    \end{adjustwidth*}
    \item Johnston, H., Leake, C., Efendiev, Y., and Mortari, D. Selected Applications of the Theory of Connections: A Technique for Analytical Constraint Embedding. \emph{Mathematics} \textbf{2019}, 7(6), 537; \href{https://doi.org/10.3390/math7060537}{[\textbf{Link}]}
    \begin{adjustwidth*}{1cm}{1cm}
     This paper highlights the utility of TFC by introducing various problems that can be solved using this framework, including (1) analytical linear constraint optimization, (2) the brachistochrone problem, (3) over-constrained differential equations; (4) inequality constraints; and (5) triangular domains.
    \end{adjustwidth*}
    
    \item Mortari, D. and Leake, C. The Multivariate Theory of Connections. \emph{Mathematics} \textbf{2019}, 7(3), 296; \href{https://doi.org/10.3390/math7030296}{[\textbf{Link}]}
    \begin{adjustwidth*}{1cm}{1cm}
    This paper extends the univariate TFC, introduced by Mortari in 2017, to the multivariate case on rectangular domains with detailed attention to the bivariate case. Although this article's focus is on two-dimensional spaces, the ﬁnal section introduces the multivariate TFC, validated by a mathematical proof; this section describes how to write constrained expressions on rectangular domains for an arbitrary number of constraints with arbitrary order derivatives in $n$-dimensions. In all, this last section was the first iteration of what is later presented in ``The Multivariate Theory of Functional Connections: Theory, Proofs, and Application in Partial Differential Equations.''
    \end{adjustwidth*}
    
    \item Wang, Y. and Topputo, F. A Homotopy Method Based on Theory of Functional Connections. \emph{arXiv} \textbf{2019}; \href{https://arxiv.org/abs/1911.04899}{[\textbf{Link}]}
    \begin{adjustwidth*}{1cm}{1cm}
    A method for solving zero-finding problems is developed by tracking homotopy paths, which define connecting channels between an auxiliary problem and the objective problem. Current algorithms’ success relies heavily on empirical knowledge, as the homotopy paths must be selected manually. This work introduces a homotopy method based on TFC. The TFC-based method implicitly defines infinite homotopy paths, from which the most promising ones are selected. A two-layer continuation algorithm is devised, where the first layer tracks the homotopy path by monotonously varying the continuation parameter, while the second layer recovers possible failures and resorts to a TFC representation of the homotopy function. Compared to pseudo-arclength methods, the proposed TFC-based method retains the simplicity of direct continuation while allowing for flexible path switching.
    \end{adjustwidth*}
    
    \item Leake, C., Johnston, H., and Mortari, D. The Multivariate Theory of Functional Connections: Theory, Proofs, and Application in Partial Differential Equations. \emph{Mathematics} \textbf{2020}, 8(8), 1303; \href{https://doi.org/10.3390/math8081303}{[\textbf{Link}]}
    \begin{adjustwidth*}{1cm}{1cm}
    This article exploits constrained expressions' underlying functional structure to ease their derivation and provides mathematical proofs regarding their properties. Furthermore, the extension of the technique to and proofs in $n$-dimensions is immediate through a recursive application of the univariate formulation.
    \end{adjustwidth*}
    
    \item Mortari, D. and Arnas, D. Bijective Mapping Analysis to Extend the Theory of Functional Connections to Non-Rectangular 2-Dimensional Domains. \emph{Mathematics} \textbf{2020}, 8(9), 1593; \href{https://doi.org/10.3390/math8091593}{[\textbf{Link}]}
    \begin{adjustwidth*}{1cm}{1cm}
    This work presents an initial analysis of using bijective mappings to extend TFC to non-rectangular, two-dimensional domains. Speciﬁcally, this manuscript proposes three different mapping techniques: 1) complex mapping, 2) the projection mapping, and 3) polynomial mapping. In that respect, an accurate least-squares approximated inverse mapping is also developed for those mappings with no closed-form inverse. 
    \end{adjustwidth*}
    
    \item Mortari, D. and Furfaro, R. Univariate Theory of Functional Connections Applied to Component Constraints, \emph{Math. Comput. Appl.} \textbf{2021},  26(1), 9; \href{https://doi.org/10.3390/mca26010009}{[\textbf{Link}]}
    \begin{adjustwidth*}{1cm}{1cm}
    This work presents a methodology to derive analytical functionals, with embedded linear constraints among the components of a vector (e.g., coordinates) that is a function a single variable (e.g., time). This work prepares the background necessary for the indirect solution of optimal control problems via the application of the Pontryagin Maximum Principle. The methodology presented is part of the univariate Theory of Functional Connections that has been developed to solve constrained optimization problems. To increase the clarity and practical aspects of the proposed method, the work is mostly presented via examples of applications than via rigorous mathematical deﬁnitions and proofs.
    \end{adjustwidth*}
    
\end{itemize}

\subsection{Solution of Differential Equations}
\begin{itemize}
    \item Mortari, D. Least-Squares Solution of Linear Differential Equations. \emph{Mathematics} \textbf{2017}, 5(4), 48; \href{https://doi.org/10.3390/math5040048}{[\textbf{Link}]}
    \begin{adjustwidth*}{1cm}{1cm}
    This is the first work utilizing the TFC method to solve linear ordinary differential equations. Herein, the constrained expressions from the TFC framework are used to embed the differential equation constraints, and the free function is defined by Chebyshev and Legendre polynomials. The process converts a differential equation subject to constraints to a linear system of equations that is solved via linear least-squares. The method is thus a unified way to solve initial-, boundary-, and multi-value problems.
    \end{adjustwidth*}
    
    \item Johnston, H. and Mortari, D. Linear Differential Equations Subject to Relative, Integral, and Infinite Constraints. Proceedings of the \emph{AAS/AIAA Astrodynamics Specialist Conference} \textbf{2018}, 167, AAS 18-273, pp. 3107-3121, Snowbird, UT, August 19-23, 2018; \href{https://www.researchgate.net/publication/344455474_Linear_Differential_Equations_Subject_to_Relative_Integral_and_Infinite_Constraints}{[\textbf{Link}]}
    \begin{adjustwidth*}{1cm}{1cm}
    This study looks into extending TFC to incorporate relative, integral, and inﬁnite constraints in the solution of differential equations. The results obtained by this method are then compared in terms of speed and accuracy with the solution provided by the \verb"Chebfun" toolbox and are shown to be more accurate with reduced computation time (two orders of magnitude). The new TFC switching-projection form in this dissertation updates the results of this paper.
    \end{adjustwidth*}
    
    \item Johnston, H. and Mortari, D. Weighted Least-Squares Solutions of Over-Constrained Differential Equations. Proceedings of the \emph{International Academy of Astronautics SciTech Forum} \textbf{2018}, AAS 18-812, Moscow, Russia, November 13-15, 2018; \href{https://www.researchgate.net/publication/344455463_Weighted_Least-Squares_Solutions_of_Over-Constrained_Differential_Equations}{[\textbf{Link}]}
    \begin{adjustwidth*}{1cm}{1cm}
    The main purpose of this paper was to explore the ability to derive over-constrained expressions. These constrained expressions satisfy the constraints subject to some relative weighting. They can be used to solve over-constrained differential equations, i.e., it is desired to incorporate more measurements than the order of the differential equation. The contents of this have been refreshed and are included in this dissertation.
    \end{adjustwidth*}
    
    \item Mortari, D., Johnston, H., and Smith, L. High accuracy least-squares solutions of nonlinear differential equations, \emph{Journal of Computational and Applied Mathematics} \textbf{2019}, Vol. 352, pp. 293-307; \href{https://doi.org/10.1016/j.cam.2018.12.007}{[\textbf{Link}]}
    \begin{adjustwidth*}{1cm}{1cm}
    The techniques developed in Mortari's ``Least-Squares Solution of Linear Differential Equations'' are extended to nonlinear differential equations by implementing a nonlinear least-squares method. This technique is compared to MATLAB's \verb"ode45" and the \verb"Chebfun" package. Additionally, the paper provides the initial scheme to handle long propagation times and is tested on the simple and duffing oscillator.
    \end{adjustwidth*}

    \item Leake, C., Johnston, H., Smith, L., and Mortari, D. Analytically Embedding Differential Equation Constraints into Least Squares Support Vector Machines Using the Theory of Functional Connections. \emph{Mach. Learn. Knowl. Extr.} \textbf{2019}, 1(4), 1058-1083; \href{https://doi.org/10.3390/make1040060}{[\textbf{Link}]}
    \begin{adjustwidth*}{1cm}{1cm}
    This work merges least-squares support vector machines (LS-SVM) with TFC to produced a technique called constrained SVMs (CSVM). In general, TFC is shown to be slightly faster (by an order of magnitude or less) and more accurate (by multiple orders of magnitude) than the LS-SVM and CSVM approaches. Therefore, this technique is not recommended for use. However, this article was an important step towards integrating TFC with machine learning algorithms.
    \end{adjustwidth*}
    
    \item  Johnston, H., Leake, C., and Mortari. D. An Analysis of the Theory of Functional Connections Subject to Inequality Constraints. Proceedings of the \emph{AAS/AIAA Astrodynamics Specialist Conference} \textbf{2019}, AAS 19-732, Portland, ME, August 11-15, 2019; \href{https://www.researchgate.net/publication/344455395_An_Analysis_of_the_Theory_of_Functional_Connections_Subject_to_Inequality_Constraints}{[\textbf{Link}]}
    \begin{adjustwidth*}{1cm}{1cm}
    This paper is the first work that incorporates inequality constraints into the TFC framework. The work shows how to extend the original theory to problems subject to equality and inequality constraints for one- and two-dimensions. All of the work in this paper has been updated in this dissertation.
    \end{adjustwidth*}
    
   \item Johnston, H. and Mortari, D. Least-squares solutions of boundary-value problems in hybrid systems. \emph{arXiv} \textbf{2019}; \href{https://arxiv.org/abs/1911.04390}{[\textbf{Link}]}
    \begin{adjustwidth*}{1cm}{1cm}
    This paper looks to apply the mathematical framework of TFC to the solution of boundary-value problems arising from hybrid systems (or a sequence of different differential equations). The approach developed in this work derives an analytical constrained expression for the entire range of a hybrid system, enforcing both the boundary conditions and the continuity conditions across the sequence of differential equations. This reduces the solution space of the hybrid system to only admissible solutions. This technique is widely used throughout this dissertation and enables the solution of problems such as fuel-optimal landing.
    \end{adjustwidth*}
    
   \item Leake, C. and Mortari, D. Deep Theory of Functional Connections: A New Method for Estimating the Solutions of Partial Differential Equations. \emph{Mach. Learn. Knowl. Extr.} \textbf{2020}, 2(1), 37-55; \href{https://doi.org/10.3390/make2010004}{[\textbf{Link}]}
    \begin{adjustwidth*}{1cm}{1cm}
    This article uses neural networks as the free function in TFC \ces\ to estimate the solutions of PDEs. Neural networks are not plagued by the same computational curse-of-dimensionality that occurs when using a linear expansion of basis functions as the free function. Neither are they typically trained via least-squares, which is also memory intensive. This new methodology, called Deep-TFC, is advantageous when estimating the solutions of complex PDEs, such as Navier-Stokes, and has broader impacts outside of differential equation solutions: the article's contents can be used to apply constraints to neural networks, which has multiple applications throughout the machine learning community.
    \end{adjustwidth*}
   
   %
   \item Johnston, H., Leake, C., and Mortari, D. Least-Squares Solutions of Eighth-Order Boundary Value Problems Using the Theory of Functional Connections. \emph{Mathematics} \textbf{2020}, 8(3), 397; \href{https://doi.org/10.3390/math8030397}{[\textbf{Link}]}
    \begin{adjustwidth*}{1cm}{1cm}
    This paper shows how to obtain highly accurate solutions of eighth-order boundary-value problems of linear and nonlinear ordinary differential equations. The results highlight that the TFC approach does not lose accuracy based on the order of the differential equation and all problems were solved with error on the order of $\mathcal{O}(10^{-13} - 10^{-16})$. In all problems, TFC outperformed current literature by at least four orders of magnitude. 
    \end{adjustwidth*}
\end{itemize}

\subsection{Optimization and Optimal Control}
\begin{itemize}
    \item Mai, T. and Mortari, D. Theory of functional connections applied to nonlinear programming under equality constraints. \emph{arXiv} \textbf{2019}; \href{https://arxiv.org/abs/1910.04917}{[\textbf{Link}]}
    \begin{adjustwidth*}{1cm}{1cm}
    This paper introduces an efficient approach to solve quadratic programming problems subject to equality constraints via TFC. This is done without using the traditional Lagrange multipliers approach, and the solution is provided in closed-form for two distinct constrained expressions (satisfying the equality constraints). The unknown optimization variable is then the free vector $\B{g}$ introduced by TFC. The solution to the general nonlinear programming problem is obtained by Newton’s method. Each iteration involves the second-order Taylor approximation, starting from an initial vector $\B{x}_0$, which is a solution of the equality constraint. Numerical results are provided, which compare the speed and accuracy of this approach to MATLAB’s \verb"quadprog". Finally, a convergence analysis of NLP using TFC is provided.
    \end{adjustwidth*}

    \item Drozd, K., Furfaro, R., and Mortari, D. Constrained Energy-Optimal Guidance in Relative Motion via Theory of Functional Connections and Rapidly-Explored Random Trees. Proceedings of the \emph{AAS/AIAA Astrodynamics Specialist Conference} \textbf{2019}, AAS 19-662, Portland, ME, August 11-15, 2019; \href{https://www.researchgate.net/publication/335842052_Constrained_Energy-Optimal_Guidance_in_Relative_Motion_via_Theory_of_Functional_Connections_and_Rapidly-Explored_Random_Trees}{[\textbf{Link}]}
    \begin{adjustwidth*}{1cm}{1cm}
    This is a preliminary study that explores using TFC as a fast and reliable TPBVP solver for kinodynamic sample-based motion planners, like RRTs. A trajectory for a deputy satellite that is energy-optimal, successfully rendezvous with a chief satellite, and is governed by the Clohessy-Wiltshire equations of motion (relative motion) is computed. Within the RRT process, multiple solutions from the many TPBVPs solved via TFC are strung together to form a trajectory that also avoids keep-out-zones.
    \end{adjustwidth*}

    \item Furfaro, R. and Mortari, D. Least-squares Solution of a Class of Optimal Guidance Problems via Theory of Connections, \emph{ACTA Astronautica}, \textbf{2020}, Vol. 168, pp. 92-103; \href{https://doi.org/10.1016/j.actaastro.2019.05.050}{[\textbf{Link}]}
    \begin{adjustwidth*}{1cm}{1cm}
    This paper is the first application of TFC to solve the TPBVPs derived from the indirect method of optimal control. The examples solved in this work include a class of optimal guidance problems, including energy-optimal landing on planetary bodies (where time is fixed for the TFC loop) and fixed-time optimal intercept for a target-interceptor scenario.
    \end{adjustwidth*}

    \item Johnston, H., Schiassi, E., Furfaro, R. and  Mortari, D. Fuel-Efficient Powered Descent Guidance on Large Planetary Bodies via Theory of Functional Connections. \emph{J Astronaut Sci}  \textbf{2020}; \href{https://doi.org/10.1007/s40295-020-00228-x}{[\textbf{Link}]}
    \begin{adjustwidth*}{1cm}{1cm}
    This paper presents a new approach to solve the fuel-efficient powered descent guidance problem on large planetary bodies with no atmosphere (e.g., Moon or Mars). The problem is formulated using the indirect method, which casts the optimal guidance problem as a system of nonlinear two-point boundary value problems that are solved with TFC. In general, the technique produces solutions with error on the order of $\mathcal{O}\left(10^{-10}\right)$. The results of this paper are contained in Chapter \ref{chap:fol} of this dissertation.
    \end{adjustwidth*}
        
    \item Schiassi, E., D’Ambrosio, A., Johnston, H., Furfaro, R., Curti, F., and Mortari, D. Complete Energy Optimal Landing on Small and Large Planetary Bodies via Theory of Functional Connections. Proceedings of the \emph{AAS/AIAA Astrodynamics Specialist Conference} \textbf{2020}, AAS 20-557, Lake Tahoe, CA, August 9-13, 2020; \href{https://www.researchgate.net/publication/343628030_Complete_Energy_Optimal_Landing_on_Small_and_Large_Planetary_Bodies_via_Theory_of_Functional_Connections}{[\textbf{Link}]}
    \begin{adjustwidth*}{1cm}{1cm}
    This paper proposes a unified approach to solve the energy optimal landing on a planetary body (e.g., planet, asteroid, comet, etc.). The method accurately computes the energy optimal landing trajectories, including the optimal time of flight, with a computation time on the order of 10-100 milliseconds, using MATLAB. The algorithms developed from this theory are validated for the landing ﬁnal descent phase in Gaspra and Bennu asteroids and Mars.
    \end{adjustwidth*}
        
    \item Schiassi, E., D’Ambrosio, A., Johnston, H., De Florio, M., Drozd, K., Furfaro, R., Curti, F., and Mortari, D. Physics-Informed Extreme Theory of Functional Connections Applied to Optimal Orbit Transfer. Proceedings of the \emph{AAS/AIAA Astrodynamics Specialist Conference} \textbf{2020}, AAS 20-524, Lake Tahoe, CA, August 9-13, 2020; \href{https://www.researchgate.net/publication/343627850_Physics-Informed_Extreme_Theory_of_Functional_Connections_Applied_to_Optimal_Orbit_Transfer}{[\textbf{Link}]}
    \begin{adjustwidth*}{1cm}{1cm}
    This paper looks to solve a class of trajectory optimization problems using the TFC framework with the free function defined as a single-layer NN. This technique, referred to as X-TFC, is used to solve the system of differential equations derived through the indirect method of optimal control. The problems studied include the Feldbaum problem, minimum time orbit transfer, and maximum radius orbit transfer. 
    \end{adjustwidth*}
\end{itemize}

\subsection{Astrodynamics}
\begin{itemize}
    \item  Johnston, H. and Mortari. D. The Theory of Connections Applied to Perturbed Lambert’s Problem. Proceedings of the \emph{AAS/AIAA Astrodynamics Specialist Conference} \textbf{2018}, AAS 18-282, Snowbird, UT, August 19-23, 2018; \href{https://www.researchgate.net/publication/344455627_The_Theory_of_Connections_Applied_to_Perturbed_Lambert's_Problem}{[\textbf{Link}]}
    \begin{adjustwidth*}{1cm}{1cm}
    This paper formulates the perturbed Lambert's problem, a boundary-value problem, in the TFC framework such that the method uses an unperturbed solution as the baseline (or initial guess) and looks to add all perturbations simultaneously with the constrained expression. The results and theory of this paper are dated, and the major issue with this work is that the constrained expressions capturing the perturbations are added to the numerical solution of the unperturbed Lambert's solver. This causes numerical issues and is remedied by only using the unperturbed Lambert's solution as an initial guess to a constrained expression describing the full solution. The updated approach to solve this problem is provided in ``Evaluation of transfer costs in the Earth-Moon system using the Theory of Functional Connections.''
    \end{adjustwidth*}
    
    \item  Johnston, H. and Mortari. D. Orbit Propagation via the Theory of Functional Connections. Proceedings of the \emph{AAS/AIAA Astrodynamics Specialist Conference} \textbf{2019}, AAS 19-736, Portland, ME, August 11-15, 2019; \href{https://www.researchgate.net/publication/344455392_Orbit_Propagation_via_the_Theory_of_Functional_Connections}{[\textbf{Link}]}
    \begin{adjustwidth*}{1cm}{1cm}
    Spurring from the study of Lambert's problem, this paper investigates the accuracy of TFC applied to the perturbed orbit propagation (initial-value) problem. The method is analyzed for accuracy and convergence behavior and is compared with the \verb"ode113" propagator and the F \& G method. This paper shows that TFC is comparable to other techniques but is better suited for boundary-value problems.
    \end{adjustwidth*}
    
    \item  de Almeida Jr., A. K., Johnston, H., Leake, C., and Mortari. D. Evaluation of transfer costs in the Earth-Moon system using the Theory of Functional Connections. Proceedings of the \emph{AAS/AIAA Astrodynamics Specialist Conference} \textbf{2020}, AAS 20-596, Lake Tahoe, CA, August 9-13, 2020; \href{https://www.researchgate.net/publication/346679594_Evaluation_of_Transfer_Costs_in_the_Earth-Moon_System_using_the_Theory_of_Functional_Connections}{[\textbf{Link}]}
    \begin{adjustwidth*}{1cm}{1cm}
    This paper uses TFC to analyze the mission design space of the two-impulse maneuver Earth-Moon orbit transfer problem by evaluating $\Delta V$ as a function of time of flight and other parameters, like the points of application of the thrusts. 
    Transfers from low-Earth orbit to the L1 Lagrange point and near-Earth orbit to a near-Moon orbit are analyzed as functions of the departure position and the time of ﬂight. Furthermore, the inﬂuence of perturbations due to the gravitational attraction of the Sun is also investigated.
    \end{adjustwidth*}
    
    \item Johnston, H., Lo, M., and Mortari, D. A Functional Interpolation Method to Compute Period Orbits in the Circular Restricted Three-Body Problem. Proceedings of the \emph{31st AAS/AIAA Space Flight Mechanics Meeting} \textbf{2021}, AAS 21-257, Virtual, February 1-4, 2021; \href{https://www.researchgate.net/publication/349151277_A_Functional_Interpolation_Approach_to_Compute_Periodic_Orbits_in_the_Circular_Restricted_Three-Body_Problem}{[\textbf{Link}]}
    \begin{adjustwidth*}{1cm}{1cm}
    In this paper, we develop a method to solve for periodic orbits, i.e. Lyapunov and Halo orbits, using a functional interpolation scheme called the Theory of Functional Connections (TFC). Using this technique, a periodic constraint is analytically embedded into the TFC constrained expression. By doing this, the system of differential equations governing the three-body problem is transformed into an unconstrained optimization problem where simple numerical schemes can be used to find a solution, e.g. nonlinear least-squares. This allows for a simpler numerical implementation with comparable accuracy and speed to the traditional differential corrector method.
    \end{adjustwidth*}
\end{itemize}

\subsection{Transport Theory}

\begin{itemize}
    \item De Florio, M. Accurate Solutions of the Radiative Transfer Problem via Theory of Connections. Thesis for: \emph{MSc in Energy and Nuclear Engineering} \textbf{2019}; \href{https://www.researchgate.net/publication/341786871_Accurate_Solutions_of_the_Radiative_Transfer_Problem_via_Theory_of_Connections?channel=doi&linkId=5ed48ae24585152945279695&showFulltext=true}{[\textbf{Link}]}
    \begin{adjustwidth*}{1cm}{1cm}
    In this thesis, a new approach to solve a class of radiative transfer problems is presented using TFC to solve the linear one-point boundary-value problem derived from the Boltzmann integrodifferential equation for radiative transfer. The proposed algorithm resides in the category of numerical methods for the solution of transport equations and is accurate and suitable for applications in atmospheric science and remote sensing.
    \end{adjustwidth*}

    \item De Florio, M., Schiassi, E.,  Furfaro, R., Ganapol, B.D., and Mostacci, D. Solutions of Chandrasekhar's Basic Problem in Radiative Transfer via  Theory of Functional Connections. \emph{Journal of Quantitative Spectroscopy and Radiative Transfer}, p.107384. \textbf{2020};
    \href{https://doi.org/10.1016/j.jqsrt.2020.107384}{[\textbf{Link}]}
    \begin{adjustwidth*}{1cm}{1cm}
    In this paper, Chandrasekhar's problem in radiative transfer is solved using TFC. The method is designed to efficiently and accurately solve the linear boundary-value problem arising from the angular discretization of the integrodifferential Boltzmann equation for radiative transfer. The proposed algorithm falls under the category of numerical methods for the solution of radiative transfer equations. The accuracy of this new method is tested by benchmark comparison for Mie and Haze L scattering laws.
    \end{adjustwidth*}
\end{itemize}

\subsection{Physics-Informed Neural Networks}

\begin{itemize}
   \item Schiassi, E., Leake, C., De Florio, M., Johnston, H., Furfaro, R., and Mortari, D. Extreme Theory of Functional Connections: A Physics-Informed Neural Network Method for Solving Parametric Differential Equations. \emph{arXiv} \textbf{2020}; \href{https://arxiv.org/abs/2005.10632}{[\textbf{Link}]}
   \begin{adjustwidth*}{1cm}{1cm}
    This article uses a single layer neural network (NN), or more precisely an Extreme Learning Machine (ELM), as the free function in TFC constrained expressions to estimate the solutions of DEs. The results show that X-TFC achieves high accuracy with low computational time but is never more accurate than the original TFC formulation with orthogonal polynomials for simple problems, nor more accurate than Deep-TFC for complex problems.
    \end{adjustwidth*}

    \item Schiassi, E.,  D'Ambrosio, A., De Florio, M., Furfaro, R., and Curti, F. Physics-Informed Extreme Theory of Functional Connections Applied to Data-Driven Parameters Discovery of Epidemiological Compartmental Models. \emph{arXiv} \textbf{2020}; \href{https://arxiv.org/abs/2008.05554}{[\textbf{Link}]}
    \begin{adjustwidth*}{1cm}{1cm}
    This paper utilizes the X-TFC framework, which combines TFC with the Physics-Informed Neural Networks (PINN) framework for data-driven parameters discovery of problems modeled via ordinary differential equations (ODEs). In particular, this work focuses on the capability of X-TFC in solving inverse problems to estimate the parameters governing the epidemiological compartmental models via a deterministic approach. The epidemiological compartmental models treated in this work are Susceptible Infectious Recovered (SIR), Susceptible Exposed Infectious Recovered (SEIR), and Susceptible Exposed Infectious Recovered Susceptible (SEIRS). The results show that these problems can be accurately solved with low computational times under the influence of unperturbed and perturbed data.
    \end{adjustwidth*}

\end{itemize}

\let\oldbibitem\bibitem
\renewcommand{\bibitem}{\setlength{\itemsep}{0pt}\oldbibitem}
\bibliographystyle{ieeetr}

\phantomsection
\addcontentsline{toc}{chapter}{REFERENCES}

\renewcommand{\bibname}{{\normalsize\rm REFERENCES}}

\bibliography{main}

%

\begin{appendices}
\titleformat{\chapter}{\centering\normalsize}{APPENDIX \thechapter}{0em}{\vskip .5\baselineskip\centering}
\renewcommand{\appendixname}{APPENDIX}

%


\phantomsection

\chapter{ORTHOGONAL BASIS FUNCTIONS}

Since the proposed method uses a set of basis functions, a summary of the candidate orthogonal polynomial basis functions is provided.

\section{Chebyshev} 

Chebyshev Orthogonal Polynomials (CP) of the first kind, $T_k (z)$, are defined on the domain $z\in[-1,+1]$ and are generated using the recursive function,
\begin{equation}\label{A1}
    T_{k + 1} = 2 \, z \, T_k - T_{k - 1} \qquad \text{starting from:} \; \begin{cases} T_0 =& 1 \\ T_1 =& z\end{cases}
\end{equation}
All derivatives of CP can be computed recursively, starting from
\begin{equation*}
    \dfrac{\dd T_0}{\dd z} = 0, \quad \dfrac{\dd T_1}{\dd z} = 1 \qquad \text{and} \qquad \dfrac{\dd^d T_0}{\dd z^d} = \dfrac{\dd^d T_1}{\dd z^d} = 0 \quad (\forall \; d > 1),
\end{equation*}
while the subsequent derivatives of Equation (\ref{A1}) are given for $k \ge 1$,
\begin{equation*}
    \begin{array}{ccccc}
        \dfrac{\dd T_{k+1}}{\dd z} &=& 2 \, \left(T_k + z \, \dfrac{\dd T_k}{\dd z}\right) &- \dfrac{\dd T_{k-1}}{\dd z} \\ [8pt]
        \dfrac{\dd^2 T_{k+1}}{\dd z^2} &=& 2 \left(2 \, \dfrac{\dd T_k}{\dd z} + z \, \dfrac{\dd^2 T_k}{\dd z^2}\right) &- \dfrac{\dd^2 T_{k-1}}{\dd z^2} \\ [4pt]
        \vdots & ~ & \vdots & \vdots \\ [4pt]
        \dfrac{\dd^d T_{k+1}}{\dd z^d} &=& 2 \left( d \, \dfrac{\dd^{d-1} T_k}{\dd z^{d-1}} + z \, \dfrac{\dd^d T_k}{\dd z^d}\right) &- \dfrac{\dd^d T_{k-1}}{\dd z^d}; & (\forall \; d \ge 1).
    \end{array}
\end{equation*}
In particular,
\begin{equation*}
    T_k (-1) = (-1)^k, \quad \left.\dfrac{\dd T_k}{\dd z}\right|_{z=-1} = (-1)^{k+1} \, k^2, \quad \left.\dfrac{\dd^2 T_k}{\dd z^2}\right|_{z=-1} = (-1)^k \, \dfrac{k^2 \, (k^2 - 1)}{3}
\end{equation*}
and
\begin{equation*}
    T_k (1) = 1, \qquad \left.\dfrac{\dd T_k}{\dd z}\right|_{z = 1} = k^2, \qquad \left.\dfrac{\dd^2 T_k}{\dd z^2}\right|_{z = 1} = \dfrac{k^2 \, (k^2 - 1)}{3}.
\end{equation*}

\section{Legendre} 

Legendre Orthogonal Polynomials (LeP), $L_k (z)$, are defined on the domain $z\in[-1,+1]$ and are generated using the recursive function,
\begin{equation}\label{B1}
    L_{k+1} = \dfrac{2k+1}{k+1} \, z \, L_k - \dfrac{k}{k+1} \, L_{k-1} \qquad \text{starting:} \; \begin{cases} L_0 =& 1 \\ L_1 =& z\end{cases}
\end{equation}
All derivatives of LeP can be computed recursively, starting from
\begin{equation*}
    \dfrac{\dd L_0}{\dd z} = 0, \quad \dfrac{\dd L_1}{\dd z} = 1 \qquad \text{and} \qquad \dfrac{\dd^d L_0}{\dd z^d} = \dfrac{\dd^d L_1}{\dd z^d} = 0 \quad (\forall \; d > 1),
\end{equation*}
while the subsequent derivatives of Equation (\ref{B1}) for $k \ge 1$, can be computed in cascade,
\begin{equation*}
    \begin{array}{ccccc}
        \dfrac{\dd L_{k+1}}{\dd z} &=& \dfrac{2k+1}{k+1} \left(L_k + z \dfrac{\dd L_k}{\dd z}\right) & - \dfrac{k}{k+1} \dfrac{\dd L_{k-1}}{\dd z} \\ [8pt]
        \dfrac{\dd^2 L_{k+1}}{\dd z^2} &=& \dfrac{2k+1}{k+1} \left(2\dfrac{\dd L_k}{\dd z} + z \dfrac{\dd^2 L_k}{\dd z^2}\right) & - \dfrac{k}{k+1} \dfrac{\dd^2 L_{k-1}}{\dd z^2} \\ [4pt]
        \vdots & ~ & \vdots & \vdots \\ [4pt]
        \dfrac{\dd^d L_{k+1}}{\dd z^d} &=& \dfrac{2k+1}{k+1} \left(d\dfrac{\dd^{d-1} L_k}{\dd z^{d-1}} + z \dfrac{\dd^d L_k}{\dd z^d}\right) & - \dfrac{k}{k+1} \dfrac{\dd^d L_{k-1}}{\dd z^d}; & (\forall \; d \ge 1).
    \end{array}
\end{equation*}

\section{Laguerre} 

Laguerre Orthogonal Polynomials (LaP), $L_k (z)$, are defined on the domain $z\in[0,\infty)$ and are generated using the recursive function,
\begin{equation*}
    L_{k+1} (z) = \dfrac{2k + 1 - z}{k + 1} \, L_k (z) - \dfrac{k}{k + 1} \, L_{k-1} (z) \qquad \text{starting:} \; \begin{cases} L_0 =& 1 \\ L_1 =& 1 - z\end{cases}
\end{equation*}
All derivatives of LaP can be computed recursively, starting from
\begin{equation*}
    \dfrac{\dd L_0}{\dd z} = 0, \quad \dfrac{\dd L_1}{\dd z} =-1 \qquad \text{and} \qquad \dfrac{\dd^d L_0}{\dd z^d} = \dfrac{\dd^d L_1}{\dd z^d} = 0 \quad (\forall \; d > 1),
\end{equation*}
then
\begin{equation*}
    \begin{array}{ccc}
        \dfrac{\dd L_{k+1}}{\dd z} &=& \dfrac{2k + 1 - z}{k + 1} \dfrac{\dd L_k}{\dd z} - \dfrac{1}{k + 1} L_k - \dfrac{k}{k + 1} \dfrac{\dd L_{k-1}}{\dd z} \\ [8pt]
        \dfrac{\dd^2 L_{k+1}}{\dd z^2} &=& \dfrac{2k + 1 - z}{k + 1} \dfrac{\dd^2 L_k}{\dd z^2} - \dfrac{2}{k + 1} \dfrac{\dd L_k}{\dd z} - \dfrac{k}{k + 1} \dfrac{\dd^2 L_{k-1}}{\dd z^2} \\ [8pt]
        \vdots & ~ & \vdots \\ [4pt]
        \dfrac{\dd^d L_{k+1}}{\dd z^d} &=& \dfrac{2k + 1 - z}{k + 1} \dfrac{\dd^d L_k}{\dd z^d} - \dfrac{d}{k + 1} \dfrac{\dd^{d-1} L_k}{\dd z^{d-1}} - \dfrac{k}{k + 1} \dfrac{\dd^d L_{k-1}}{\dd z^d}
    \end{array}
\end{equation*}

\section{Hermite} 

There are two Hermite Orthogonal Polynomials (HP), the probabilists, indicated by $E_k (z)$ defined on the domain $z\in(-\infty,\infty)$, and the physicists, indicated by $H_k (z)$ also defined on the domain $z\in(-\infty,\infty)$. They both are generated using recursive functions.

The probabilistists are defined as
\begin{equation*}
    E_{k + 1} (z) = z \, E_k (z) - k E_{k - 1} (z) \qquad \text{starting:} \; \begin{cases} E_0 (z) =& 1 \\ E_1 (z) =& z\end{cases}
\end{equation*}
All derivatives can be computed recursively, starting from
\begin{equation*}
    \dfrac{\dd E_0}{\dd z} = 0, \quad \dfrac{\dd E_1}{\dd z} = 1 \qquad \text{and} \qquad \dfrac{\dd^d E_0}{\dd z^d} = \dfrac{\dd^d E_1}{\dd z^d} = 0 \quad (\forall \; d > 1),
\end{equation*}
then
\begin{equation*}
    \begin{array}{ccl}
        \dfrac{\dd E_{k+1}}{\dd z} &=& E_k + z \dfrac{\dd E_k}{\dd z} - k \dfrac{\dd E_{k-1}}{\dd z} \\ [8pt]
        \dfrac{\dd^2 E_{k+1}}{\dd z^2} &=& 2 \dfrac{\dd E_k}{\dd z} + z \dfrac{\dd^2 E_k}{\dd z^2} - k \dfrac{\dd^2 E_{k-1}}{\dd z^2} \\ [8pt]
        \vdots & ~ & \vdots \\ [4pt]
        \dfrac{\dd^d E_{k+1}}{\dd z^d} &=& d \dfrac{\dd^{d-1} E_k}{\dd z^{d-1}} + z \dfrac{\dd^d E_k}{\dd z^d} - k \dfrac{\dd^d E_{k-1}}{\dd z^d}
    \end{array}
\end{equation*}

The physicists are defined as
\begin{equation*}
    H_{k + 1} (z) = 2 z \, H_k (z) - 2 k \, H_{k - 1} (z) \qquad \text{starting:} \; \begin{cases} H_0 (z) =& 1 \\ H_1 (z) =& 2z\end{cases}
\end{equation*}
All derivatives can be computed recursively, starting from
\begin{equation*}
    \dfrac{\dd H_0}{\dd z} = 0, \quad \dfrac{\dd H_1}{\dd z} = 2 \qquad \text{and} \qquad \dfrac{\dd^d H_0}{\dd z^d} = \dfrac{\dd^d H_1}{\dd z^d} = 0 \quad (\forall \; d > 1),
\end{equation*}
then
\begin{equation*}
    \begin{array}{ccl}
        \dfrac{\dd H_{k+1}}{\dd z} &=& 2 H_k + 2 z \dfrac{\dd H_k}{\dd z} - 2k \dfrac{\dd H_{k-1}}{\dd z} \\ [8pt]
        \dfrac{\dd^2 H_{k+1}}{\dd z^2} &=& 4 \dfrac{\dd H_k}{\dd z} + 2 z \dfrac{\dd^2 H_k}{\dd z^2} - 2k \dfrac{\dd^2 H_{k-1}}{\dd z^2} \\ [8pt]
        \vdots & ~ & \vdots \\ [4pt]
        \dfrac{\dd^d H_{k+1}}{\dd z^d} &=& 2 d \dfrac{\dd^{d-1} H_k}{\dd z^{d-1}} + 2 z \dfrac{\dd^d H_k}{\dd z^d} - 2 k \dfrac{\dd^d H_{k-1}}{\dd z^d}
    \end{array}
\end{equation*}

\section{Fourier Basis} 
The Fourier Series (FS) is defined on the domain $z\in[-\pi,\pi]$; however, it does not have a recursive generating function like the other basis sets. In general, the FS can be written as
\begin{equation*}
    g (z) = \frac{1}{2} a_0 + \ds\sum_{k = 1}^m \Big(a_k \, \cos(k z) + b_k \, \sin(k z)\Big)
\end{equation*}
The derivatives are of the following based on the order $d$, where $d > 0$
\begin{align*}
    \dfrac{\dd^d g (z)}{\dd z^d} = 
    \begin{cases} 
    k^d \, \ds\sum_{k = 1}^m \Big(a_k \, \cos(k z) + b_k \, \sin(k z)\Big) &\mod(d,4) = 0\\
    k^d \, \ds\sum_{k = 1}^m \Big(-a_k \, \sin(k z) + b_k \, \cos(k z)\Big)  &\mod(d,4) = 1\\
    k^d \, \ds\sum_{k = 1}^m \Big(-a_k \, \cos(k z) - b_k \, \sin(k z)\Big)  &\mod(d,4) = 2\\
    k^d \, \ds\sum_{k = 1}^m \Big(a_k \, \sin(k z) - b_k \, \cos(k z)\Big)  &\mod(d,4) = 3
    \end{cases}
\end{align*}

%


\chapter{LINEAR LEAST-SQUARES METHODS}\label{chap:LS}
There are different numerical techniques to compute the linear least-squares (LS) solution of $A \, \B{\xi} = \B{b}$. These are:
\begin{itemize}
\item The Moore-Penrose inverse,
    \begin{equation*}
        \B{\xi} = (A\T \, A)^{-1} \, A\T \, \B{b}.
    \end{equation*}
\item QR decomposition,
    \begin{equation*}
         A = Q \, R \qquad \to \qquad \B{\xi} = R^{-1} \, Q\T \, \B{b},
    \end{equation*}
    where $Q$ is an orthogonal matrix and $R$ an upper triangular matrix.
\item SVD decomposition,
    \begin{equation*}
         A = U \,  \Sigma \, V\T \qquad \to \qquad \B{\xi} = A^+ \, \B{b} = V \,  \Sigma^+ \, U\T \, \B{b}
    \end{equation*}
     where $U$ and $V$ are two orthogonal matrices, and where $\Sigma^+$ is the pseudo-inverse of $\Sigma$, which is formed by replacing every non-zero diagonal entry by its reciprocal and transposing the resulting matrix.
\item Cholesky decomposition,
    \begin{equation*}
        A\T A \, \B{\xi} = U\T U \B{\xi} = A\T \, \B{b} \qquad \to \qquad \B{\xi} = U^{-1} \left(U^{-\mbox{\tiny T}} A\T \, \B{b}\right),
    \end{equation*}
    where $U$ is a upper triangular, and consequently, $U^{-1}$ and $U^{-\mbox{\tiny T}}$ are easy to compute.
\end{itemize}

One can reduce the condition number of the matrix to be inverted by scaling the columns of $A$,
\begin{equation*}
    A \left(S S^{-1}\right) \B{\xi} = \left(A S\right) \left(S^{-1} \B{\xi}\right) = B \, \B{\eta} = \B{b} \; \to \; \B{\xi} = S \, \B{\eta} = S \, (B\T B)^{-1} B\T \B{b},
\end{equation*}
where $S$ is the $m\times m$ scaling diagonal matrix whose diagonal elements are the inverse of the norms of the corresponding columns of $A$: $s_{kk} = |\B{a}_k|^{-1}$ or the maximum absolute value, $s_{kk} = \max\limits_{i} |a_{ki}|$. 

In this dissertation, the least-squares problem is solved using two methods: (1) the SVD decomposition introduced above (2) a combination of QR decomposition and the previously mentioned scaling, called the scaled QR approach. This approach performs the QR decomposition of the scaled matrix,
\begin{equation*}
    B = A \, S = Q \, R \qquad \to \qquad \B{\xi} = S \, R^{-1} \, Q\T \, \B{b}.
\end{equation*}

A weighted LS solution can be obtained by introducing an $n\times n$ diagonal matrix of weights, $W$. This technique exactly follows the Moore-Penrose inverse, however, the weight matrix $W$ allows for unequal emphasis given to the fitting of the solution,
\begin{equation*}
    W \, A \, \B{\xi} =W \, \B{b} \qquad \to \qquad \B{\xi} = (A\T \,W^2 \, A)^{-1} \, A\T \,W \, \B{b}.
\end{equation*}
Furthermore, it can also be shown that a simple scaling of the rows of $A$ is equivalent to weighted LS.

\pagebreak{}
%


\chapter{SOME COMMON CONSTRAINED EXPRESSIONS}\label{chap:commonSwithcingFunctions}

\begin{wblankBox}{Point and derivative}
\noindent \textbf{Constraints:}
\begin{equation*}
    y(x_0) = \kappa_1 \quad \text{and} \quad y_x(x_0) = \kappa_2
\end{equation*}
\textbf{Projection functionals:}
\begin{equation*}
\rho_1(x,g(x)) = \kappa_1 - g(x_0) \quad \text{and} \quad \rho_2(x,g(x)) = \kappa_2 - g_x(x_0)
\end{equation*}
\textbf{Switching functions:}
\begin{equation*}
    \phi_1(x) = 1 \quad \text{and} \quad \phi_2(x) = x - x_0
\end{equation*}
\end{wblankBox}

\begin{wblankBox}{Initial and final point}
\noindent \textbf{Constraints:}
\begin{equation*}
   y(x_0) = \kappa_1 
   \quad \text{and} \quad 
   y(x_f) = \kappa_2 
\end{equation*}
\textbf{Projection functionals:}
\begin{equation*}
\rho_1(x,g(x)) = \kappa_1 - g(x_0) 
\quad \text{and} \quad 
\rho_2(x,g(x)) = \kappa_2 - g(x_f)
\end{equation*}
\textbf{Switching functions:}
\begin{equation*}
    \phi_1 = \frac{x_f - x}{x_f - x_0}
    \quad \text{and} \quad 
    \phi_2 = \frac{x - x_0}{x_f - x_0}
\end{equation*}
\end{wblankBox}

\begin{wblankBox}{Initial point and final point/derivative}
\noindent \textbf{Constraints:}
\begin{equation*}
   y(x_0) = \kappa_1,
   \quad
   y(x_f) = \kappa_2,
   \quad \text{and} \quad 
   y_x(x_f) = \kappa_3 
\end{equation*}
\textbf{Projection functionals:}
\begin{equation*}
\rho_1(x,g(x)) = \kappa_1 - g(x_0),
\quad
\rho_2(x,g(x)) = \kappa_2 - g(x_f)
\quad \text{and} \quad 
\rho_3(x,g(x)) = \kappa_3 - g_x(x_f)
\end{equation*}
\textbf{Switching functions:}
\begin{align*}
    \phi_1(x) &= \frac{1}{(x_f - x_0)^2}\Big( x_f^2 - 2 x_f x + x^2 \Big) \\ 
    \phi_2(x) &= \frac{1}{(x_f - x_0)^2}\Big(x_0(x_0-2x_f) + 2x_f x - x^2\Big) \\ 
    \phi_3(x) &= \frac{1}{x_f - x_0}\Big( x_0 x_f - (x_0 + x_f)x + x^2\Big)
\end{align*}
\end{wblankBox}

\begin{wblankBox}{Initial point/derivative and final point}
\noindent \textbf{Constraints:}
\begin{equation*}
   y(x_0) = \kappa_1,
   \quad
   y_x(x_0) = \kappa_2,
   \quad \text{and} \quad 
   y(x_f) = \kappa_3 
\end{equation*}
\textbf{Projection functionals:}
\begin{equation*}
\rho_1(x,g(x)) = \kappa_1 - g(x_0),
\quad
\rho_2(x,g(x)) = \kappa_2 - g_x(x_0)
\quad \text{and} \quad 
\rho_3(x,g(x)) = \kappa_3 - g(x_f)
\end{equation*}
\textbf{Switching functions:}
\begin{align*}
    \phi_1(x) &= \frac{1}{(x_f - x_0)^2} \Big(x_f(x_f-2x_0) + 2 x_0 x - x^2 \Big)\\ 
    \phi_2(x) &= \frac{1}{x_f - x_0} \Big(-x_f x_0 + (x_f + x_0) x - x^2 \Big)\\ 
    \phi_3(x) &= \frac{1}{(x_f - x_0)^2} \Big(x_0^2 - 2 x_0 x + x^2 \Big)
\end{align*}
\end{wblankBox}

\begin{wblankBox}{Initial point/derivative and final point/derivative}
\noindent \textbf{Constraints:}
\begin{equation*}
   y(x_0) = \kappa_1,
   \quad
   y(x_f) = \kappa_2,
   \quad
   y_x(x_0) = \kappa_3,
   \quad \text{and} \quad 
   y_x(x_f) = \kappa_4
\end{equation*}
\textbf{Projection functionals:}
\begin{align*}
\rho_1(x,g(x)) &= \kappa_1 - g(x_0) \qquad \rho_3(x,g(x)) = \kappa_3 - g_x(x_0) \\
\rho_2(x,g(x)) &= \kappa_2 - g(x_f) \qquad \rho_4(x,g(x)) = \kappa_4 - g_x(x_f)
\end{align*}
\textbf{Switching functions:}
\begin{align*}
    \phi_1(x) &= \frac{1}{(x_f - x_0)^3} \Big( -x_f^2 (3  x_0- x_f) + 6  x_0  x_f x  -3 ( x_0+ x_f) x^2 + 2x^3 \Big)\\
	\phi_2(x) &= \frac{1}{(x_f - x_0)^3} \Big( - x_0^2 ( x_0-3  x_f) -6  x_0  x_f x  + 3 ( x_0+ x_f)x^2 -2x^3 \Big)\\
	\phi_3(x) &= \frac{1}{(x_f - x_0)^2} \Big(-x_0  x_f^2 +  x_f (2  x_0+ x_f)x -(x_0+2  x_f)x^2 + x^3 \Big)\\
	\phi_4(x) &= \frac{1}{(x_f - x_0)^2} \Big( - x_0^2  x_f + x_0 ( x_0+2  x_f) x - (2  x_0+ x_f) x^2 + x^3 \Big)
\end{align*}
\end{wblankBox}
 
\pagebreak{}
%

\chapter{ANALYTICAL TERMS FOR SELECTED PROBLEMS}\label{chap:app_jacob}

The analytical terms of this section are provided for completeness; however, in code, these terms are handled through JAX \cite{JaxGithub,JaxOriginalPaper} and the TFC toolbox (\href{https://github.com/leakec/tfc}{\textcolor{blue}{\underline{TFC GitHub}}}) \cite{TfcGithub} where all of the partial derivatives are taken by automatic differentiation.

\section{Linear-Nonlinear differential equation Jacobian terms from Section \ref{sect:s3_linearNonlinearDE}}\label{sect:app_linearNonlinearDE}
\begin{equation}\label{eq:Jacobian}
\mathbb{J}(\Xi) = \begin{bmatrix} \dfrac{\partial \p{1}{\tilde{F}(x_0,\Xi)}}{\partial \p{1}{\B{\xi}}} & \B{0}_{1 \times m} & \dfrac{\partial \p{1}{\tilde{F}(x_0,\Xi)}}{\partial y_1} & \dfrac{\partial \p{1}{\tilde{F}(x_0,\Xi)}}{\partial y_{1_x}} \\ 
\vdots & \vdots & \vdots & \vdots \\
\dfrac{\partial \p{1}{\tilde{F}(x_1,\Xi)}}{\partial \p{1}{\B{\xi}}} & \B{0}_{1 \times m} & \dfrac{\partial \p{1}{\tilde{F}(x_1,\Xi)}}{\partial y_1} & \dfrac{\partial \p{1}{\tilde{F}(x_1,\Xi)}}{\partial y_{1_x}}\\
\B{0}_{1 \times m} & \dfrac{\partial \p{2}{\tilde{F}(x_1,\Xi)}}{\partial \p{2}{\B{\xi}}} & \dfrac{\partial \p{2}{\tilde{F}(x_1,\Xi)}}{\partial y_1} & \dfrac{\partial \p{2}{\tilde{F}(x_1,\Xi)}}{\partial y_{1_x}} \\
\vdots & \vdots & \vdots & \vdots \\
\B{0}_{1 \times m} & \dfrac{\partial \p{2}{\tilde{F}(x_f,\Xi)}}{\partial \p{2}{\B{\xi}}} & \dfrac{\partial \p{2}{\tilde{F}(x_f,\Xi)}}{\partial y_1} & \dfrac{\partial \p{2}{\tilde{F}(x_f,\Xi)}}{\partial y_{1_x}}
\end{bmatrix}
\end{equation}
For this problem all terms of Equation \eqref{eq:Jacobian} are provided below:
\begin{align*}
    \dfrac{\partial \p{1}{\tilde{F}}}{\partial \p{1}{\B{\xi}}} &= \Big[c^2 \,\B{h}_{zz}(z)- \p{1}{\phi}_{1_{xx}} \B{h}(z_0) - \p{1}{\phi}_{2_{xx}} \B{h}(z_1) - \p{1}{\phi}_{3_{xx}} c \, \B{h}_{z}(z_1) \\ & \qquad + \B{h}(z)- \p{1}{\phi}_1\B{h}(z_0) -\p{1}{\phi}_2 \B{h}(z_1) - \p{1}{\phi}_3 c \, \B{h}_{z}(z_1)\Big]\T \nonumber\\
    \dfrac{\partial \p{1}{\tilde{F}}}{\partial y_1} &= \p{1}{\phi}_{2_{xx}}(x) + \p{1}{\phi}_{2}(x) \\
    \dfrac{\partial \p{1}{\tilde{F}}}{\partial y_{1_x}} &= \p{1}{\phi}_{3_{xx}}(x) + \p{1}{\phi}_{3}(x) \\
    \dfrac{\partial \p{2}{\tilde{F}}}{\partial \p{2}{\B{\xi}}} &= \Big[c^2 \, \B{h}_{zz}(z)- \p{2}{\phi}_{1_{xx}} \B{h}(z_1) - \p{2}{\phi}_{2_{xx}} c \, \B{h}(z_1) - \p{2}{\phi}_{3_{xx}} \B{h}(z_f) \\ & \qquad + \p{2}{y}_x \Big(\B{h}(z)- \p{2}{\phi}_1 \B{h}(z_1) - \p{2}{\phi}_2  c \, \B{h}(z_1) - \p{2}{\phi}_3 \B{h}(z_f) \Big) \nonumber\\ & \qquad + \p{2}{y}\Big(c \, \B{h}(z)- \p{2}{\phi}_{1_x}\B{h}(z_1) - \p{2}{\phi}_{2_x} c \, \B{h}(z_1) - \p{2}{\phi}_{3_x}\B{h}(z_f) \Big)\Big]\T \nonumber\\
    \dfrac{\partial \p{2}{\tilde{F}}}{\partial y_1} &= \p{2}{\phi}_{1_{xx}}(x) + \p{2}{y} \p{2}{\phi}_{1_{x}}(x) + \p{2}{y}_x \p{2}{\phi}_1(x)  \\
    \dfrac{\partial \p{2}{\tilde{F}}}{\partial y_{1_x}} &= \p{2}{\phi}_{2_{xx}}(x)  + \p{2}{y} \p{2}{\phi}_{2_{x}}(x) + \p{2}{y}_x \p{2}{\phi}_2(x)
\end{align*}
where $\Xi$ is the vector of unknown coefficients such that,
\begin{equation*}
    \Xi = \begin{Bmatrix} \p{1}{\B{\xi}}\T & \p{2}{\B{\xi}}\T & y_1 & y_{1_x} \end{Bmatrix}\T,
\end{equation*}

\section{Convection-diffusion equation from Section \ref{sec:s3_condif} }\label{sect:app_condif}
The Jacobian is of the form,
\begin{equation*}
    \mathbb{J}(\Xi) = \begin{bmatrix} \frac{\partial \p{1}{\mathbb{L}}}{\partial \p{1}{\B{\xi}}} & \B{0}_{N \times m} & \frac{\partial \p{1}{\mathbb{L}}}{\partial y_1} &\frac{\partial \p{1}{\mathbb{L}}}{\partial y_{1_x}} &\frac{\partial \p{1}{\mathbb{L}}}{\partial \bar{c}} \\ \B{0}_{N \times m} & \frac{\partial \p{2}{\mathbb{L}}}{\partial \p{1}{\B{\xi}}} & \frac{\partial \p{2}{\mathbb{L}}}{\partial y_1} &\frac{\partial \p{2}{\mathbb{L}}}{\partial y_{1_x}} &\frac{\partial \p{2}{\mathbb{L}}}{\partial \bar{c}}  \end{bmatrix}
\end{equation*}
where the following equations are the detailed Jacobian terms from the convection-diffusion equation from Section \ref{sec:s3_condif}. For clarity, the constrained expressions are,
\begin{align*}
    \p{1}{y}(z,\p{1}{\B{\xi}}) = \Big(\B{h}(z)- \p{1}{\phi}_1(z) \B{h}(z_0) &+ \p{1}{\phi}_2(z) \B{h}(z_f) + \p{1}{\phi}_3(z) \B{h}_z(z_f) \Big)\T \p{1}{\B{\xi}} \\ &+ \p{1}{\phi}_1(z) y_0 + \p{1}{\phi}_2(z) y_1 + \p{1}{\phi}_3(z) \frac{y_{1_x}}{\p{1}{c}} \\
    \p{2}{y}(z,\p{2}{\B{\xi}}) = \Big(\B{h}(z)- \p{2}{\phi}_1(z) \B{h}(z_0) &+ \p{2}{\phi}_2(z) \B{h}_z(z_0) + \p{2}{\phi}_3(z) \B{h}(z_f) \Big)\T \p{2}{\B{\xi}} \\&+ \p{2}{\phi}_1(z) y_1 + \p{2}{\phi}_2(z) \frac{y_{1_x}}{\p{2}{c}} + \p{2}{\phi}_3(z) y_f
\end{align*}
where the loss vectors of each segment are,
\begin{equation*}
    \p{1}{\mathbb{L}}(\Xi) =  \begin{Bmatrix} \p{1}{\tilde{F}}(z_0,\Xi) \\ \vdots \\ \p{1}{\tilde{F}}(z_f,\Xi) \end{Bmatrix} = \begin{Bmatrix}  \bar{c}^2 \p{1}{y}_{xx}(z_0,\Xi) - \text{Pe}\, \bar{c} \, \p{1}{y}_x(z_0,\Xi) \\ \vdots \\ \bar{c}^2 \p{1}{y}_{xx}(z_f,\Xi) - \text{Pe}\, \bar{c} \, \p{1}{y}_x(z_f,\Xi) \end{Bmatrix} 
\end{equation*}
and
\begin{equation*}
    \p{2}{\mathbb{L}}(\Xi) =  \begin{Bmatrix} \p{2}{\tilde{F}}(z_0,\Xi) \\ \vdots \\ \p{2}{\tilde{F}}(z_f,\Xi)   \end{Bmatrix} = \begin{Bmatrix}  \Big(\dfrac{\bar{c} \Delta z}{\bar{c} - \Delta z}\Big)^2 \p{2}{y}_{xx}(z_0,\Xi) - \text{Pe}\, \Big(\dfrac{\bar{c} \Delta z}{\bar{c} - \Delta z}\Big) \, \p{2}{y}_x(z_0,\Xi) \\ \vdots \\ \Big(\dfrac{\bar{c} \Delta z}{\bar{c} - \Delta z}\Big)^2 \p{2}{y}_{xx}(z_f,\Xi) - \text{Pe}\, \Big(\dfrac{\bar{c} \Delta z}{\bar{c} - \Delta z}\Big) \, \p{2}{y}_x(z_f,\Xi)   \end{Bmatrix}. 
\end{equation*}
The following equations are the Jacobians of the loss vectors with respect to the unknowns:

\begin{equation*}
    \frac{ \p{1}{\mathbb{L}}(\Xi)}{\partial \p{1}{\B{\xi}}} = \begin{bmatrix} \Big[\bar{c}^2 \dfrac{\partial \p{1}{y}_{zz}}{\partial \p{1}{\B{\xi}}}(z_0) - \text{Pe}\, \bar{c} \, \dfrac{\partial \p{1}{y}_{z}}{\partial \p{1}{\B{\xi}}}(z_0)\Big]\T
    \\ \vdots \\ 
    \Big[\bar{c}^2 \dfrac{\partial \p{1}{y}_{zz}}{\partial \p{1}{\B{\xi}}}(z_f) - \text{Pe}\, \bar{c} \, \dfrac{\partial \p{1}{y}_{z}}{\partial \p{1}{\B{\xi}}}(z_f)\Big]\T \end{bmatrix}
\end{equation*}
\begin{equation*}
    \frac{ \p{1}{\mathbb{L}}(\Xi)}{\partial y_1} = \begin{Bmatrix} \bar{c}^2 \p{1}{\phi}_{2_{zz}}(z_0) -  \text{Pe}\, \bar{c} \p{1}{\phi}_{2_{z}}(z_0) \\ \vdots \\ \bar{c}^2 \p{1}{\phi}_{2_{zz}}(z_f) -  \text{Pe}\, \bar{c} \p{1}{\phi}_{2_{z}}(z_f)\end{Bmatrix}
\end{equation*}
\begin{equation*}
    \frac{ \p{1}{\mathbb{L}}(\Xi)}{\partial y_{1_x}} = \begin{Bmatrix} \bar{c} \p{1}{\phi}_{3_{zz}}(z_0) -  \text{Pe}\, \p{1}{\phi}_{3_{z}}(z_0) \\ \vdots \\ \bar{c} \p{1}{\phi}_{3_{zz}}(z_f) -  \text{Pe}\, \p{1}{\phi}_{3_{z}}(z_f)\end{Bmatrix}
\end{equation*}
\begin{equation*}
    \frac{ \p{1}{\mathbb{L}}(\Xi)}{\partial \bar{c}} = \begin{Bmatrix} 2 \bar{c} \p{1}{y}_{zz}(z_0)- \p{1}{\phi}_{3_{zz}}(z_0) y_{1_x} -\text{Pe} \,  \p{1}{y}_{z}(z_0) +  \text{Pe} \dfrac{\p{1}{\phi}_{3_{z}}(z_0) y_{1_x}}{\bar{c}} \\ \vdots \\ 2 \bar{c} \p{1}{y}_{zz}(z_f)- \p{1}{\phi}_{3_{zz}}(z_f) y_{1_x} -\text{Pe} \,  \p{1}{y}_{z}(z_f) +  \text{Pe} \dfrac{\p{1}{\phi}_{3_{z}}(z_f) y_{1_x}}{\bar{c}}  \end{Bmatrix}
\end{equation*}

\begin{equation*}
    \frac{ \p{2}{\mathbb{L}}(\Xi)}{\partial \p{2}{\B{\xi}}} = \begin{bmatrix} \Big[ \Big(\dfrac{\bar{c} \Delta z}{\bar{c} - \Delta z}\Big)^2 \dfrac{\partial \p{2}{y}_{zz}}{\partial \p{2}{\B{\xi}}}(z_0) - \text{Pe}\, \Big(\dfrac{\bar{c} \Delta z}{\bar{c} - \Delta z}\Big) \, \dfrac{\partial \p{2}{y}_{z}}{\partial \p{2}{\B{\xi}}}(z_0)  \Big]\T
    \\ \vdots \\ 
    \Big[ \Big(\dfrac{\bar{c} \Delta z}{\bar{c} - \Delta z}\Big)^2 \dfrac{\partial \p{2}{y}_{zz}}{\partial \p{2}{\B{\xi}}}(z_f) - \text{Pe}\, \Big(\dfrac{\bar{c} \Delta z}{\bar{c} - \Delta z}\Big) \, \dfrac{\partial \p{2}{y}_{z}}{\partial \p{2}{\B{\xi}}}(z_f)  \Big]\T \end{bmatrix}
\end{equation*}
\normalsize
\begin{equation*}
    \frac{ \p{2}{\mathbb{L}}(\Xi)}{\partial y_1} = \begin{Bmatrix} \Big(\dfrac{\bar{c} \Delta z}{\bar{c} - \Delta z}\Big)^2 \p{2}{\phi}_{1_{zz}}(z_0) -  \text{Pe}\, \Big(\dfrac{\bar{c} \Delta z}{\bar{c} - \Delta z}\Big) \p{2}{\phi}_{1_{z}}(z_0) \\ \vdots \\ \Big(\dfrac{\bar{c} \Delta z}{\bar{c} - \Delta z}\Big)^2 \p{2}{\phi}_{1_{zz}}(z_f) -  \text{Pe}\, \Big(\dfrac{\bar{c} \Delta z}{\bar{c} - \Delta z}\Big) \p{2}{\phi}_{1_{z}}(z_f)\end{Bmatrix}
\end{equation*}
\begin{equation*}
    \frac{ \p{2}{\mathbb{L}}(\Xi)}{\partial y_{1_x}} = \begin{Bmatrix} \Big(\dfrac{\bar{c} \Delta z}{\bar{c} - \Delta z}\Big) \p{2}{\phi}_{2_{zz}}(z_0) -  \text{Pe}\,  \p{2}{\phi}_{2_{z}}(z_0) \\ \vdots \\ \Big(\dfrac{\bar{c} \Delta z}{\bar{c} - \Delta z}\Big) \p{2}{\phi}_{2_{zz}}(z_f) -  \text{Pe}\,  \p{2}{\phi}_{2_{z}}(z_f) \end{Bmatrix}
\end{equation*}
\footnotesize
\begin{equation*}
    \frac{ \p{2}{\mathbb{L}}(\Xi)}{\partial \bar{c}} = \begin{Bmatrix} -\dfrac{\Delta z^2}{(\bar{c} - \Delta z)^2} \left[2 \Big(\dfrac{\bar{c} \Delta z}{\bar{c} - \Delta z}\Big) \p{2}{y}_{zz}(z_0) - \p{2}{\phi}_{2_{zz}}(z_0) y_{1_x} -\text{Pe} \,  \p{2}{y}_{z}(z_0) +  \text{Pe} \dfrac{\p{2}{\phi}_{2_{z}}(z_0) y_{1_x}}{\Big(\dfrac{\bar{c} \Delta z}{\bar{c} - \Delta z}\Big)}\right] \\ \vdots \\ -\dfrac{\Delta z^2}{(\bar{c} - \Delta z)^2} \left[ 2 \Big(\dfrac{\bar{c} \Delta z}{\bar{c} - \Delta z}\Big) \p{2}{y}_{zz}(z_f) - \p{2}{\phi}_{2_{zz}}(z_f) y_{1_x} -\text{Pe} \,  \p{2}{y}_{z}(z_f) +  \text{Pe} \dfrac{\p{2}{\phi}_{2_{z}}(z_f) y_{1_x}}{\Big(\dfrac{\bar{c} \Delta z}{\bar{c} - \Delta z}\Big)} \right] \end{Bmatrix}
\end{equation*}
\normalsize

\section{Terms for Outer-loop approach in the energy optimal landing problem from Section \ref{sec:s5_outLoop}}\label{sect:app_sec:s5_outLoop}
By discretizing the domain the linear system becomes,
\begin{equation*}
    \begin{bmatrix} \mathbb{A} & \B{0}_{N\times m} & \B{0}_{N\times m} & -\mathbb{C} & \B{0}_{N\times 2} & \B{0}_{N\times 2} \\
    \B{0}_{N\times m} & \mathbb{A} & \B{0}_{N\times m} & \B{0}_{N\times 2} & -\mathbb{C} & \B{0}_{N\times 2} \\
    \B{0}_{N\times m} & \B{0}_{N\times m} & \mathbb{A} & \B{0}_{N\times 2} & \B{0}_{N\times 2} & -\mathbb{C}\\
    \end{bmatrix}  \begin{Bmatrix} \B{\xi}_1 \\ \B{\xi}_2 \\ \B{\xi}_3 \\ \B{\xi}_{u_1} \\ \B{\xi}_{u_2} \\ \B{\xi}_{u_3} \end{Bmatrix} =  - \begin{Bmatrix} \mathbb{B}_1 \\ \mathbb{B}_2 \\ \mathbb{B}_3 \end{Bmatrix}
\end{equation*}
where $\mathbb{A}$, $\mathbb{B}_i$, $\mathbb{C}$ are defined as,
\begin{align*}
    \mathbb{A} = \begin{bmatrix} \Big(c^2 \B{h}_{zz}(z_0) - \ddot{\phi}_1(t_0) \B{h}(z_0) - \ddot{\phi}_2(t_0) \B{h}(z_f) - \ddot{\phi}_3(t_0) c\B{h}_z(z_0) - \ddot{\phi}_4(t_0) c\B{h}_z(z_f) \Big)\T \\ \vdots \\ \Big(c^2 \B{h}_{zz}(z_k) - \ddot{\phi}_1(t_k) \B{h}(z_0) - \ddot{\phi}_2(t_k) \B{h}(z_f) - \ddot{\phi}_3(t_k) c\B{h}_z(z_0) - \ddot{\phi}_4(t_k) c\B{h}_z(z_f) \Big)\T \\ \vdots \\ \Big(c^2 \B{h}_{zz}(z_f) - \ddot{\phi}_1(t_f) \B{h}(z_0) - \ddot{\phi}_2(t_f) \B{h}(z_f) - \ddot{\phi}_3(t_f) c\B{h}_z(z_0) - \ddot{\phi}_4(t_f) c\B{h}_z(z_f) \Big)\T \end{bmatrix}
\end{align*}
\begin{align*}
    \mathbb{B}_i = \begin{bmatrix} \ddot{\phi}_1(t_0) r_{0_i} + \ddot{\phi}_2(t_0) r_{f_i} + \ddot{\phi}_3(t_0) v_{0_i} + \ddot{\phi}_4(t_0) v_{f_i} - a_{g_i} \\ \vdots \\ \ddot{\phi}_1(t_k) r_{0_i} + \ddot{\phi}_2(t_k) r_{f_i} + \ddot{\phi}_3(t_k) v_{0_i} + \ddot{\phi}_4(t_k) v_{f_i} - a_{g_i} \\ \vdots \\ \ddot{\phi}_1(t_f) r_{0_i} + \ddot{\phi}_2(t_f) r_{f_i} + \ddot{\phi}_3(t_f) v_{0_i} + \ddot{\phi}_4(t_f) v_{f_i} - a_{g_i} \end{bmatrix} \quad \mathbb{C}= \begin{bmatrix} \B{h}_{u}\T(z_0) \\ \vdots \\ \B{h}_{u}\T(z_k) \\ \vdots \\ \B{h}_{u}\T(z_f) \end{bmatrix}.
\end{align*}

\section{Single-loop approach Jacobian terms in the energy optimal landing problem from Section \ref{sec:s5_singLoop}}\label{sect:app_sec:s5_singLoop}
The partial derivatives for the state loss function, $\mathbb{L}_i$, when $i=j$ are,
\begin{eqnarray*}
    \frac{\partial \mathbb{L}_i}{\partial \B{\xi}_j} &=& b^4 \Big(\B{h}_{zz}(z) - \phiz{z}{\phi}_{1_{zz}} \B{h}_0 - \phiz{z}{\phi}_{2_{zz}} \B{h}_f - \phiz{z}{\phi}_{3_{zz}} c\B{h}_{z}(z_0) - \phiz{z}{\phi}_{4_{zz}} c\B{h}_{z}(z_f) \Big)\T \\
    \frac{\partial \mathbb{L}_i}{\partial \B{\xi}_{u_j}} &=& \B{h}_{u}\T.
\end{eqnarray*}
If $i \neq j$
\begin{eqnarray*}
    \frac{\partial \mathbb{L}_i}{\partial \B{\xi}_j} &=& \B{0}_{N \times m} \\
    \frac{\partial \mathbb{L}_i}{\partial \B{\xi}_{u_j}} &=& \B{0}_{N \times 2}
\end{eqnarray*}
and
\begin{align*}
    \frac{\partial \mathbb{L}_i}{\partial b} = 4 b^3\Big[ \Big(\B{h}_{zz}(z) - \phiz{z}{\phi}_{1_{zz}} \B{h}_0 - \phiz{z}{\phi}_{2_{zz}} \B{h}_f &- \phiz{z}{\phi}_{3_{zz}} c\B{h}_{z}(z_0) - \phiz{z}{\phi}_{4_{zz}} c\B{h}_{z}(z_f) \Big)\T \\&+ \phiz{z}{\phi}_{1_{zz}} r_{0_i} + \phiz{z}{\phi}_{2_{zz}} r_{f_i}\Big] \\ &+ 2b\Big[ \phiz{z}{\phi}_{3_{zz}} v_{0_i} + \phiz{z}{\phi}_{4_{zz}} v_{f_i} \Big].
\end{align*}
Similarly, the partial derivatives for $\mathbb{L}_H$ are,
\begin{eqnarray*}
    \frac{\partial \mathbb{L}_H}{\partial \B{\xi}_j} &=& \B{0}_{1 \times m} \\
    \frac{\partial \mathbb{L}_H}{\partial \B{\xi}_{u_j}} = \frac{\partial \mathbb{L}_H}{\partial u_j} \frac{\partial u_j}{\partial \B{\xi}_{u_j}} &=&  \B{h}_{u}\T \Big(a_{g_j} -u_j(z_f)  \Big) \\ 
    \frac{\partial \mathbb{L}_H}{\partial b} &=& 0.
\end{eqnarray*}
Combining these into a single Jacobian term leads to,
\begin{equation*}
    \mathbb{J} = \begin{bmatrix}\frac{\partial \mathbb{L}_1}{\partial \B{\xi}_1} & \B{0}_{N \times m} & \B{0}_{N \times m} & \frac{\partial \mathbb{L}_1}{\partial \B{\xi}_{u_1}} & \B{0}_{N \times 2} & \B{0}_{N \times 2} &  \frac{\partial \mathbb{L}_1}{\partial b} \\  \B{0}_{N \times m} & \frac{\partial \mathbb{L}_2}{\partial \B{\xi}_2} & \B{0}_{N \times m} & \B{0}_{N \times 2} & \frac{\partial \mathbb{L}_2}{\partial \B{\xi}_{u_2}} & \B{0}_{N \times 2} &  \frac{\partial \mathbb{L}_2}{\partial b} \\ \B{0}_{N \times m} & \B{0}_{N \times m} & \frac{\partial \mathbb{L}_3}{\partial \B{\xi}_3} & \B{0}_{N \times 2} & \B{0}_{N \times 2} & \frac{\partial \mathbb{L}_3}{\partial \B{\xi}_{u_3}} &  \frac{\partial \mathbb{L}_1}{\partial b} \\ \B{0}_{1 \times m} & \B{0}_{1 \times m} & \B{0}_{1 \times m} & \frac{\partial \mathbb{L}_H}{\partial \B{\xi}_1} & \frac{\partial \mathbb{L}_H}{\partial \B{\xi}_2} & \frac{\partial \mathbb{L}_H}{\partial \B{\xi}_3} & \frac{\partial \mathbb{L}_H}{\partial b} \end{bmatrix}_{(3N + 1) \times (3m + 7)}
\end{equation*}
with the augmented loss function and unknown vector defined as
\begin{align*}
    \mathbb{L} &= \begin{Bmatrix} \mathbb{L}_1\T & \mathbb{L}_2\T &\mathbb{L}_3\T & \mathbb{L}_H\end{Bmatrix}\T_{(3N + 1) \times 1} \\ 
    \Xi &= \begin{Bmatrix} \B{\xi}_1\T & \B{\xi}_2\T & \B{\xi}_3\T & \B{\xi}_{u_1}\T & \B{\xi}_{u_2}\T & \B{\xi}_{u_3}\T &  b \end{Bmatrix}\T _{(3m + 7) \times 1} .
\end{align*}

\section{Fuel-Optimal Landing from Section \ref{sect:s6_fol_jac}}\label{sect:app_FOL}
In the fuel-optimal landing problem the analytical partial derivatives of the state loss function are:
\begin{equation*}
\begin{aligned}
    \frac{\partial \p{s}{\mathbb{L}}_i}{\partial \p{s}{\B{\xi}_i}} &= \p{s}{\left(c^2\B{h}_{zz} - \ddot{\phi}_1(t)\B{h}(z_0) - \ddot{\phi}_2(t)\B{h}(z_f) - \ddot{\phi}_3(t)c \B{h}_{z}(z_0) - \ddot{\phi}_4(t)c \B{h}_{z}(z_f) \right)\T} \\ 
    \frac{\partial \p{1}{\mathbb{L}}_i}{\partial r_{1_i}} &= \p{1}{\ddot{\phi}}_2(t) \\
    \frac{\partial \p{1}{\mathbb{L}}_i}{\partial v_{1_i}} &= \p{1}{\ddot{\phi}}_4(t) \\  \frac{\partial \p{2}{\mathbb{L}}_i}{\partial r_{1_i}} &= \p{2}{\ddot{\phi}}_1(t), \qquad  \frac{\partial \p{2}{\mathbb{L}}_i}{\partial v_{1_i}} = \p{2}{\ddot{\phi}}_3(t) \\ 
    \frac{\partial \p{2}{\mathbb{L}}_i}{\partial r_{2_i}} &= \p{2}{\ddot{\phi}}_2(t), \qquad \frac{\partial \p{2}{\mathbb{L}}_i}{\partial v_{2_i}} = \p{2}{\ddot{\phi}}_4(t) \\
    \frac{\partial \p{3}{\mathbb{L}}_i}{\partial r_{2_i}} &= \p{3}{\ddot{\phi}}_1(t) \\ \frac{\partial \p{3}{\mathbb{L}}_i}{\partial v_{2_i}} &= \p{3}{\ddot{\phi}}_3(t).
\end{aligned}
\end{equation*}
For the costate portion, if $i = j$
\begin{equation*}
   \frac{\partial \mathbb{L}_i}{\partial \B{\xi_{\lambda_i}}} = \beta(t) \left[ \left(\displaystyle\sum_{j=1}^3 \lambda^2_{v_j}\right)^{-1/2} - \lambda^2_{v_i} \, \left(\displaystyle\sum_{j=1}^3 \lambda^2_{v_j}\right)^{-3/2}\right] \B{h}_{\lambda}\T 
\end{equation*}
if $i \neq j$
\begin{equation*}
    \frac{\partial \mathbb{L}_i}{\partial \B{\xi_{\lambda_j}}} = \beta(t) \left[ - \lambda_{v_i} \, \lambda_{v_j} \, \left(\displaystyle\sum_{j=1}^3 \lambda^2_{v_j}\right)^{-3/2}\right] \B{h}_{\lambda}\T.
\end{equation*}

For the loss function associated with the transversality conditions for the Hamiltonian, $\mathbb{L}_H$, the only non-zero partial is with respect to $\B{\xi}_\lambda$, which is defined by
\begin{equation*}
   \frac{\partial \mathbb{L}_H}{\partial \B{\xi_{\lambda_i}}} = \left[a_{g_i} - \beta(t_f) \, \lambda_{v_i}(t_f) \, \left(\sum_{j=1}^3 \lambda^2_{v_j}(t_f)\right)^{-1/2}\right] \B{h}_{\lambda}\T(t_f).
\end{equation*}
The augmented loss functions for the discretized points become
\begin{equation*}
    \mathbb{L} = \begin{Bmatrix} \p{1}{\mathbb{L}}_1\T & \p{1}{\mathbb{L}}_2\T & \p{1}{\mathbb{L}}_3\T & \p{2}{\mathbb{L}}_1\T & \p{2}{\mathbb{L}}_2\T & \p{2}{\mathbb{L}}_3\T & \p{3}{\mathbb{L}}_1\T & \p{3}{\mathbb{L}}_2\T & \p{3}{\mathbb{L}}_3\T & \mathbb{L}_H\end{Bmatrix}_{(\{9N+1\}\times 1)}\T
\end{equation*}
with the unknown vector
\begin{align*}
    \Xi = \Big\{ &\p{1}{\B{\xi}}_1\T \quad \p{1}{\B{\xi}}_2\T \quad \p{1}{\B{\xi}}_3\T \quad \p{2}{\B{\xi}}_1\T \quad \p{2}{\B{\xi}}_2\T \quad \p{2}{\B{\xi}}_3\T \quad \p{3}{\B{\xi}}_1\T \quad \p{3}{\B{\xi}}_2\T \quad \p{3}{\B{\xi}}_3\T \\ &\B{\xi}_{\lambda_1}\T \quad \B{\xi}_{\lambda_2}\T \quad \B{\xi}_{\lambda_3}\T  \quad \B{r}_1\T \quad \B{v}_1\T \quad  \B{r}_2\T \quad \B{v}_2\T  \Big\}_{(9m + 18)}\T.
\end{align*}
All partials can be combined into one augmented matrix,
\begin{equation}\label{eq:app_aug_J}
    \mathbb{J} = \begin{bmatrix} 
    \p{1}{J}_{\B{\xi}} & \B{0}_{(3N\times 3m)} & \B{0}_{(3N\times 3m)} &  \p{1}{J}_{\B{\xi}_\lambda} & \p{1}{J}_{r_1, v_1} & \B{0}_{(3N\times 6)} \\ 
    \B{0}_{(3N\times 3m)}  &   \p{2}{J}_{\B{\xi}}  & \B{0}_{(3N\times 3m)} &  \p{2}{J}_{\B{\xi}_\lambda} & \p{2}{J}_{r_1, v_1} & \p{2}{J}_{r_2, v_2}\\ 
    \B{0}_{(3N\times 3m)} &  \B{0}_{(3N\times 3m)}  &  \p{3}{J}_{\B{\xi}}  & \p{3}{J}_{\B{\xi}_\lambda} & \B{0}_{(3N\times 6)} & \p{3}{J}_{r_2, v_2} \\ \B{0}_{(1\times 3m)}  & \B{0}_{(1\times 3m)}  & \B{0}_{(1\times 3m)} & J_H & \B{0}_{(1\times 6)} & \B{0}_{(1\times 6)} \end{bmatrix}_{(\{9N+1\}\times \{9m + 18\})}
\end{equation}
The terms of Equation (\ref{eq:app_aug_J}) are defined by the following equations:
\begin{equation*}
    \p{s}{J}_{\B{\xi}} = \begin{bmatrix} \frac{\partial \p{s}{\mathbb{L}}_1}{\partial \p{s}{\B{\xi}_1}} & \B{0} & \B{0} \\ \B{0} & \frac{\partial \p{s}{\mathbb{L}}_2}{\partial \p{s}{\B{\xi}_2}} & \B{0} \\ \B{0} & \B{0} & \frac{\partial \p{s}{\mathbb{L}}_3}{\partial \p{s}{\B{\xi}_3}} \end{bmatrix}_{(3N\times 3m)},
    \quad  \p{s}{J}_{\B{\xi}_\lambda} = \begin{bmatrix} J_{\B{\xi}_{\lambda_{11}}} & J_{\B{\xi}_{\lambda_{12}}} & J_{\B{\xi}_{\lambda_{13}}} \\ J_{\B{\xi}_{\lambda_{21}}} & J_{\B{\xi}_{\lambda_{22}}} & J_{\B{\xi}_{\lambda_{23}}} \\ J_{\B{\xi}_{\lambda_{31}}} & J_{\B{\xi}_{\lambda_{32}}} & J_{\B{\xi}_{\lambda_{33}}}  \end{bmatrix}_{(3N\times 6)}
\end{equation*}

\begin{eqnarray*}
    \p{1}{J}_{r_1,v_1} = \begin{bmatrix} \p{1}{ \ddot{\phi}}_2 & \B{0} & \B{0} & \p{1}{ \ddot{\phi}}_4 & \B{0} & \B{0} \\ \B{0} & \p{1}{ \ddot{\phi}}_2 & \B{0} & \B{0} & \p{1}{ \ddot{\phi}}_4 & \B{0} \\ \B{0} & \B{0} & \p{1}{ \ddot{\phi}}_2 & \B{0} & \B{0} & \p{1}{ \ddot{\phi}}_4 \end{bmatrix}_{(3N\times 6)} \\ \p{2}{J}_{r_1,v_1} = \begin{bmatrix} \p{1}{ \ddot{\phi}}_1 & \B{0} & \B{0} & \p{1}{ \ddot{\phi}}_3 & \B{0} & \B{0} \\ \B{0} & \p{1}{ \ddot{\phi}}_1 & \B{0} & \B{0} & \p{1}{ \ddot{\phi}}_3 & \B{0} \\ \B{0} & \B{0} & \p{1}{ \ddot{\phi}}_1 & \B{0} & \B{0} & \p{1}{ \ddot{\phi}}_3 \end{bmatrix}_{(3N\times 6)}\\
    \p{2}{J}_{r_2,v_2} = \begin{bmatrix} \p{2}{ \ddot{\phi}}_2 & \B{0} & \B{0} & \p{2}{ \ddot{\phi}}_4 & \B{0} & \B{0} \\ \B{0} & \p{2}{ \ddot{\phi}}_2 & \B{0} & \B{0} & \p{2}{ \ddot{\phi}}_4 & \B{0} \\ \B{0} & \B{0} & \p{2}{ \ddot{\phi}}_2 & \B{0} & \B{0} & \p{2}{ \ddot{\phi}}_4 \end{bmatrix}_{(3N\times 6)} \\
    \p{3}{J}_{r_2,v_2} = \begin{bmatrix} \p{3}{ \ddot{\phi}}_1 & \B{0} & \B{0} & \p{3}{ \ddot{\phi}}_3 & \B{0} & \B{0} \\ \B{0} & \p{3}{ \ddot{\phi}}_1 & \B{0} & \B{0} & \p{3}{ \ddot{\phi}}_3 & \B{0} \\ \B{0} & \B{0} & \p{3}{ \ddot{\phi}}_1 & \B{0} & \B{0} & \p{3}{ \ddot{\phi}}_3 \end{bmatrix}_{(3N\times 6)}
\end{eqnarray*}
\begin{equation*}
    J_H = \begin{bmatrix} \dfrac{\partial \mathbb{L}_H}{\partial \B{\xi_{\lambda_1}}}, & \quad \dfrac{\partial \mathbb{L}_H}{\partial \B{\xi_{\lambda_2}}}, & \quad \dfrac{\partial \mathbb{L}_H}{\partial \B{\xi_{\lambda_3}}}\end{bmatrix}_{(1\times 6)}.
\end{equation*}

\end{appendices}

\end{document}